\documentclass[3p,times]{elsarticle}

\usepackage{latexsym}

\usepackage{indentfirst}
\usepackage[english]{babel}
\usepackage{cancel}

\usepackage{lipsum}
\usepackage{array}
\usepackage{color} 
\usepackage{tabularx}
\usepackage{graphicx} 
\usepackage{amsmath}
\usepackage{amssymb}
\usepackage{amsfonts}	
\usepackage{moreverb}
\usepackage{dsfont}
\usepackage{tipa}
\usepackage{upgreek}
\usepackage{bm}
\usepackage{multirow}
\usepackage{soul}
\usepackage{ textcomp }
\usepackage[dvipsnames]{xcolor}
\usepackage{mathtools}
\usepackage{setspace}
\usepackage{algorithm}
\usepackage{algorithmic}
\allowdisplaybreaks

\newcommand\BibTeX{{\rmfamily B\kern-.05em \textsc{i\kern-.025em b}\kern-.08em
T\kern-.1667em\lower.7ex\hbox{E}\kern-.125emX}}

\usepackage{hyperref} 
\hypersetup{
    colorlinks=true,                          
    linkcolor=blue, 
    citecolor=red, 
    urlcolor=blue  } 

\newcommand{\x}{\mbf{x}}

\newcommand{\xcn}{\mbf{x}_{\mbf{c}}^n}               
\newcommand{\xcin}{\mbf{x}_{\mbf{c}_i}^n}    
\newcommand{\xcinp}{\mbf{x}_{\mbf{c}_i}^{n+1}}    

\newcommand{\phitilde}{\tilde{\phi}}

\newcommand{\xbin}{\mbf{x}_{\mbf{b}_i}^n}    
\newcommand{\xbinp}{\mbf{x}_{\mbf{b}_i}^{n+1}}    
\newcommand{\xbixn}{{x}_{{b}_i}^n}    
\newcommand{\xbiyn}{{y}_{{b}_i}^n}

\newcommand{\mbf}[1]{\mathbf{#1}}			%

\newcommand{\Q}{\mathbf{Q}}
\renewcommand{\S}{\mathbf{S}}
\renewcommand{\u}{\mathbf{u}}

\newcommand{\q}{\mathbf{q}}
\newcommand{\F}{\mathbf{F}}
\newcommand{\f}{\mathbf{f}}
\newcommand{\g}{\mathbf{g}}

\newcommand{\de}[2]{\frac {\partial #1}{\partial#2}}

\renewcommand{\v}{\mathbf{v}}
\newcommand{\B}{\mathbf{B}}

\newcommand{\RIIcolor}[1]{{\leavevmode\color{black} #1}}

\newcommand{\numberset}{\mathbb}
\newcommand{\N}{\numberset{N}}
\newcommand{\R}{\numberset{R}}

\newcommand{\halb}{\frac{1}{2}}

\newcommand{\bdm}{\begin{displaymath}}
\newcommand{\edm}{\end{displaymath}}

\newcommand{\bea}{\begin{eqnarray} }
\newcommand{\eea}{\end{eqnarray} }

\newfont{\numerikEleven}{ecrm1000}
\newfont{\numerikTen}{cmss10}
\newfont{\numerikNine}{cmss9}
\newfont{\numerikEight}{cmss8}

\journal{Computers \& Fluids}

\begin{document} 
\begin{frontmatter}
\title{High order well-balanced Arbitrary-Lagrangian-Eulerian {ADER} discontinuous Galerkin schemes on {general} polygonal {moving} meshes} 
\author[UniVR]{Elena Gaburro\corref{cor1}}
\ead{elena.gaburro@univr.it}
\cortext[cor1]{Corresponding author}

\address[UniVR]{Department of Computer Science, University of Verona, Strada le Grazie 15, Verona, 37134, Italy}

\begin{abstract} \color[rgb]{0,0,0}
In this work, we present a novel family of high order accurate numerical schemes for the solution of
hyperbolic partial differential equations (PDEs) which combines several geometrical and physical
structure preserving properties. First, we settle our methods in the Lagrangian framework, where each
element of the mesh evolves following as close as possible the local fluid flow, so to
\textit{reduce} the \textit{numerical dissipation} at contact waves and moving interfaces and to
satisfy the \textit{Galilean and rotational invariance} properties of the studied PDEs system. 
In particular, we choose the direct Arbitrary-Lagrangian-Eulerian (ALE) approach which, 
in order to always guarantee the \textit{high quality of the moving mesh}, 
allows to combine the Lagrangian motion with mesh optimization techniques. 
The employed polygonal tessellation is thus regenerated at each time step, the previous one is
connected with the new one by spacetime control volumes, including \textit{hole-like sliver}
elements in correspondence of topology changes, over which we integrate a spacetime divergence form
of the original PDEs through a \textit{high order} accurate ADER discontinuous Galerkin (DG) scheme.
Mass \textit{conservation} and adherence to the \textit{GCL condition} are guaranteed by
construction thanks to the integration over closed control volumes, and \textit{robustness} over
shock discontinuities is ensured by the use of an \textit{a posteriori} subcell finite volume (FV)
limiting technique. 
%
On top of this effective moving mesh framework, we have also modified the full ADER-DG scheme
with \textit{a posteriori} subcell FV limiter 
to be, for the first time in literature, \textit{well-balanced}. 
This is achieved by ensuring that any projection, reconstruction and integration procedures were always performed by
summing up the exact value of a given equilibrium plus the high order accurate evolution of the
fluctuations w.r.t. said equilibrium. %
The paper is closed by a wide set of numerical results, including simulations of Keplerian disks,
which demonstrate all the claimed properties and the increased accuracy of our novel family of
schemes, in particular for the evolution of small perturbations arising over moving equilibrium
profiles. \end{abstract}

\end{frontmatter}


%

\section{Introduction} 
\label{sec.introduction}

The objective of this paper is to combine the accuracy of high order schemes with the structure 
preserving properties made available from the Lagrangian and well-balanced techniques 
in order to develop a robust and effective numerical method 
for the solution of hyperbolic systems of partial differential equations~(PDEs),
in particular one that is well suited for the study of vortical phenomena over long simulation times.

\paragraph{High order ADER discontinuous Galerkin schemes} 
With this aim in mind, we chose to work within the framework of 
high order ADER discontinuous Galerkin (DG) schemes.
DG finite element schemes for first order hyperbolic equations were introduced by Cockburn and Shu in
the 1990s~\cite{cockburn1990runge,cockburn2000development,cockburn2001runge}: 
they base their discretization technique on the use 
of a piecewise high order \textit{spatial} polynomial representation, of the unknown conserved variables, 
in each cell of the computational domain. 
This representation is then inserted in the weak 
formulation of the PDE and most commonly a semi-discrete approach is 
adopted which evolves the data in time following the method of lines,
for example using a single-step multistage scheme like a 
high order Runge-Kutta time integrator~\cite{butcher1996history}. 
In this work, we adopt a different technique for the time evolution, 
the ADER approach (Arbitrary high order DErivative Riemann problem), 
introduced in~\cite{mill,toro1,toro3}, then reworked in its modern formulation in the 
seminal paper~\cite{dumbser2008unified}, and widely used in literature 
(we cite here just a few recent works that span a 
wide range of technical improvements, analytical results, and 
applications~\cite{busto2020high,han2021dec,busto2021high,boscheriAFE2022,ciallella2023shifted,popov2023space,micalizzi2023efficient,dumbser2023WBGR}).
ADER methods make use of a predictor-corrector technique to obtain uniform high order of accuracy 
in space and in time through a \textit{one-step fully discrete} procedure 
which works on data in the form of \textit{spacetime} high order polynomials. 
In addition, this time integration technique proves to be particularly well suited for 
constructing Lagrangian schemes based on a spacetime approach. 

\paragraph{Lagrangian methods} Our aim is the use of a Lagrangian approach, 
because Lagrangian methods~\cite{Neumann1950, Wilkins1964,munz94,LoubereShashkov2004,Maire2010,chengshu2,scovazzi2,
despres2017numerical,morgan2021origins,cremonesi2010lagrangian,cremonesi2017explicit,Dobrev3,plessier2023implicit,del2023triangular,despres2024lagrangian}, 
thanks to the flow driven motion of the computational domain,
allow to significantly reduce the errors due to the convection terms, 
to sharply capture moving interfaces and contact discontinuities, 
and they can be made automatically entropy stable, rotationally invariant, and discretely Galilean invariant. 
However, in their pure form, where the mesh is forced to move exactly with the 
fluid velocity, they are commonly affected by issues related to  
mesh distortion or mesh tangling that may slow down the computation or cause it to halt entirely,
in particular if strong shear flows or 
vortical flows are simulated for long times.

A first possible approach here is provided by the wide class of \textit{meshless methods} and SPH methods~\cite{liu2010smoothed, lind2020review,wang2021improved,kincl2023unified}, 
which are remarkably flexible from a geometrical point of view, 
but generally less accurate than their mesh-based counterparts and therefore not corresponding to our needs.

\RIIcolor{However, already in the '70s, novel attracting \textit{mesh-based approaches}, 
able to relax the constraint of an exact match between fluid flow and mesh motion,  
while at the same time being capable of maintaining as much as possible the benefits of Lagrangian schemes, 
have been developed.
The first papers about this relaxed approach, that named it as Arbitrary-Lagrangian-Eulerian~(ALE)
technique, are~\cite{wilkins1964methods,trulio1966air} and~\cite{hirt1974arbitrary} which all address fluid-dynamics problems studied with a finite difference discretization.

Then, this ALE moving mesh strategy was also adopted in the finite element community for the description of fluid-structure interfaces, for example in~\cite{donea1977lagrangian, hughes1981lagrangian, donea1982arbitrary}, and also in~\cite{belytschko1978computer} where  this approach was addressed as a \textit{quasi-Eulerian} scheme.  
An interesting review paper about the original literature on ALE methods can be found in~\cite{benson1992computational}.
Later, the development of ALE schemes continued both linked to \textit{finite element} simulations~\cite{ghosh1991arbitrary,takashi1992arbitrary,nobile1999stability, sarrate2001arbitrary,saksono2007adaptive,khoei2008extended,guardone2011arbitrary,ortega2011geometrically,zeng2014frame,cremonesi2016lagrangian,di20243d}, for studying interfaces, large media deformations, and used as a powerful tool for adaptive mesh techniques, 
and to \textit{finite volume} methods for hyperbolic equations, on which we will concentrate in the next paragraphs. 

In both the communities ALE techniques can generally be distinguished into two categories: 
the first is the one of \textit{indirect} ALE schemes, characterized by a rezoning procedure (where the mesh quality is optimized) and then a remapping phase, 
where the numerical solution defined on the old Lagrangian mesh is transferred onto the new optimized grid; 
here, a list of references particularly relevant for the simulation of hyperbolic equations is given by~\cite{ReALE2010,ReALE2011,ReALE2015,ShashkovCellCentered,ShashkovRemap1,MaireMM2,indALE-AWE2016,wu2021cell,kenamond2021positivity,lei2023high}.

The second category is that of \textit{direct} ALE schemes, where the discretization method responsible for the PDE evolution directly provides the updated solution on the new mesh without the need of a projection-reconstruction procedure (as the method object of this paper).
Because of the complexity of developing high order accurate projection-reconstruction techniques, 
the {direct} ALE framework represents a convenient option when interested in high order accurate numerical methods.
In particular, the direct ALE approach object of this paper finds its foundation in~\cite{Lagrange2D,Lagrange3D,ALELTS2D,springel2010pur} 
and it is based on the idea of integrating a \textit{spacetime conservation formulation} 
of the governing PDE system over closed, non-overlapping \textit{spacetime} control volumes 
that directly connect the mesh at time $t^n$ with the new one at time $t^{n+1}$, 
where the new computational grid is obtained from the previous one 
by combining the motion along the Lagrangian trajectories with mesh optimization flow-congruent techniques. 
}

Nevertheless, if complex flow characteristics are present and the mesh is allowed to deform only on an 
element-by-element basis, that is, under a constrained fixed connectivity between cells, 
the mesh elements can quickly reach undesirable shapes, even if stretching and skewing are controlled in an 
optimized way,
see for example~\cite{gaburro2020high,gaburro2021bookchapter} and also~\cite{Dobrev1, Dobrev2, Dobrev3, dobrev2020simulation}, 
where such issues are strongly mitigated by the use of extremely curved elements, rather than straight-edge
polygons like in the present work.
Thus, additional freedom should be granted to the mesh optimization procedures, 
in such a way that the elements do not only \textit{move} and deform,
but their vertex count and/or the mesh connectivity may change from time to time. 
This flexibility has been introduced through sliding lines techniques already in~\cite{Caramana2009, DelPino2010, LoubereSL2013},
it is an intrinsic feature of the indirect ALE methods,
and appeared in the context of direct ALE with a \textit{nonconforming} 
approach in~\cite{gaburro2017direct,gaburro2021unified} and in the form of \textit{topology changes} in~\cite{springel2010pur,gaburro2020high}.
The main novelty of these recent works has been the introduction of \textit{hole-like sliver spacetime elements}
to deal with the high order PDE integration around a topology change in meshes made by general polygons or Voronoi elements. 
The characteristics and the necessity of this new type of control volumes have been 
fully described in a previous work of the author~\cite{gaburro2020high} which represent the starting point for the algorithm presented here. 
Indeed, the work presented in~\cite{gaburro2020high} has demonstrated optimal robustness and accuracy in the long time simulations of complex flows,
because there
\textit{i)} we can truly optimize the mesh motion, while following the fluid flow closely,
having the freedom of introducing topology changes;
\textit{ii)} we maintain the high order of accuracy of the underlying ALE ADER-DG method, 
having naturally extended the \textit{direct} PDE integration to the \textit{hole-like} elements
thus avoiding accuracy penalties due to low order projection-reconstruction techniques;
\textit{iii)} we enforce exact conservativity and the satisfaction of the GCL condition everywhere
by redistributing the fluxes around sliver elements among their neighbors
through PDE-based information, being this last flux-rescaling procedure further improved in the present work.

\paragraph{Well-balanced techniques}
Even when combined with Lagrangian techniques, 
high order methods and fine meshes 
are not always enough to provide accurate numerical results. 
A typical situation where such sophisticated machinery tends to fall short
is the long time simulation of equilibrium solutions, both in presence and in absence of small physical perturbations.
In this case we need to endow our numerical scheme with so-called \textit{well-balanced} (WB) techniques, 
i.e. techniques that, exploiting in different ways some available information on the equilibrium profiles, are able to guarantee 
\textit{i)} machine precision accuracy in the simulation of the equilibria themselves and
\textit{ii)} increased resolution of small deviations from equilibrium, not comparable to that of non well-balanced schemes, even when using
fine meshes and with high order of accuracy. 

Well-balanced methods were originally introduced for atmosphere models and shallow water 
equations in~\cite{cargo1994schema,bermudez1994upwind,leveque1998balancing,gosse2001well,bouchut2004nonlinear,
audusse2004fast,Castro2006,Pares2006,castro2007well,VMDansac,Noelle1,Noelle2} 
and then they have been used in many other context: as an illustration, we cite~\cite{Castro2008,arpaia2020well,castro2020well,Pimentel,abgrall2022hyperbolic}
and for the extension to discontinuous Galerkin schemes we mention~\cite{guerrero2021well,mantri2021well,wu2021uniformly,dumbser2023WBGR,tabernero2024high,caballero2024semi}.
In particular, nowadays well-balancing is of interest in the framework of astrophysical applications.
Indeed, as in our manuscript, we find the use of WB techniques for 
the Newtonian Euler equations with gravity and the magnetohydrodynamics equations 
in~\cite{BottaKlein,kappeli2014well,KM15_630,Klingenberg2015,bermudez2016numerical, gaburro2018well,desveaux2016well,
klingenberg2019arbitrary,thomann2019second,Thomann2020,grosheintz2019high,berberich2021high,KlingenbergWBMHD1,BirkeBoscheri,FAMBRI2023112493},  
and they have been applied also to the more complex Einstein field equations in~\cite{gaburro2021well,dumbser2023WBGR}.

\paragraph{Main novelty of the paper}
In our present work, 
we {follow} in particular the WB methodology introduced in~\cite{berberich2021high} for finite volume schemes,
and our previous experience on taking advantage of the combination of WB and ALE~\cite{gaburro2017direct,gaburro2018well}
in order to \textit{significantly extend} the existing framework to 
\textit{i)} arbitrary high order accurate ADER discontinuous Galerkin schemes 
with \textit{ii)} \textit{a posteriori} subcell FV limiter on 
\textit{iii)} arbitrarily moving ALE meshes with topology changes. 

We remark that, this is the first time in literature where well-balanced techniques are applied to high order DG schemes on \textit{moving} meshes, being all the previously cited references applied mainly to low order methods on Eulerian meshes,
and in a few cases, or only to high order methods on Cartesian meshes~\cite{dumbser2023WBGR} or to simpler Lagrangian-type techniques. 
Existing well-balanced works thus take advantage of simplifications 
that, when combining high order DG and moving meshes, as in this work, are not employable. 

In particular, it is worthy to underline the structural difference between this work and a recent, apparently similar, contribution co-authored by the author~\cite{dumbser2023WBGR}.
In~\cite{dumbser2023WBGR}, it is assumed that the equilibrium (supposed to be known \textit{a priori}) is not only constant in time but also \textit{discretely} constant in time. 
This assumption allows to construct the WB scheme simply by 
i) appending, to the standard set of conserved variables, the equilibrium variables also discretized and evolved through (high order) polynomials, and then 
ii) subtracting, from the main scheme for the evolutionary variables, 
each term evaluated on the corresponding \textit{discrete} equilibrium variable. 
However, when the mesh moves, the evolution of the discrete polynomial equilibrium, due to the presence of the advection terms integrated in space and time on a element which is changing, 
does not coincide with the continuum equilibrium, thus the ALE scheme must be really re-written, as explained in details in this paper, in order to correctly introduce the subtraction of the continuous equilibrium terms.  
This difference makes the present approach substantially different 
(and unfortunately significantly more intrusive to be implemented) 
w.r.t. the only other high order WB DG scheme with a posteriori FV limiter 
existing in literature~\cite{dumbser2023WBGR}, which is strictly limited to the Eulerian formalism.

\paragraph{Structure of the paper}
The rest of the paper is organized as follows. 
Section~\ref{sec.numerical} is devoted to the description of our novel numerical method which combines, 
in an effective way, simple well-balanced techniques with our direct ALE approach over moving 
meshes with topology changes. Then, in Section~\ref{sec.tests} we present a large number of 
numerical results that show the increased accuracy provided by our approach that automatically allows 
to simulate very small perturbations of equilibrium solutions, not achievable with standard (non well-balanced) numerical methods.
Finally, in Section~\ref{sec.conclusions} we give some conclusive remarks and an outlook to future developments.

\section{Numerical method}
\label{sec.numerical}

The numerical method developed in this work is applicable to any first order partial differential
equation of hyperbolic type that can be cast into the following form \begin{equation}
\label{eq.generalform}
\partial_t{\Q} + \nabla \cdot \F(\Q) = \mathbf{S}(\Q), \qquad \x \in \Omega(t) \subset \mathbb{R}^d, 
\qquad \Q \in \Omega_{\Q} \subset \mathbb{R}^{\nu},
\end{equation}
where we work in dimension $d=2$, $\x = (x,y)$ is the spatial position vector and $t$ represents the
time. Moreover, $\Q = (q_1,q_2, \dots, q_{\nu})$ is the vector of the $\nu$ conserved variables
defined in the space of the admissible states $\Omega_{\Q} \subset \mathbb{R}^{\nu}$, $ \F(\Q) =
(\,\f(\Q), \g(\Q)\,) $ is the non linear flux tensor and $\mathbf{S}(\Q)$ represents a non linear
algebraic source term. In particular, the PDEs considered in our benchmarks are the 
Euler equations of gasdynamics with and without gravity source term and the magnetohydrodynamics
equations; they are precisely described in Section~\ref{ssec.test_equations}. 

To discretize our moving two-dimensional domain $\Omega(t)$ we employ a \textit{moving tessellation}
made of $N_P$ non overlapping \textit{general polygons} $P_i, i=1, \dots N_P$, first built at time
$t=0$, then regenerated at each time step (according to the motion of a set of 
so-called generator points) and connected in spacetime between each time level via
spacetime control volumes, as specified in Section~\ref{sec.nm-domain}. 

Next, the data, used to discretize the numerical solution and the physical fluctuations, are
represented via polynomials \RIIcolor{of degree up to $N$} in space in each polygon ($\u_h$ and $\u_f$), as detailed in
Section~\ref{sec.nm-quantities}, and are then evolved via our novel well-balanced
Arbitrary-Lagrangian-Eulerian~(ALE) discontinuous Galerkin~(DG) scheme with \textit{a posteriori}
subcell finite volume (FV) limiter. 

In Section~\ref{sec.nm-ader} and~\ref{sec.nm-limiter} we provide a description of 
our high order algorithm on moving regenerated
meshes~\cite{gaburro2020high,gaburro2021unified,gaburro2021bookchapter} and we explain 
in detail how
to make sure that the well-balanced property is preserved throughout all of the stages of 
the algorithm. 

We also remark that the presented method evolves the solution in an \textit{explicit} way, 
timestep by timestep i.e from $t^n$ to $t^{n+1} = t^n + \Delta t$, 
so it is stable under the following CFL stability condition on the timestep size 
\begin{equation}
	\Delta t < \textnormal{CFL} \left( 
	\frac{ |P_i^n| }{ (2N+1) \, |\lambda_{\max,i}| \, \sum_{\partial P_{i_j}^n} |\ell_{i_j}| } 
	\right), \qquad \forall P_i^n \in \Omega^n. 
	\label{eq:timestep}
\end{equation}
In the above formula, $\ell_{i_j}$ is the length of the edge $j$ of $P_i^n$, 
$|\lambda_{\max,i}|$ is the spectral radius of the Jacobian of the flux $\mathbf{F}$ 
\RIIcolor{and $N$ is the degree of the polynomials used for representing the data in space, which leads to a method of formal order of accuracy $N+1$.} 
On unstructured meshes the CFL stability condition requires the 
inequality $\textnormal{CFL} < \frac{1}{d}$ to be satisfied, see~\cite{dumbser2008unified}.

\subsection{Discretization of the moving domain}
\label{sec.nm-domain}

At time $t^n=0$ we cover the computational domain $\Omega$ and its boundary with $N_P$ points
\begin{equation}
	\xcin, \ i=1,\dots, N_P,
\end{equation}
that we call \textit{generator} points (the orange points in Figure~\ref{fig.mesh}). 
The position of these points evolves at each timestep according to 
\begin{equation}
\xcinp = \xcin + \Delta t \, \mbf{v}(\xcin),
\label{eqn.xcnew} 
\end{equation}
where the velocity $\mbf{v}(\xcin)$ is chosen in such a way as to balance between 
two often contrasting requirements: first, closely follow the fluid flow, 
and second, retain advantageous element shapes from the point of view of 
numerical discretization. The motion is obtained via a high order 
integration of the trajectories of the generator points, in order to 
exploit the advantages of Lagrangian schemes, then followed by an
optimization of the quality of the mesh to be built at the next timestep.
The latter step aims at reducing the numerical errors due to excessive element distortion, 
see~\cite{gaburro2020high} (Sections 2.4 and 2.5) for more details. 

Once the position of the generators at $t^n$ and $t^{n+1}$ has been fixed, 
we need to construct both the \textit{spatial} mesh at the new timestep $t^{n+1}$ 
(and, only at the first timestep, also the initial mesh at time $t^n=0$) and also 
the \textit{spacetime} mesh that connects $t^n$ with $t^{n+1}$ by completely 
filling the spacetime between the two levels while respecting the time-slicing. 

\subsubsection{Delaunay triangulation and polygonal tessellation} 

\begin{figure}[!bp] \centering 
	{\includegraphics[width=0.24\linewidth]{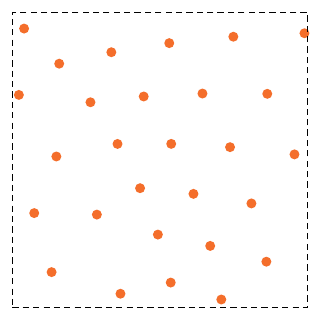}}\
	{\includegraphics[width=0.24\linewidth]{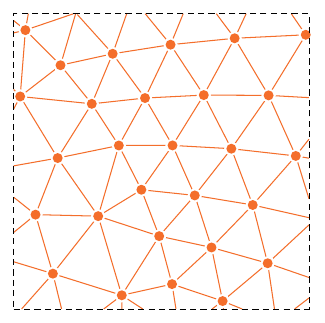}}\
	{\includegraphics[width=0.24\linewidth]{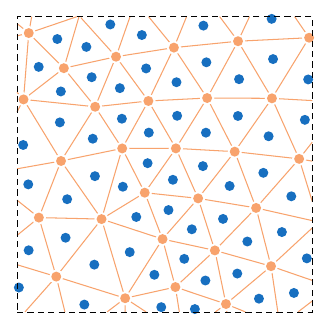}}\		
	{\includegraphics[width=0.24\linewidth]{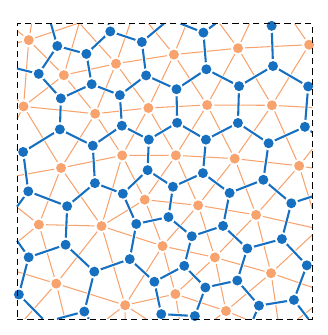}}
	\caption{The domain $\Omega$ is covered with a set of generator points (orange points) 
	that are connected via a Delaunay triangulation. Then, around each generator we construct a 
	polygon by connecting the centroids of the attached triangles (blue points). }
	\label{fig.mesh}
\end{figure}

To build the \textit{spatial} mesh at a generic time $t^n$, 
we first connect the generators via a Delaunay triangulation in such a way that the 
generators $\xcn$ are the vertexes of the Dealunay triangles, 
following standard algorithms as~\cite{bowyer1981computing,watson1981computing,mucke1999fast} (see the second panel of Figure~\ref{fig.mesh}). 
Then, around each generator $\xcin$ we construct a polygon $P_i^n$ by connecting in 
counterclockwise order the \textit{barycenters} of all the Delaunay triangles sharing 
this $\xcin$ as a vertex (see the third and fourth panel of Figure~\ref{fig.mesh}). 
Note that the use of the barycenters (instead of circumcenters) to obtain the polygonal 
elements (instead of Voronoi elements) avoids short edges and in particular zero-lengths ones.

Given a polygon $P_i^n$ we introduce the following notation: 
we call $\mathcal{V}(P_i^n) = \{v_{i_1}^n, \dots, v_{i_j}^n, \dots, v_{i_{N_{V_i}^n}}^n \}$ the set of its $N_{V_i}^n$ neighbors, 
$\mathcal{E}(P_i^n) = \{e_{i_1}^n, \dots, e_{i_j}^n, \dots, e_{i_{N_{V_i}^n}}^n \}$ the set of its $N_{V_i}^n$ edges, 
and 
$\mathcal{D}(P_i^n) = \{d_{i_1}^n, \dots, d_{i_j}^n, \dots, d_{i_{N_{V_i}^n}}^n \}$ the set of its $N_{V_i}^n$ 
vertexes, consistently \textit{ordered counterclockwise} \RIIcolor{(the subscript $i$ may be omitted when the notation is non-ambiguous)}.
Next, we denote the barycenter of $P_i^n$, that in general does not coincide with the generator point, 
as $\xbin = (\xbixn, \xbiyn)$ and, for computational reasons, for example for fixing adequate quadrature points, 
we subdivide $P_i^n$ in $N_{V_i}^n$ subtriangles denoted as $\mathcal{T}(P_i^n) = \{T_{i_1}^n, 
\dots, T_{i_j}^n, \dots, T_{i_{N_{V_i}^n}}^n \}$ by connecting $\xbin$ with each vertex of $\mathcal{D}(P_i)$.

When, at a new timestep, we generate a new mesh, the \textit{only guaranteed invariant} between the 
tessellations at time $t^n$ and $t^{n+1}$ is the number $N_P$ of generator points (i.e. of total 
polygons) and their \textit{global numbering}. Instead, the shape of each polygon is allowed to 
change, i.e. $N_{V_i}^n \ne N_{V_i}^{n+1}$, and consequently also the connectivities, i.e. for 
example $\mathcal{V}(P_i^n) \ne \mathcal{V}(P_i^{n+1})$.
This possibility to change the topology of the grid is actually the \textit{strength} of the present ALE algorithm, 
since it allows us to highly optimize the mesh construction at each timestep: for example, 
we can adapt the triangulation connecting the generator points any time they move by
enforcing the Delaunay property of the underlying triangulation and thus improve the 
quality of the corresponding polygonal tessellation, avoiding distorted or tangled elements.

\subsubsection{Spacetime connectivity: partitioning the spacetime volumetric slices with discrete control volumes}

\begin{figure}[!bp] \centering
	{\includegraphics[width=0.28\linewidth,height=0.19\textheight]{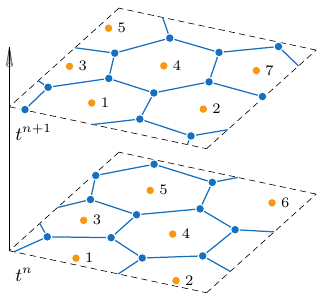}}\qquad 
	{\includegraphics[width=0.28\linewidth,height=0.19\textheight]{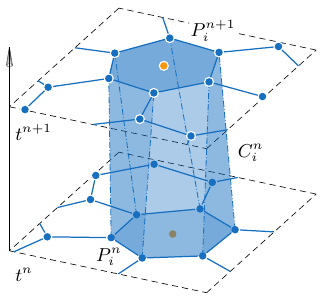}}\qquad 
	{\includegraphics[width=0.28\linewidth,height=0.19\textheight]{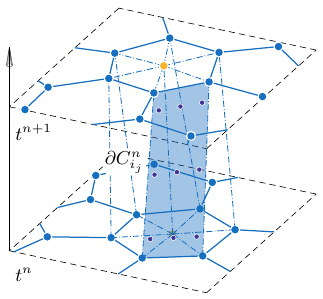}}\\
	{\includegraphics[width=0.28\linewidth,height=0.19\textheight]{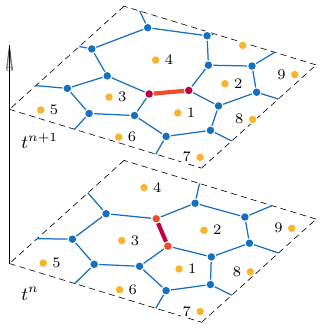}}\qquad	
	{\includegraphics[width=0.28\linewidth,height=0.19\textheight]{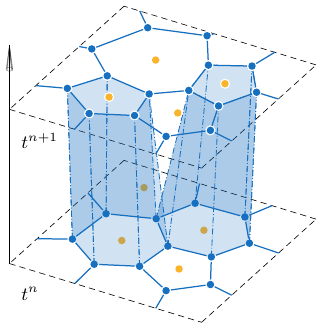}}\qquad
	{\includegraphics[width=0.28\linewidth,height=0.19\textheight]{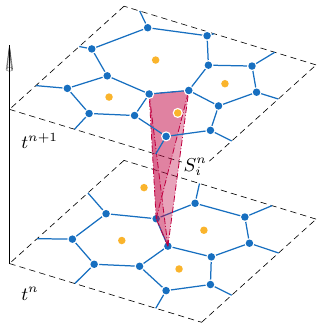}}%
	\caption{Spacetime connectivity \textit{without} (top) and \textit{with} (bottom) topology changes, 
	with corresponding standard spacetime control volumes (blue, top), degenerate spacetime control volumes 
	(blue, bottom) and a hole-like sliver control volume (pink, bottom).}
	\label{fig.controlvolumes} \vspace{-10pt}
\end{figure}

The topology changes really help Lagrangian simulations in case of of complex {flow fields}, 
like strong shear flows or vortical flows, 
however they make it more complex to evolve the solutions in their vicinity while retaining 
formal high order of accuracy.
In particular, to obtain a \textit{direct} ALE scheme {of order greater than two} we need a
\textit{complete knowledge} of the \textit{spacetime} structure between the two time levels, so we
are required to connect them via spacetime control volumes and to keep track of 
a complete description of their geometry and their spacetime neighbors. 

Wherever $N_{V_i}^n = N_{V_i}^{n+1}$ and $\mathcal{V}(P_i^n) =
\mathcal{V}(P_i^{n+1})$, i.e. the polygon $P_i^n$ is not affected by topology changes (see
the first line of Figure~\ref{fig.controlvolumes}), in order to construct the spacetime control
volume $C_i^n$ it is enough to connect, via straight line segments, each node of $P_i^n$ 
with the corresponding one in $P_i^{n+1}$. Note that the correspondence between nodes can be
established by inspecting the neighbors numbering, since the element number is a fixed label between
the timesteps. In this way we obtain a closed oblique frustum (in the case of uniform 
translational motion, a prism) with polygonal bottom and top faces, $P_i^n$ and $P_i^{n+1}$, 
and $N_{V_i}^{n,st}=N_{V_i}^n$ quadrilateral faces $ \partial C_{i_j}^n$ composing its total spacetime lateral surface 
\begin{equation}
	\label{eqn.dC}
	\partial C_i^n = P_i^n \cup P_i^{n+1} \bigcup_{j=1}^{N_{V_i}^{n,st}} \partial C_{i_j}^n .
\end{equation} 

Wherever instead we find $N_{V_i}^n \ne N_{V_i}^{n+1}$, we first have to fix the set of the
spacetime neighbors of $C_i^n$, i.e. $\mathcal{V}(C_i^n)$ which does not coincide with either 
$\mathcal{V}(P_i^n)$ or $\mathcal{V}(P_i^{n+1})$ and counts $N_{V_i}^{n,st}$ elements. 
$\mathcal{V}(C_i^n)$ contains all the polygons of $\mathcal{V}(P_i^n)$ and $\mathcal{V}(P_i^{n+1})$
counted once (i.e. without multiple entries) and  ordered counterclockwise following the order of 
both $\mathcal{V}(P_i^n)$ and $\mathcal{V}(P_i^{n+1})$. This spacetime ordering allows also to
establish a correspondence between the nodes of $P_i^n$ and $P_i^{n+1}$,
which should then be connected again via straight line segments; this process induces the 
appearance of degenerate elements of two types.

First, \textit{degenerate} spacetime control volumes (see the fifth panel of
Figure~\ref{fig.controlvolumes}) for which the top and the bottom faces are polygons with a
different number of nodes. We call them \textit{degenerate} because they 
feature a \textit{triangular} spacetime lateral face
(instead of a \textit{quadrilateral} one) due to the fact that two distinct   
nodes at one level correspond to
the same node at the other level. However, their numerical treatment is identical to {that of} 
standard spacetime control volumes, because this degeneracy simply changes the spacetime 
distributions of the quadrature points.

Second, the \textit{hole-like sliver elements} $S_i^n$ (see the sixth panel of
Figure~\ref{fig.controlvolumes}), that are stretched tetrahedra allowing to fill the
empty spacetime volume left by the connection of existing spatial elements around a topology
change. These additional control volumes can be obtained by connecting the endpoints of the edges that at
a certain time level are shared between two neighboring elements, and at the other level run between
the same two elements which however, because of a topology change, are no longer neighbors. 
We remark that they are additional control volumes which does not exist at time 
$t^n$ or at time $t^{n+1}$, since they coincide with edges of the tessellation, and, as such, they
have zero area in space at $t^n$ and $t^{n+1}$. These two characteristics are responsible of the two
main numerical difficulties connected with their treatment: 
\textit{i)} the numerical solution is not clearly defined for them at time $t^n$, because they lie at the 
boundary between two elements with in principle discontinuous information at $t^n$, and \textit{ii)} 
the contributions across them, computed with an explicit numerical scheme, since the  
area of $S_i^n$ at $t^{n+1}$ is zero, should be redistributed among the neighbors, in order to ensure conservation. 

For additional details on the construction of standard, degenerate and sliver spacetime control volumes, 
and a complete justification of the robustness of their construction procedure and their necessity 
for long time simulations, we refer to~\cite{gaburro2020high,gaburro2021bookchapter}. Instead, 
their numerical treatment and the strategies to overcome the above described difficulties 
are clearly explained in Section~\ref{ssec.nm-sliver}.


\subsection{Outline of the employed discretization quantities and basis functions}
\label{sec.nm-quantities}

In this work, we are proposing a well-balanced method, i.e. a numerical strategy capable of simulating 
equilibrium solutions of the studied PDE 
system (in our case known \textit{a priori}) with errors of the order of machine precision, thus
it is particularly appropriate to study small physical \textit{fluctuations} arising around the equilibria.
Note that this does not mean that the scheme is a perturbation method 
based on linearized solutions around a certain state:
instead it exhibits the standard nonlinear behaviour far from the equilibrium, 
capturing the solution as conventional DG techniques, but it also implements some improvements that allow
much better behaviour wherever indeed the solution is close to an equilibrium state.
In addition, the method is high order accurate in space and time thanks to the ADER 
predictor-corrector approach~\cite{dumbser2008unified} and it belongs to the family of the \textit{direct} ALE methods.

\subsubsection{Discretization quantities: equilibrium, complete solution, fluctuations}

Directly from the above short summary, 
we can justify the use, along all our algorithm, of \textit{three main groups of quantities}.

First, we consider the \textit{equilibrium solution}, that we assume  
to be the same throughout the simulation and it is supposed to be
analytically known anywhere or exactly available at any quadrature point whenever required. 
We denote by \RIIcolor{$\Q_E: \R^d \times \R \rightarrow \R^\nu$}
the equilibrium profile \RIIcolor{which} exactly satisfies the PDE, i.e.
\begin{equation}
	\label{eq.equilibrium}
	\nabla \cdot \F(\Q_E) - \mathbf{S}(\Q_E) = \mathbf{0}.
\end{equation} 
We highlight that $\Q_E$ is \textit{not} the polynomial \RIIcolor{representation} of the equilibrium 
and it is \textit{not} the result of the numerical evolution of the equilibrium (as it is in other works as~\cite{dumbser2023WBGR}), it is the \textit{exact} equilibrium value.
We also recall that a numerical scheme, similar to the one presented here, 
could be developed also for equilibrium profiles not completely known \textit{a priori}, 
but recoverable at the beginning of each timestep from the numerical data plus additional information
(such as for example ordinary differential equations describing the equilibrium), 
see~\cite{gomez2021high,guerrero2020second}.
In this case as well, however, the equilibrium value $\Q_E$
would not come from the numerical evolution of the equilibrium, 
but from additional/external procedures. So, our scheme could in principle be 
extended to more general situations.

Second, we have the \textit{complete numerical solution} that describes entirely the physical 
situation accounting both for the equilibrium and its oscillations (or also a non-equilibrium situation); 
classically, this would be the unique result of a standard non well-balanced numerical scheme.
Here, we denote by $\RIIcolor{\Q: \R^d \times \R \rightarrow \R^\nu} $ the conserved variables, $ \RIIcolor{\Q_{\text{IC}}: \R^d \rightarrow \R^\nu} $ the initial condition, 
\RIIcolor{$\u_h^n: \R^d \rightarrow \R^\nu$}
the numerical \textit{spatial} solution at time $t^n$ discretized in the 
space of polynomials up to degree $N$ spanned by the basis functions given in~\eqref{eq.Dubiner_phi_spatial},
and by \RIIcolor{$\q_h^n: \R^d \times \R \rightarrow \R^\nu$}
the spacetime predictor (see Section~\ref{ssec.nm-predictor}) 
discretized in the space of polynomials depending both on \RIIcolor{\textit{space} and \textit{time}} and spanned by~\eqref{eq.Dubiner_phi}.

Third, we consider the so-called \textit{fluctuations}.
At the beginning of the simulation we initialize \RIIcolor{$\u_f^0$} as the spatial 
polynomial expressed via~\eqref{eq.Dubiner_phi_spatial} representing $\Q_\text{IC} - \Q_E$. 
Then, $\u_f^n$ is the result of our numerical well-balanced simulation at each 
timestep $t^n$, see~\eqref{eqn.corrector-f},
and $\q_f^n$ the corresponding predictor~\eqref{eqn.predictor-f}.
We highlight that not all the operations of a standard high order DG scheme can be 
performed directly on the fluctuations (deviations from equilibrium), 
because they, considered without the corresponding equilibrium, do not give raise to physically valid state vectors. 
So, during the algorithm, we will need several times 
to recover the values of $\u_h$ and $\q_h$ at certain \RIIcolor{spacetime points}
(e.g. for evaluating fluxes/sources, the transformation from conservative to primitive 
variables and viceversa) and then go immediately back to $\u_f$ by exploiting the 
knowledge and the properties of the equilibrium.

Note also that in the non well-balanced version of our numerical scheme, the polynomial $\u_h$ 
represents the direct result of the numerical simulation (see~\eqref{eqn.corrector} and 
refer to~\cite{gaburro2020high} for the description of the pure non well-balanced scheme).
Instead, in the \textit{well-balanced} version $\u_h$ is always recovered 
in \RIIcolor{each spacetime point $(\x,t)$ where it is needed (e.g. on quadrature points)}
as the sum between the exact equilibrium $\Q_E$ and the 
polynomial representing the fluctuations $\u_f$ \RIIcolor{evaluated at that point}, i.e.
\begin{equation}
	\label{eqn.uhfromuf}
	\u_h(\x,t) = \Q_E(\x, t) + \u_f(\x,t),
\end{equation} 
and, in the same way we obtain $\q_h$ in each needed \RIIcolor{point} as 
\begin{equation}
	\label{eqn.qhfromqf}
	\q_h(\x,t) = \Q_E(\x, t) + \q_f(\x,t);
\end{equation} 
here, we never reconstruct the polynomial interpolants of these values because we do not need them.\\ 

\subsubsection{Modal basis functions in space, in spacetime and moving basis}

In order to define a high order direct ALE discontinuous Galerkin scheme on moving meshes which 
\textit{i)}~updates the solution in time via a \textit{one-step} procedure and \textit{ii)} \textit{directly} evolves 
the solution from one mesh configuration to the new one, we need to introduce three sets of basis functions. 
\RIIcolor{Note that, throughout this paper, to lighten the notation regarding the basis functions and the scheme formulation, we employ the Einstein summation convention, which implies summation over two repeated indices.}

First, our discrete \textit{spatial} data, namely $\u_h$ and $\u_f$, are represented inside each 
polygon $P_i^n$ ($\x \in P_i^n$) via a cell-centered approach through {piecewise polynomials} of degree up to $N \geq 0$ as 
\begin{equation}
	\mathbf{u}_h^n(\x,t^n) = \phi_\ell(\x,t^n) \, \hat{\mathbf{u}}^{n}_{\ell}, \qquad
	\mathbf{u}_f^n(\x,t^n) = \phi_\ell(\x,t^n) \, \hat{\mathbf{u}}^{n}_{f,\ell}, \qquad  \RIIcolor{\forall} \ell \in \{0, \dots, {\mathcal{L}(N,d)}-1\},
	\label{eqn.uh}
\end{equation}
\RIIcolor{where $\phi_\ell$ are the basis functions described in the next paragraph and $\hat{\mathbf{u}}^{n}_{\ell}$, $\hat{\mathbf{u}}^{n}_{f,\ell} \in \R$ are the so-called degrees of freedom;
the number of degrees of freedom and basis functions used to span the space of spatial 
polynomials up to degree $N$ is $\mathcal{L}(N,d)$ 
\begin{equation}
	\mathcal{L}(N,d) = \frac{1}{d!} \prod \limits_{m=1}^{d} (N+m),
	\label{eqn.nDOF}
\end{equation}
where $d=2$ for two-dimensional domains.}

\RIIcolor{
The chosen functions \RIIcolor{$ \x \mapsto \phi_\ell (\x, t^n)$} are \textit{modal} spatial basis functions obtained as 
rescaled monomials of degree up to $N$ in the variables $\mathbf{x}=(x,y)$ directly 
defined on the \textit{physical element} $P_i^n$, 
expanded about its current barycenter $\xbin$ at time $t^n$ and normalized by its current characteristic length $h_i$.
To express them we employ a multi-index notation $\ell \rightarrow (\ell_1,\ell_2)$ with $\ell_1, \ell_2 \in \N$
and 
we have
\begin{equation} 
	\begin{aligned} 
	\label{eq.Dubiner_phi_spatial}
	\phi_\ell(\x,t^n) |_{P_i^n} = \frac{(x - \xbixn)^{\ell_1}}{h_i^{\ell_1}} \, \frac{(y - \xbiyn)^{\ell_2}}{h_i^{\ell_2}}
	\qquad  \forall \ell \in \{0, \dots, {\mathcal{L}(N,d)}-1\} \ \text{ and } \ \ 0 \leq \ell_1 + \ell_2 \leq N. 
	\end{aligned} 
\end{equation} 
Here $h_i$ is a length scale estimated by setting 
\begin{equation}
h_i = \sqrt{\frac{J^{xx}_i + J^{yy}_i}{|P_i^n|}}, 
\end{equation}
where $J^{xx}_i$ and $J^{yy}_i$ are
\begin{equation} 
	\begin{aligned} 
		\label{eq.Jxx_Jyy}	
		J^{xx}_i = \frac{1}{12} \, \sum_{j=1}^{N_{V_i}^{n}} ( y_{d_{j}}^2 + y_{d_{j}} y_{d_{j+1}} + x_{d_{j+1}}^2     ) 
		(\,  x_{d_{j}} y_{d_{j+1}} - x_{d_{j+1}} y_{d_{j}}  \, ),\\
		J^{yy}_i = \frac{1}{12} \, \sum_{j=1}^{N_{V_i}^{n}} ( x_{d_{j}}^2 + x_{d_{j}} x_{d_{j+1}} + y_{d_{j+1}}^2     ) 
		(\,  x_{d_{j}} y_{d_{j+1}} -  x_{d_{j+1}} y_{d_{j}}  \, ),	
	\end{aligned} 
\end{equation} 
i.e. the two second moments of area of $P_i^n$, and their sum thus represents
the polar moment of area of $P_i^n$ (the integral of the squared distance from the centroid of $P_i^n$).}
\RIIcolor{We want to clarify that this choice of $h$ has been found to optimize the condition number associated to the basis functions in case of deformed elements in a few numerical experiments we did; however, any other classical choice for $h$, as the incircle diameter or the square root of the element volume, is also possible.}
%

We also introduce our \textit{spacetime} basis functions $(\x,t) \mapsto \theta(\x,t)$ that will be used to 
represent the \textit{predictor} $\q_h$ and $\q_f$ (see Section~\ref{ssec.nm-predictor}) inside each control volume $C_i^n$.
The  $\theta_\ell$ are simple monomial-type \textit{modal spacetime} basis of the 
polynomials of degree up to $N$ in $d+1$ dimensions ($d$ space dimensions plus time), 
expanded around $(\xbin,t^n)$, which, \RIIcolor{using again a multi-index notations $\ell \rightarrow (\ell_1,\ell_2,\ell_3)$ with $\ell_1, \ell_2, \ell_3 \in \N$,  read
\begin{equation}
	\label{eq.Dubiner_phi}
	\theta_\ell(\x,t)|_{C_i^n} = \frac{(x - \xbixn)^{\ell_1}}{h_i^{\ell_1}} \, \frac{(y - \xbiyn)^{\ell_2}}{h_i^{\ell_2}}
	\, \frac{(t - t^n)^{\ell_3}}{h_i^{\ell_3}}, 
	\ \RIIcolor{\forall} \ell \in \{0, \dots, \mathcal{L}(N,d+1)-1\},
	\, \text{ and } \ 0 \leq \ell_1 + \ell_2 + \ell_3 \leq N. \\
\end{equation}  
}
Using them, we can write the predictors as 
\begin{equation}
\q_h^n(\x, t) = \theta_\ell (\x, t) \hat{\q}_{\ell}^n, \qquad
\q_f^n(\x, t) = \theta_\ell (\x, t) \hat{\q}_{f,\ell}^n, \qquad 
\quad \RIIcolor{\forall}\ell \in \{ 0, \dots, \mathcal{L}(N,d+1) \},
\label{eqn.qh}
\end{equation}
\RIIcolor{where $\hat{\q}_{\ell}^n$ and $\hat{\q}_{f,\ell}^n$ $\in \R$ are the degrees of freedom characterizing the polynomials}.
We remark that, these basis functions $\theta$ are redefined at the beginning of each time step 
in function of the current position $\xbin$, thus they are directly linked to the current mesh 
configuration and used to represent information only valid \textit{locally} inside each $C_i^n$.

Next, we introduce an \textit{essential} basis for our \textit{direct} ALE algorithm:  
the so-called \textit{moving} spatial modal test functions \RIIcolor{$(\x,t) \mapsto \phitilde(\x,t)$}, 
which will be used as test functions in the corrector step, see~\eqref{eqn.corrector} and~\eqref{eqn.corrector-f}. 
They coincide with~\eqref{eq.Dubiner_phi_spatial} at $t=t^n$ and at $t=t^{n+1}$, 
i.e. $\phitilde_k(\x,t^n) = \phi_k(\x,t^n)$ and 
$\phitilde_k(\x,t^{n+1}) = \phi_k(\x,t^{n+1})$, 
but unlike $\phi$, the $\phitilde$ can be used at both the time levels 
because they automatically adapt to the mesh evolution and thus they act as test 
functions in a moving spacetime mesh method.
Indeed, they are tied to the motion of the barycenter $\x_{\mathbf{b}_i}(t)$ and move together 
with $P_i(t)$ in such a way that at time $t = t^{n+1}$ they automatically refer to the new barycenter $\xbinp$. They read 
\begin{eqnarray} 
	\label{eq.Dubiner_phi_movingspatial}
	&& \phitilde_\ell(\x,t)|_{C_i^{n}} = 
	\frac{(x - x_{b_i}(t))^{\ell_1}}{h_i^{\ell_1}} \, \frac{(y - y_{b_i}(t))^{\ell_2}}{h_i^{\ell_2}}, \quad 
	\text{ with } \ \x_{\mbf{b}_i} (t) = \frac{t-t^n}{\Delta t} \xbin + \left(1-\frac{t-t^n}{\Delta t}\right) \xbinp, \\[2pt] 
	&&\RIIcolor{ \forall \ell \in \{0, \dots, {\mathcal{L}(N,d)}-1\}, \ \,  \ell \rightarrow (\ell_1,\ell_2), \ \text{ and } \ \ 0 \leq \ell_1 + \ell_2 \leq N}.  \nonumber 
\end{eqnarray}

\subsection{Well-balanced direct Arbitrary-Lagrangian-Eulerian ADER discontinuous Galerkin method}
\label{sec.nm-ader}

The governing equations~\eqref{eq.generalform} are now evolved, in a succession of timesteps, 
with our high order fully-discrete \textit{one-step predictor-corrector} ADER-DG method applied 
between the generic timesteps $t^n$ and $t^{n+1}$,
by \textit{directly} updating the solution from one mesh configuration to the following one 
simply through spacetime integration, without needing any projection-remapping steps.

To achieve this result, following~\cite{Lagrange2D, Lagrange3D, LagrangeISO, gaburro2020high}, we 
first rewrite the governing PDE system~\eqref{eq.generalform} in a spacetime divergence form as
\begin{equation}
	\tilde \nabla \cdot \tilde{\F} (\Q) = \mathbf{S} (\Q), 
	\label{eqn.st.pde}
\end{equation}
with $\tilde \nabla = \left( \partial_x, \, \partial_y, \, \partial_t \right)$ denoting the 
spacetime divergence operator and \RIIcolor{$\tilde{\F}(\Q) = \left( \mathbf{f}(\Q), \, \mathbf{g}(\Q), \, \Q \right)$} 
being the corresponding spacetime flux.
Then, we multiply~\eqref{eqn.st.pde} by the set of \textit{moving} spatial modal test 
functions \RIIcolor{$\phitilde_k$} introduced in~\eqref{eq.Dubiner_phi_movingspatial} and we integrate 
over the closed spacetime control volume $C_i^n$ obtaining 
\begin{equation}
	\int_{C_i^n} \phitilde_k \tilde \nabla \cdot \tilde{\F}(\Q) \, d\mathbf{x} dt = \int_{C_i^n} \phitilde_k \mathbf{S} (\Q)\, d\x dt, \qquad \RIIcolor{ \forall k \in \{0, \dots, {\mathcal{L}(N,d)}-1\} }.
	\label{eqn.intPDE}
\end{equation}
Next, by applying the Gauss theorem, we get 
\begin{equation}
	\int_{\partial C_{i}^n} 
	\phitilde_k  \tilde{\F}(\Q) \cdot \ \mathbf{\tilde n} \, dS   
	- \int_{C_i^n} 
	\tilde \nabla \phitilde_k \cdot \tilde{\F}(\Q) \, d\mathbf{x} dt = \int_{C_i^n} \phitilde_k\mathbf{S}(\Q) \, d\x dt, \qquad \RIIcolor{ \forall k \in \{0, \dots, {\mathcal{L}(N,d)}-1\} },
	\label{eqn.GaussPDE}
\end{equation}
where $\mathbf{\tilde n} = (\tilde n_x,\tilde n_y,\tilde n_t)$ denotes the outward pointing 
spacetime unit normal vector on the spacetime faces composing the boundary $\partial C_{i}^n$. 

Now, we decompose the surface integral over the faces of $\partial C_{i}^n$ given in~\eqref{eqn.dC},
we introduce the value $\u_h^n$, discretized by~\eqref{eqn.uh}, representing the solution at $t^n$, 
$\u_h^{n+1}$ representing the solution at $t^{n+1}$ and the predictor $\q_h^{n}$ that will be a 
high order approximation of the solution valid \textit{inside} $C_i^n$ (see Section~\ref{ssec.nm-predictor}). 
Thus, we obtain our explicit direct ALE ADER-DG scheme
\begin{equation}
	\begin{aligned} 
		\left( \,  \int_{P_i^{n+1}} \phitilde_k \phi_\ell \, d\x \,  \right )
		 \hat{\u}_{\ell}^{n+1}   = & \phantom{+}
		\left( \,  \int_{P_{i}^n} 
		\phitilde_k  \phi_\ell \, d\x  \, \right) \hat{\u}_{\ell}^n  \\
		& - \sum_{j=1}^{N_{V_i}^{n,st}} \int_{\partial C_{ij}^n} 
		\phitilde_k  \mathcal{F}(\q_h^{n,-},\q_h^{n,+}) \cdot \mathbf{\tilde n} \, dS  
		+ \int_{C_i^n} 
		\tilde \nabla \phitilde_k \cdot \tilde{\F} (\q_h)\, d\mathbf{x} dt + \int_{C_i^n} \phitilde_k \mathbf{S}(\q_h) \, d\x dt,\\
		& \quad \RIIcolor{ \forall k \in \{0, \dots, {\mathcal{L}(N,d)}-1\} }
	\end{aligned} 
	\label{eqn.corrector}
\end{equation}
where the DOFs of the unknown solution at the new time step $\hat{\u}_\ell^{n+1}$ can be computed
\textit{directly} from the previous ones $\hat{\u}_\ell^n$ through the integration of the fluxes and
the source terms over $C_i^n$ thanks to the use of the predictor $\q_h^{n}$ and the numerical flux
function $ \mathcal{F}(\q_h^{n,-},\q_h^{n,+}) \cdot \mathbf{\tilde n} $ which couples the interactions between neighbors.
Hence, the above formula is termed the \textit{corrector} step of our ADER-DG scheme.

The two-point numerical flux function $ \mathcal{F}(\q_h^{n,-},\q_h^{n,+}) \cdot \mathbf{\tilde n} $ 
is computed via an ALE Riemann solver applied to   
$\q_h^{n,-}$ and $\q_h^{n,+}$ which are respectively the inner and outer boundary-extrapolated data, 
i.e. the values assumed by the predictors of the two neighbor elements ($C_i^n$ and $C_j^n$)
 at a given quadrature point on the shared spacetime lateral surface. 
Here, the simplest choice consists in adopting a Rusanov-type~\cite{Rusanov:1961a} ALE flux, 
\begin{equation}
	\mathcal{F}(\q_h^{n,-},\q_h^{n,+}) \cdot \mathbf{\tilde n} =  
	\frac{1}{2} \left( \tilde{\F}(\q_h^{n,+}) + \tilde{\F}(\q_h^{n,-})  \right) \cdot \mathbf{\tilde n}_{ij}  - 
	\frac{1}{2} s_{\max} \left( \q_h^{n,+} - \q_h^{n,-} \right),  
	\label{eq.rusanov} 
\end{equation} 
where $s_{\max}$ is the maximum \RIIcolor{of the spectral radii} of 
$\mathbf{A}^{\hspace{-1pt}\mathbf{V}}_{\mathbf{n}}(\q_h^{n,+})$ 
and $\mathbf{A}^{\hspace{-1pt} \mathbf{V}}_{\mathbf{n}}(\q_h^{n,-})$, 
and $\mathbf{A}^{\hspace{-1pt}\mathbf{V}}$ is the ALE Jacobian matrix w.r.t. the normal direction in space 
\begin{equation} 
	\label{eq.ALEjacobianMatrix}
	\mathbf{A}^{\hspace{-1pt} \mathbf{V}}_{\mathbf{n}}(\Q)=\left(\sqrt{\tilde n_x^2 + 
	\tilde n_y^2}\right)\left[\frac{\partial \mathbf{F}}{\partial \Q} \cdot \mathbf{n}  - 
	(\mathbf{V} \cdot \mathbf{n}) \,  \mathbf{I}\right], \qquad    
	\mathbf{n} = \frac{(\tilde n_x, \tilde n_y)^T}{\sqrt{\tilde n_x^2 + \tilde n_y^2}},  
\end{equation} 
with $\mathbf{I}$ representing the identity matrix and $\mathbf{V} \cdot \mathbf{n}$ denoting the local normal mesh velocity. 
Also note that $\mathbf{n}$ is the spatial normalized normal vector, which is different from the spacetime normal vector $\mathbf{\tilde{n}}$. 
In this work, we will also adopt the 
the
HLL flux~\cite{hll}, which can be written  
\begin{equation} \label{eq:fluxhll}
    \mathcal{F}\left(\q_h^{n,-}, \q_h^{n,+}\right) = 
    \dfrac{S_r\, \tilde{\F}(\q_h^{n,-}) - S_\ell\,\tilde{\F}(\q_h^{n,+})}{S_r - S_\ell} + 
    \dfrac{S_{r}\,S_\ell}{S_r - S_\ell}\,\left(\q_h^{n,+} - \q_h^{n,-}\right),
\end{equation}
and we compute the necessary estimates of the minimum and maximum wave speeds as
\begin{equation}
    S_\ell = \min{\left(0, \lambda_{min}\left(\q_h^{n,-}\right), \lambda_{min}
        \left(\overline{\q}\right)\right)},\ \
    S_r = \max{\left(0, \lambda_{max}\left(\q_h^{n,+}\right), \lambda_{max}
        \left(\overline{\q}\right)\right)},\ \ \text{with }
     \ \overline{\q} = \dfrac{1}{2}\left(\q_h^{n,-} + \q_h^{n,+}\right),
\end{equation}
where $\lambda_{min}(\q)$ and $\lambda_{max}(\q)$ 
have been obtained in a similar manner to $s_{\max}$ for the Rusanov flux, 
and they represent the minimum and maximum
eigenvalues of the ALE Jacobian matrix in normal direction.

Alternatively, we also employ the less dissipative Osher-type~\cite{osherandsolomon,OsherUniversal} ALE flux
\begin{equation}
	\mathcal{F}(\q_h^{n,-},\q_h^{n,+}) \cdot \mathbf{\tilde n} = 
	\frac{1}{2} \left( \tilde{\F}(\q_h^{n,+}) + \tilde{\F}(\q_h^{n,-})  \right) \cdot \mathbf{\tilde n}_{ij}  - 
	\frac{1}{2} \left( \int_0^1 \left| \mathbf{A}^{\hspace{-1pt} \mathbf{V}}_{\mathbf{n}}(\boldsymbol{\Psi}(s)) 
	\right| ds \right) \left( \q_h^{n,+} - \q_h^{n,-} \right),  
	\label{eqn.osher} 
\end{equation} 
where we choose to connect the left and the right state across the discontinuity using a simple straight--line segment path  
\begin{equation}
	\boldsymbol{\Psi}(s) = \q_h^{n,-} + s \left( \q_h^{n,+} - \q_h^{n,-} \right), \qquad 0 \leq s \leq 1.  
	\label{eqn.path} 
\end{equation} 
The absolute value of $\mathbf{A}^{\hspace{-1pt} \mathbf{V}}_{\mathbf{n}}$ is evaluated as usual as 
$ \mathbf{R} |\boldsymbol{\Lambda}| \mathbf{R}^{-1}$,  
where $\mathbf{R}$, $\mathbf{R}^{-1}$ and $\boldsymbol{\Lambda}$ denote, respectively, the right eigenvector matrix, 
its inverse and the diagonal matrix of the eigenvalues of $\mathbf{A}^{\hspace{-1pt} \mathbf{V}}_{\mathbf{n}}$.

Next, the volume integrals in the above expression~\eqref{eqn.corrector} can be computed directly on the 
physical control volume $C_i^n$ by summing up the contributions on each subtriangular prism extruded 
from $\mathcal{T}(P_i^n)$, where Gaussian quadrature rules are easily available, see~\cite{stroud}. 
Instead, for the lateral spacetime surface $\partial C_{ij}^n$ we introduce a bi-linear parametrization 
that maps it to a reference square (even when it has a triangular shape, which is simply seen as 
a quadrilateral where one edge has been collapsed), see~\cite{gaburro2020high,Lagrange2D}.

Finally, we remark that the PDE integration over the spacetime volumes $C_i^n$ automatically satisfies
 the geometric conservation law (GCL) for all test functions $\phitilde_k$.
This follows from the application of the Gauss theorem to our \textit{closed} control volumes; a complete 
proof can be find in~\cite{Lagrange3D}.

\subsubsection{Well-balanced formulation of the ALE ADER-DG corrector step}
\label{ssec.nm-corrector-wb}

Now, we need to modify~\eqref{eqn.corrector} to make the scheme \textit{well-balanced}. 
To reach this objective, with an approach similar to the one introduced 
in~\cite{Ghosh2016,berberich2021high, gaburro2021well} for finite volume schemes, 
we start by writing~\eqref{eqn.corrector} for a given equilibrium solution $\Q_E$, 
\begin{equation}
	\begin{aligned} 
		\int_{P_i^{n+1}} \phitilde_k \Q_E(\x, t^{n+1}) \, d\x \,  = & \phantom{+}
		\int_{P_{i}^n} \phitilde_k  \Q_E(\x, t^n) \, d\x  
		\, - \, \sum_{j=1}^{N_{V_i}^{n,st}} \int_{\partial C_{ij}^n} 
		\phitilde_k  \mathcal{F}\left(\Q_E(\x^-, t),\Q_E(\x^+,t)\right) \cdot \mathbf{\tilde n} \, dS  \\
		& + \int_{C_i^n} 
		\tilde \nabla \phitilde_k \cdot \tilde{\F} (\Q_E)\, d\mathbf{x} dt \, + \, \int_{C_i^n} \phitilde_k \mathbf{S}(\Q_E) \, d\x dt, \quad \RIIcolor{ \forall k \in \{0, \dots, {\mathcal{L}(N,d)}-1\} },
	\end{aligned} 
	\label{eqn.corrector-eq}
\end{equation}
where for smooth equilibria, as those considered in this work, the boundary extrapolated value 
of $\Q_E$ coincides, i.e. $\Q_E(\x^-, t) = \Q_E(\x^+,t)$.
{Note that since $\Q_E$ is a solution of the PDE, see~\eqref{eq.equilibrium}, 
it is also a discrete equilibrium for the well-balanced scheme.} 

Here, to avoid any misunderstanding, we underline that integrals in~\eqref{eqn.corrector-eq}, 
as well as all the integrals of this paper, 
are computed numerically via Guassian quadrature formulas of the same order of accuracy of the DG method, 
thus relation~\eqref{eqn.corrector-eq} \textit{a priori} is not exact up to machine precision but only up to the order of the method.
However, the point of the WB reasoning is that the subtraction performed in the next lines will reduce the final discretization errors up to machine precision for equilibrium benchmarks.
 
So, we now subtract~\eqref{eqn.corrector-eq} from~\eqref{eqn.corrector} obtaining the 
\textit{update formula} written directly for the \textit{fluctuations} $\u_f$ w.r.t. $\Q_E$
\begin{equation}
	\begin{aligned} 
		\left( \,  \int_{P_i^{n+1}} \phitilde_k \phi_\ell \, d\x \,  \right )
		\hat{\u}_{f,\ell}^{n+1}   = &
		\phantom{+} \left( \,  \int_{P_{i}^n} 	\phitilde_k  \phi_\ell \, d\x  \, \right) \hat{\u}_{f,\ell}^n  \\
		& - \sum_{j=1}^{N_{V_i}^{n,st}} \int_{\partial C_{ij}^n} 
		\phitilde_k  \mathcal{F}(\q_h^{n,-},\q_h^{n,+}) \cdot \mathbf{\tilde n} \, dS  
		  \, + \,  \sum_{j=1}^{N_{V_i}^{n,st}} \int_{\partial C_{ij}^n} 
		\phitilde_k  \mathcal{F}(\Q_E,\Q_E) \cdot \mathbf{\tilde n} \, dS  \\
		& + \int_{C_i^n} 
		\tilde \nabla \phitilde_k \cdot \tilde{\F} (\q_h)\, d\mathbf{x} dt 
		  \, - \, \int_{C_i^n} 
		\tilde \nabla \phitilde_k \cdot \tilde{\F} (\Q_E)\, d\mathbf{x} dt  \\
		& + \int_{C_i^n} \phitilde_k \mathbf{S}(\q_h) \, d\x dt
		  \, - \, \int_{C_i^n} \phitilde_k \mathbf{S}(\Q_E) \, d\x dt, \quad \RIIcolor{ \forall k \in \{0, \dots, {\mathcal{L}(N,d)}-1\} }.
	\end{aligned} 
	\label{eqn.corrector-f}
\end{equation}
We remark that in the above well-balanced update formula the employed predictor $\q_h$ will be 
obtained, at each necessary \RIIcolor{spacetime} point $(\x,t)$ as
\begin{equation}
	\q_h(\x,t) = \Q_E(\x,t) + \q_f(\x,t),
\end{equation}
where $\q_f$ is the predictor approximating with high order of accuracy the evolution of 
the perturbations inside the spacetime control volumes, see~\eqref{eqn.predictor-f}.

\paragraph{WB properties} Moreover, we highlight that \textit{three fundamental properties} of a 
well-balanced scheme hold true for the update formula~\eqref{eqn.corrector-f}, provided that the 
same can be said for the predictor $\q_f$, as will be proved in the next Section~\ref{ssec.nm-predictor}.

First, whenever $\Q_{\text{IC}} = \Q_E$, then 
$\u_f^0 = \mathbf{0}$ and $\q_f^n=\mathbf{0}$, 
thus $\q_h(\x,t) = \Q_E(\x,t)$ in each \RIIcolor{spacetime} point and all the terms in~\eqref{eqn.corrector-f} vanish, so that and $\u_f^n = \mathbf{0}$ $\forall n$, i.e. the scheme is \textit{exact} on the equilibria. 
This will be also numerically substantiated 
in several of our benchmarks, see in particular Section~\ref{ssec.test_wb}.

Second, in~\eqref{eqn.corrector-f} all the computations requiring a physically valid state vector (fluxes, 
sources, Riemann solvers) are performed on the 
sum between the equilibrium and the fluctuations, not on the sole fluctuations.
This assures a safe use of the scheme in any situation and will be numerically shown in Section~\ref{ssec.nonequilibrium}.

The third one is related to the way the scheme is written: note indeed that 
in~\eqref{eqn.corrector-f} each operation performed on $\q_h$ is then repeated, with the opposite sign, on $\Q_E$. 
In this way, we drastically reduce the numerical errors arising in the case 
where $\u_h$ and $\q_h$ are just small perturbations of $\Q_E$, i.e. in the case where $\u_f$ and $\q_f$ are small. 
The fundamental benefits due to this feature will be apparent 
in the numerical results presented in Section~\ref{ssec.test_wb}.

\subsubsection{High order spacetime predictor and its well-balanced formulation}
\label{ssec.nm-predictor}

This subsection is dedicated to the so-called \textit{predictor} step.
It consists in \textit{locally}, iteratively, solving a weak form of the governing 
PDE~\eqref{eq.generalform} \textit{in the small}, as written in~\cite{eno}, \textit{inside} 
each spacetime control volume $C_i^n$ and $S_i^n$, 
starting from the information available at $t^n$, enforcing the causality principle, 
but without coupling interactions between different elements. 
The predictor provides, for each spacetime element, a spacetime polynomial, $\q_h$ or $\q_f$
of high order of accuracy both in space and time, which serves as the preliminary 
solution in~\eqref{eqn.corrector} and~\eqref{eqn.corrector-f} and is used for evaluating the numerical 
fluxes and the source terms at each spacetime quadrature point. 

To obtain the final predictor formula, we start again by multiplying the 
governing PDE~\eqref{eq.generalform} by a test function, 
this time the spacetime functions $\theta$ given in~\eqref{eq.Dubiner_phi}, 
next we integrate over $C_i^n$ and insert the discrete solution $\q_h^n$ in place of $\mathbf{Q}$, 
obtaining 
\begin{equation}
	\int_{C_i^n} \theta_k(\x,t) \de{\q_h^n}{t} \, d\x dt  + \int_{C_i^n} \theta_k(\x,t) 
	\nabla \cdot \F(\q_h^n) \, d\x dt  - \int_{C_i^n} \theta_k(\x,t) \, \S(\q_h^n) \, d\x dt 
	= \mathbf{0}, \quad \forall k \in \RIIcolor{\{}0, \dots, \mathcal{L}(N,d+1)\RIIcolor{\}} .
	\label{eqn.pred-step1}
\end{equation}
Then, we rewrite the first integral taking into account potential jumps of $\q_h$ at the boundary of $C_i^n$
via a simplified path-conservative approach~\cite{Pares2006,Castro2006,Castro2008},
with the test functions only taken from within $C_i^n$, combined with an upwinding of the fluxes in time
(thus, for prismatic control volumes, we apply it only on $P_i^n$ and  
for sliver elements see~\eqref{eqn.predictor-jumpformula-sliver}) obtaining 

\begin{equation}
	\begin{aligned} 
	  \int_{C_i^n \backslash P_i^n} \theta_k(\x,t) \de{\q_h^n}{t} \, d\x dt + 
	  \int_{P_i^n} \theta_k(\x,t^n) \, \left( \q_h^{n}(\x,t^n) - \u_h^n(\x,t^n) \right) \, d\x  & \\[1pt]
	   + \int_{C_i^n} \theta_k(\x,t) \nabla \cdot \F(\q_h^n) \, d\x dt  - 
	   \int_{C_i^n} \theta_k(\x,t) \, \S(\q_h^n) \, d\x dt & \, = \,  \mathbf{0},
	   \qquad \forall k \in \RIIcolor{\{}0, \dots, \mathcal{L}(N,d+1)\RIIcolor{\}} .
	\label{eqn.pred-step2}
	\end{aligned} 
\end{equation}
Finally, we insert the polynomial expansion of $\q_h^n$, see~\eqref{eqn.qh}, and we reorder the terms of the equation
\begin{equation}
	\begin{aligned} 
        &\phantom{\, -\,} \left ( \int_{C_i^n \backslash P_i^n} \theta_k(\x,t) \de{\theta_k(\x,t)}{t} \, d\x dt  + 
		\int_{P_i^n} \theta_k(\x,t^n) \theta_\ell(\x,t^n) 	\right ) \hat{\q}_\ell^n  		
		 = \left ( \int_{P_i^n} \theta_k(\x,t^n) \phi_\ell(\x,t^n) \, d\x  \right ) \hat{\u}_\ell^n  	\\
	    &\, - \, \left ( \int_{C_i^n} \theta_k(\x,t) \de{\theta_\ell(\x,t)}{x}   \, d\x dt    \right ) \hat{\f}_\ell^n
		\, - \, \left ( \int_{C_i^n} \theta_k(\x,t) \de{\theta_\ell(\x,t)}{y}   \, d\x dt    \right ) \hat{\g}_\ell^n		
		 \, + \, \left( \int_{C_i^n} \theta_k(\x,t) \theta_\ell(\x,t)  \, d\x dt \right ) \hat{\S}_\ell^n, \\
		&\hspace{9.85cm}  \forall k \in \RIIcolor{\{}0, \dots, \mathcal{L}(N,d+1)\RIIcolor{\}}.
		\label{eqn.pred-step3}
	\end{aligned} 
\end{equation}
With this last formulation it is easy to use a fixed point discrete Picard iteration (over the index $r$), 
as detailed in~\cite{dumbser2008unified,hidalgo2011ader,busto2020high}, which, at convergence, allows to recover the coefficients $\hat{\q}_\ell^n$
\begin{equation}
	\begin{aligned} 
		&\phantom{\, -\,} \left ( \int_{C_i^n \backslash P_i^n} \theta_k(\x,t) \de{\theta_k(\x,t)}{t} \, d\x dt  + 
		\int_{P_i^n} \theta_k(\x,t^n) \theta_\ell(\x,t^n) 	\right ) \,  \hat{\q}_\ell^{r+1}  		
		\ = \ \left ( \int_{P_i^n} \theta_k(\x,t^n) \phi_\ell(\x,t^n) \, d\x  \right ) \, \hat{\u}_\ell^n  	\\
		&\, - \, \left ( \int_{C_i^n} \theta_k(\x,t) \de{\theta_\ell(\x,t)}{x}   \, d\x dt    \right ) \, \hat{\f}_\ell^r
		\, - \, \left ( \int_{C_i^n} \theta_k(\x,t) \de{\theta_\ell(\x,t)}{y}   \, d\x dt    \right ) \,  \hat{\g}_\ell^r		
		\, + \, \left( \int_{C_i^n} \theta_k(\x,t) \theta_\ell(\x,t)  \, d\x dt \right ) \, \hat{\S}_\ell^r \\
		&\hspace{7.7cm}  \forall k \in \RIIcolor{\{}0, \dots, \mathcal{L}(N,d+1)\RIIcolor{\}}, \quad r = 1, \dots, 10.
		\label{eqn.predictor}
	\end{aligned} 
\end{equation}
where $\hat{\f}_\ell^{n/r}$, $\hat{\g}_\ell^{n/r}$ and $\hat{\S}_\ell^{n/r}$ are 
the projection of the nonlinear fluxes and source terms in the chosen spacetime polinomial space.

In the above iterative procedure~\eqref{eqn.predictor} the initial guess for $\hat{\q}_\ell^{r=0}$ 
is taken equal to $\hat{\u}_\ell^n$
for the common spatial degrees of freedom and zero for the other ones. We also remark that the 
procedures~\eqref{eqn.predictor} has been proved to be convergent and to yield the formal order of 
accuracy $N$ in $N+1$ iterations, see~\cite{busto2020high,han2021dec} for more details.

\paragraph{WB formulation} Now, in order to make~\eqref{eqn.predictor} \textit{well-balanced}, we 
start again by writing the predictor formula for the equilibrium~$\Q_E$
\begin{equation}
	\begin{aligned} 
		&\phantom{\, -\,}  \int_{C_i^n \backslash P_i^n} \theta_k(\x,t) \de{\,\Q_E(\x,t)}{t} \, d\x dt  + 
		\int_{P_i^n} \theta_k(\x,t^n)\, \Q_E(\x,t^n) 	
		\ = \ \int_{P_i^n} \theta_k(\x,t^n) \, \Q_E(\x,t^n) \, d\x  	\\
		&\, - \,  \int_{C_i^n} \theta_k(\x,t) \de{\, \f (\Q_E (\x,t) ) }{x}   \, d\x dt  
		\, - \,  \int_{C_i^n} \theta_k(\x,t) \de{\, \g ( \Q_E (\x,t) ) }{y}   \, d\x dt   
		\, + \, \int_{C_i^n} \theta_k(\x,t) \, \S( \Q_E (\x,t) ) \, d\x dt,  	\\
		&  \quad \RIIcolor{ \forall k \in \{0, \dots, {\mathcal{L}(N,d)}-1\} },
		\label{eqn.predictor-eq}
	\end{aligned} 
\end{equation}
and we subtract~\eqref{eqn.predictor-eq} from~\eqref{eqn.predictor} obtaining 
\begin{equation}
	\begin{aligned} 
		\left ( \int_{C_i^n \backslash P_i^n} \theta_k(\x,t) \de{\theta_k(\x,t)}{t} \, d\x dt  \right . & \left. + 
		\int_{P_i^n} \theta_k(\x,t^n) \theta_\ell(\x,t^n) 	\right ) \, \hat{\q}_{f, \ell}^{r+1} \  = 	\\[1pt]
		& \phantom{\, + \,} \left ( \int_{P_i^n} \theta_k(\x,t^n) \phi_\ell(\x,t^n) \, d\x  \right ) \, \hat{\u}_\ell^n 
		- \int_{P_i^n} \theta_k(\x,t^n) \, \Q_E(\x,t^n) \, d\x 	\\[1pt]
		&\, - \, \left ( \int_{C_i^n} \theta_k(\x,t) \de{\theta_\ell(\x,t)}{x}   \, d\x dt    \right ) \, \hat{\f}_\ell^r 
		\, + \,  \int_{C_i^n} \theta_k(\x,t) \de{\, \f (\Q_E (\x,t) ) }{x}   \, d\x dt  \\[1pt]
		&\, - \, \left ( \int_{C_i^n} \theta_k(\x,t) \de{\theta_\ell(\x,t)}{y}   \, d\x dt    \right ) \, \hat{\g}_\ell^r		
		\, + \,  \int_{C_i^n} \theta_k(\x,t) \de{\, \g ( \Q_E (\x,t) ) }{y}   \, d\x dt   \\[1pt]
		&\, + \, \left( \int_{C_i^n} \theta_k(\x,t) \theta_\ell(\x,t)  \, d\x dt \right ) \, \hat{\S}_\ell^r 
		\, - \, \int_{C_i^n} \theta_k(\x,t) \, \S( \Q_E (\x,t) ) \, d\x dt  \\[2pt]
		&\hspace{0.6cm}  \forall k \in \RIIcolor{\{}0, \dots, \mathcal{L}(N,d+1)\RIIcolor{\}}, \quad r = 1, \dots, 10,
		\label{eqn.predictor-f}
	\end{aligned} 
\end{equation}
where an additional convergence criterion is evaluated by measuring the difference between two subsequent 
iteration values of $\q_f^r$, so that the fixed point procedure can be stopped earlier than the preset maximum of 
10 iterations. Note that the 
equilibrium part (the right terms of the equations) can be computed only once (per element and timestep) 
and re-used throughout all iterations of a given predictor step.
We also highlight that the result of the procedure~\eqref{eqn.predictor-f} at each iteration $r$ 
is the sole value of $\hat{\q}_{f, \ell}^{r}$ i.e. the coefficients of the predictor of the fluctuations. 
So to compute $\hat{\f}_\ell^r$, $\hat{\g}_\ell^r$, $\hat{\S}_\ell^r$ we sum, in each needed quadrature point, 
the exact value of the equilibrium with the value of the predictor of the fluctuations $\q_f^{r}$ available from the previous iteration $r-1$.

\paragraph{WB properties}  By inspecting this well-balanced formulation of the predictor it can be seen that 
the fundamental properties highlighted for the corrector step in Section~\ref{ssec.nm-corrector-wb} also hold true for $\q_f$. 
Indeed, first, when $\u_f=\mathbf{0}$ also $\q_f=\mathbf{0}$,
because in each line of the formula~\eqref{eqn.predictor-f}
the left and the right terms assume exactly the same value,
being the {integrands} equal in each quadrature point (due to the way we obtain $\hat{\f}_\ell^r$, $\hat{\g}_\ell^r$, $\hat{\S}_\ell^r$). 
So the predictor is \textit{exact} on equilibria. 
Second, the flux and source evaluations are again performed on physically meaningful state vectors 
(the sum between equilibrium solution and fluctuations, not on the fluctuations alone) so the predictor is valid also 
in non-equilibrium situations. 
And, third, the continuous subtraction between operations performed on $\q_h^r$ and $\Q_E$ 
drastically reduces the numerical errors obtained in the case where $\u_h$ is just a small perturbation of $\Q_E$.

\subsubsection{Notes on the numerical treatment of hole-like sliver elements}
\label{ssec.nm-sliver}

We address in this section the strategies we have adopted 
to overcome the main difficulties, or unusual situations, arising from the application of our numerical 
scheme over the \textit{hole-like sliver} elements.

\paragraph{Quadrature points} First, the most straightforward issue regards the choice of \textit{quadrature points} to employ over them. Our 
hole-like elements are simply stretched tetrahedra with triangular lateral faces, so for them we can 
use standard quadrature points both for the volume integrals and the surface integrals.

\paragraph{Initial condition for the predictor step} 
Then, we need to handle the fact that our sliver elements at time $t^n$ just coincide with an edge of the tessellation, 
where in principle discontinuities are located, and 
hence no \textit{uniquely} defined value for $\u_h$ and $\u_f$ at~$t^n$ is present.
This poses a problem in the \textit{predictor} step because 
\textit{i)} we cannot use the trivial initial guess for the predictor $\hat{\q}_{h/f}^{r=0} = \hat{\u}_{h/f}$, and 
\textit{ii)} we do not have an easy initial condition to feed the weak formulation~\eqref{eqn.predictor} with 
information coming from the past (enforcing the causality principle), 
namely we cannot compute the term
\begin{equation}
	\left ( \int_{P_i^n} \theta_k(\x,t^n) \phi_\ell(\x,t^n) \, d\x  \right ) \, \hat{\u}_\ell^n
	\label{eqn.predicotr-t0info}
\end{equation}
which is for example the last term of the first line of~\eqref{eqn.predictor} 
and the first term of the second line of~\eqref{eqn.predictor-f}.

For the initial guess ${\q}_{h/f}^{r=0}$, 
we simply average the spatial information found in its four (degenerate) neighbors.
Note that since this is used only as an initial \textit{guess} and will be corrected 
during the predictor iterations, 
no accuracy is lost due to this simple choice.

Instead, for the term~\eqref{eqn.predicotr-t0info}, which is missing in the case of a sliver element, 
we recall that it comes from this choice 
\begin{equation}
	\begin{aligned} 
	\label{eqn.predictor-jumpformula-standard}
	\int_{C_i^n} \theta_k(\x,t) \de{\q_h^n}{t} \, d\x dt  & = 	\int_{C_i^n \backslash P_i^n} 
	\theta_k(\x,t) \de{\q_h^n}{t} \, d\x dt + \int_{P_i^n} \theta_k(\x,t^n) \, \left( \q_h^{n}(\x,t^n) - \u_h^n(\x,t^n) \right) \, d\x, \\[1pt]
	& \quad \RIIcolor{ \forall k \in \{0, \dots, {\mathcal{L}(N,d)}-1\} },
	\end{aligned} 
\end{equation}
where we have used a path conservative approach to treat the boundary of $C_i^n$, 
but considering only the \textit{easiest} information (because already explicitly available) 
coming from the past, i.e. from $P_i^n$.
Since, the fact of considering the bottom part of a sliver element has no meaning, 
we extend~\eqref{eqn.predictor-jumpformula-standard} by including all the neighbors of $S_i^n$
\begin{equation}
	\begin{aligned} 
	\int_{S_i^n} \theta_k(\x,t) \de{\q_h^n}{t} \, d\x dt & = \int_{S_i^n \backslash \partial S_i^n} \theta_k(\x,t) \de{\q_h^n}{t} \, d\x dt +
	 \sum_{j=1}^{4}
	 \int_{\partial S_{ij}^n} \theta_k(\x,t) \, \left( \q_h^{n,+} - \q_h^{n,-} \right) \, \tilde{\mbf{n}}_t^- \, dS,\\
	 &  \quad \RIIcolor{ \forall k \in \{0, \dots, {\mathcal{L}(N,d)}-1\} },
	 \end{aligned} 
	\label{eqn.predictor-jumpformula-sliver}
\end{equation}
where, by writing $\tilde{\mbf{n}}_t^-$ we mean that we only consider the spacetime neighbors 
$C_j^n$ whose common surface $\partial S_{ij}^n = S_i^n \cap C_j^n$ exhibits a \textit{negative} 
time component of the outward pointing spacetime normal vector, namely $ \tilde{\mbf{n}}_t^- < 0$, 
so that the causality principle is still satisfied. 
We also remark that, since in our algorithm a sliver element is always surrounded by four standard neighbors, 
we can first compute their predictor and then use their values for the 
formula~\eqref{eqn.predictor-jumpformula-sliver} so that the predictor remains an explicit formula for each element. 
With the choice of~\eqref{eqn.predictor-jumpformula-sliver}, 
the term~\eqref{eqn.predicotr-t0info} can be easily substituted both in~\eqref{eqn.predictor} and in~\eqref{eqn.predictor-f} with
\begin{equation}
	 \sum_{j=1}^{4}
	\int_{\partial S_{ij}^n} \theta_k(\x,t) \, \left( \q_h^{n,+} - \q_h^{n,-} \right) \, \tilde{\mbf{n}}_t^- \, dS,
	\label{eqn.predicotr-sliver-t0info}
\end{equation}
so that a high order accurate predictor polynomial is computable also for sliver 
elements both in the standard and in the well-balanced version of our method.

\paragraph{Flux-redistribution around a topology change} 
Finally, with regards to sliver elements, we have a last 
matter to address
with the final update formula~\eqref{eqn.corrector} and~\eqref{eqn.corrector-f};
in both the cases, the two terms in the first line are zero by construction,
since zero are the area values of the slivers at $t^n$ and $t^{n+1}$,
but the sum of the other terms does not vanish at the discrete level,
due to the time integration being explicit.

\newcommand{\areanew}{|\Omega_i^{n+1}|}
\newcommand{\areaold}{|\Omega_i^{n}|}
\newcommand{\sliverfluxdof}{\mathbf{F}_{ij}}
\newcommand{\sliverfluxavg}{\bar{\mathbf{F}}_{ij}}
\newcommand{\sliverfluxavgk}{\bar{\mathbf{F}}_{ij}^\nu}
\newcommand{\slivercellaveragenew}{\Q_i^{n+1}}
\newcommand{\slivercellaverageold}{\Q_i^{n}}
\newcommand{\sliverspacetimefacearea}{\partial S_{ij}^n}
\newcommand{\sliverspacetimevolume}{|S_{i}^n|}
\newcommand{\sliveraveragedsource}{\bar{\mathbf{S}}_i}
\newcommand{\sliveraveragedsourcek}{\bar{\mathbf{S}}_i^\nu}

To restore discrete (machine precision-accurate) conservativity, we introduce a simple \textit{flux-rescaling procedure}, which
is based on enforcing that the integral form of the governing equations be satisfied exactly
over the hole-like sliver spacetime volume.

The flux-rescaling method assumes that the update rule for the (never computed) cell average value for the
sliver can be cast as
\begin{equation}
\areanew\,\slivercellaveragenew = \areaold\,\slivercellaverageold - \sum_{j=1}^4 \sliverfluxavg\,\sliverspacetimefacearea + 
\sliverspacetimevolume\,\sliveraveragedsource, 
\end{equation}
where $\areaold = 0$ and $\areanew = 0$ are the areas relative to the hole-like sliver element
at times $t^n$ and $t^{n+1}$ respectively, while $\slivercellaverageold$ and $\slivercellaveragenew$
are the corresponding cell averages, which are never used explicitly in the algorithm, 
since they would pertain to zero-area spatial control volumes at every discrete time level.
We indicate with $\sliverfluxavg$ the spacetime face-averaged numerical flux across one of the four faces 
of the sliver, and with $\sliverspacetimefacearea$ the extension of such a spacetime surface.
Analogously $\sliveraveragedsource$ is the integral average of the algebraic source term over the sliver spacetime control volume.
Since by definition the area associated with the sliver control volume vanishes at times $t^n$ and $t^{n+1}$, 
the constraint 
\begin{equation}
0 = 0 - \sum_{j=1}^4 \sliverfluxavg\,\sliverspacetimefacearea + 
\sliverspacetimevolume\,\sliveraveragedsource
\end{equation}
has to be explicitly imposed, being in principle verified only up to the accuracy of the numerical method
instead of machine precision.

In previous works~\cite{gaburro2020high, gaburro2021bookchapter} the constraint was enforced by merging 
the hole-like sliver element with one of its standard spacetime neighbors, while here we introduce 
a simple procedure that allows for a more fine-grained redistribution of the constraint violation.
In particular, what we want to impose is not only that the scheme be indeed discretely conservative, 
but also that the magnitude of the required flux correction be as limited as possible.
A straightforward way of achieving both of these goals is to rescale each one of the four average fluxes, 
conserved variable by conserved variable (indexed through $\nu$), with a scalar 
coefficient $\alpha_{j}^\nu$ (local to each spacetime face indexed by $j$).
Formally, we can then write the discrete constraint with respect to the rescaled fluxes $\alpha_{j}^\nu\, \sliverfluxavgk$ as 
\begin{equation} \label{eq:sliverdiscreteconstraint}
0 = -\sum_{j=1}^4 \alpha_{j}^\nu\,\sliverfluxavgk\,\sliverspacetimefacearea + \sliverspacetimevolume\,\sliveraveragedsourcek.
\end{equation}
Subsequently we seek a set of $\alpha_j^\nu$ such that, for each conserved variable $\nu$, equation~\eqref{eq:sliverdiscreteconstraint} is
satisfied exactly. Moreover, we also minimize the deviation from unity of all coefficients $\alpha_j^\nu$.
We remark that $\alpha_j^\nu = 1$ would simply mean that the constraint was already satisfied exactly by the 
preliminary fluxes generated by the high order direct ALE scheme, and that the closer $\alpha_j^\nu$ is to unity, 
the smaller the entity of the flux correction will be.
Then the scale factors $\alpha_j^\nu$ are computed by minimizing $f^\nu = \sum_{j=1}^4\left(\alpha_j^\nu - 1\right)^2$ under
the constraint~\eqref{eq:sliverdiscreteconstraint}.
As a physical consequence of the structure of such a minimization problem, we have that each flux 
will tend to be modified proportionally to its magnitude or the magnitude of the face through which 
it flows. This also implies that fluxes across faces with zero or small magnitude (due to either a small face or
stagnant flux) will remain essentially unchanged
and that inversion of the sign of fluxes is highly unlikely, hence preserving the qualitative 
physical behavior naturally captured by the numerical scheme.

Finally, the correction can be applied to 
update the data in all spacetime neighbors of the sliver elements, by multiplying 
 each of the numerical fluxes $\mathcal{F}$ evaluated
in the general DG update formulas~\eqref{eqn.corrector-f} and~\eqref{eqn.corrector}, 
by the face-local rescaling factor $\alpha_{j}^\nu$
and using these rescaled update formulas in place of the standard ones.

\subsection{\textit{A posteriori} subcell finite volume limiter}
\label{sec.nm-limiter}

The high order scheme that we have presented up to now is linear in the sense of Godunov~\cite{godunov},
that is,  
their update rule is not data-dependent and the state vectors at the new time levels 
are linear functions of the old data and the corresponding fluxes.
Hence, as proven by the Godunov theorem~\cite{godunov}, starting from its second order version, 
the method may produce dangerous oscillations in presence of discontinuities; 
to avoid them we need to introduce a last ingredient in our high order methodology: a \textit{nonlinear limiting} strategy. 

Among the different approaches available in literature, as those inspired to Cockburn and Shu~\cite{cockburn1990runge,cockburn2000development},
based on the use of a total variation bounded limiter, 
or the moment limiters~\cite{krivodonova2004shock}, 
the artificial viscosity procedures~\cite{persson2006sub} 
WENO-type limiters~\cite{qiu2004hermite,qiu2005runge}, 
or gradient-based limiters~\cite{kuzmin2020gradient,kuzmin2021new},
we have selected the so-called \textit{a posteriori} subcell finite volume (FV) limiter. 
This type of limiter is based on the MOOD approach~\cite{CDL1,ADERMOOD}, which has already 
been successfully applied in the ALE finite volume framework in~\cite{ALEMOOD1,ALEMOOD2} and in 
the discontinuous Galerkin case in~\cite{SonntagDG,Sonntag2,DeLaRosaMunzDGMHD,hajduk2019new,vilar2019,markert2021sub,rueda2022subcell,popov2023space}
and, with a notation similar to the one used here, 
in~\cite{DGLimiter1,DGLimiter2,DGLimiter3,Zanotti2015d,kemm2020simple,gaburro2020high,gaburro2021posteriori}.
We finally remark that shock-capturing techniques, based on subcell finite volume schemes, can also be
applied in a predictive 
(\textit{a priori}) fashion, for example as in~\cite{SonntagDG, Sonntag2, beck2020, mossier2022, gao2024}.
While referring to the cited literature for more details, 
and in particular to~\cite{gaburro2020high} for the complete description of our limiter on moving 
polygonal meshes, here we just briefly recall the key passages and illustrate the small details necessary
to make the resulting scheme well-balanced. 

Since we choose to work with a limiter which acts \textit{a posteriori}, 
we first run our unlimited DG scheme (as presented in the above sections) everywhere on the domain 
and we obtain 
the updated solution at time $t^{n+1}$ that we consider as a \textit{candidate} solution: $\u_{h/f}^{n+1,*}$. 

Also, since the chosen limiter acts at the level of \textit{subcells}, 
for each polygon we consider its subtriangulation 
$\mathcal{T}(P_i^n)$
and we further subdivide each triangle in~$N^2$ smaller subtriangles that we call $s_{i,\alpha}^n$ 
with $\alpha = 1, \dots, N_{V_i}^{n,st} \cdot N^2 $,  with $|s_{i,\alpha}^n|$ their area. 

We then check the admissibility of the {candidate} solution $\u_h^{n+1,*}$, 
i.e. we verify that both the cell averages of $\u_h^{n+1,*}$ on each $s_{i,\alpha}^{n+1}$,
obtained via the \RIIcolor{projection operator}
\begin{equation}
	\RIIcolor{\mathbf{v}_{i,\alpha}^{n+1} = \frac{1}{|s_{i,\alpha}^{n+1}|} \int_{s_{i,\alpha}^{n+1}} 
	\mathbf{u}_{h}^{n+1}(\x,t^{n+1}) \, d\x :=\mathcal{P}(\mathbf{u}_h^{n+1}),}
	\label{eqn.projection}
\end{equation}
and the values of $\u_h^{n+1,*}$ on each vertex of $s_{i,\alpha}^{n+1}$ 
satisfy certain conditions.
They \textit{i)} should be acceptable from a numerical point of view, i.e.
avoid \textit{not-a-number} or \textit{infinity} values potentially found in degenerate solutions of the 
unlimited DG scheme;
\textit{ii)} should be valid from a physical point of view, i.e. densities and pressures should be positive 
numbers and eventually other physical criteria can be checked, see~\cite{guermond2018second,guermond2019invariant} 
and also~\cite{kuzmin2021entropy,kuzmin2022limiter,hajduk2019new} for entropy based limiter;
and \textit{iii)} they should satisfy a relaxed discrete maximum principle, to explicitly 
enforce absence of overshoots and undershoots in the solution, see~\cite{DGLimiter1,DGLimiter2}.

In particular, we remark that in the \textit{well-balanced} case we check the admissibility 
of $\u_{h}^{n+1,*}$ not of $\u_{f}^{n+1,*}$. 
To recover the values of $\u_{h}^{n+1,*}$ at a vertex $(\x, t)$ of $s_{i,\alpha}^n$,
we simply sum $\Q_E(\x,t)$ with $\u_f^{n+1}(\x,t)$; 
instead, to obtain the cell average on $s_{i,\alpha}^{n+1}$, we project $\u_f^{n+1}$ on the 
triangle $s_{i,\alpha}^{n+1}$ via~\eqref{eqn.projection} and we sum the obtained value 
with the cell average of the equilibrium $\Q_E$.

For all the polygons where the {candidate} solution satisfies the above criteria we simply take 
$\u_{h/f}^{n+1} = \u_{h/f}^{n+1,*}$.
Instead, when the {candidate} solution in $P_i^{n+1}$ is found to be \textit{troubled}, 
we reject it and we \textit{locally} recompute the solution inside $P_i^{n+1}$ with a more robust 
scheme, namely a finite volume scheme. 
Moreover, in order to preserve as much as possible the resolution properties of 
the high order DG scheme, 
we apply the FV method at the subtriangulation level.
This means that via the projection operator~\eqref{eqn.projection} we compute the cell averages 
$\mathbf{v}_{i,\alpha}^{n}$ coming from the available solution at time $t^n$.
We then evolve $\mathbf{v}_{i,\alpha}^{n}$ through a finite volume method and we recover the cell 
averages values $\mathbf{v}_{i,\alpha}^{n+1}$.
From those robust values we then reconstruct the DG polynomial by means of the following reconstruction operator 
\begin{equation}
	\int_{s_{i,\alpha}^n} \mathbf{u}_{h}^{n+1}(\x,t^{n+1}) \, d\x = \int_{s_{i,\alpha}^n} 
	\RIIcolor{\mathbf{v}_{i,\alpha}^{n+1} \, d\x :=\mathcal{R}(\mathbf{v}_{i,\alpha}^{n+1})}  \qquad \forall \alpha, 
	\label{eqn.intRec}
\end{equation}
that we enforce conservation on the main cell $P_i^{n+1}$ via the additional linear constraint
\begin{equation}
	\int_{P_i^{{n+1}}} \mathbf{u}_{h}^{n+1}(\x,t^{n+1}) \, d\x = \int_{P_i^{{n+1}}} \RIIcolor{\mathbf{v}_{h}^{n+1}} \, d\x.
	\label{eqn.LSQ}
\end{equation}

To ensure that this last part of the algorithm does in fact maintain the
\textit{well-balanced} property, we need to 
split also the cell averages $\mathbf{v}_{i,\alpha}^{n}$ in the equilibrium part and the 
fluctuations part $\mathbf{v}_{f,i,\alpha}^{n}$ 
and then apply a well-balanced ALE finite volume methodology for the evolution of the fluctuations $\mathbf{v}_{f,i,\alpha}^{n}$.
For example, a first order ALE FV method simply consists in taking the formula~\eqref{eqn.corrector-f} with $N=0$ and $\q_h^n=\u_h^n$.
For a higher order FV technique instead one can combine a 
\textit{i)} high order well-balanced \textit{spatial} reconstruction technique, as those presented 
in~\cite{xing2006high,castro2007well,castro2020well,grosheintz2019high,gaburro2018wellb,berberich2021high}, 
(that, being applied at a single spatial level, is independent on the ALE framework),
with \textit{ii)} the well-balanced ALE predictor step~\eqref{eqn.predictor-f} 
(where, instead of $\u_h^n$, one could use the high order polynomial value coming from the 
well-balanced reconstruction). 

Finally, we remark that the reconstruction operator~\eqref{eqn.intRec}-\eqref{eqn.LSQ} might 
still lead to oscillations, since it reconstructs $\u_{h/f}^{n+1}$ from 
the subcell finite volume data via a \textit{linear} procedure. 
When this happens, the cell $P_i^{n+1}$ will be 
detected as \textit{troubled} also at the next time level $t^{n+2}$, 
and thus will be treated again with a finite volume method which will evolve the cell averages 
$\v_{i,\alpha}^{n+1}$ and obtain $\v_{i,\alpha}^{n+2}$.
This means that whenever the limiter activates on the same cell for two consecutive times, 
to ensure the robustness of our scheme, the cell averages to start the FV scheme at the second 
timestep, $\v_{i,\alpha}^{n+1}$, 
will not come from the projection~\eqref{eqn.projection} of the reconstructed solution obtained 
at the previous timestep, but will be directly the cell averages available from the subcell FV 
scheme at the previous timestep. For this reason, we always store in memory both 
the polynomial data and the computed cell 
averages until a troubled polygon is eventually found to be acceptable at the next timestep.

\section{Numerical results}
\label{sec.tests}

In order to validate our novel algorithm and demonstrate its capabilities we present a series of test cases 
aimed at showing all the fundamental properties of a functional and effective Lagrangian \textit{well-balanced} high order scheme.
First, we verify some basic properties: we show numerically the high order of convergence and
that the scheme is able to handle correctly non-equilibrium situations. 
Next, we show the exact preservation of equilibrium solutions for very long simulation times and large mesh deformations. 
Then, we showcase the capabilities of our novel scheme in describing small 
perturbations arising around moving equilibrium profiles, which would be not properly captured by 
standard numerical methods in a reasonable computational time. 
In particular, we highlight the role of both the Lagrangian framework and the well-balancing in reaching this last objective. 

We verify these key properties of our scheme on three sets of hyperbolic governing equations, 
namely the Euler equations of gasdynamics with and without gravity and the equations of ideal magnetohydrodynamics,
which are briefly recalled below. 

\subsection{The considered governing equations and the general setting of our benchmarks}
\label{ssec.test_equations}

In this section we recall the expression of the hyperbolic systems considered in this work and we 
describe the general setting of our benchmarks.

All our numerical tests are carried out on unstructured two-dimensional polygonal tessellations 
whose generators move in a Arbitrary-Lagrangian-Eulerian fashion, i.e. 
following as close as possible the local fluid velocity but also applying carefully designed smoothing techniques 
that guarantee a high quality of the moving mesh while preserving the Lagrangian character of the algorithm. 

The CFL safety factor in~\eqref{eq:timestep} is taken to be $0.1$ for the first $10$ time steps 
of each simulations (as suggested in~\cite{torobook}) and is always set to CFL~$ = 0.5$ afterwards. We
verify the convergence order and the WB property of our approach by testing all the schemes from the
$\mathbb{P}_1$, order $2$, to the $\mathbb{P}_4$, order $5$; then, we
select the $\mathbb{P}_2$ scheme of order $3$ to present all the comparisons between
the standard Eulerian DG framework and our non well-balanced and well-balanced ALE DG approach.

\subsubsection{Euler equations of gasdynamics with and without gravity source term}
\label{ssec_test_eq_euler}

The Euler equations of compressible gas dynamics represent a well-known system of hyperbolic 
equations that can be cast in the form~\eqref{eq.generalform} by taking $\Q$, $\F$ and, optionally, 
the gravity source term $\S$, as follows
\begin{equation}
	\label{eq.eulerTerms}
	\Q = \left( \begin{array}{c} \rho   \\ \rho u  \\ \rho v  \\ \rho E \end{array} \right), \quad
	\mathbf{F} = \left( \begin{array}{ccc}  \rho u       & \rho v        \\ 
		\rho u^2 + p & \rho u v          \\
		\rho u v     & \rho v^2 + p      \\ 
		u(\rho E + p) & v(\rho E + p)   
	\end{array} \right), 
	\qquad \S = \, 0  \ \text{ or } \ \,
	\left( \begin{array}{c} 0 \\[3pt] -\cos(\phi) \rho \, \frac{G \, m_s}{r^2} \\[3pt] -\sin(\phi) 
	\rho \, \frac{G \, m_s}{r^2} \\[3pt] -  (u_x \cos(\phi)+ u_y \sin(\phi)) \rho \frac{G  m_s}{r^2}  \end{array} \right).
\end{equation}
The vector \RIIcolor{$\Q = \Q(\x,t)$} of the conserved variables includes the fluid density $\rho$, the momentum  
vector $\rho \v=(\rho u, \rho v)$ and the total energy density $\rho E$. The corresponding set of 
primitive variables is $\mathbf{V} = (\rho, u, v, p)$.

In this work we always discretize the equations as they are written in Cartesian coordinates; 
however, sometimes it is convenient to employ the cylindrical coordinates $(r,\phi)$ according to 
the conventional relations $x = r \cos(\phi), \ y = r \sin(\phi)$ to express the gravity source term or some equilibrium profiles. 
We will also denote by $u_r$ and $u_\phi$ the radial and the angular component of the velocity, linked to $u$ and $v$ by 
\begin{equation} 
\label{eq.PolarVelocities}
u =  \cos(\phi) u_r - \sin(\phi)  u_\phi,   \qquad 
v =  \sin(\phi) u_r + \cos(\phi) u_\phi. 
\end{equation}
The Euler system is closed with the ideal gas equation of state which relates the 
fluid pressure $p$ to $\Q$ with
\begin{equation}
	\label{eqn.eos} 
	p = (\gamma-1)\left(\rho E - \frac{1}{2} \rho \mathbf{v}^2 \right), 
\end{equation}
where $\gamma$ is the ratio of specific heats so that the speed of sound takes the form $c=\sqrt{\frac{\gamma p}{\rho}}$.
To characterize parametrically the gravity source term, we take $G=1$ and $m_s=1$. 

We recall that the Euler equations, with $\S=0$ are a fundamental system used for simulations in 
gas and fluid-dynamics and, with the (externally given) gravity source term written above, they 
can be used to study prototype problems in computational astrophysics connected for example with 
the rotation of gas clouds around central objects as in the case of Keplerian 
disks~\cite{klingenberg2019arbitrary,gaburro2018well,kappeli2014well}.

For the test case concerning these sets of equations, we set $\gamma=1.4$ and, 
as Riemann solver, we use either the 
HLL~\eqref{eq:fluxhll} or the Osher Riemann solvers~\eqref{eqn.osher}.

\subsubsection{Equations of ideal magnetohydrodynamics}
\label{ssec_test_eq_mhd}

Furthermore, we consider the equations of classical ideal magnetohydrodynamics (MHD) that also 
include the modeling of the magnetic field $\B=(B_x,B_y,B_z)$. 
The state vector $\Q$ and the flux tensor $\F$ for writing the MHD equations in the general form~\eqref{eq.generalform} are 
\begin{equation}
	\label{eq.MHDTerms}
	\Q=\left( \begin{array}{c} \rho \\ \rho \mathbf{v} \\ \rho E \\ \B \\ \psi \end{array} \right) \ \text{ and } \ \quad \F = 
	\left( \begin{array}{c}  \rho \v  \\ 
		\rho \v \otimes \v + p_t \mathbf{I} - \frac{1}{4 \pi} \B \otimes \B \\ 
		\v (\rho E + p_t ) - \frac{1}{4 \pi} \B ( \v \cdot \B ) \\ 
		\v \otimes \B - \B \otimes \v + \psi \mathbf{I} \\
		c_h^2 \B  \end{array} \right),
\end{equation}
where $p_t=p+\frac{1}{8\pi}\mathbf{B}^2$ is the total pressure; 
the system is then closed by the following equation of state 
\begin{equation}
	p = \left(\gamma - 1 \right) \left(\rho E - \frac{1}{2}\mathbf{v}^2 - \frac{\mathbf{B}^2}{8\pi}\right).
	\label{MHDeos}
\end{equation}
We remark that the MHD system requires an additional constraint, on the divergence of the magnetic 
field, to be satisfied, i.e.
$ \nabla \cdot \mathbf{B} = 0$; for this reason, following~\cite{Dedneretal}, we have included one 
additional scalar PDE for the evolution of a so-called cleaning variable $\psi$, which is used to 
transport divergence errors outside the computational domain with an artificial divergence cleaning speed $c_h$. 

For all the test case concerning this set of equations, we set $\gamma=1.4$, and
as Riemann solver we use the HLL flux~\eqref{eq:fluxhll}.

\subsection{Order of convergence of our WB ALE ADER-DG scheme}
\label{ssec.test_order}

We consider two translating vortical solutions of the considered hyperbolic systems 
and we numerically verify the order of convergence of our well-balanced scheme on moving meshes. 
Here, it is important to remark that our well-balanced scheme is able to preserve any known 
stationary solution with machine precision, 
if the chosen initial condition coincides with the equilibrium profile selected to be preserved 
(this fact will be also shown in the next Section~\ref{ssec.test_wb}). 
For this reason, in order to show the order of convergence of our scheme,
we prescribe an initial condition that, despite being a stationary solution, is chosen to be 
different from the equilibrium preserved by the scheme.
As a clarification, note that such a translating solution can be seen as stationary if 
a scheme is properly capturing the Lagrangian motion of the bulk flow, while for an Euleran 
observer the term stationary would be in
contrast with the fact that a radially symmetric vortex is indeed translating in space.

\subsubsection{A moving vortical solution of the Euler equations: the Shu vortex} 
\label{ssec.test_order_euler}

\begin{table}[!tp] 
	\centering
	\numerikNine
	\begin{tabular}{|l|ccccc|cc|}
		\hline
		\multicolumn{8}{|c|}{ {Moving vortex solution of the Euler equations 
		    with $\epsilon=5$ and $\epsilon_E=0.5$} }\\
		\hline
		\hline
		& $\overline{h}(\Omega(t_f))$  & No. timestep  & No. sliver  & $L_2(\rho)$ error & $L_2(p)$ error &  order $\rho$ & order $p$ \\
		\hline	
		\hline
		\multirow{4}{*}{\rotatebox{90}{{DG $\mathbb{P}_1$}}}
		& 4.67e-02 &   449 &   260 & 1.18e-03 & 3.85e-03 & - & - \\     
		& 3.04e-02 &   688 &   678 & 4.98e-04 & 1.62e-03 & 2.0 & 2.0 \\ 
		& 2.33e-02 &   899 &  1052 & 2.90e-04 & 9.38e-04 & 2.0 & 2.1 \\ 
		& 1.83e-02 &  1144 &  1772 & 1.78e-04 & 5.83e-04 & 2.0 & 2.0 \\ 
		\cline{2-6}
		\hline
		\hline
		\multirow{4}{*}{\rotatebox{90}{{DG $\mathbb{P}_2$}}}
		& 8.53e-02 &   470 &    66 & 2.47e-04 & 9.73e-04 & - & - \\     
		& 4.67e-02 &   867 &   242 & 3.98e-05 & 1.57e-04 & 3.0 & 3.0 \\ 
		& 3.04e-02 &  1330 &   624 & 1.08e-05 & 4.26e-05 & 3.1 & 3.0 \\ 
		& 2.33e-02 &  1739 &  1016 & 4.87e-06 & 1.91e-05 & 3.0 & 3.0 \\ 
		\cline{2-6}
		\hline
		\hline
		\multirow{4}{*}{\rotatebox{90}{{DG $\mathbb{P}_3$}}}
		& 1.49e-01 &   446 &    18 & 1.18e-04 & 5.10e-04 & - & - \\     
		& 8.52e-02 &   794 &    60 & 1.17e-05 & 5.47e-05 & 4.1 & 4.0 \\ 
		& 4.66e-02 &  1469 &   230 & 1.03e-06 & 4.87e-06 & 4.0 & 4.0 \\ 
		& 3.04e-02 &  2257 &   562 & 1.79e-07 & 8.74e-07 & 4.1 & 4.0 \\ 
		\cline{2-6}
		\hline
		\hline
		\multirow{4}{*}{\rotatebox{90}{{DG $\mathbb{P}_4$}}}
		& 1.49e-01 &   644 &    18 & 1.21e-05 & 5.43e-05 & - & - \\     
		& 1.11e-01 &   881 &    24 & 2.68e-06 & 1.19e-05 & 5.1 & 5.1 \\ 
		& 8.52e-02 &  1148 &    56 & 6.96e-07 & 3.55e-06 & 5.1 & 4.6 \\ 
		& 6.09e-02 &  1622 &   130 & 1.40e-07 & 6.85e-07 & 4.8 & 4.9 \\ 
		\cline{2-6}
		\hline
	\end{tabular}
	\caption{   
		Convergence results for the Shu-type stationary vortex solution of the Euler equations 
		(with $\epsilon=5$), solved with our well-balanced ALE DG scheme set to preserve an equilibrium
		solution (with $\epsilon_E=0.5$) which differs from the imposed initial condition. We show the $L_2$
		error norms for $\rho$ and $p$ and the corresponding order of accuracy at time $t=1$. We remark that
		these results have been obtained after hundreds of sliver elements have been originated 
		thus showing that the order of convergence is maintained also in presence of large mesh deformations.}
	\label{table.order_euler_shu}
\end{table}

We open our set of benchmarks with a translating smooth isentropic vortex, as the one proposed in~\cite{HuShuVortex1999}, 
which represents a moving vortex solution of the Euler equations~\eqref{eq.eulerTerms}.

As computational domain, we take the square $\Omega(t)=[0+t,10+t]\times[0+t,10+t]$ with periodic boundary conditions, 
and we cover it with an increasingly refined set of polygonal tessellations, whose averaged mesh size is
 denoted by $\overline{h}$ and which moves together with the fluid flow.
The initial conditions and the equilibrium profile are given in terms of the primitive variables 
\begin{equation}
	\label{eq.ic_shu}
	\mathbf{V} = 
	\left \{
	\begin{aligned} 
	\rho  & = 1 + (1+\delta T)^{\frac{1}{\gamma-1}}-1, \\
	 u & = u_t - (y-5-t) \frac{\epsilon}{2\pi}e^{\frac{1-r^2}{2}},  \\
     v & = v_t + (x-5-t) \frac{\epsilon}{2\pi}e^{\frac{1-r^2}{2}},   \\
	 p & =  (1+\delta T)^{\frac{\gamma}{\gamma-1}}, \\
	\end{aligned} 
	\right .
\end{equation}
with the temperature fluctuation $\delta T = -\frac{(\gamma-1)\epsilon^2}{8\gamma\pi^2}e^{1-r^2}$,
$r=\sqrt{(x-5-t)^ 2+ (y-5-t)^2)}$ and the translation velocity $u_t=v_t=1$.
The vortex strength is set to $\epsilon=5$ for what concerns the initial condition $\Q_{\text{IC}}$ and 
to $\epsilon=\epsilon_E$ when setting the equilibrium profile to be preserved {$\Q_E$}.

We report in Table~\ref{table.order_euler_shu} the $L_2$ error norms of the numerical results $\u_h$, 
obtained for the density and pressure profile with our WB ALE DG scheme at time $t=1$, w.r.t. the 
imposed initial stationary conditions $\Q_{\text{IC}}$. The final time has been chosen in such a way that
during the simulation, due to the vortical rotation and in order to always maintain a good mesh
quality, hundreds of sliver elements originate; in this way, we can highlight that the expected
theoretical convergence order is always obtained, even in presence of large mesh 
deformations that trigger the presence of complex sliver elements in 
our WB ALE DG scheme.

\subsubsection{Moving vortical solution of the MHD equations}
\label{ssec.test_order_mhd}

\begin{table}[!tp] 
	\centering
	\numerikNine
	\begin{tabular}{|l|ccc|ccc|ccc|}
		\hline
		\multicolumn{10}{|c|}{{Moving vortical solution of the MHD equations with $\epsilon=5$ and $\epsilon_E=1$}}\\
		\hline
		\hline
		& $\overline{h}(\Omega(t_f))$  & timestep  &  sliver & $L_2(\rho)$ error & $L_2(p)$ error & $L_2(B_y)$ error &  order $\rho$ & order $p$ & order $B_y$ \\
		\hline	
		\hline
		\multirow{4}{*}{\rotatebox{90}{{DG-$\mathbb{P}_1$}}}
		& 8.52e-02 &   612 &    60 & 6.13e-02 & 1.08e-01 & 1.58e-02 & - & - & -\\     
		& 6.09e-02 &   863 &   142 & 3.84e-02 &  5.47e-02 & 1.03e-02 & 1.4 & 2.0 & 1.3 \\ 
		& 4.67e-02 &  1130 &   234 & 1.38e-02 &  2.96e-02 & 5.78e-03 & 3.8 & 2.3 & 2.0 \\ 
		& 3.04e-02 &  1735 &   584 & 3.55e-03 &  9.69e-03 & 2.84e-03 & 3.2 & 2.6 & 1.4 \\ 	
		\cline{2-6}
		\hline
		\hline
		\multirow{4}{*}{\rotatebox{90}{{DG-$\mathbb{P}_2$}}}
		& 1.11e-01 &   907 &    24 & 1.44e-03 & 4.04e-03 & 4.67e-04 & - & - & -\\     
		& 8.52e-02 &  1182 &    56 & 9.80e-04 &  1.93e-03 & 2.50e-04 & 1.5 & 2.8 & 2.4 \\ 
		& 6.09e-02 &  1669 &   130 & 5.77e-04 &  7.24e-04 & 9.43e-05 & 1.6 & 2.9 & 2.9 \\ 
		& 4.66e-02 &  2188 &   222 & 3.52e-04 &  3.24e-04 & 4.53e-05 & 1.8 & 3.0 & 2.6 \\ 	
		\cline{2-6}
		\hline
		\hline
		\multirow{4}{*}{\rotatebox{90}{{DG-$\mathbb{P}_3$}}}
		& 1.49e-01 &  1123 &    12 & 4.66e-04 & 8.65e-04 & 5.81e-04 & - & - & - \\     
		& 1.11e-01 &  1537 &    18 & 1.46e-04 &  2.53e-04 & 2.46e-04 & 3.9 & 4.1 & 2.9 \\ 
		& 8.52e-02 &  2005 &    54 & 5.68e-05 &  9.05e-05 & 1.02e-04 & 3.6 & 3.9 & 3.4 \\ 
		& 6.10e-02 &  2833 &   102 & 1.61e-05 &  2.41e-05 & 2.91e-05 & 3.8 & 4.0 & 3.7 \\ 	
		\cline{2-6}
		\hline
		\hline
		\multirow{4}{*}{\rotatebox{90}{{DG-$\mathbb{P}_4$}}}
		& 1.49e-01 &  1624 &     2 & 7.13e-05 & 7.72e-05 & 3.89e-05 & - & - & - \\     
		& 1.11e-01 &  2224 &    18 & 2.03e-05 &  1.89e-05 & 1.21e-05 & 4.2 & 4.7 & 3.9 \\ 
		& 8.51e-02 &  2902 &    46 & 6.74e-06 &  5.53e-06 & 4.31e-06 & 4.2 & 4.7 & 4.0 \\ 
		& 6.10e-02 &  4102 &    86 & 1.65e-06 &  1.10e-06 & 1.05e-06 & 4.2 & 4.8 & 4.2 \\ 
		\cline{2-6}
		\hline
	\end{tabular}
	\caption{
		Convergence results for the MHD vortex (with $\epsilon=5$), solved with our well-balanced 
		ALE DG scheme set to preserve an equilibrium solution (with $\epsilon_E=1$) which differs from
		the imposed initial condition. We show the $L_2$ error norms for $\rho$, $p$ and $B_y$ and the
		corresponding order of accuracy at time $t=1$. We remark that these results have been obtained after
		hundreds of sliver elements have been originated thus showing that the order of convergence is 
		maintained also in presence of large mesh deformations.	
	}
	\label{table.order_mhd_vortex}
\end{table}

We repeat now the convergence study on a moving vortical solution of the MHD equations taken from~\cite{Balsara2004}.  
The computational domain is $\Omega=[0,10]\times[0,10]$ with wall boundary conditions imposed everywhere. 
The initial condition is given in terms of the vector of primitive variables as 
\begin{equation}
	\mathbf{V}_{\text{IC}} =  ( \rho, u, v, w, p, B_x, B_y,  B_z, \Psi )^T = 
	( 1, \delta u, \delta v, 0, 1+\delta p, \delta B_x, \delta B_y,  0, 0 )^T, 
\end{equation}
with $\delta \mathbf{v} = (\delta u, \delta v, 0)^T$, $\ \delta \mathbf{B} = ( \delta B_x, \delta B_y, 0 )^T$,
$\mathbf{e}_z = (0,0,1)$, $\mathbf{r} = (x-5,y-5,0)$, $r = | \mathbf{r} | $,  $\mu=\sqrt{4 \pi}$ and 
\begin{equation}
	\left \{ 
	\begin{aligned} 
		\label{eqn.mhd3d.ic1}
		&\delta \mathbf{v} = \frac{\epsilon}{2\pi} e^{ \halb(1-r^2)} \mathbf{e}_z \times \mathbf{r}  \\ 
		&\delta \mathbf{B} = \frac{\mu}{2\pi}    e^{ \halb(1-r^2)} \mathbf{e}_z \times \mathbf{r},  \\ 
		&\delta p = \frac{1}{32 \pi^3} \left( \mu^2 (1 - r^2) - 4 \epsilon^2 \pi \right) e^{(1-r^2)}.
	\end{aligned} 
	\right .
\end{equation}
The divergence cleaning speed is chosen as $c_h=3$.
The vortex strength is set to $\epsilon=5$ for what concerns the initial condition $\Q_{\text{IC}}$ and 
to $\epsilon_E=1$ when setting the equilibrium profile to be preserved {$\Q_E$}.

The obtained convergence results are reported in Table~\ref{table.order_mhd_vortex}.

\subsection{Non equilibrium benchmarks}
\label{ssec.nonequilibrium}

In this section we verify that our well-balanced scheme is able to produce the correct solution 
also in situations far from the equilibrium, independently from the choice of initial conditions.
Hence, 
our method does not fall into the class of perturbations methods making it more general widely applicable. 
In addition, we take the occasion to highlight the positive effects of the Lagrangian framework 
when heavily convection dominated phenomena are considered.

\subsubsection{Travelling Sod shock tube problem}
\label{ssec.Sod}
\begin{figure}[!bp]
	\includegraphics[width=0.25\linewidth, trim=8 8 8 8, clip]{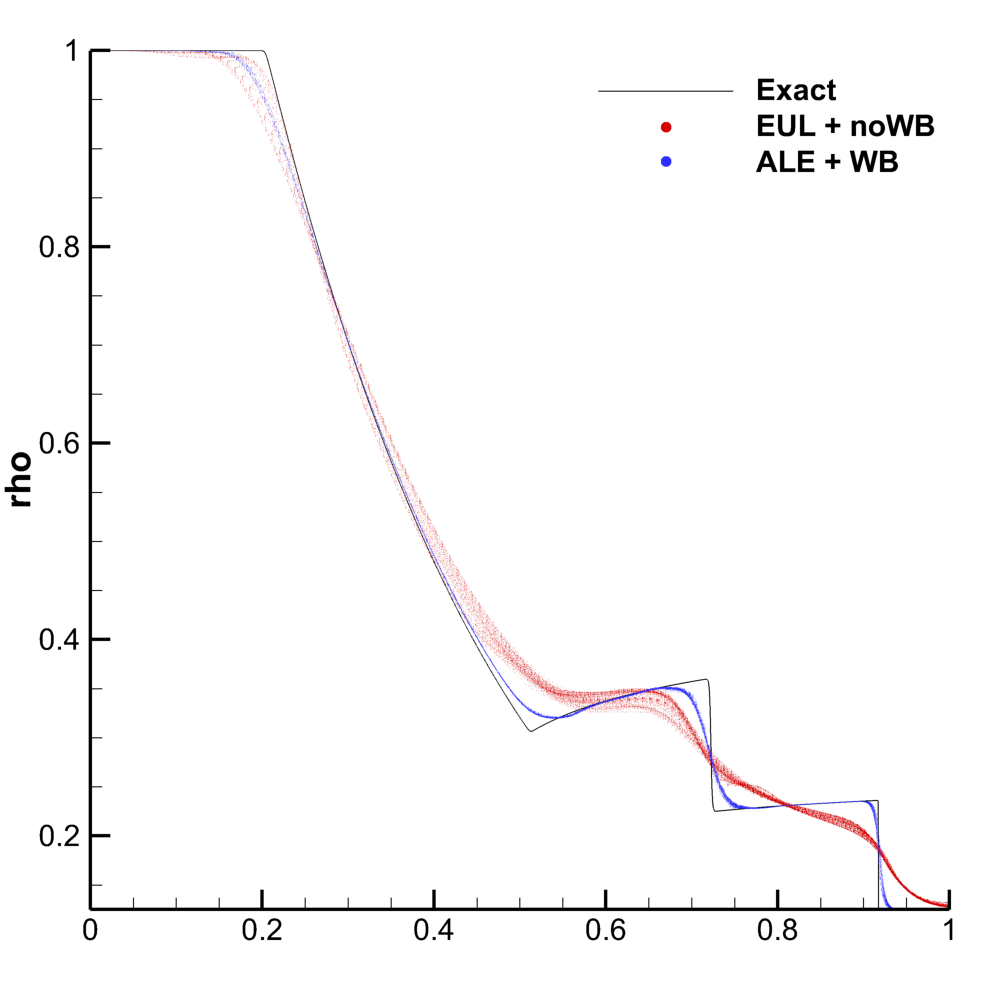}\quad\quad\ %
	\includegraphics[width=0.35\linewidth, trim=8 8 8 8, clip]{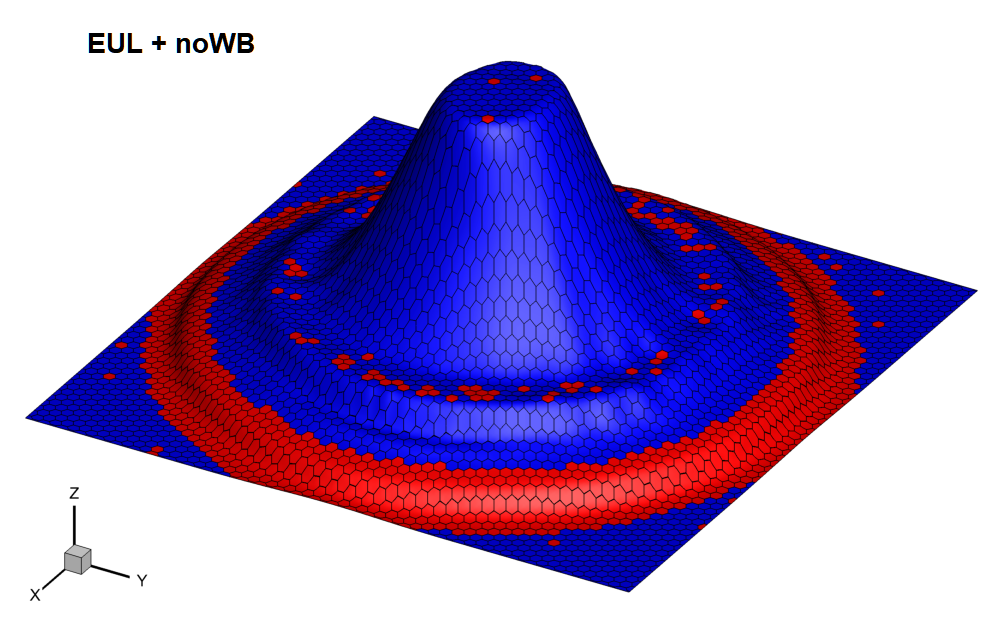}%
	\includegraphics[width=0.35\linewidth, trim=8 8 8 8, clip]{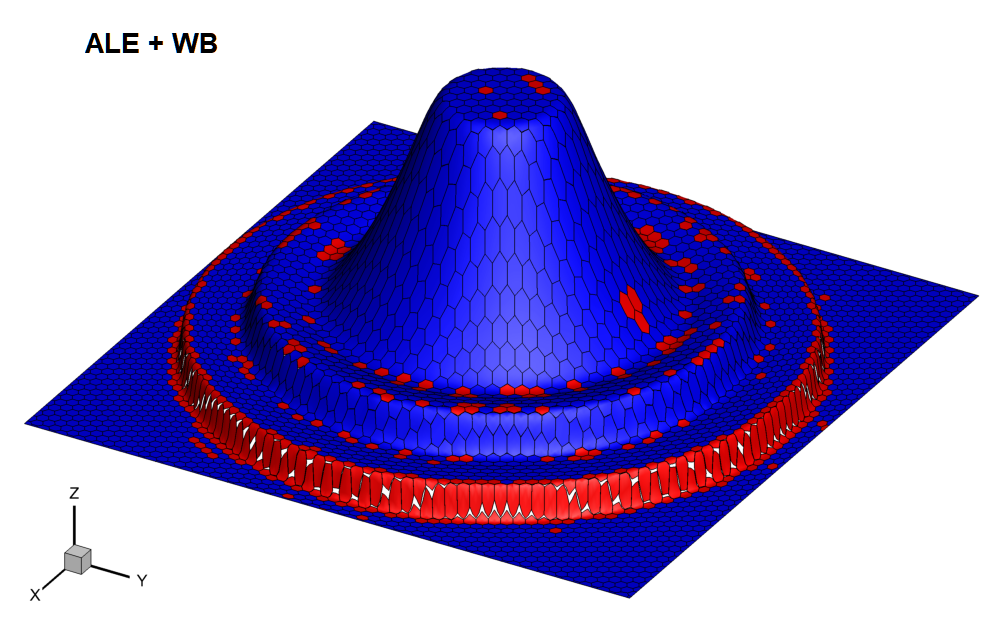}%
	\caption{Traveling Sod-type explosion problem solved with our DG scheme of order $3$ in 
	both its non well-balanced Eulerian version (middle) and the well-balanced ALE version (right). 
	In particular, we compare the density $\rho$ profile at time $t=0.25$ (left) and we show the 
	limiter activation (middle and right). 
		We can note that: \textit{i)} the use of our Lagrangian techniques allows to sharply capture the 
		discontinuities even on a fast moving background and \textit{ii)} the use of the well-balancing 
		does not affect the accuracy and robustness of the overall scheme.
	}
	\label{fig.Sodtrasl}
\end{figure}

Here, we solve a moving multidimensional explosion problem constructed as an extension of the classical Sod test case.
We set as equilibrium $\Q_E$ the isentropic vortex~\eqref{eq.ic_shu} with $u_t=v_t=0$.
As computational domain we take a square of dimension $\Omega(t) = [-1+40t;1+40t]\times[-1;1]$ covered 
with a mesh counting $4105$ polygonal elements, and the initial condition is composed of two 
different states, separated by a discontinuity delimited by a circle radius $r_d=0.5$ centered about the origin, or formally
\begin{equation}
	\mathbf{V}_{\text{IC}} = 
	\begin{cases} 
		\rho = 1, \hspace{0.98cm} u=40, \quad v= 0, \quad p = 1, \quad  &  r \leq r_d, \\
		\rho = 0.125, \quad u=40, \quad v=0, \quad  p = 0.1, \quad & r > r_d.
	\end{cases} 
\end{equation}
Note, in particular, the fast background speed imposed in the $x$-direction chosen such that at 
the final simulation time $t_f=0.25$ the initial square $\Omega(t=0)$ 
will have been displaced by $5$ times its initial size.

We report the results obtained with a standard Eulerian DG scheme and with our well-balanced 
ALE DG approach, both of order $3$, in Figure~\ref{fig.Sodtrasl}. We can notice that 
\textit{i)} the use of WB does not interfere with the resolution of the scheme, and specifically the rarefaction wave, 
the contact wave and the shock discontinuity are correctly captured and \textit{ii)} that the ALE algorithm 
allows to sharply reduce the errors due to the convection which adds a significant amount of 
artificial diffusion when the Eulerian scheme is used instead. 
We finally remark that analogous observations could be made in~\cite{gaburro2021bookchapter}, 
where the same problem was 
solved without the use of WB techniques, just exploiting the Lagrangian features of the scheme, and 
that the novelties here presented in the context of well-balancing do not negatively affect other 
aspects of the method.

\subsubsection{MHD rotor problem}
\label{ssec.MHDrotor}

\begin{figure}[!bp]
	\includegraphics[width=0.333\linewidth, trim=8 8 8 8, clip]{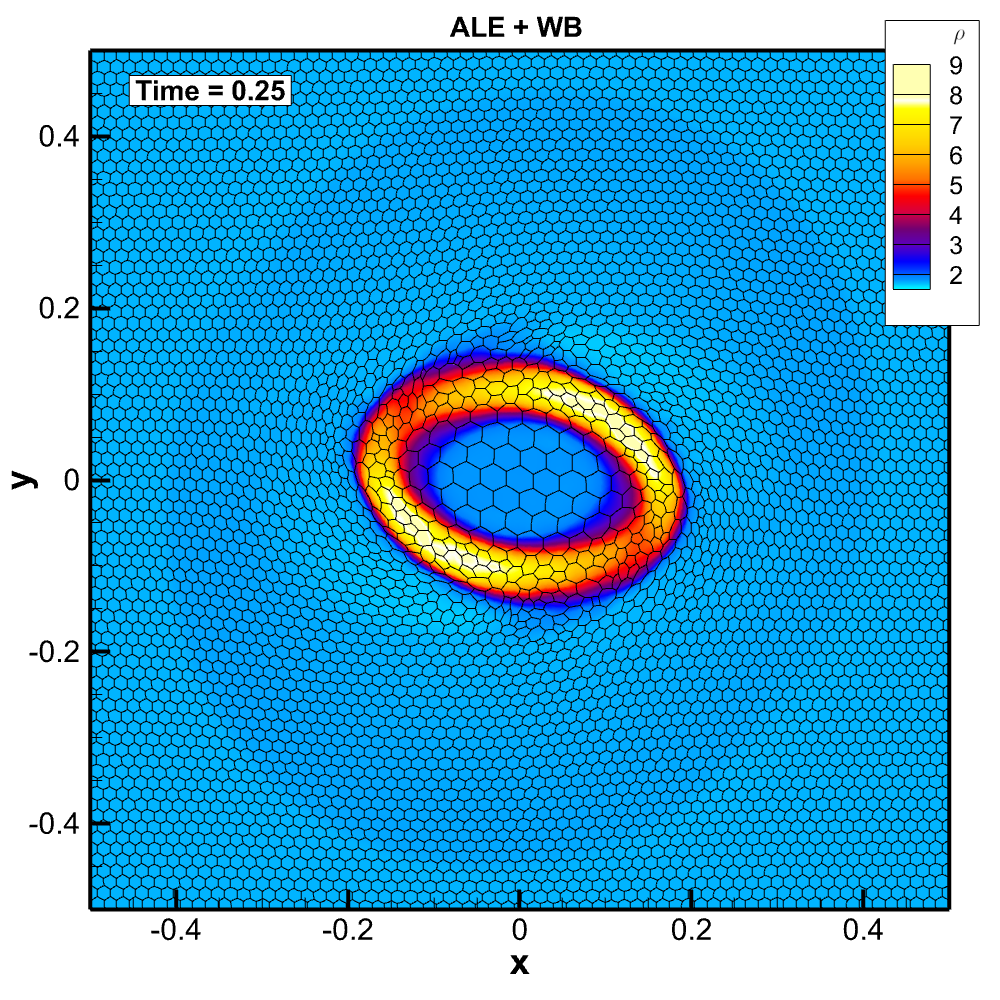}%
	\includegraphics[width=0.333\linewidth, trim=8 8 8 8, clip]{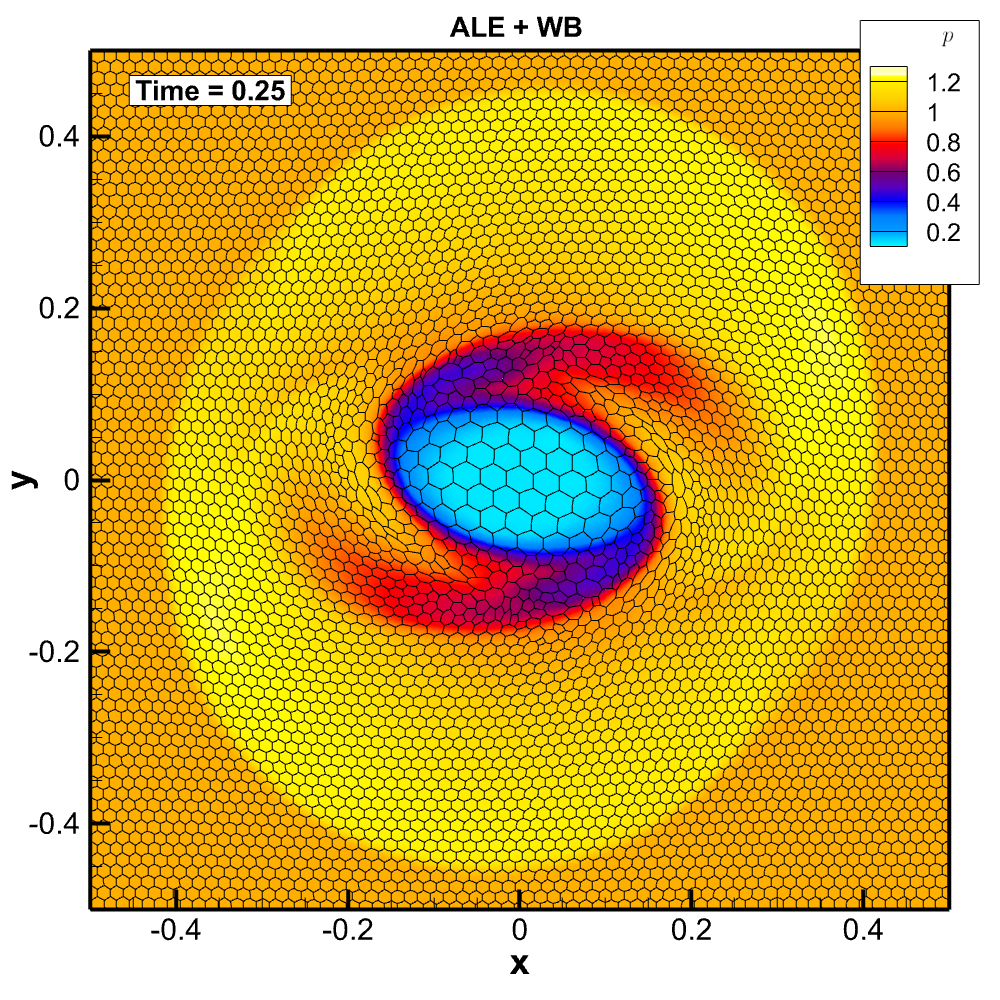}%
	\includegraphics[width=0.333\linewidth, trim=8 8 8 8, clip]{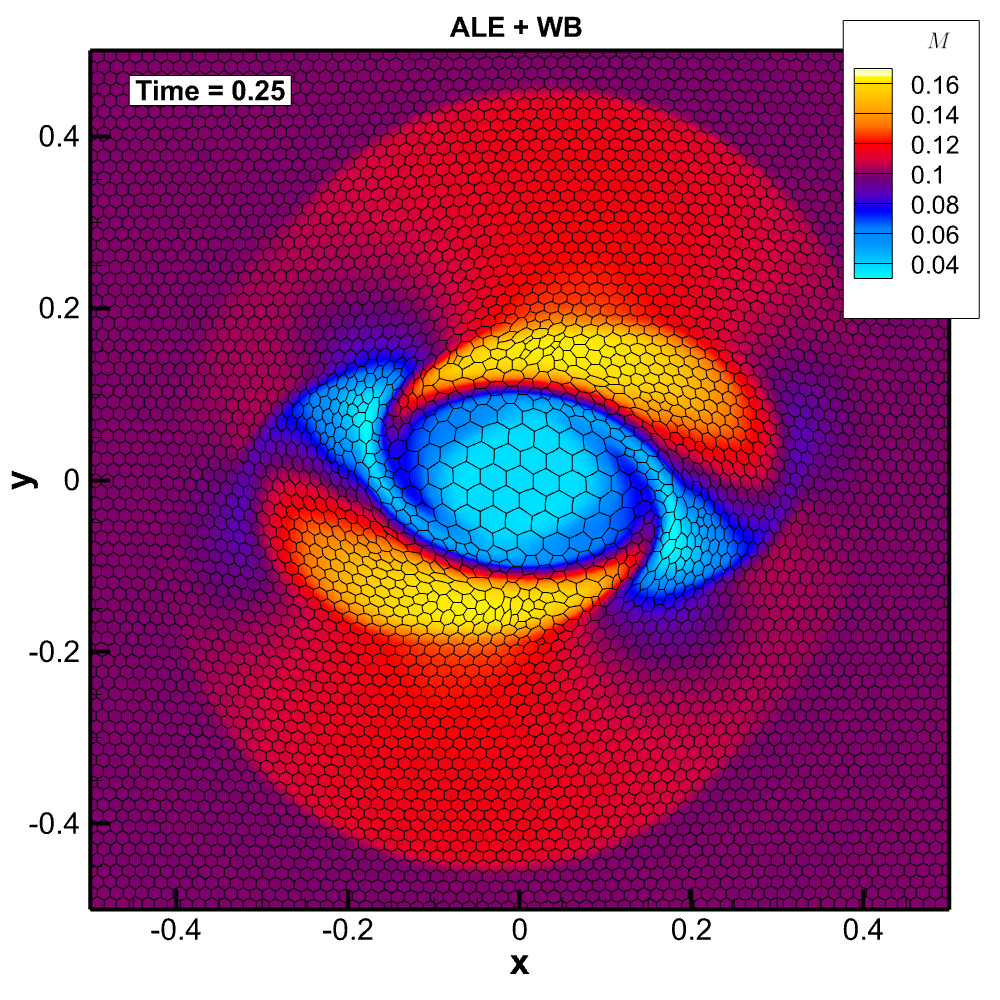}\\
	\includegraphics[width=0.333\linewidth, trim=8 8 8 8, clip]{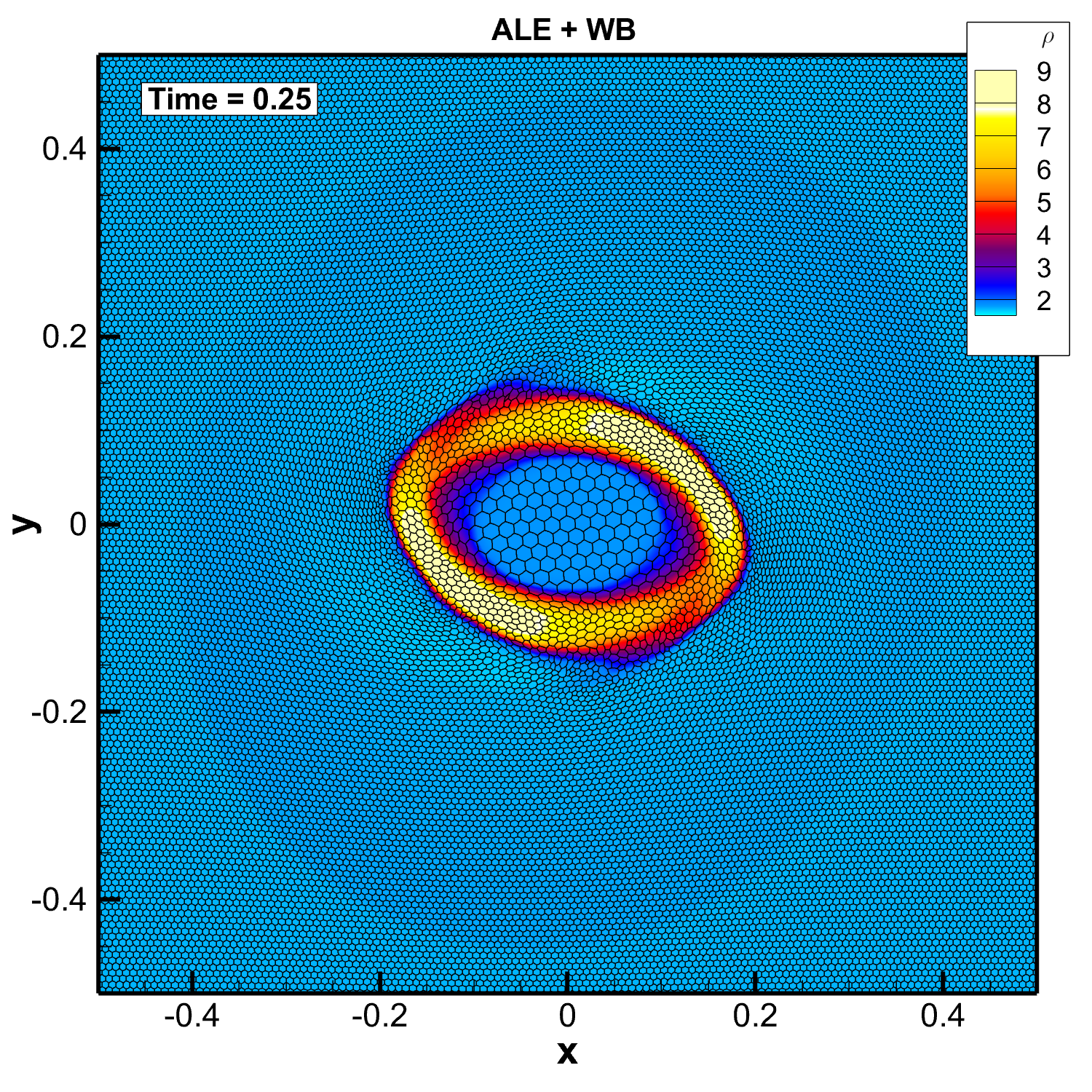}%
	\includegraphics[width=0.333\linewidth, trim=8 8 8 8, clip]{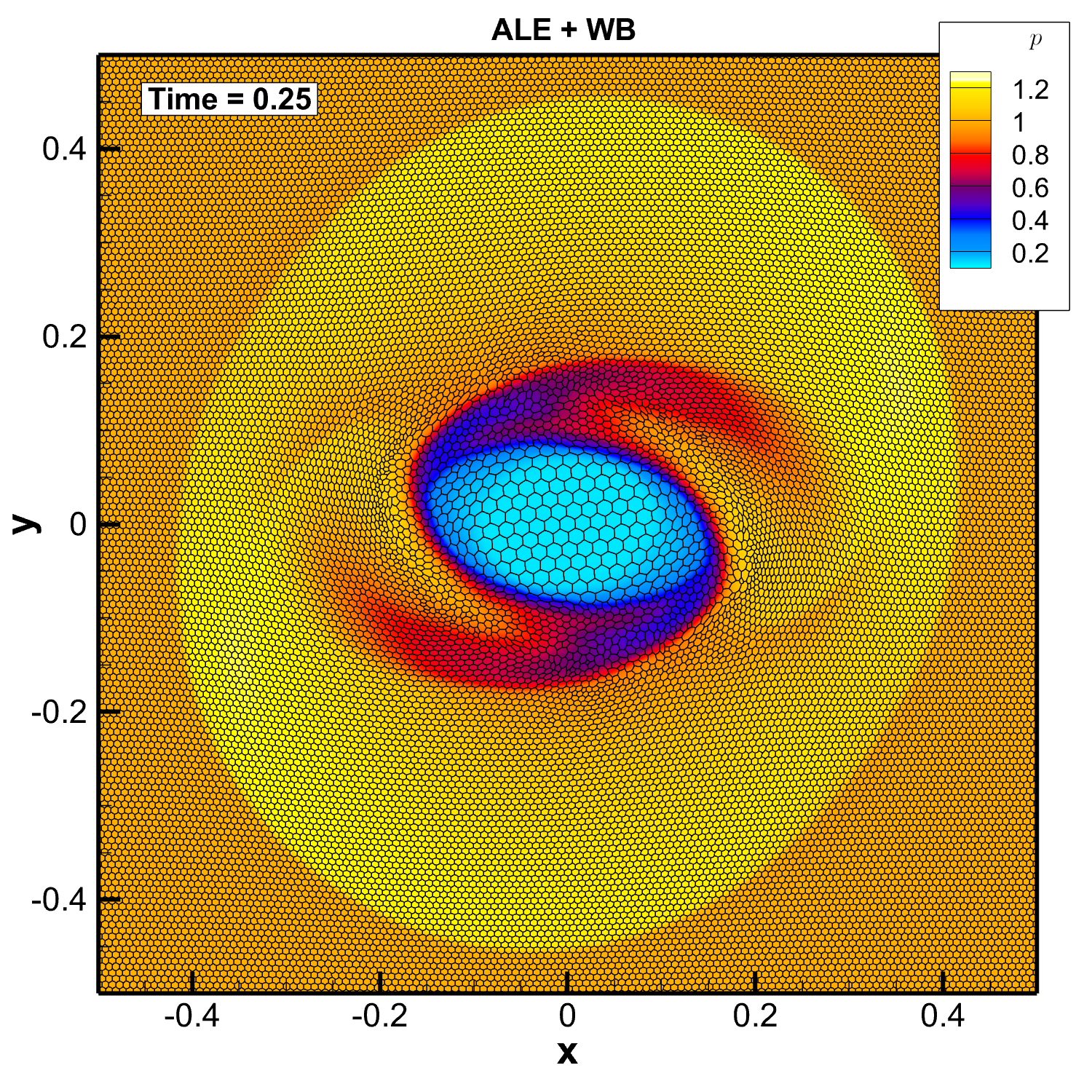}%
	\includegraphics[width=0.333\linewidth, trim=8 8 8 8, clip]{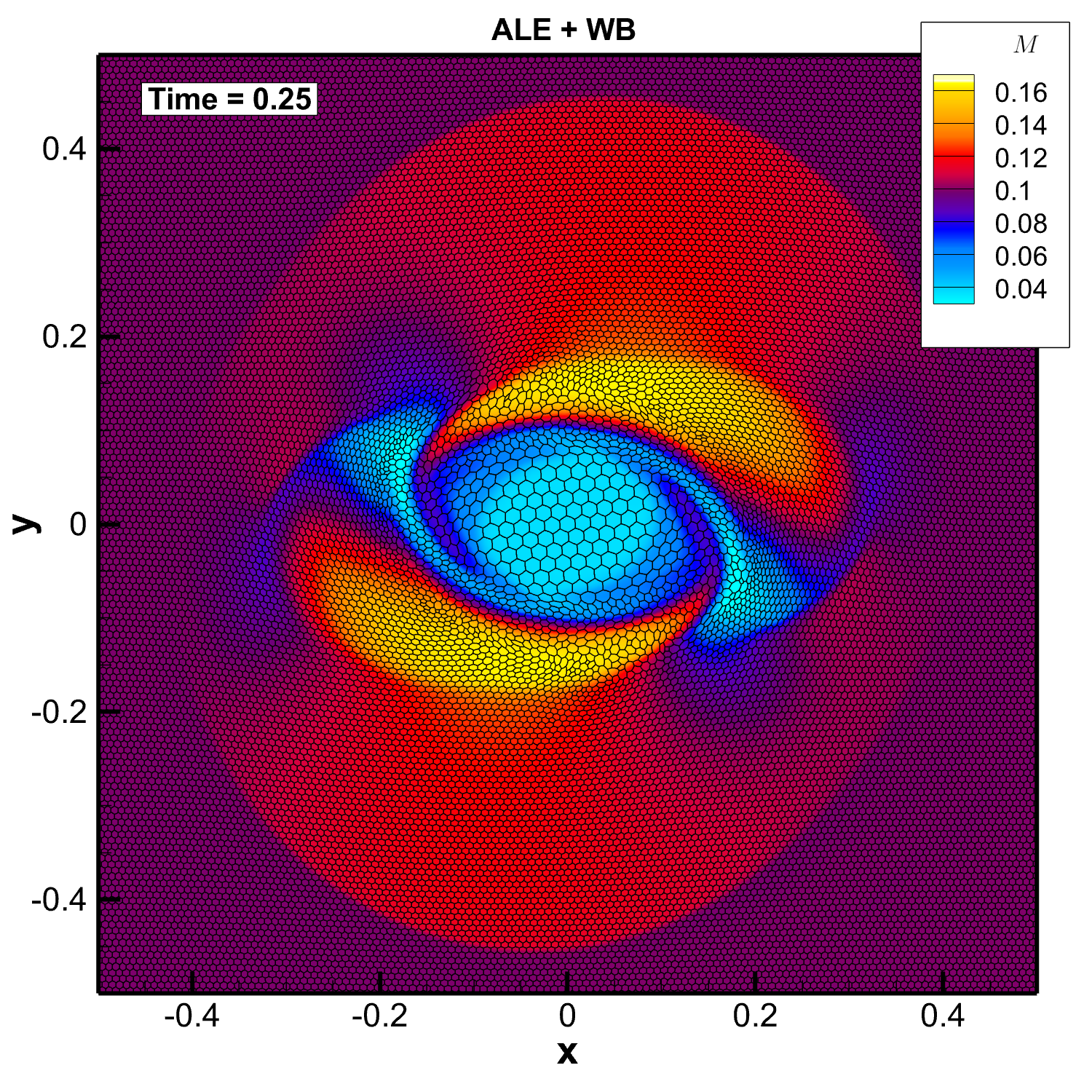}\\
	\includegraphics[width=0.333\linewidth, trim=8 8 8 8, clip]{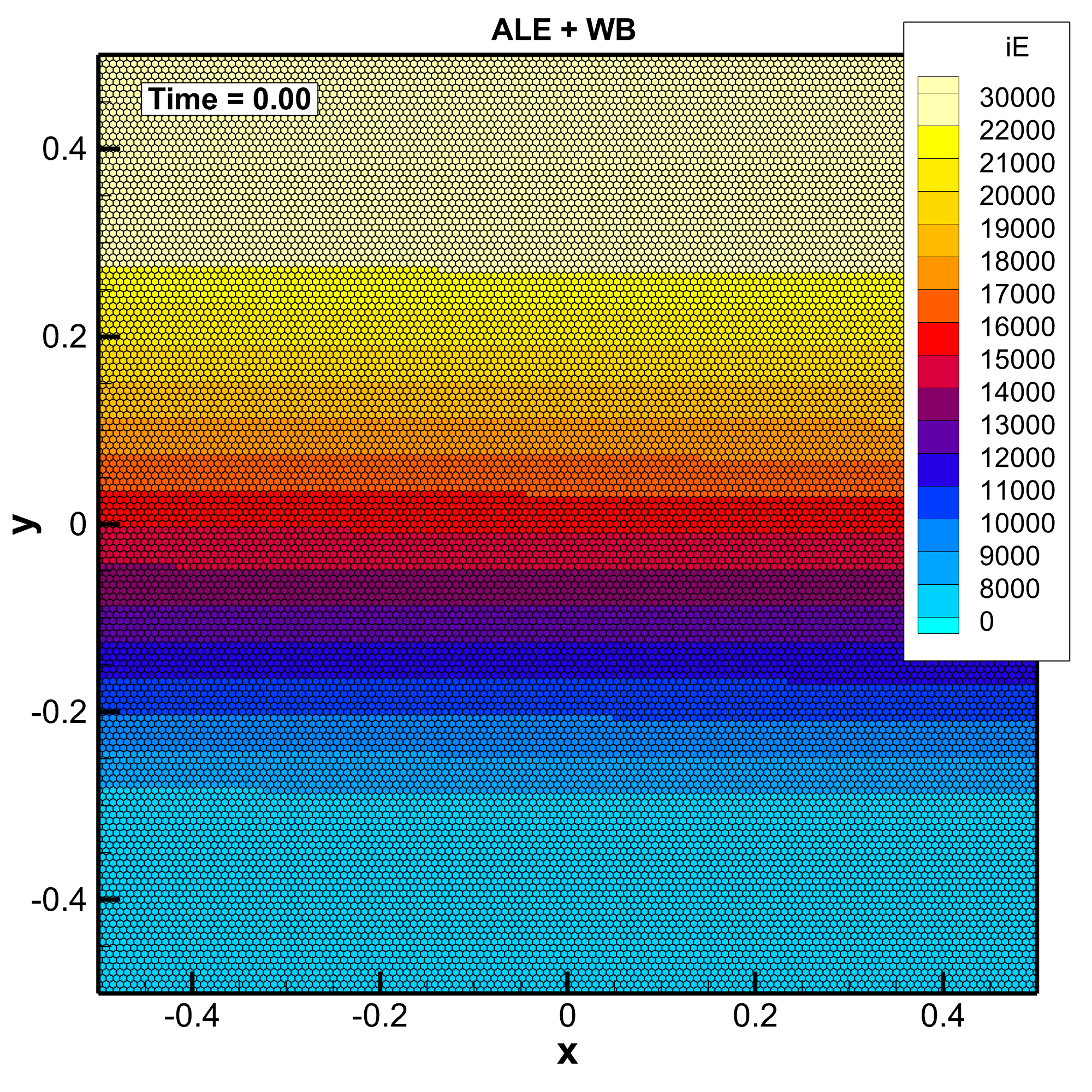}%
	\includegraphics[width=0.333\linewidth, trim=8 8 8 8, clip]{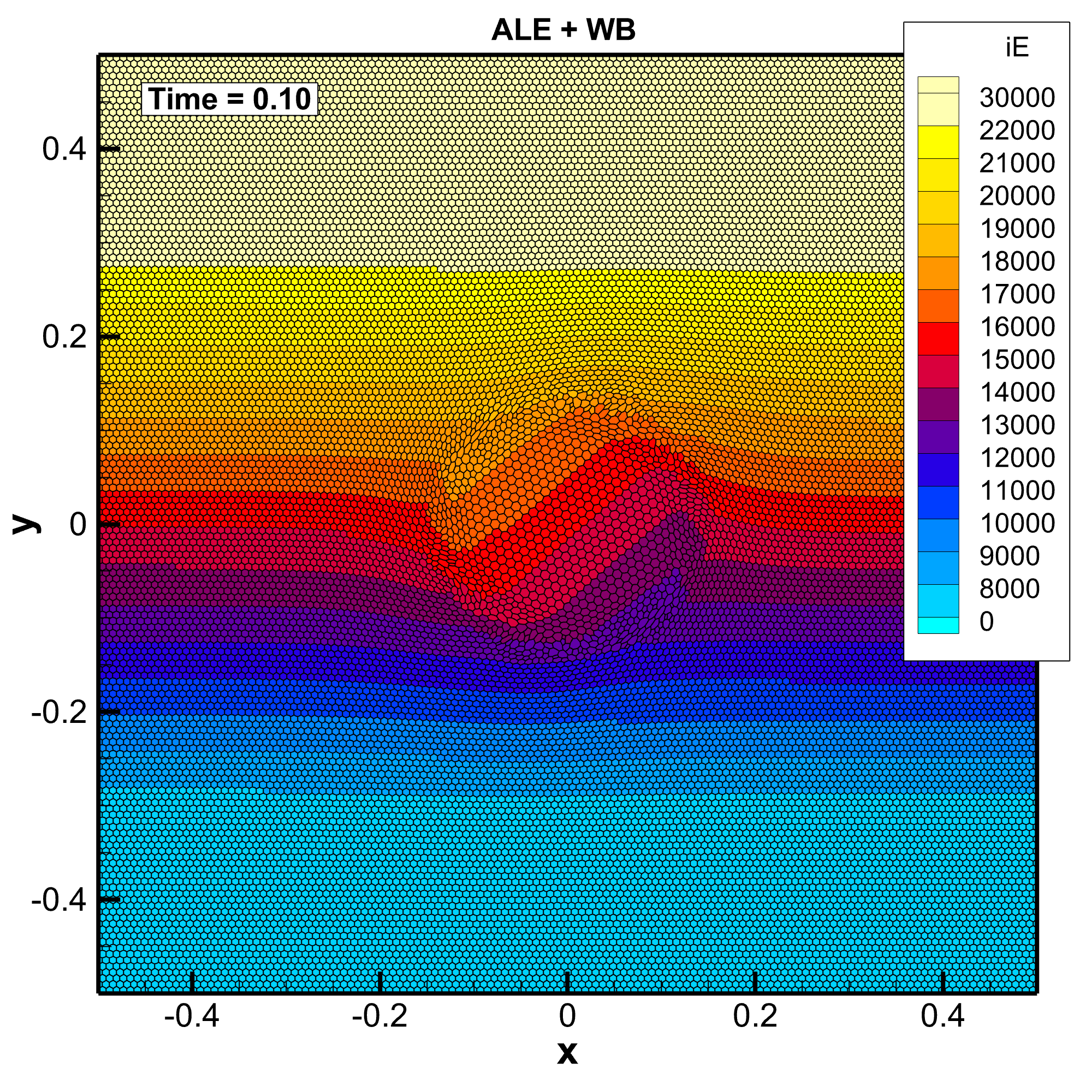}%
	\includegraphics[width=0.333\linewidth, trim=8 8 8 8, clip]{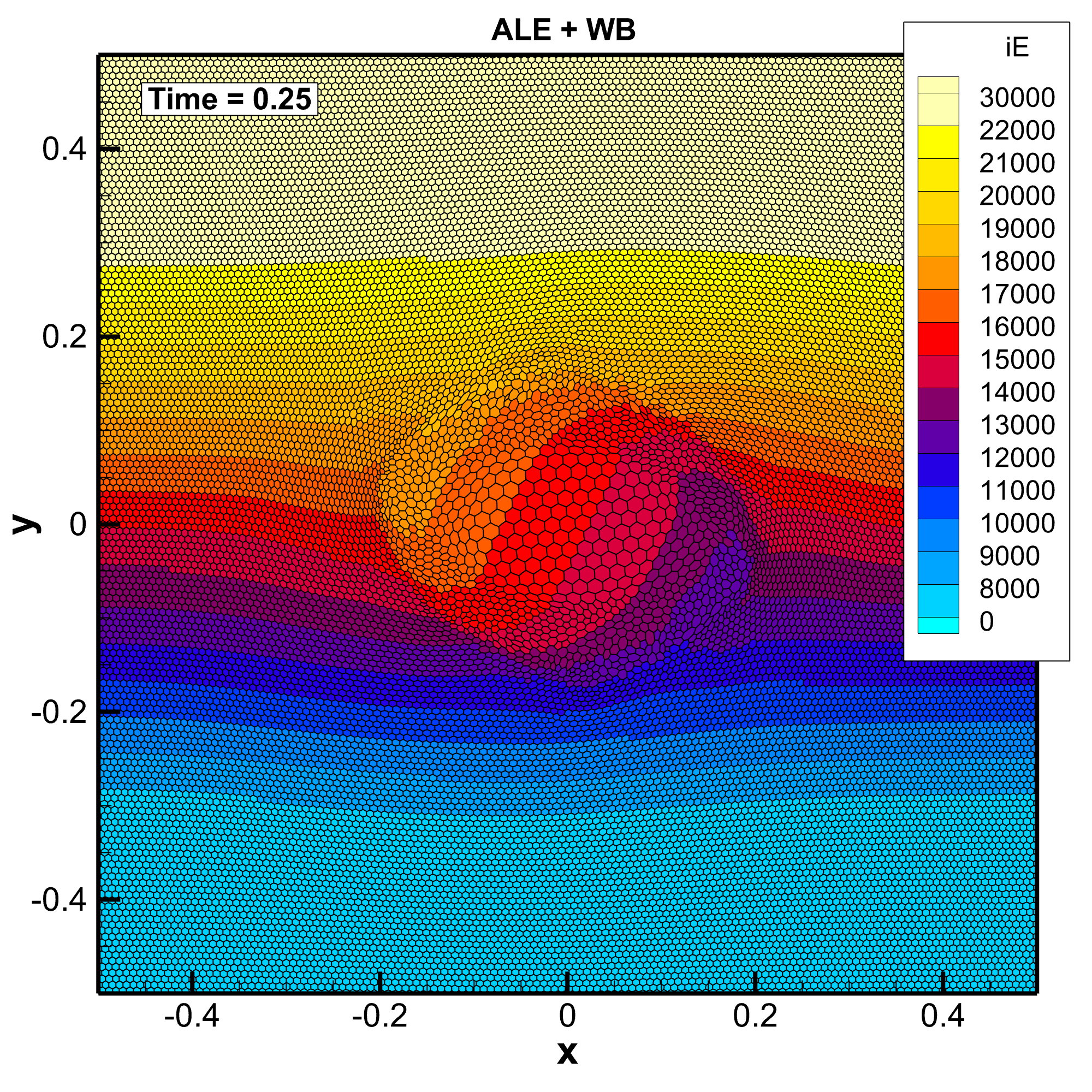}\\[-10pt]
	\caption{MHD rotor problem solved with our WB ALE DG scheme of order $3$ over a coarse 
	mesh of $7579$ polygonal elements (first line) and a finer mesh of $30387$ polygonal elements (second and third line). 
	We report the density $\rho$, pressure $p$ and the magnetic density 
	profile $M = (B_x^2 + B_y^2+B_z^2)/(8\pi)$ from which we can observe a clear mesh convergence and thus 
	verify that our scheme perfectly works also in non-equilibrium situations. 
	We also show, in the third line, the evolution of the mesh numbering to highlight the Lagrangian 
	mesh movement. 
	}
	\label{fig.MHDrotor}
\end{figure}

We also consider, as non-equilibrium test case, the MHD rotor problem proposed 
by Balsara and Spicer in~\cite{BalsaraSpicer1999}. The computational domain is taken to 
be $\Omega=[-0.55, 0.55]\times[-0.55, 0.55]$ with wall boundary conditions. 
The benchmark consists of a rapidly rotating high density fluid ($\rho=10$ for $r<0.1$) embedded in a low density one ($\rho=1$). 
Both fluids are subject to an initially constant magnetic field  $\B = (2.5, 0, 0)^T$ and the initial pressure is $p=1$ in the entire domain.  
The rotor has an angular velocity of $\omega=10$ for $r<0.1$, 
which produces torsional Alfv\'en waves that are launched into the outer fluid at rest, 
resulting in a decrease of angular momentum of the spinning rotor. 
As proposed in~\cite{BalsaraSpicer1999}, we apply a linear taper to velocity and density in the range from 
$0.1 \leq r \leq 0.12$ so that they match the ambient fluid at a radius of $r=0.12$; next, we take $c_h=2$. 

We report the results obtained with our WB ALE DG scheme of order $3$ in Figure~\ref{fig.MHDrotor} 
where we visually highlight both the convergence of our scheme and the Lagrangian motion of the mesh.

\subsection{Equilibrium solutions and their perturbations}
\label{ssec.test_wb}

In this section we show the key feature of our novel scheme, 
i.e. its capability of preserving equilibrium solutions with machine precision even in complex moving situations, 
and its superiority, w.r.t. classical non well-balanced algorithms, 
in the simulations of small perturbations around these equilibria.

In particular, we would like to remark that Lagrangian discontinuous Galerkin schemes with 
\textit{a posteriori} FV limiter, already represent one of the most accurate numerical schemes 
available for hyperbolic equations, especially on unstructured meshes, which only fail when the 
details to be simulated are of the order of the numerical errors associated with the scheme 
(which depend on the mesh size and the selected order of accuracy).
The use of well-balanced techniques allow to further enlarge the area of applicability of a 
selected scheme, without the need of excessively refining the mesh or having to increase its 
order of accuracy, making it possible to obtain reliable results even in situations that otherwise would be 
intractable due to the very high computational costs.


\subsubsection{Stationary solutions of the Euler equations: isentropic Shu vortex and Gresho vortex}
\label{ssec.test_wb_eulervortexes}
\begin{table}[!tp] 
	\centering
	\numerikNine
	\begin{tabular}{|l|l||cccccc|}
		\hline
		\multicolumn{8}{|c|}{Isentropic Shu vortex with $\epsilon = 5$} \\
		\hline
		\hline
		& Time $t = $ & $0.0$ & $0.5$ & $1.0$ & $2.5$ & $5.0$ & $10.0$ \\
		\hline
		\hline
		\multirow{4}{*}{\rotatebox{90}{{DG-$\mathbb{P}_1$}}}
		& $L_2(\rho)$ error & 5.3166E-12 & 5.2877E-12	& 5.2655E-12	& 5.2015E-12	& 4.9641E-12	& 4.9723E-12 \\
		& $L_2(p)$ error  	& 0.0        & 5.6353E-12	& 5.7406E-12	& 5.6696E-12	& 4.8272E-12	& 4.9928E-12 \\
		& No. timestep 	    & 0     	 & 116       	& 221       	& 531       	& 1043	        & 2063       \\
		& No. sliver   	    & 0     	 & 0       	    & 13        	& 69        	& 160	        & 369       \\
		\cline{2-8}
		\hline
		\hline
		\multirow{4}{*}{\rotatebox{90}{{DG-$\mathbb{P}_2$}}}
		& $L_2(\rho)$ error & 5.2246E-12 & 5.2074E-12	& 5.1984E-12	& 5.2064E-12	& 5.0257E-12	& 5.0395E-12 \\
		& $L_2(p)$ error  	& 0.0        & 5.5197E-12	& 5.6551E-12	& 5.9331E-12	& 5.2350E-12	& 5.4774E-12 \\
		& No. timestep 	    & 0     	 & 217       	& 419       	& 1017       	& 2008	        & 3985       \\
		& No. sliver   	    & 0     	 & 0       	    & 4        	    & 56       	    & 141	        & 320       \\
		\cline{2-8}
		\hline
		\hline
		\multirow{4}{*}{\rotatebox{90}{{DG-$\mathbb{P}_3$}}}
		& $L_2(\rho)$ error & 5.1698E-12 & 5.1760E-12	& 5.1915E-12	& 5.2039E-12	& 5.2153E-12	& 5.2566E-12 \\
		& $L_2(p)$ error  	& 0.0        & 5.4693E-12	& 5.4849E-12	& 5.4933E-12	& 5.4448E-12	& 5.3596E-11 \\
		& No. timestep 	    & 0     	 & 360       	& 701       	& 1716       	& 3397	        & 6759      \\
		& No. sliver   	    & 0     	 & 0         	& 1        	    & 41        	& 110	        & 259       \\
		\cline{2-8}
		\hline
		\hline
		\multirow{4}{*}{\rotatebox{90}{{DG-$\mathbb{P}_4$}}}
		& $L_2(\rho)$ error & 5.1122E-12 & 5.1284E-12	& 5.1565E-12	& 5.2007E-12	& 5.2894E-12	& 5.3559E-12 \\
		& $L_2(p)$ error  	& 0.0        & 5.4239E-12	& 5.4440E-12	& 5.4853E-12	& 5.5523E-13	& 5.5692E-11 \\
		& No. timestep 	    & 0     	 & 517      	& 1010       	& 2478      	& 4912	        & 9787       \\
		& No. sliver   	    & 0     	 & 0       	    & 0          	& 33        	& 88	        & 232       \\
		\cline{2-8}
		\hline
	\end{tabular}\vspace{-5pt}
 \caption{Verification of the well-balanced property on the Shu vortex. In this table, and the following
ones, we report the number of performed iterations, the number of sliver elements originated during
the simulation and the $L_2$ norm of the difference between the numerical solution and the
exact equilibrium solution.  We can notice that the equilibrium solution, even if initially perturbed with a random error of
1E-12 distributed everywhere on the domain, is preserved with machine accuracy for very long
simulation times and after handling thousands of sliver elements. This holds true for any
employed polynomial representation order.} 	
\label{table.wb_shu}
\end{table}
\begin{table}[!tp] \vspace{-10pt}
	\centering
	\numerikNine
	\begin{tabular}{|l|l||cccccc|}
		\hline
		\multicolumn{8}{|c|}{Gresho vortex} \\
		\hline
		\hline
		& Time $t = $ & $0.0$ & $0.5$ & $1.0$ & $2.5$ & $5.0$ & $10.0$ \\
		\hline
		\hline
		\multirow{4}{*}{\rotatebox{90}{{DG-$\mathbb{P}_1$}}}
		& $L_2(\rho)$ error & 1.0334E-12 & 1.0166E-12	& 1.0164E-12	& 1.0148E-12	& 1.0140E-12	& 1.0179E-12 \\
		& $L_2(p)$ error  	& 0.0        & 2.1747E-13	& 2.0289E-13	& 1.9675E-13	& 1.9399E-13	& 1.2650E-12 \\
		& No. timestep 	    & 0     	 & 672       	& 1338       	& 3335       	& 6663	        & 13318       \\
		& No. sliver   	    & 0     	 & 83       	& 181        	& 486       	& 1017	        & 2041       \\
		\cline{2-8}
		\hline
		\hline
		\multirow{4}{*}{\rotatebox{90}{{DG-$\mathbb{P}_2$}}}
		& $L_2(\rho)$ error & 1.0290E-12 & 1.0179E-12	& 1.0188E-12	& 1.0177E-12	& 1.0175E-12	& 1.2205E-12 \\
		& $L_2(p)$ error  	& 0.0        & 1.9886E-13	& 1.8870E-13	& 1.9020E-13	& 2.1258E-13	& 1.0530E-11 \\
		& No. timestep 	    & 0     	 & 1298       	& 2590       	& 6466       	& 12926	        & 25845       \\
		& No. sliver   	    & 0     	 & 75       	& 155        	& 422       	& 869	        & 1751       \\
		\cline{2-8}
		\hline
		\hline
		\multirow{4}{*}{\rotatebox{90}{{DG-$\mathbb{P}_3$}}}
		& $L_2(\rho)$ error & 1.0282E-12 & 1.0378E-12	& 1.0471E-12	& 1.0248E-12	& 1.1667E-12	& 1.6059E-12 \\
		& $L_2(p)$ error  	& 0.0        & 2.2900E-13	& 4.2956E-13	& 9.1580E-13	& 8.2288E-13	& 1.2231E-12 \\
		& No. timestep 	    & 0     	 & 2202         & 4399       	& 10988       	& 21969	        & 43930      \\
		& No. sliver   	    & 0     	 & 64       	& 140        	& 375        	& 776	        & 1538        \\
		\cline{2-8}
		\hline
		\hline
		\multirow{4}{*}{\rotatebox{90}{{DG-$\mathbb{P}_4$}}}
		& $L_2(\rho)$ error & 1.0226E-12 & 1.0339E-12	& 1.0344E-12	& 1.0303E-12	& 1.0358E-12	& 1.055E-12 \\
		& $L_2(p)$ error  	& 0.0        & 1.0248E-13	& 1.0461E-13	& 1.0523E-13	& 1.1391E-13	& 1.1890E-13 \\
		& No. timestep 	    & 0     	 & 3189      	& 6372       	& 15921      	& 31835	        & 63664       \\
		& No. sliver   	    & 0     	 & 64       	& 123        	& 335        	& 690	        & 1360       \\
		\cline{2-8}
		\hline
	\end{tabular}\vspace{-5pt}
 \caption{Verification of the well-balanced property on the Gresho vortex. As in the previous test case, we
can notice that the equilibrium solution, even if initially perturbed is preserved with machine
accuracy for very long simulation times and after handling thousands of sliver elements.}
\label{table.wb_gresho}
\end{table}

We start by showing that our scheme is well-balanced, i.e. that when the initial condition of a
simulation is given by a machine precision perturbation of the prescribed equilibrium solution, the
scheme is able to preserve this equilibrium maintaining the numerical errors at the level of machine
precision for very long computational time and large mesh deformation. 

We first perform this kind of test case by choosing as equilibrium profile the Shu vortex described
in the previous Section~\ref{ssec.test_order_euler} with $\epsilon=\epsilon_E=5$ both to describe
$\Q_{\text{IC}}$ and $\Q_E$, and setting $u_t=v_t=0$. We then add a random perturbation of order 1E-12 to the
initial density profile everywhere on the domain. We discretize the domain with $516$ polygonal
elements and we move the mesh together with the fluid flow. %
The obtained numerical results are reported in Table~\ref{table.wb_shu}: we can notice that, even
after a very large number of iterations and handling a large number of sliver elements, the
equilibrium solution is perfectly preserved up to machine precision. 
(In order not to overload this section we do not report every numerical test done, 
but we specify that very similar results have been obtained with initial random perturbation of order 1E-13 or 1E-14, 
for which the numerical errors after long times are of the same order of the initial perturbation. 
This holds true for all the benchmarks of this section.)
\medskip 

Next, we perform the same type of test by choosing the Gresho vortex, 
which is another moving  solution of the Euler equations, see~\cite{liska2003comparison}.
Here, the density is $\rho=0$ as well as the radial velocity $u_r=0$.
The centrifugal force is balanced by the gradient of the pressure; the angular velocity and the pressure are given by 
\begin{equation}
	\label{eq.gresho} 
	(u_\phi(r), p(r)) = 	
	\begin{cases} 
	(5r, \ 5 + \frac{25}{2}r^2) \quad &\text{ if }  \quad  0 \le r < 0.2, \\
	(2-5r, \ 9 - 4\ln(0.2) + 12.5 r^2 - 20r + 4\ln(r) ) \quad &\text{ if }  \quad 0.2 \le r < 0.4, \\
	(0, \ 3 + 4\ln(2)) \quad &\text{ if }  \quad  r < 0.4.\\
	\end{cases} 
\end{equation}
The computational domain is $\Omega = [-1,1]\times[-1,1]$ and is covered with $516$ polygonal elements. 
The obtained numerical results are reported in Table~\ref{table.wb_gresho}, where we can notice
again that the equilibrium solution is perfectly preserved for very long simulation times. 


\subsubsection{Transport on a Keplerian disk} 
\label{ssec.Transport_Kepleriandisk}

\begin{table}[!tp] 
	\centering
	\numerikNine
	\begin{tabular}{|l|l||cccccc|}
		\hline
		\multicolumn{8}{|c|}{ Keplerian disk with constant density } \\
		\hline
		\hline
		& Time $t = $ & $0.0$ & $0.5$ & $1.0$ & $2.5$ & $5.0$ & $10.0$ \\
		\hline
		\hline
		\multirow{4}{*}{\rotatebox{90}{{DG-$\mathbb{P}_1$}}}
		& $L_2(\rho)$ error & 1.6317E-12 & 1.5758E-12	& 1.5614E-12	& 1.5479E-12	& 1.5382E-12	& 1.5278E-12 \\
		& $L_2(p)$ error  	& 0.0        & 5.3637E-13	& 5.4055E-13	& 5.3275E-13	& 5.2557E-13	& 5.1567E-12 \\
		& No. timestep 	    & 0     	 & 197       	& 382       	& 913       	& 1787	        & 3533       \\
		& No. sliver   	    & 0     	 & 39       	& 171        	& 545       	& 1175	        & 2426       \\
		\cline{2-8}
		\hline
		\hline
		\multirow{4}{*}{\rotatebox{90}{{DG-$\mathbb{P}_2$}}}
		& $L_2(\rho)$ error & 1.5914E-12 & 1.59139E-12	& 1.6498E-12	& 1.6316E-12	& 1.4777E-12	& 1.4941E-12 \\
		& $L_2(p)$ error  	& 0.0        & 6.5942E-13	& 9.2098E-13	& 9.4815E-13	& 4.2645E-13	& 4.2411E-13 \\
		& No. timestep 	    & 0     	 & 376       	& 735       	& 17746       	& 3472	        & 6871       \\
		& No. sliver   	    & 0     	 & 39       	& 168        	& 555       	& 1180	        & 2426       \\
		\cline{2-8}
		\hline
		\hline
		\multirow{4}{*}{\rotatebox{90}{{DG-$\mathbb{P}_3$}}}
		& $L_2(\rho)$ error & 1.5739E-12 & 1.6077E-12	& 1.7032E-12	& 2.9009E-12	& 1.2019E-12	& 3.6240E-12 \\
		& $L_2(p)$ error  	& 0.0        & 7.2762E-13	& 1.0546E-12	& 9.1580E-12	& 3.2491E-12	& 9.6924E-12 \\
		& No. timestep 	    & 0     	 & 634       	& 1243       	& 553      	    & 5900          & 11679      \\
		& No. sliver   	    & 0     	 & 39       	& 169        	& 4440        	& 1180          & 2430       \\
		\cline{2-8}
		\hline
		\hline
		\multirow{4}{*}{\rotatebox{90}{{DG-$\mathbb{P}_4$}}}
		& $L_2(\rho)$ error & 1.5997E-12 & 1.7448E-12	& 2.0628E-12	& 2.9009E-12	    & 5.3900E-12	    & 4.9691E-11 \\
		& $L_2(p)$ error  	& 0.0        & 8.1351E-13	& 1.3706E-12	& 2.9173E-12	    & 1.3212E-11	    & 9.5823E-11 \\
		& No. timestep 	    & 0     	 & 916         	& 1799         	& 4440       		& 8653	   		    & 17019        \\
		& No. sliver   	    & 0     	 & 39       	& 170        	& 553        		& 1180	      	    & 2426         \\
		\cline{2-8}
		\hline
	\end{tabular}
	\caption{Check of the well-balanced property on a Keplerian disk with constant density. As in the previous 
	test cases, we can notice that the equilibrium solution, even if initially perturbed with a random error 
	of 1E-12 distributed everywhere on the domain, is preserved with machine accuracy for very long simulation 
	times and after handling thousands of sliver elements. This holds true for any employed polynomial representation order.}
	\label{table.KeplerianDisk_constant}
\end{table}

Next we consider the Euler equations with the gravity source term as described in~\eqref{eq.eulerTerms}. 
Here, there exist an entire class of stationary solutions characterized by the exact balance 
between the pressure gradient, the centrifugal force and the gravity force, i.e. such that 
\begin{equation} 
\label{eq.EquilibriaConstraint}
\begin{cases} 
&\rho = \rho(r), \\ 
& u_r = 0,  \qquad \partial_\phi v = 0, \\
&\partial_r (rP) = - \rho \left (  \frac{Gm_s}{r} - v^2 \right) + P.
\end{cases} 
\end{equation} 
The objective of our work is that such equilibrium solutions are preserved exactly also at the discrete level, 
so to be able to model with increased (and otherwise unreachable) precision small perturbations arising around 
those stationary profiles. 

We start by considering an equilibrium solution, belonging to the above family~\eqref{eq.EquilibriaConstraint}, 
with constant density $\rho=1$ and pressure $p=1$, angular 
velocity $u_\phi = \sqrt{ (G m_s)/r }$ and we consider as computational domain $\Omega$ the ring
$\RIIcolor{(r, \varphi)} \in [1,2] \times [0, 2\pi]$. 
(We recall that all the simulations of this paper are performed in Cartesian coordinates).

First, we verify the well-balanced property of our scheme also on this set of equations including a non zero source term.
In Table~\ref{table.KeplerianDisk_constant}, 
we can see again that our equilibrium solution, initially perturbed with a random error of 1E-12 distributed 
over all the domain, is maintained with machine precision for long simulations times on a moving 
coarse mesh of $539$ polygonal elements.

\medskip 

Then, we employ our scheme for its main purpose: the study of events happening close to equilibrium 
solutions as for example the transport of some density perturbations over the Keplerian disks. 

\begin{figure}[!bp] \vspace{0pt} \centering
	\includegraphics[width=0.31\linewidth, trim=8 8 8 8, clip]{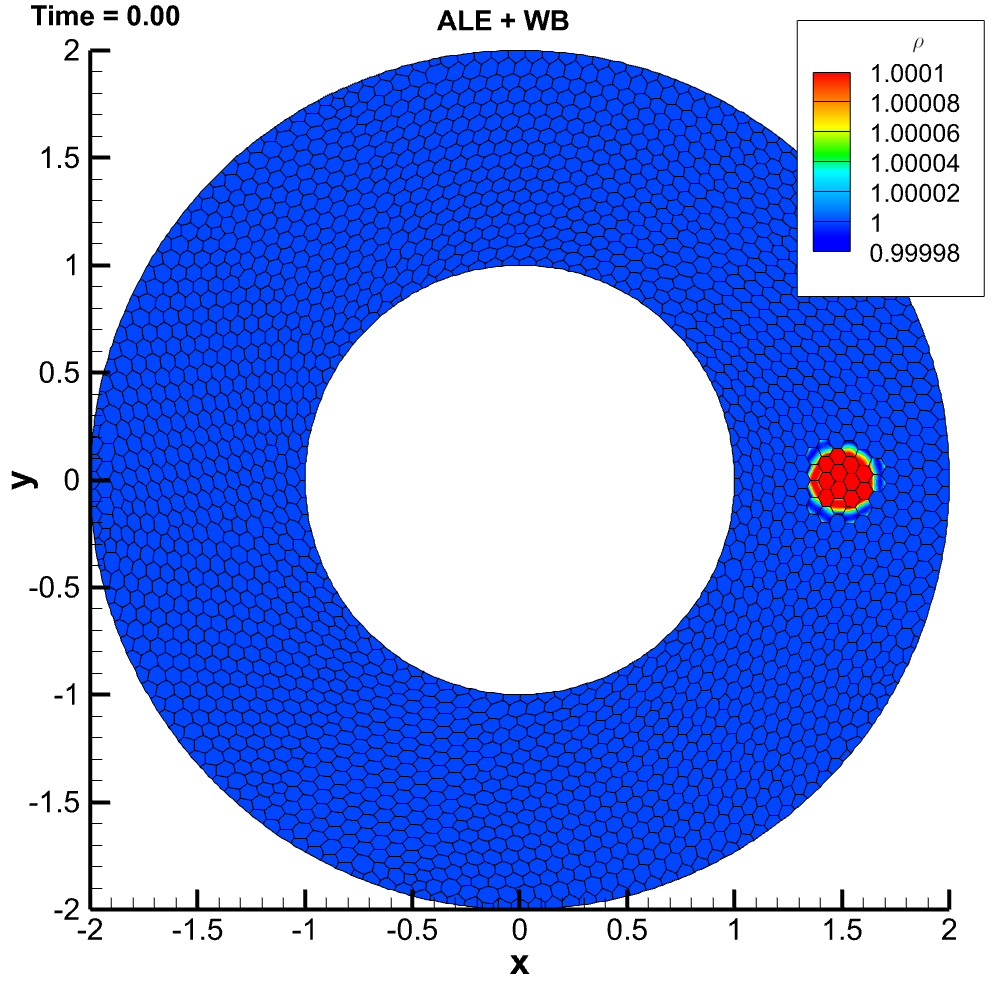}%
	\includegraphics[width=0.31\linewidth, trim=8 8 8 8, clip]{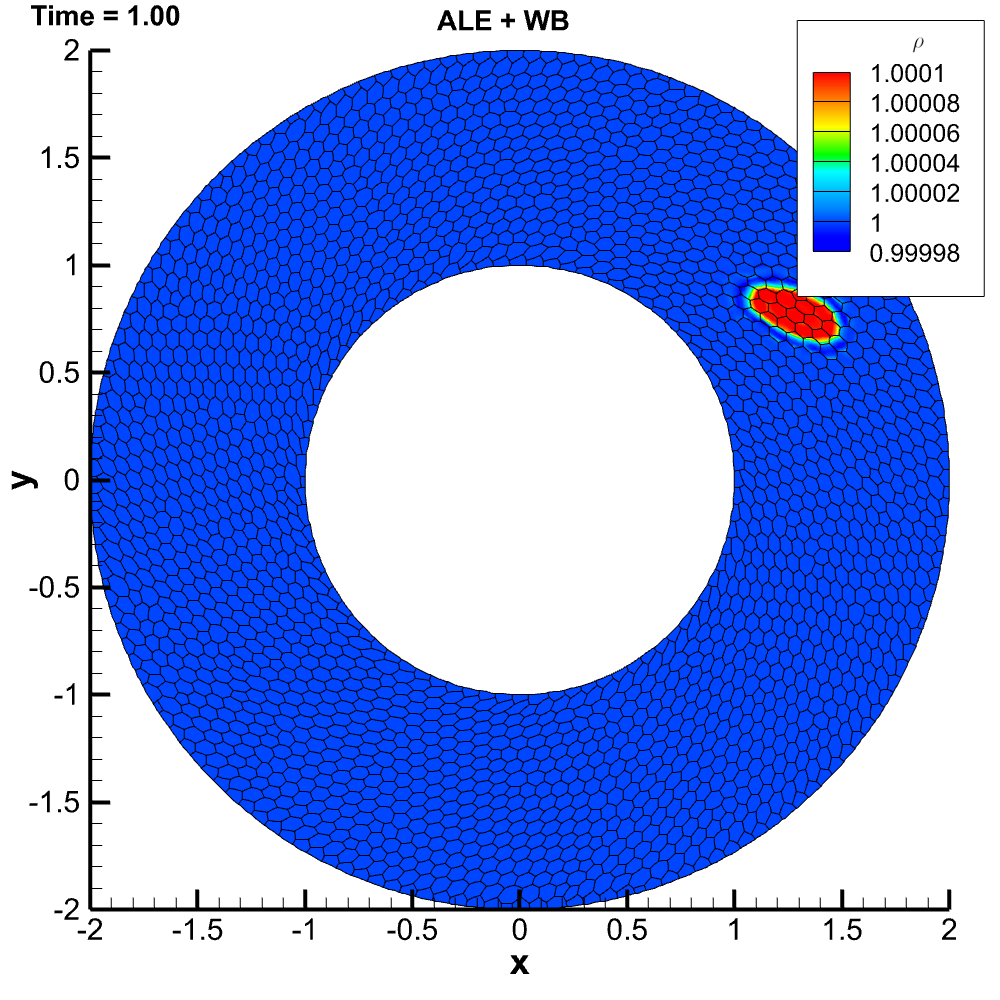}%
	\includegraphics[width=0.31\linewidth, trim=8 8 8 8, clip]{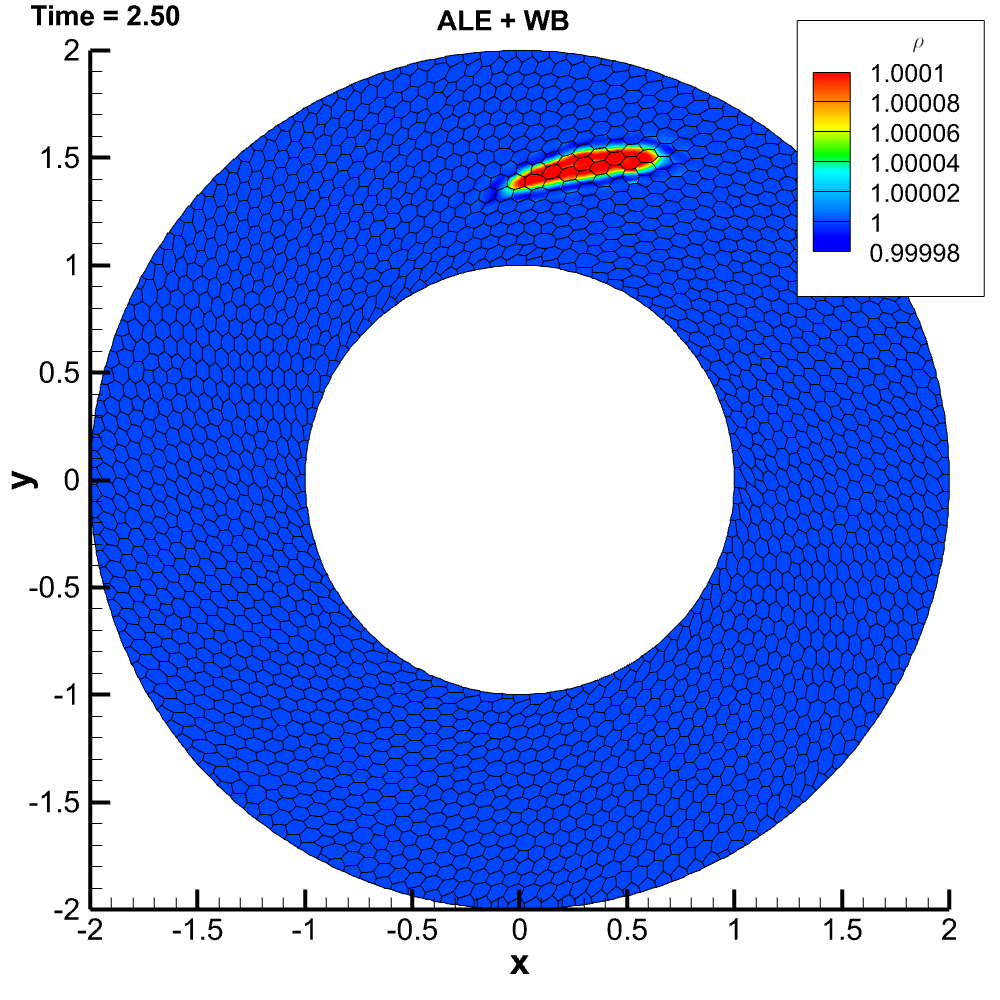}\\[0pt]
	\includegraphics[width=0.31\linewidth, trim=8 8 8 8, clip]{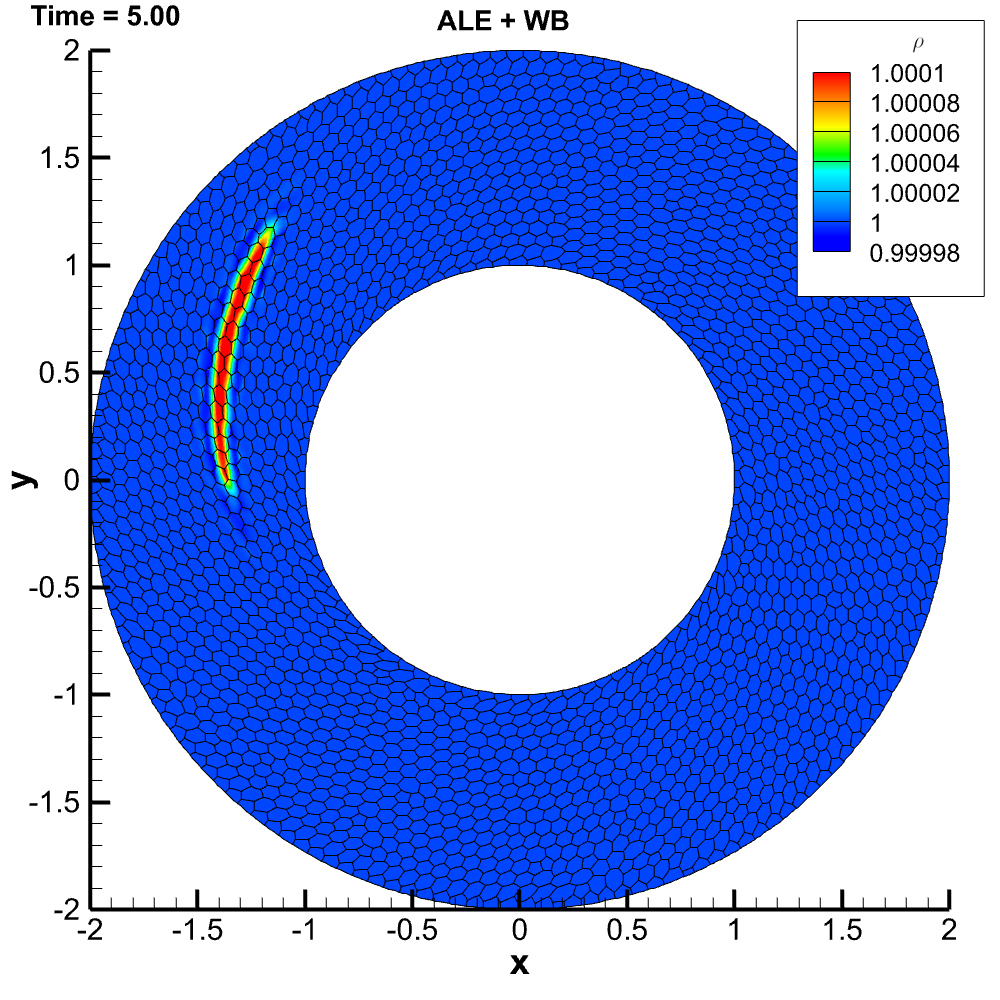}%
	\includegraphics[width=0.31\linewidth, trim=8 8 8 8, clip]{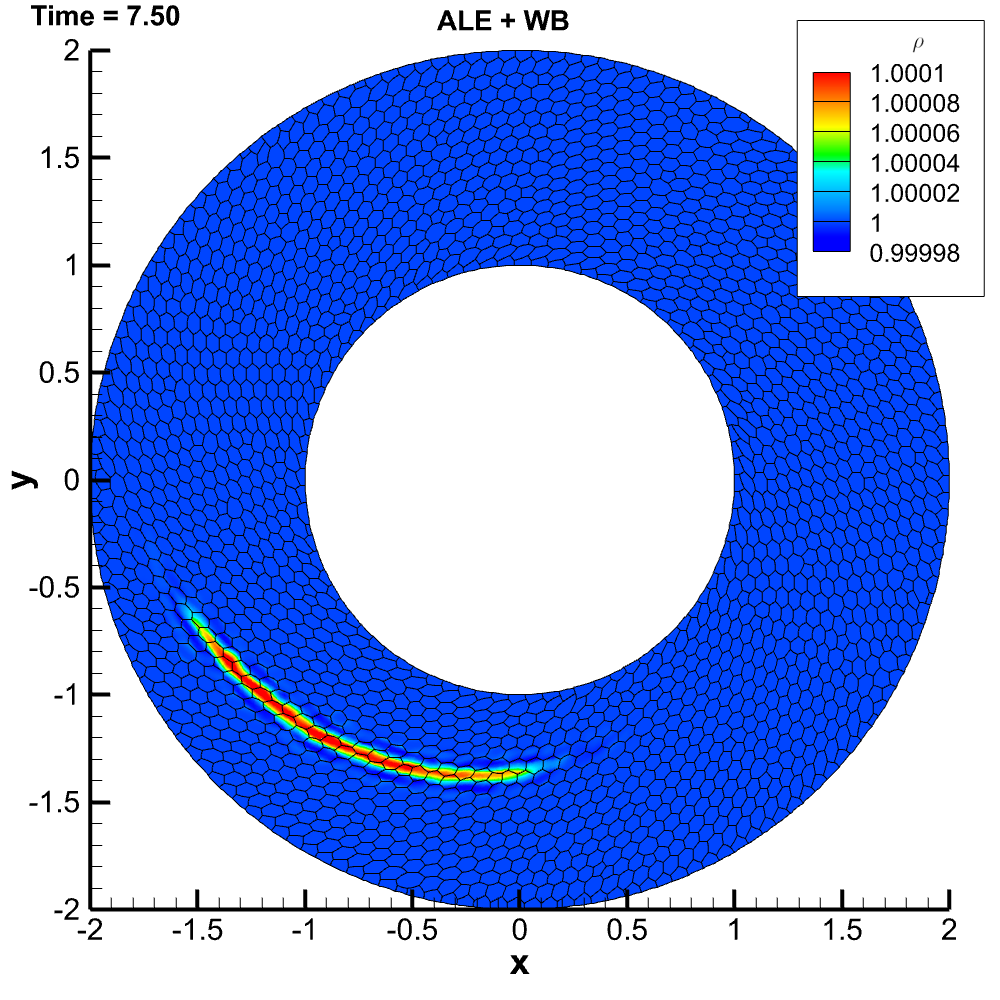}%
	\includegraphics[width=0.31\linewidth, trim=8 8 8 8, clip]{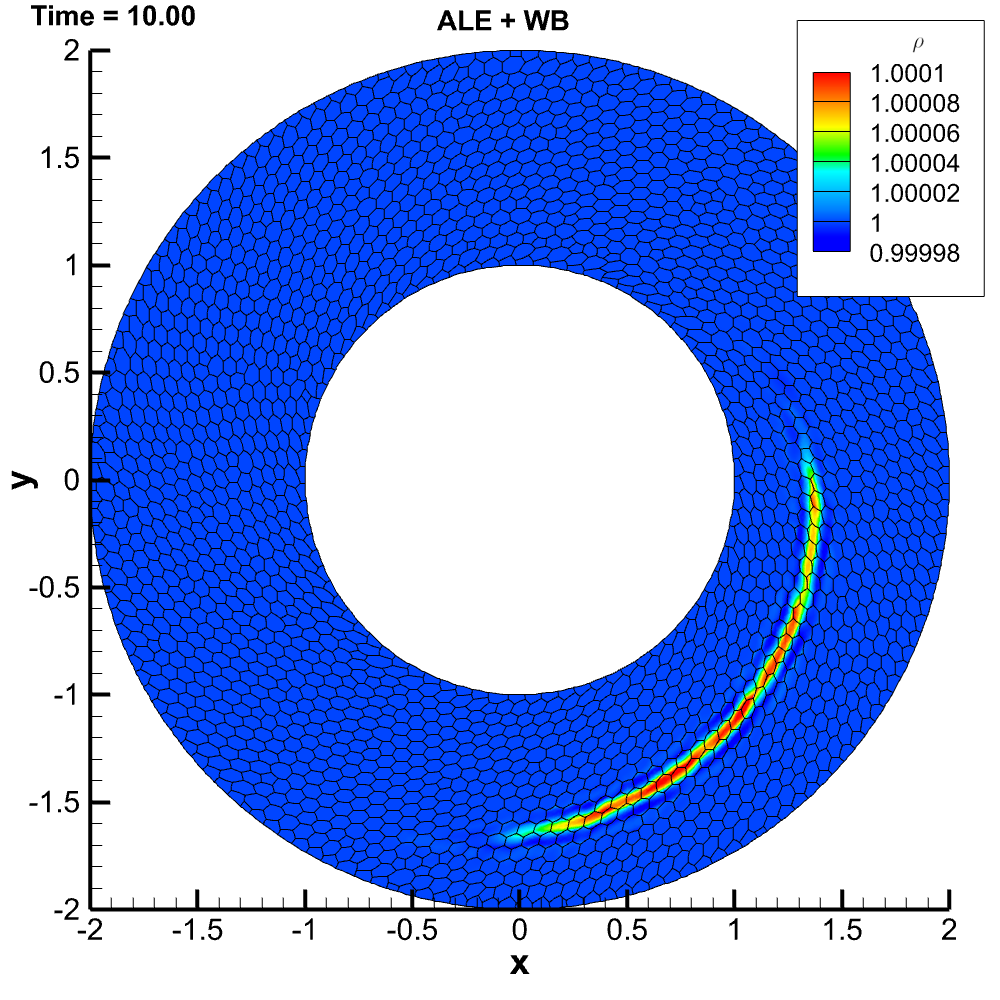}\\[-10pt]
	\caption{Transport of a small density perturbation (of magnitude 1E-4) on a Keplerian disk. 
		We plot the density profile obtained with our WB ALE DG scheme of order $3$ on a quite 
		coarse mesh of only $2152$ polygonal elements at successive times from top left to bottom right. 
		Despite the coarse mesh and the tiny amplitude of the phenomenon we want to observe, 
		thanks to the joint effect of our Lagrangian method and the well-balanced techniques, 
		the density perturbation is transported along the ring, subject to the differential vortical rotation, 
		showing only a quite low numerical dissipation and so it is clearly distinguishable from the background.}
	\label{fig.pallina4_film} \vspace{3pt}
\end{figure}

\begin{figure}[!bp] \vspace{0pt} \centering
	\includegraphics[width=0.31\linewidth, trim=8 8 8 8, clip]{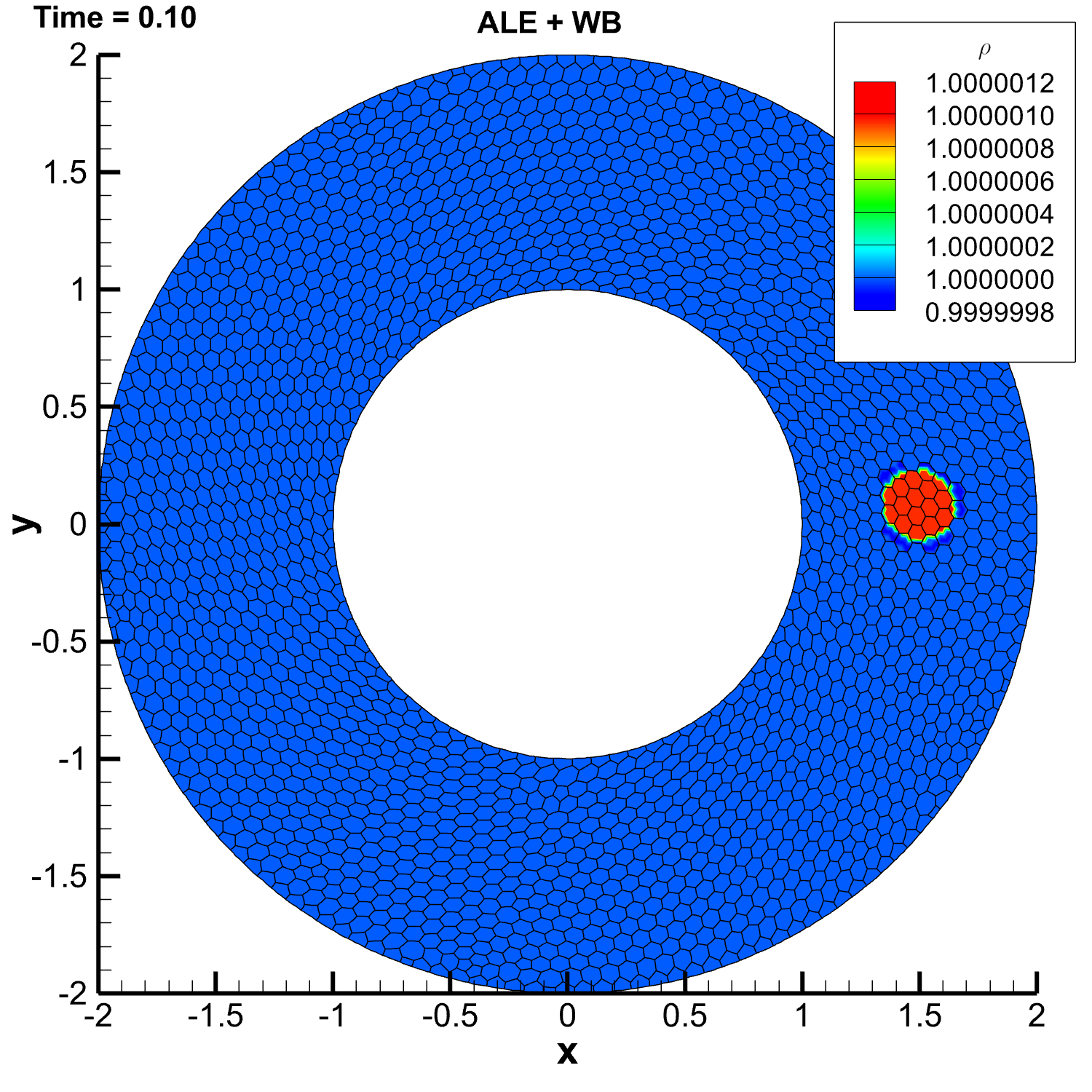}%
	\includegraphics[width=0.31\linewidth, trim=8 8 8 8, clip]{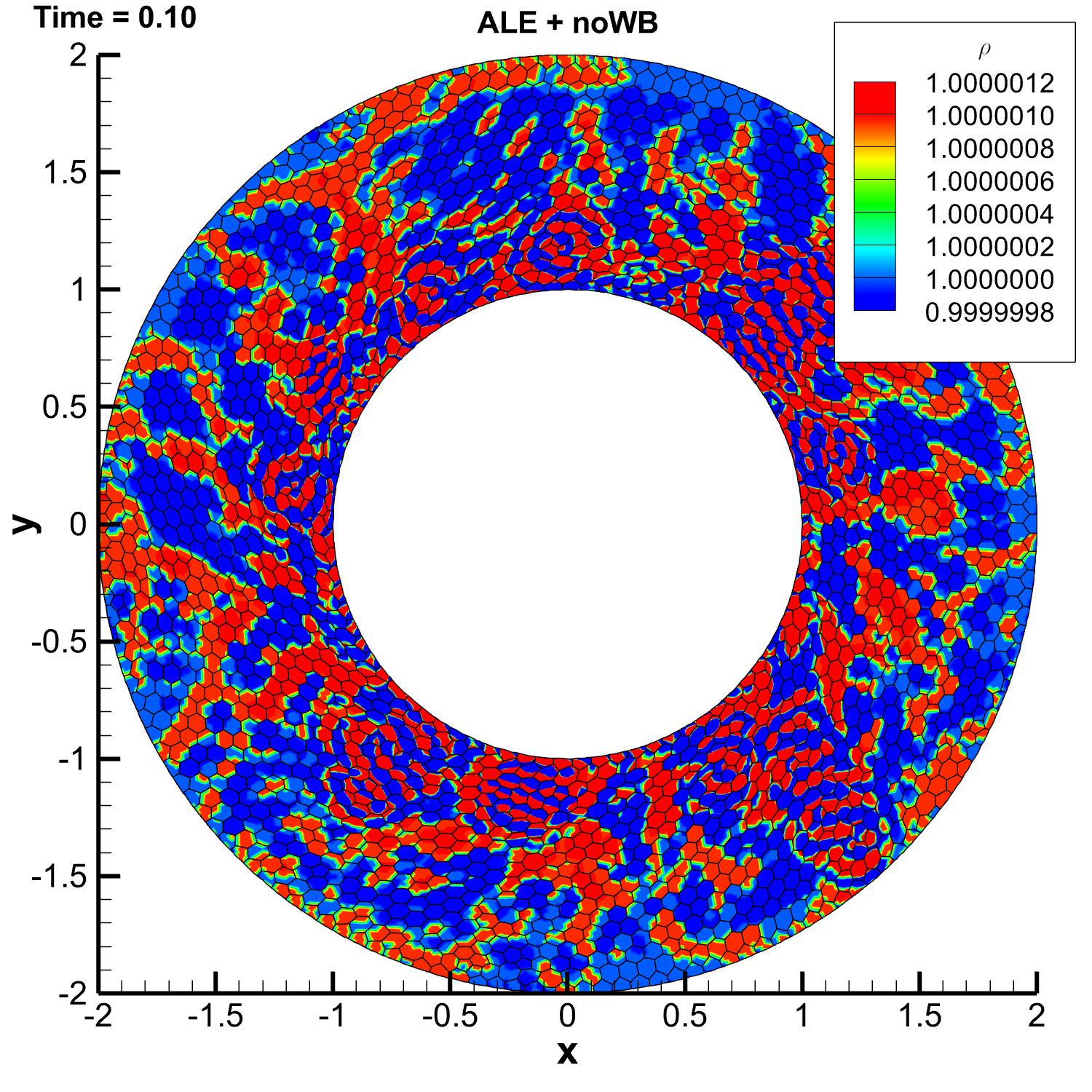}%
	\includegraphics[width=0.31\linewidth, trim=8 8 8 8, clip]{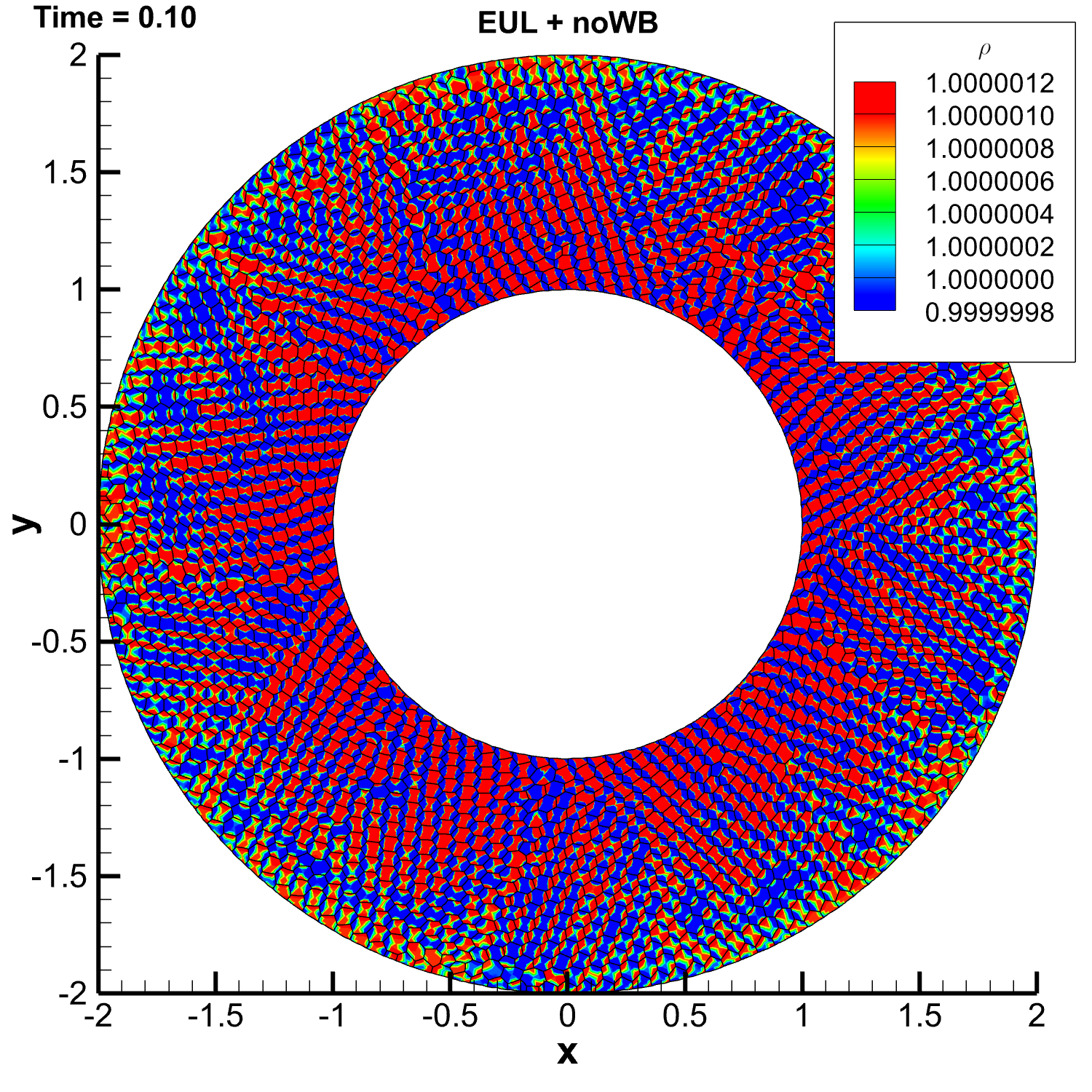}\\
	\includegraphics[width=0.31\linewidth, trim=8 8 8 8, clip]{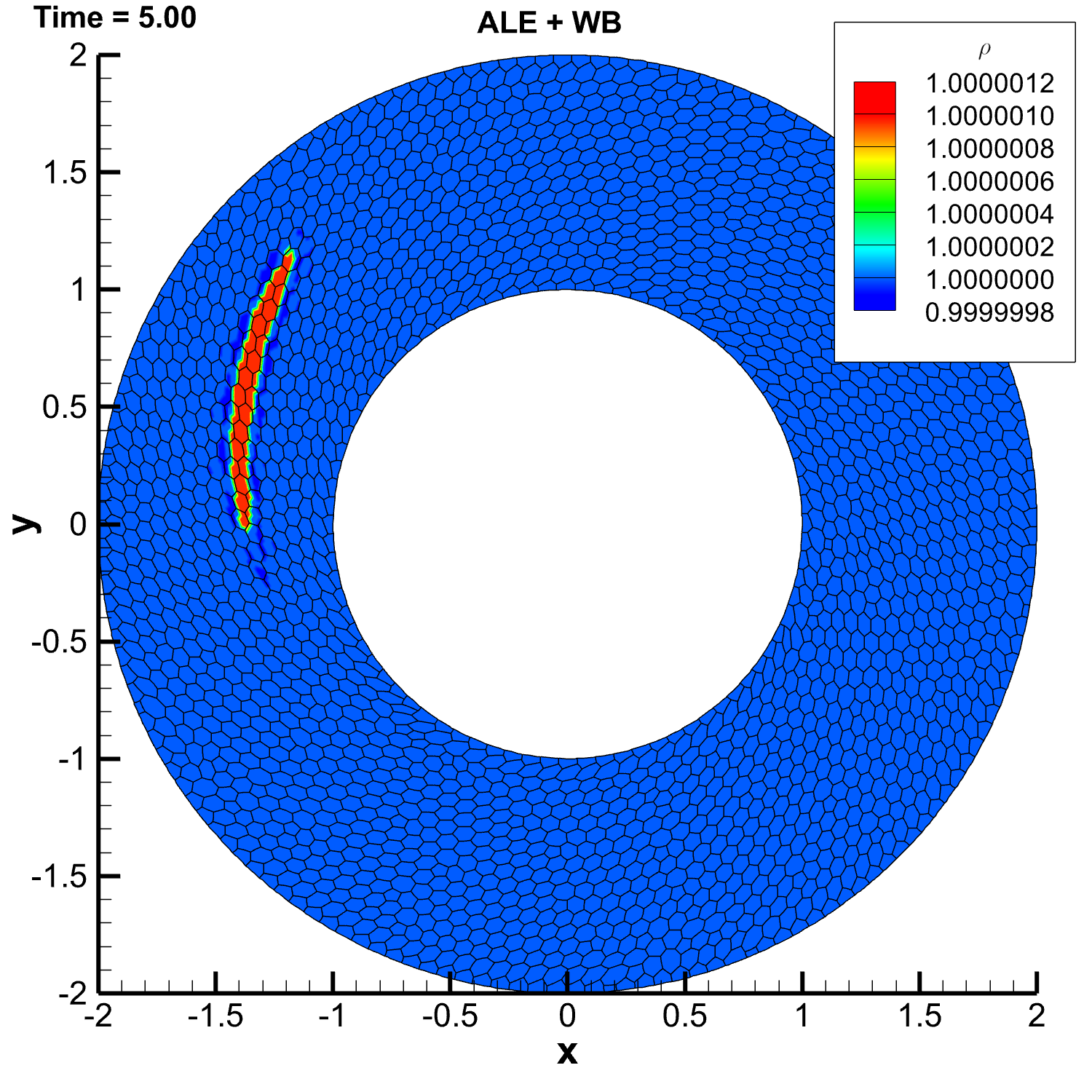}%
	\includegraphics[width=0.31\linewidth, trim=8 8 8 8, clip]{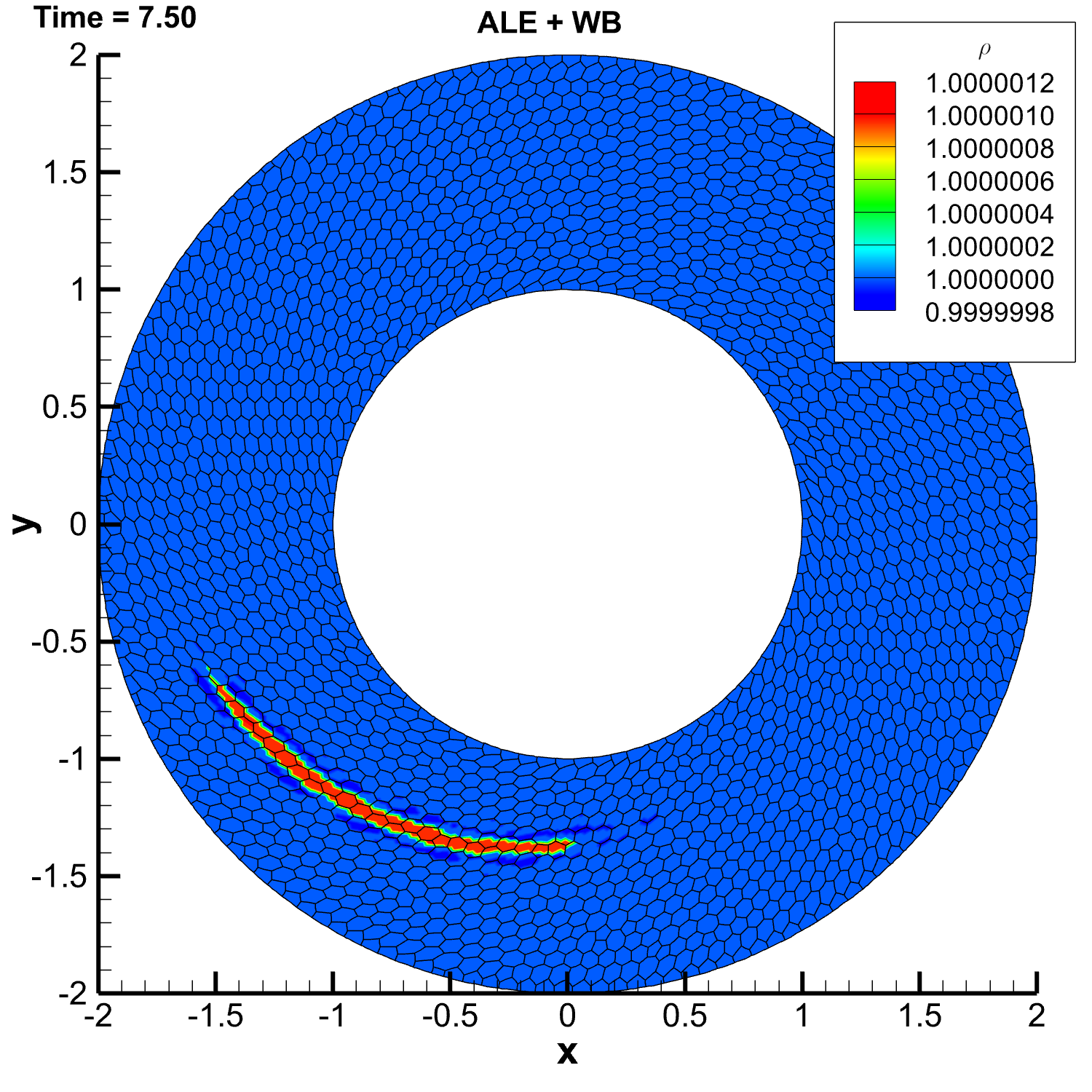}%
	\includegraphics[width=0.31\linewidth, trim=8 8 8 8, clip]{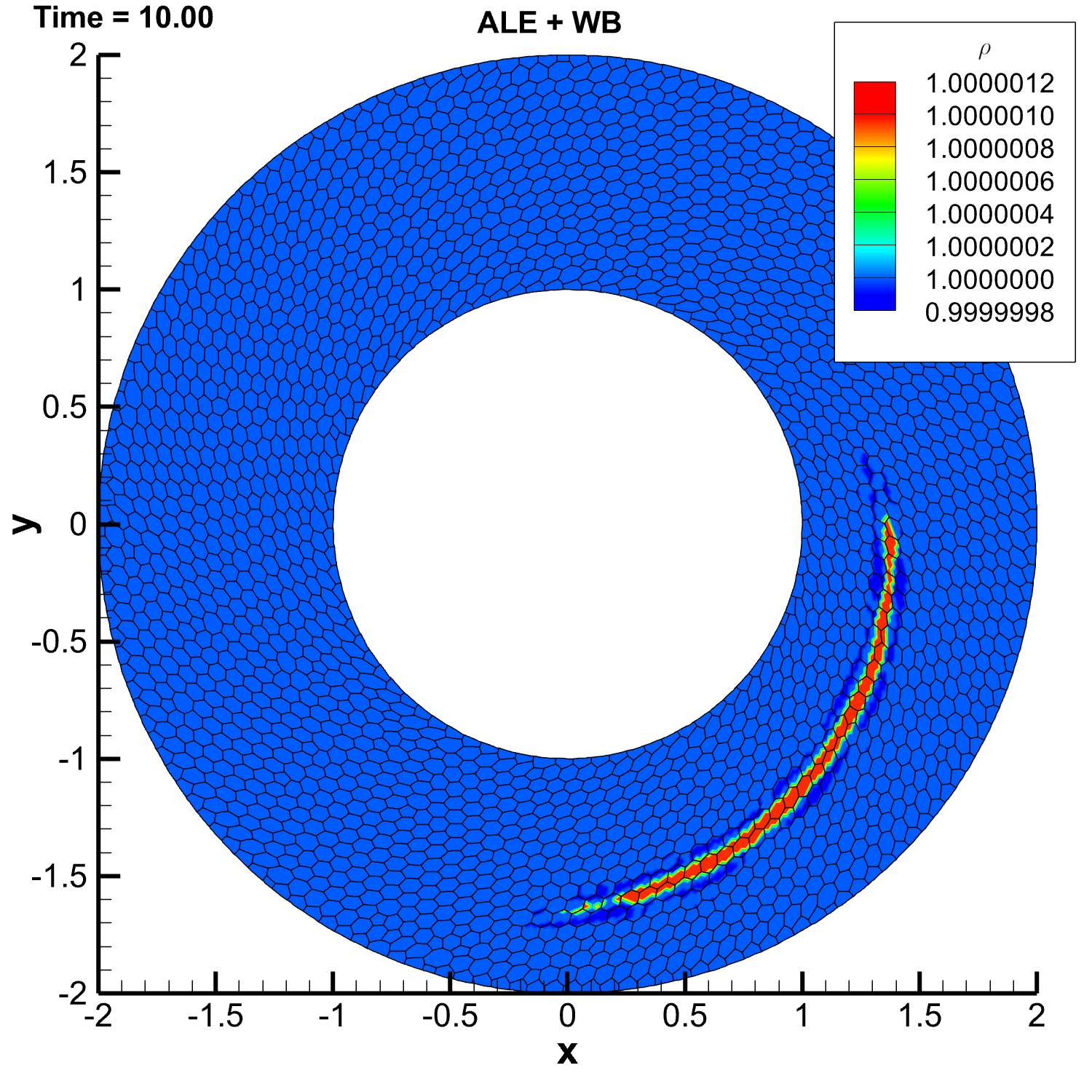}\\[-12pt]
	\caption{Transport of a very small density perturbation (of magnitude 1E-6) on a Keplerian disk. We observe that 
	our well-balanced ALE DG scheme is able to perfectly simulate this phenomenon for a long time (see pictures 
	up to $t=10$), while the non well-balanced ALE and Eulerian schemes, already at time $t=0.1$, show 
	numerical errors that strongly interfere with such a small density perturbation.} 
	\label{fig.pallina6_film_comparison}  \vspace{0pt}
\end{figure}

\begin{figure}[!bp] 	\vspace{0pt} \centering
	\includegraphics[width=0.31\linewidth, trim=8 8 8 8, clip]{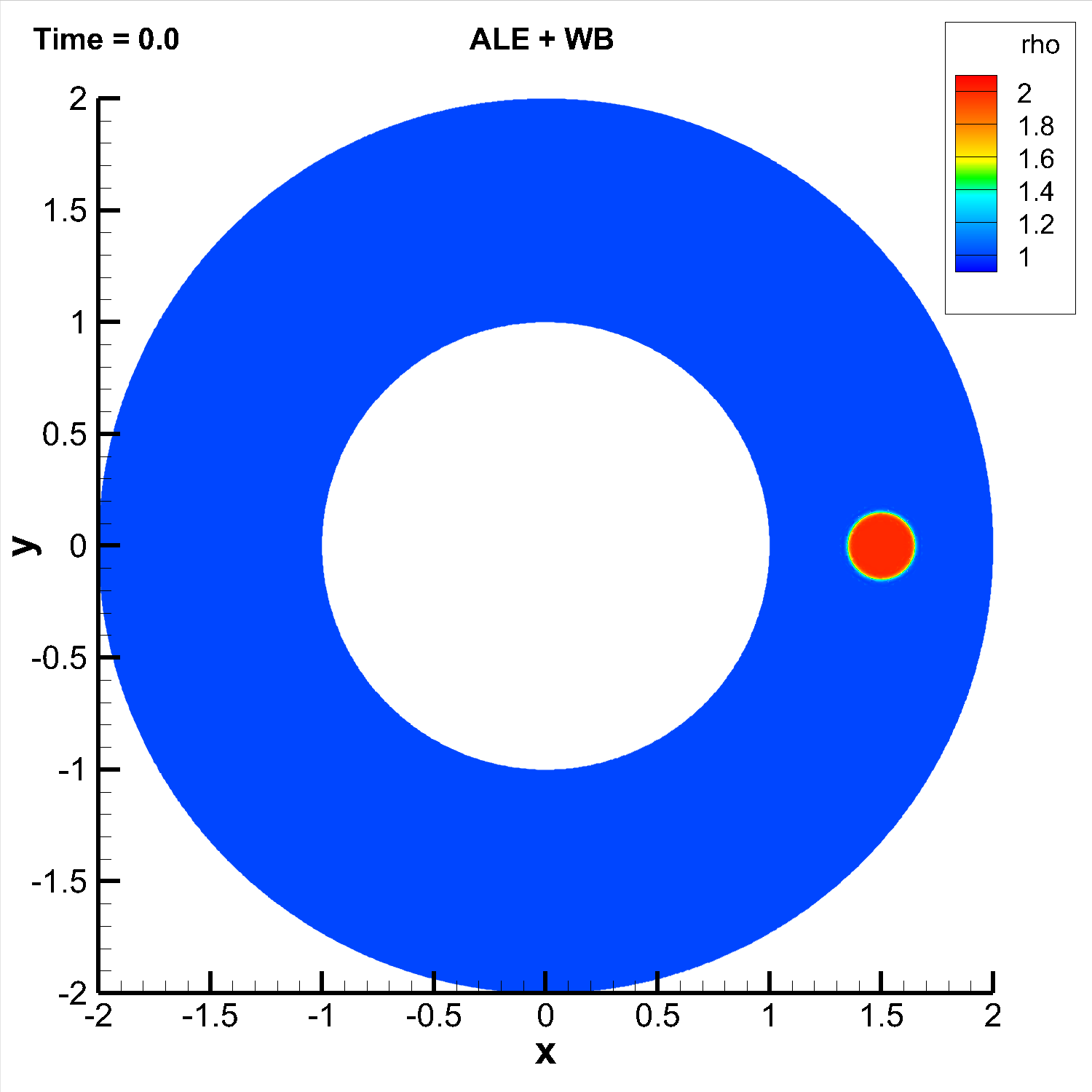}%
	\includegraphics[width=0.31\linewidth, trim=8 8 8 8, clip]{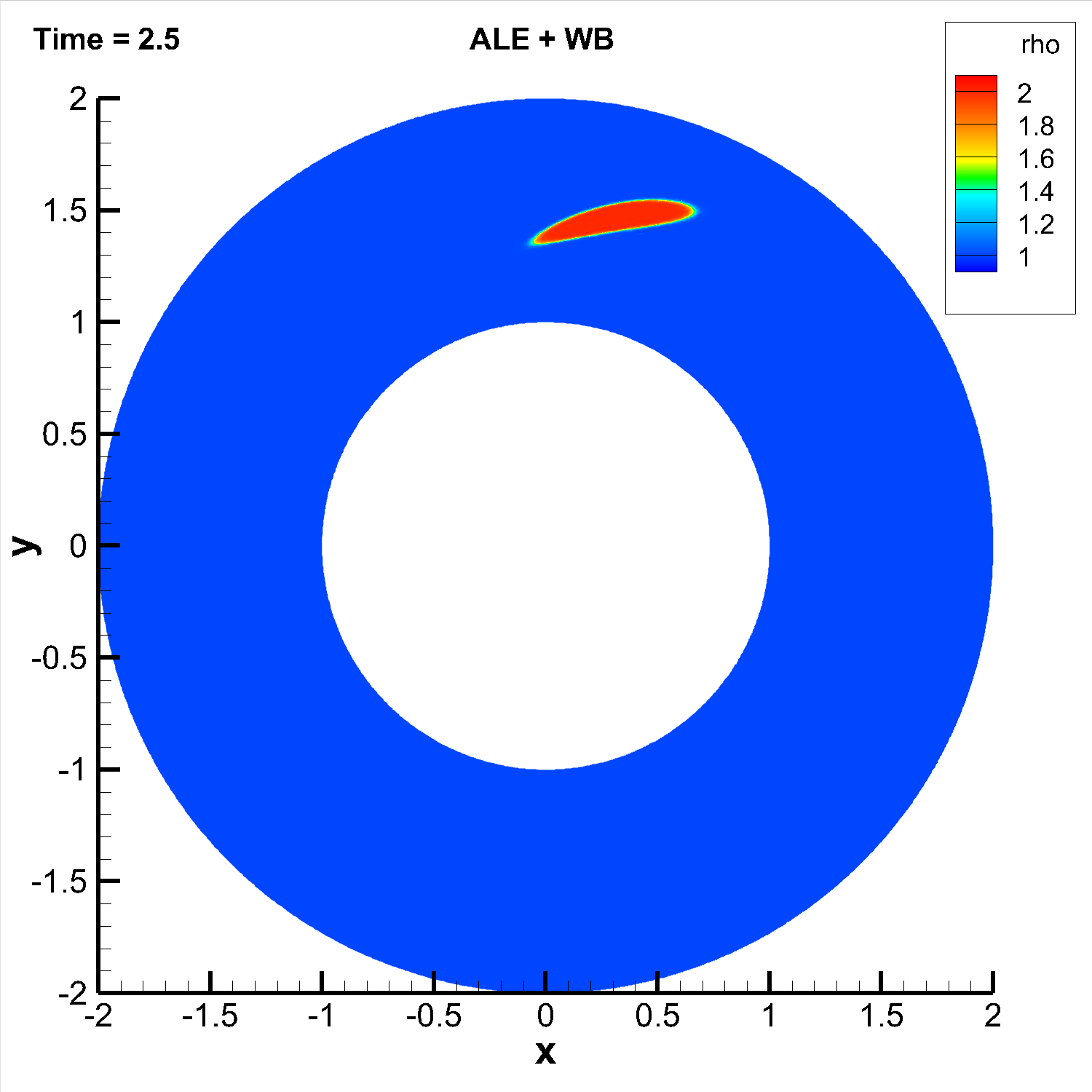}%
	\includegraphics[width=0.31\linewidth, trim=8 8 8 8, clip]{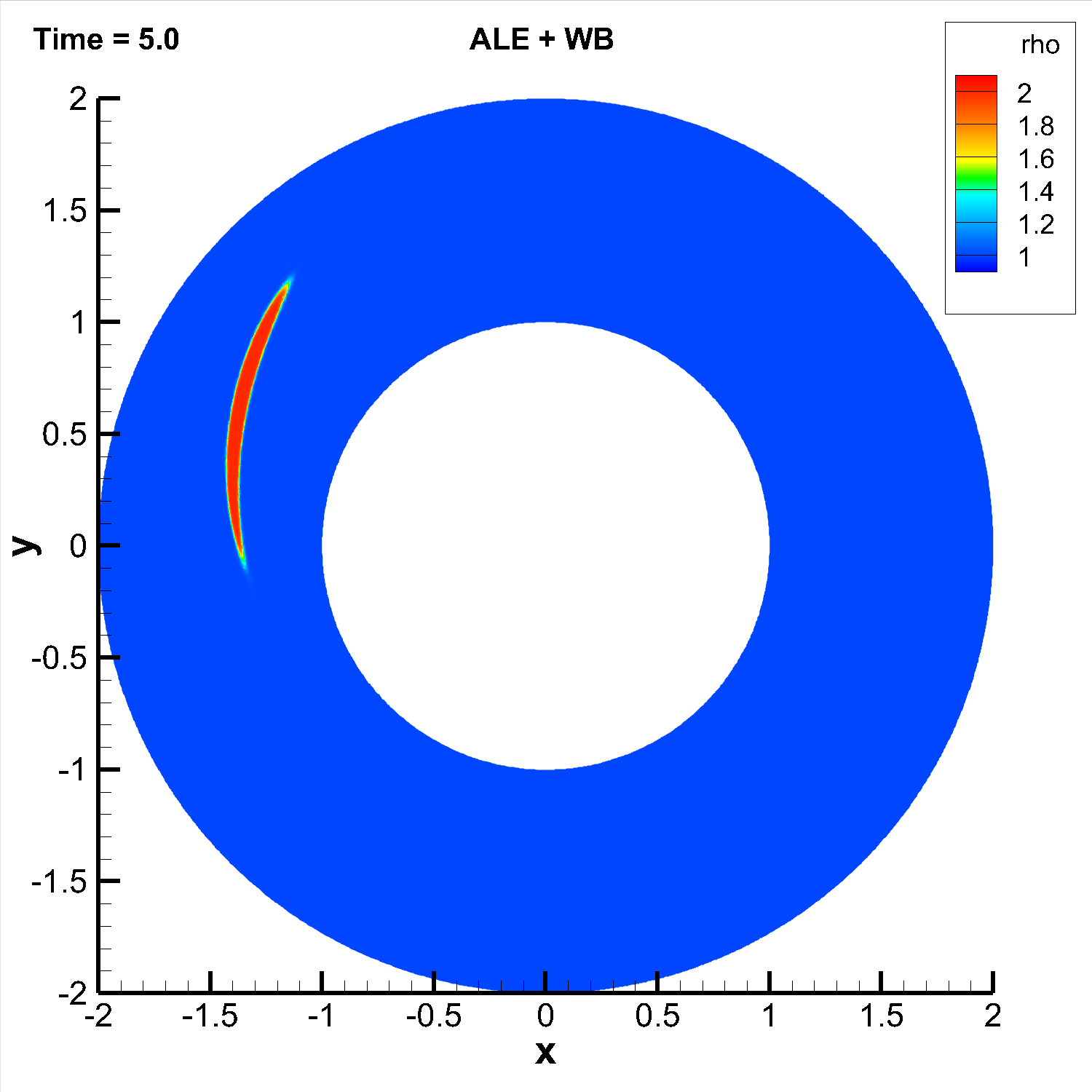}\\[2pt]%
	\includegraphics[width=0.31\linewidth, trim=8 8 8 8, clip]{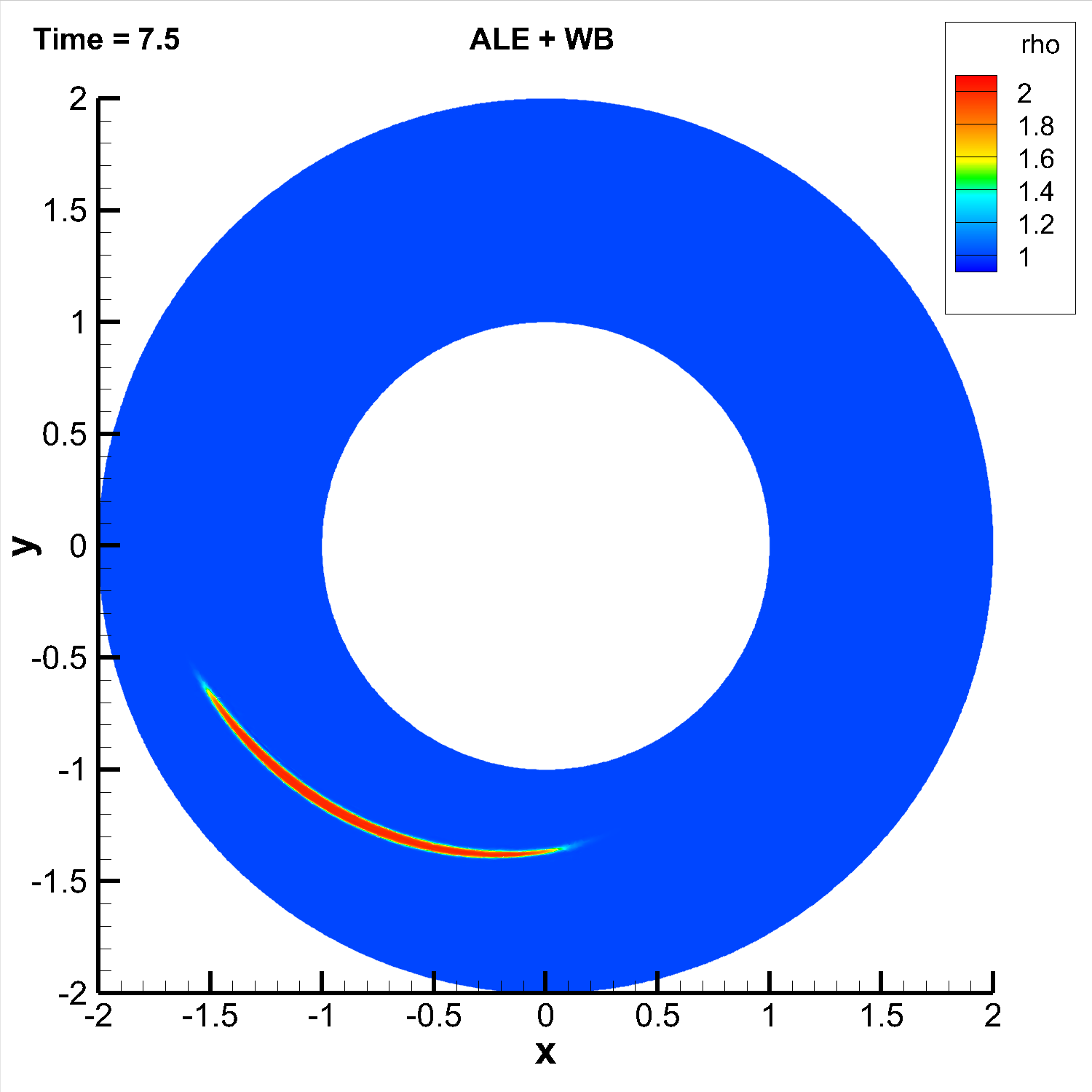}%
	\includegraphics[width=0.31\linewidth, trim=8 8 8 8, clip]{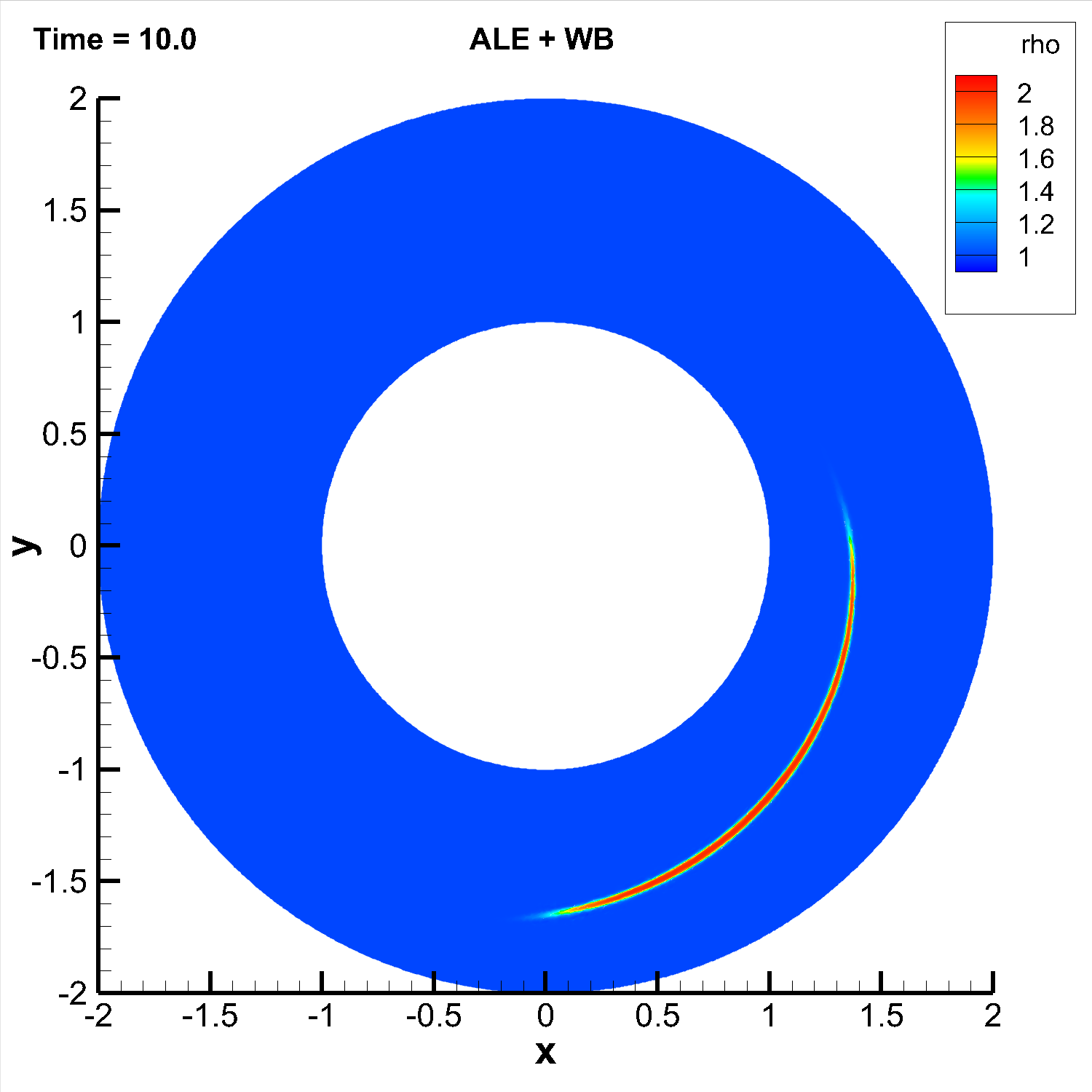}%
	\includegraphics[width=0.31\linewidth, trim=8 8 8 8, clip]{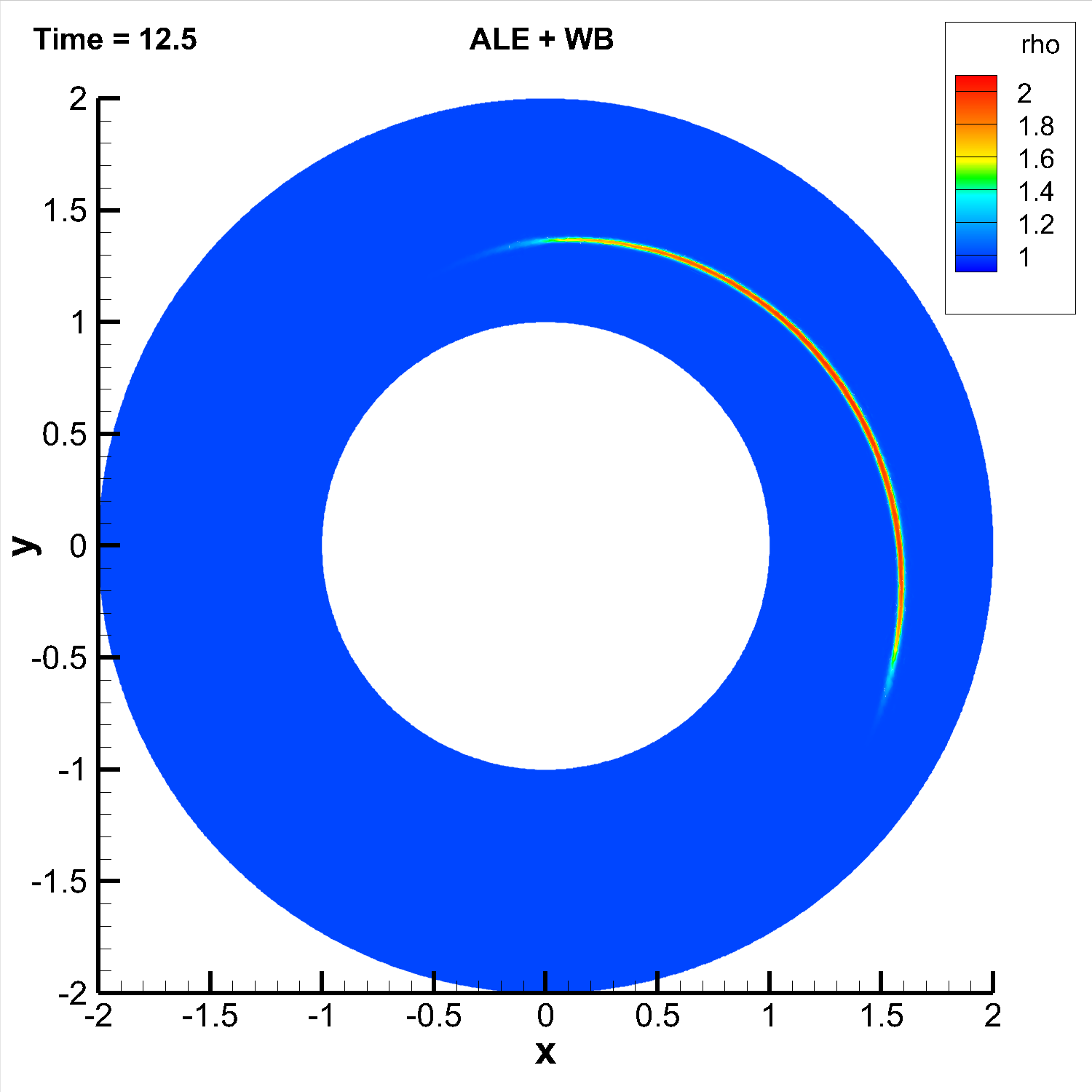}\\[-10pt]
	\caption{Transport of a heavy mass perturbation (of magnitude $1$) on a Keplerian disk. We plot the 
	density profile obtained with our WB ALE DG scheme of order $3$ on a fine mesh of $11080$ 
	polygonal elements at successive times from top left to bottom right.}	
	\label{fig.pallina1000_film}
\end{figure}

Thus, we can take the above discussed equilibrium solution and in the region 
\begin{equation}
	\Omega_\text{disk}  \ \text{ s.t. } \  r < 0.25 \ \text{ with } r = \sqrt{ (x-1.5)^2 + y^2 }
\end{equation} 
we add a quantity of amplitude $A$ to the density profile.
For a visual interpretation, one can refer to the first panel of 
Figures~\ref{fig.pallina4_film},~\ref{fig.pallina6_film_comparison} and~\ref{fig.pallina1000_film}.
This disk with a higher density is then advected along the ring, almost without any dissipation, 
following a velocity field which is strongest as $r \rightarrow 1$: hence, the disk is stretched
in accordance to the advective effects governing the expected behavior of this test problem.

\begin{figure}[!bp] \vspace{3pt} \centering
	\includegraphics[width=0.31\linewidth, trim=8 8 8 8, clip]{pallinacoarse_ale_wb_075}%
	\includegraphics[width=0.31\linewidth, trim=8 8 8 8, clip]{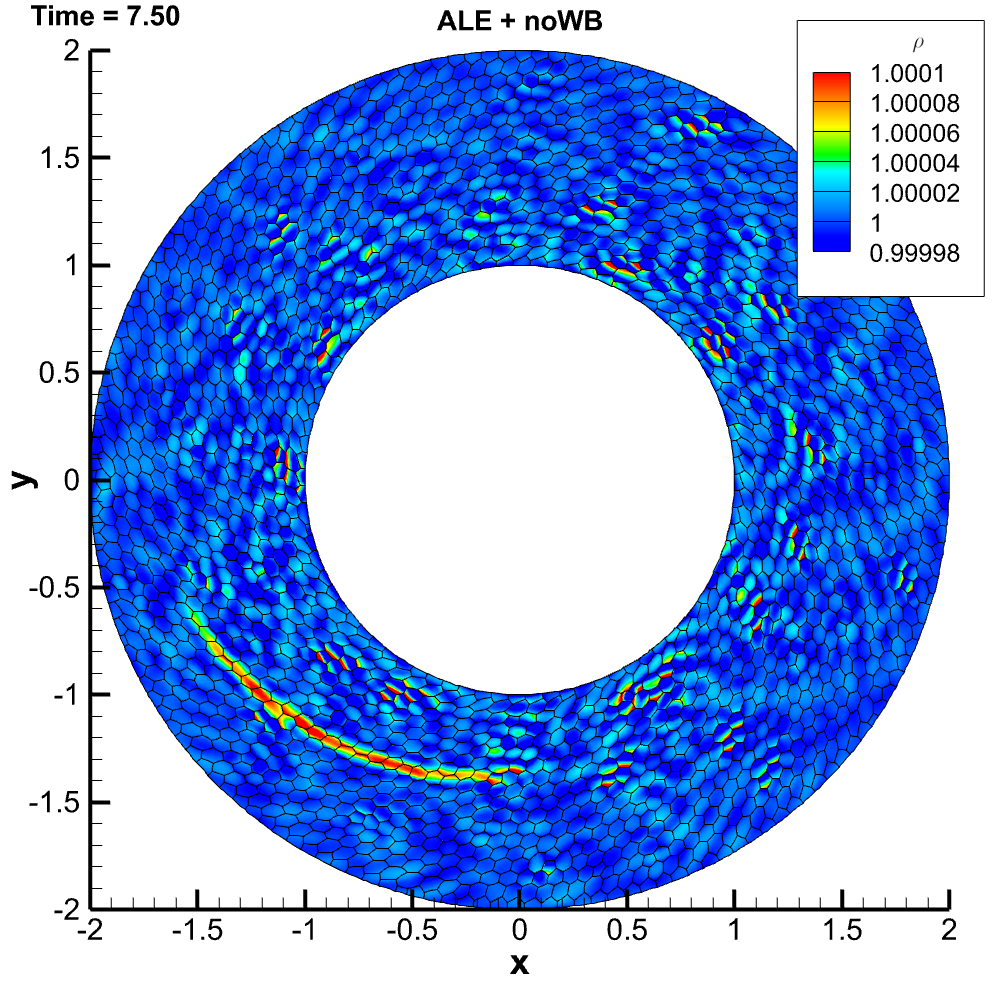}%
	\includegraphics[width=0.31\linewidth, trim=8 8 8 8, clip]{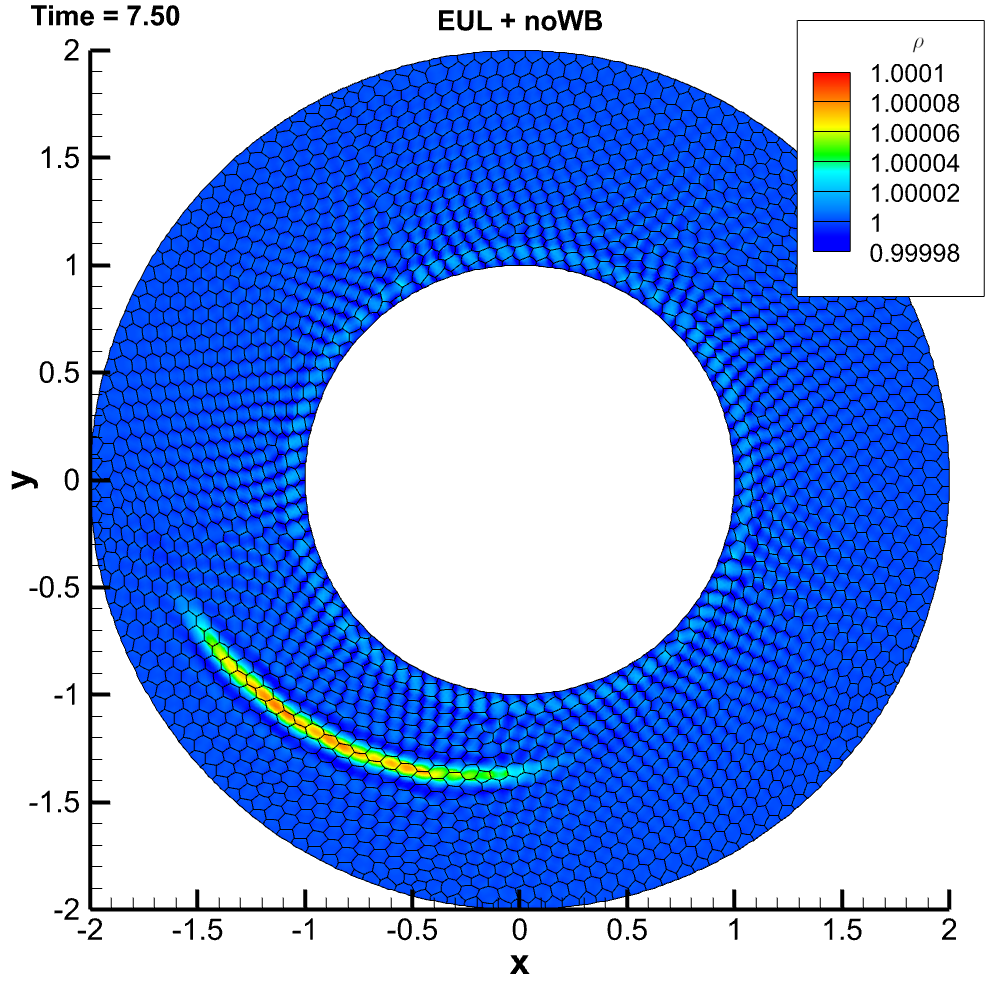}\\[-10pt]
	\caption{Transport of a small density perturbation (of magnitude 1E-4) on a Keplerian disk. In this figure, 
	we compare the results obtained at time $t= 7.5$ from our WB ALE DG scheme of order $3$ with 
	those obtained from \textit{non} well-balanced Lagrangian and Eulerian schemes of the same order, to show 
	the superior resolution and reliability of our approach.}
	\label{fig.pallina4_comparison}
\end{figure}

\begin{figure}[!bp] \vspace{-10pt} \centering
	\includegraphics[width=0.31\linewidth, trim=8 8 8 8, clip]{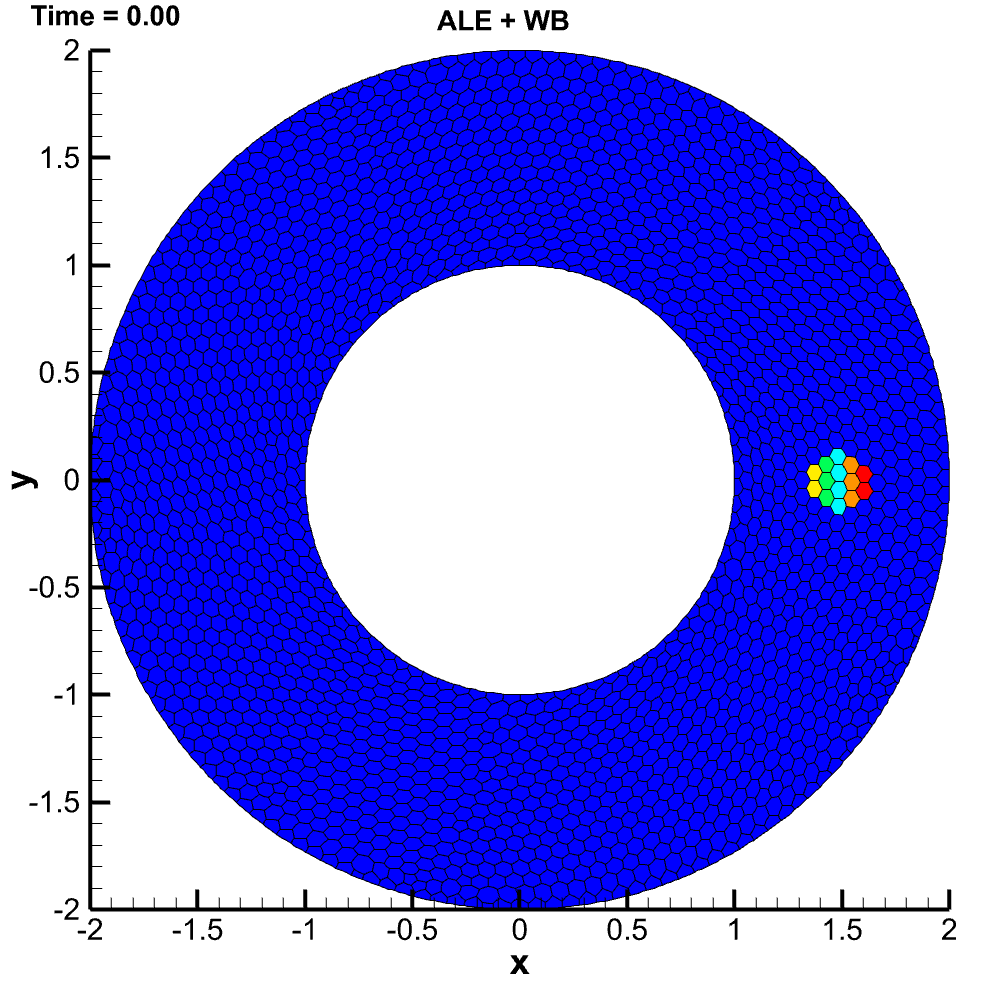}%
	\includegraphics[width=0.31\linewidth, trim=8 8 8 8, clip]{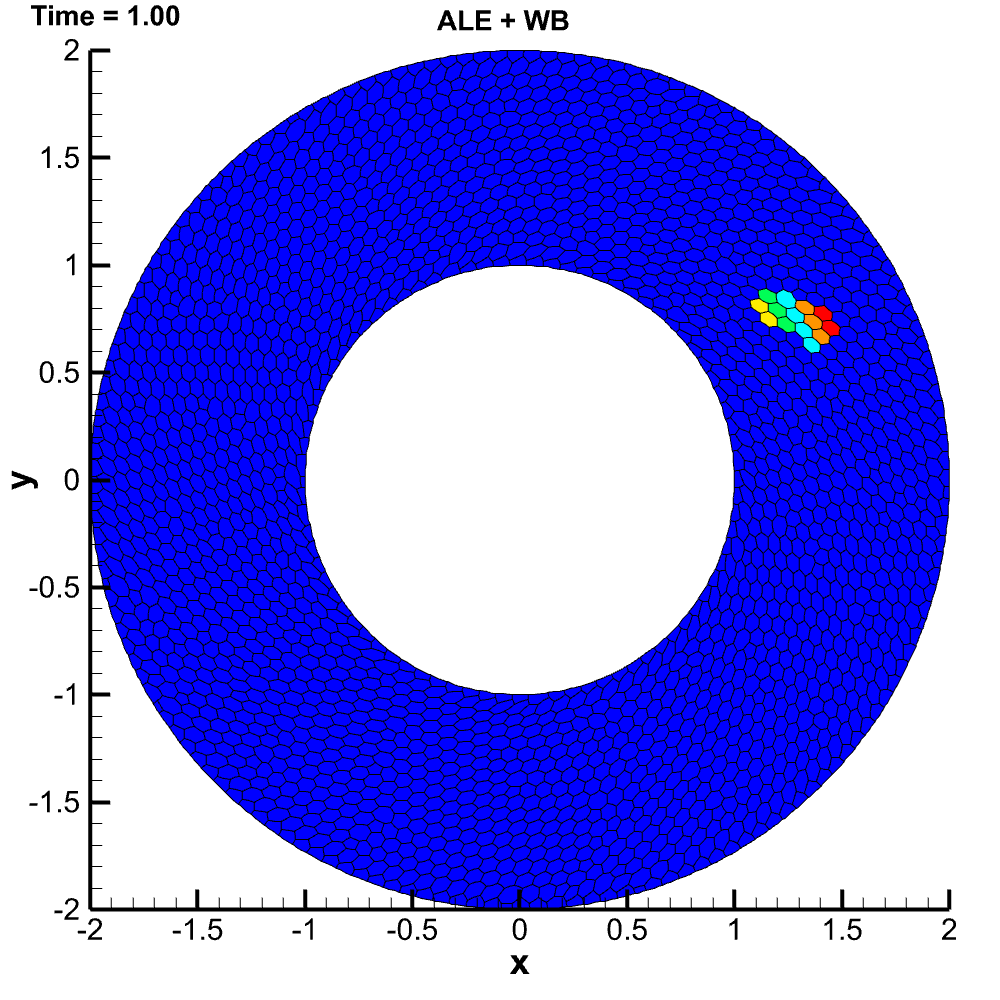}%
	\includegraphics[width=0.31\linewidth, trim=8 8 8 8, clip]{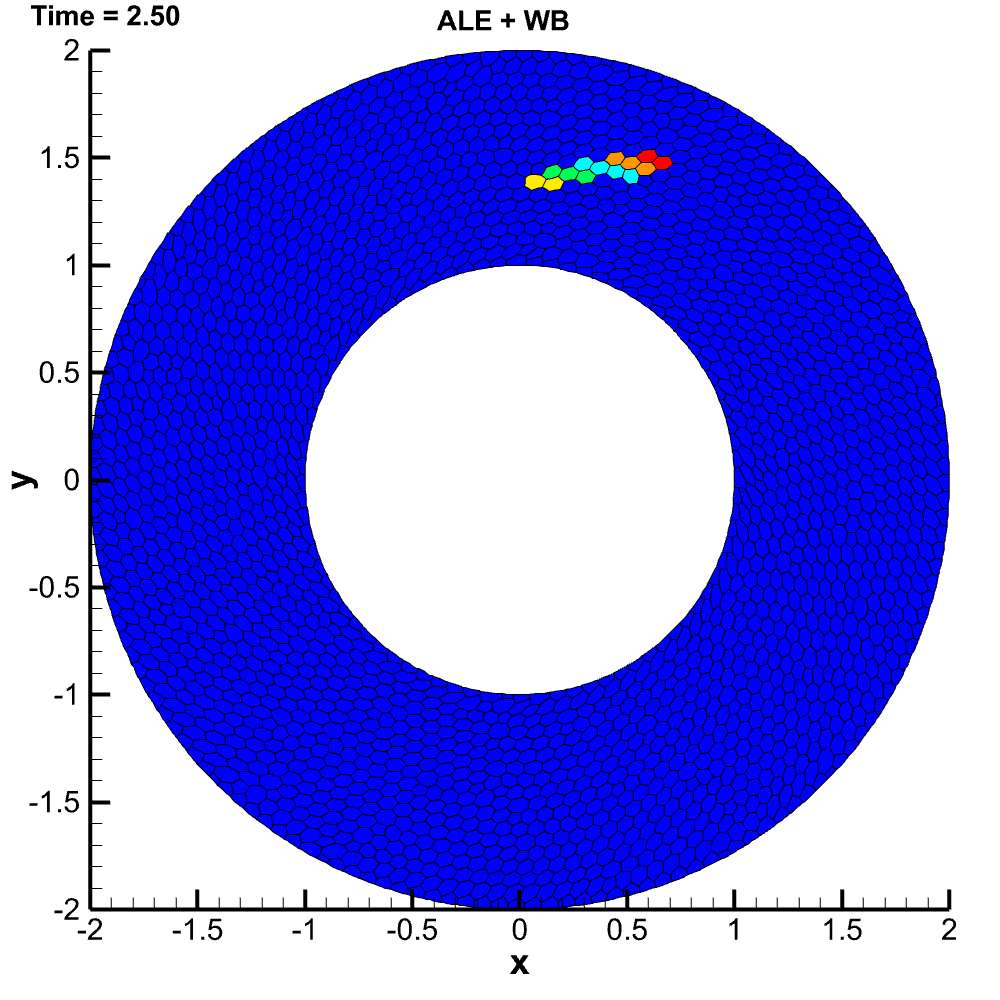}\\[0pt]
	\includegraphics[width=0.31\linewidth, trim=8 8 8 8, clip]{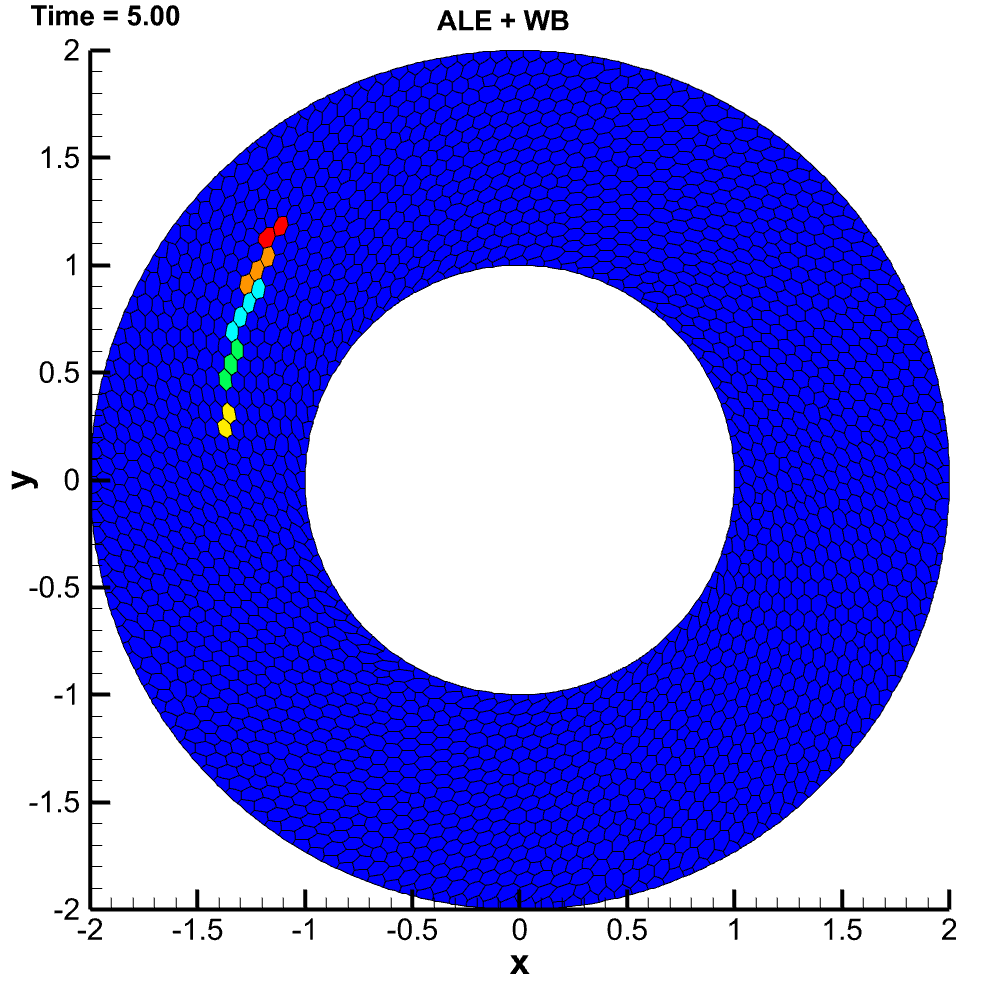}%
	\includegraphics[width=0.31\linewidth, trim=8 8 8 8, clip]{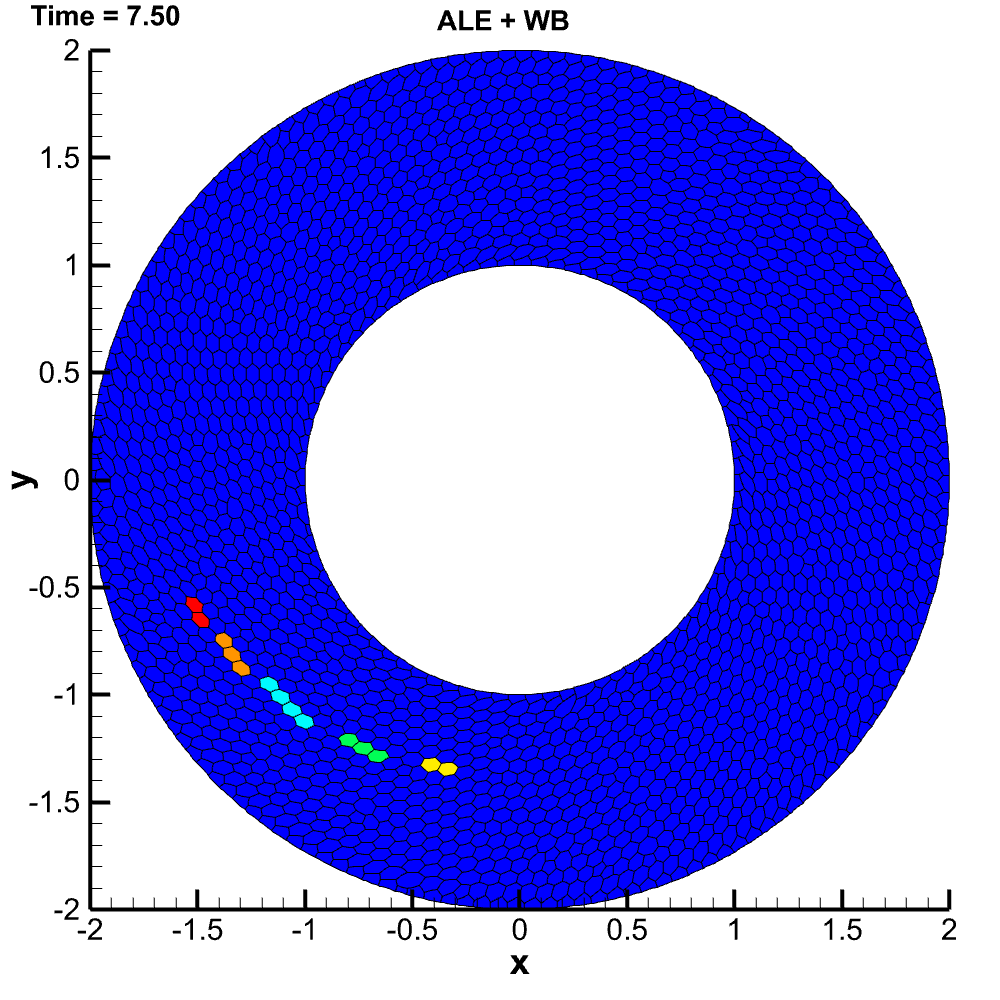}%
	\includegraphics[width=0.31\linewidth, trim=8 8 8 8, clip]{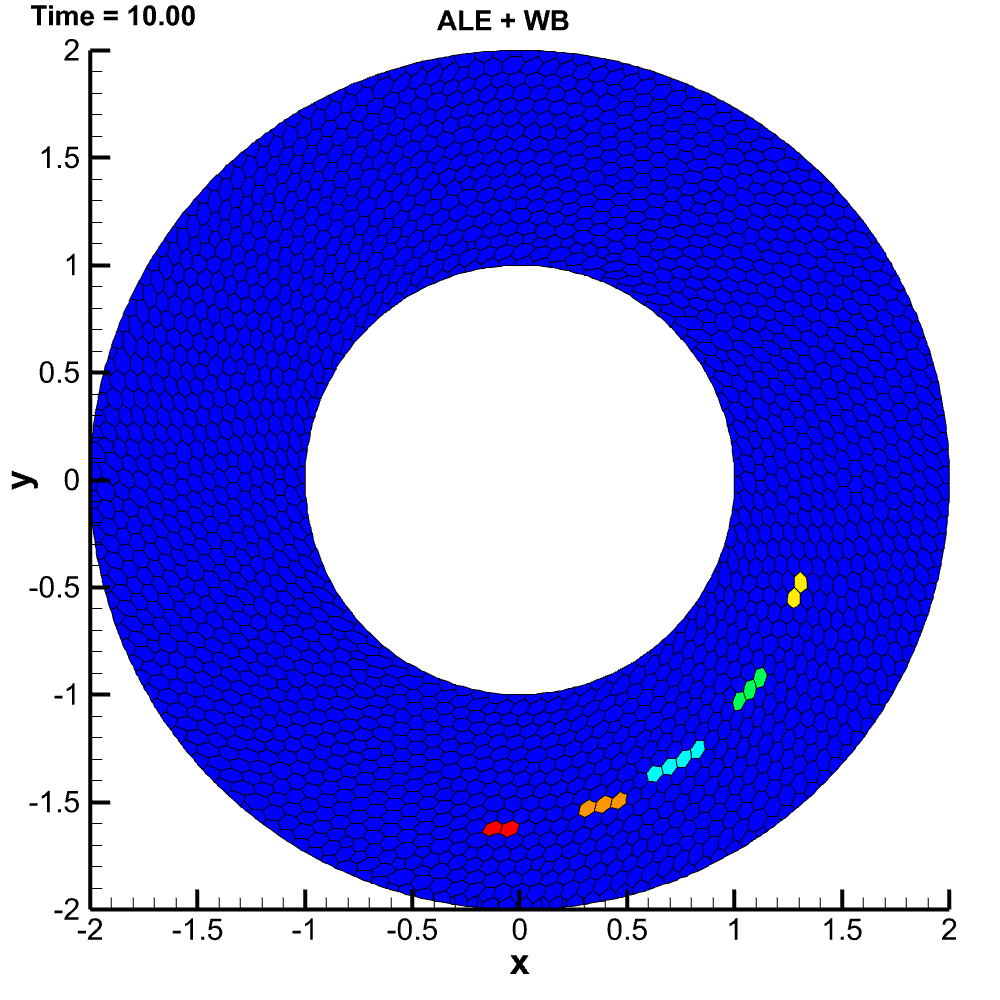}\\[-12pt]
	\caption{Transport of a small density perturbation (of magnitude 1E-4) on a Keplerian disk. 
		In this figure, we highlight a bunch of elements which at time $t=0$ are located in the same position
		of the mass perturbation, and we follow them during the simulation. 
		This image clearly shows that the mesh is moving together with the fluid flow, allowing to considerably 
		reduce the convection errors.}
	\label{fig.pallina4_numbering} \vspace{-2pt}
\end{figure}

\begin{figure}[!bp] \vspace{0pt}  \centering
	\includegraphics[width=0.31\linewidth, trim=8 8 8 8, clip]{pallina_rho_alewb_t075_1500}
	\includegraphics[width=0.31\linewidth, trim=8 8 8 8, clip]{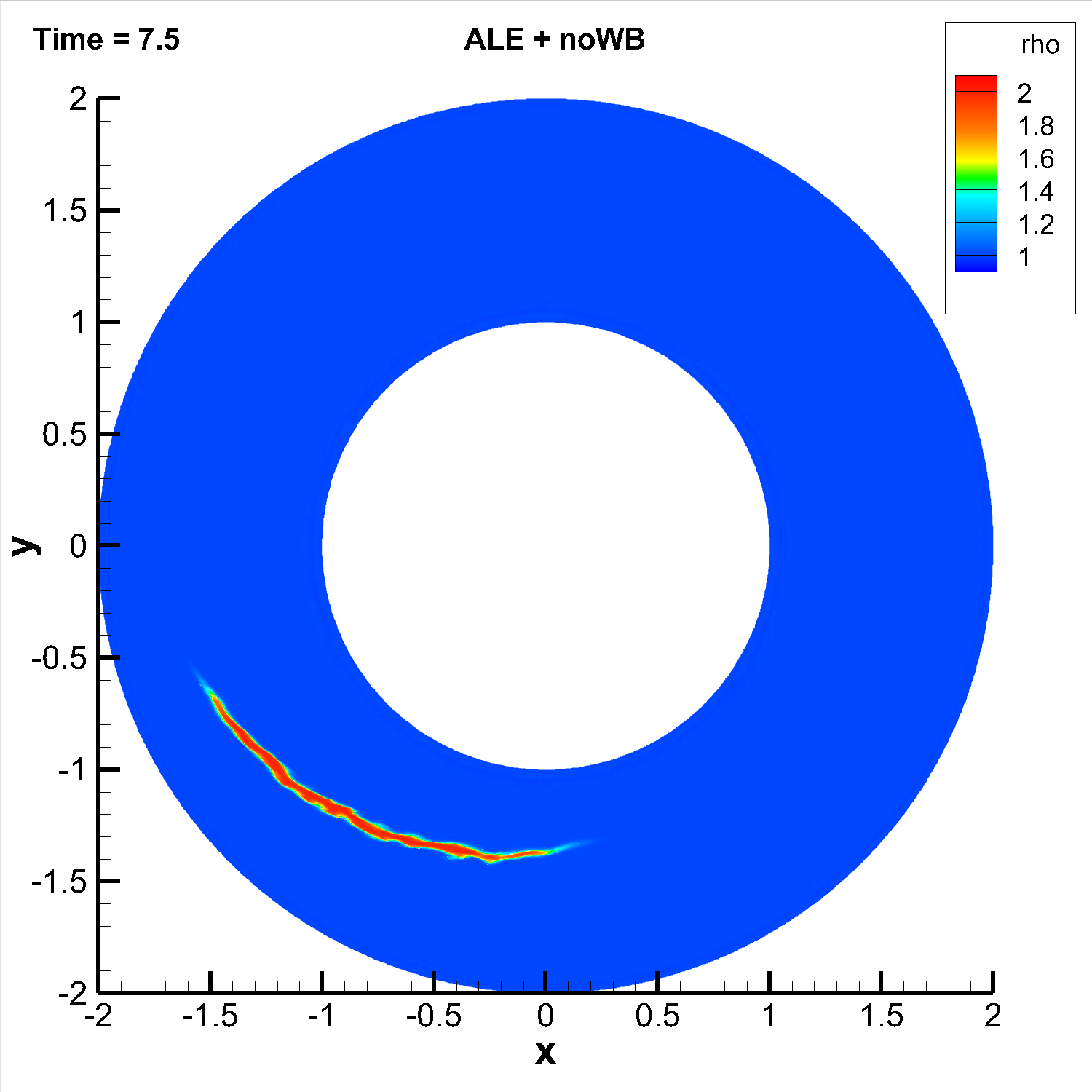}
	\includegraphics[width=0.31\linewidth, trim=8 8 8 8, clip]{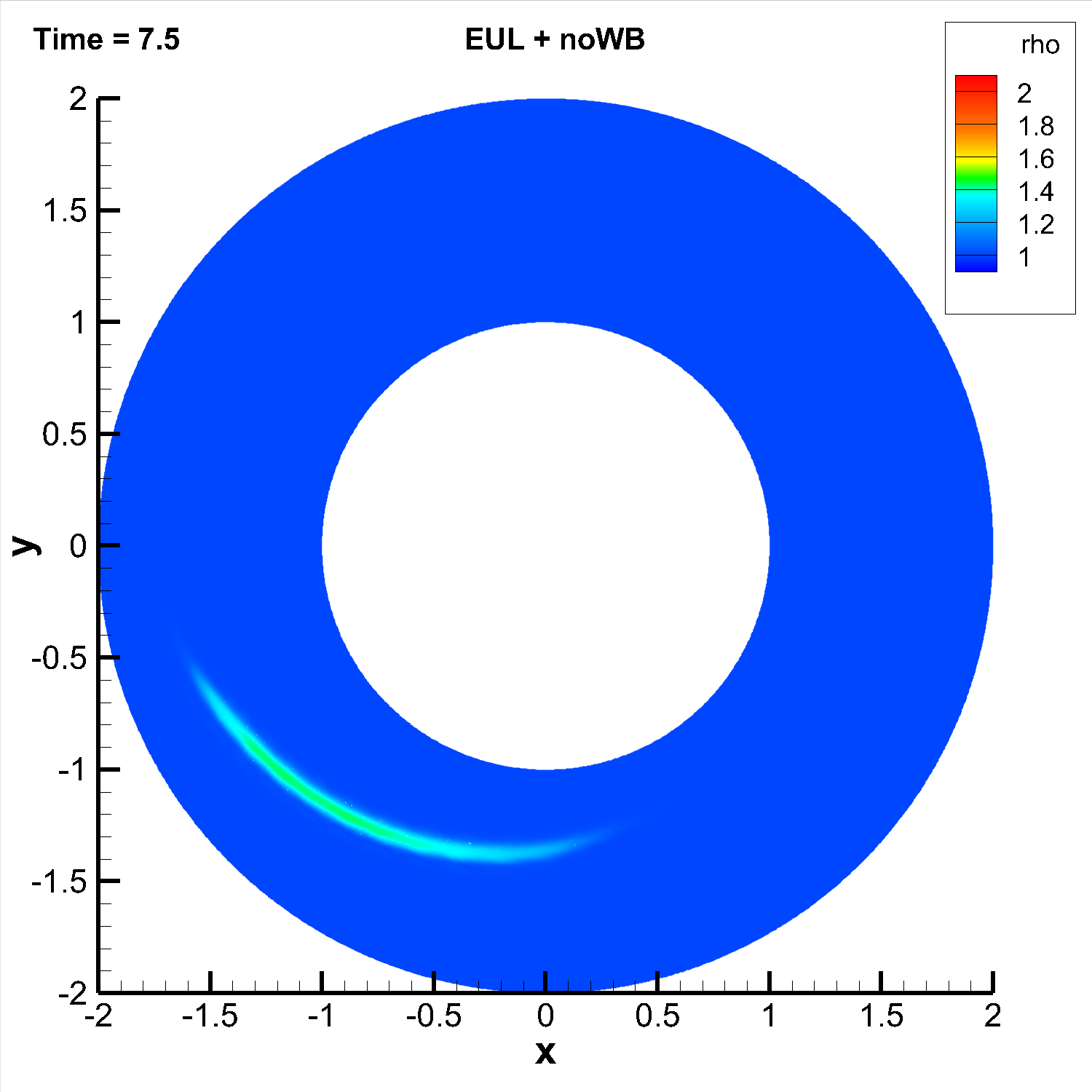}\\[-10pt]
	\caption{Transport of a large density perturbation (of magnitude $1$) on a Keplerian disk. In this figure, 
	we compare the results obtained at time $t= 7.5$ with our WB ALE DG scheme of order $3$ with 
	those obtained with \textit{non} well-balanced Lagrangian and Eulerian schemes of the same order. 
		We can notice that the Lagrangian framework greatly reduces the dissipation w.r.t the 
		Eulerian one. In addition, even on such a large perturbation, the use of well-balancing, 
		albeit not essential, allows to obtain more accurate results.} 
	\label{fig.pallina1000_comparison}
\end{figure}

In particular, we have decided to show and comment the results obtained with 
$A=$~1E-4 (see Figures~\ref{fig.pallina4_film},~\ref{fig.pallina4_numbering},~\ref{fig.pallina4_comparison}), 
$A=$~1e-6 (see Figure~\ref{fig.pallina6_film_comparison})
and $A=$~1 (see Figures~\ref{fig.pallina1000_film} and~\ref{fig.pallina1000_comparison}).


The perturbation of amplitude $A=$~1E-4 is comparable in magnitude with the numerical errors of the 
scheme we are employing (ALE DG $3$ on $2152$ polygonal elements), thus it is appropriate to 
compare the capabilities of standard Eulerian DG schemes, non well-balanced Lagrangian schemes and our approach. 
In Figure~\ref{fig.pallina4_film} we can see that, thanks to the joint beneficial effects of our ALE 
scheme and the well-balanced techniques, the position of the \textit{red} disk with higher density is 
always sharply separated from the background equilibrium and no numerical errors are visible. This 
precise results is due also to the mesh motion that closely follows the fluid flow, as shown in Figure~\ref{fig.pallina4_numbering}.


On the contrary, on such a coarse mesh, of only $2152$ polygonal elements, this simulation is quite 
inaccurate with standard schemes, see Figure~\ref{fig.pallina4_comparison}, because 
they are too dissipative and 
affected by numerical errors of the same amplitude of the 
perturbation we want to model which thus is almost hidden by spurious modes.

It is of course true, that, fixing $A=$~1E-4, we could refine the mesh and simulate at least for a while the mass 
transport with a classical method; 
however, this strategy would be excessively expensive with a smaller perturbation: refer, for example, 
to Figure~\ref{fig.pallina6_film_comparison} where we report the results obtained for \mbox{$A=$~1E-6}. 
Already at the very initial moments of the simulation the numerical errors characterizing the non well-balanced 
schemes completely hide the perturbation we would like to model, while, with our WB ALE approach, 
we easily reach the preset final simulation time. 

Next, we show the results obtained for $A=$~1 in Figures~\ref{fig.pallina1000_film} and~\ref{fig.pallina1000_comparison}.
Here, the effects of the well-balancing are less evident because the perturbation has a higher 
magnitude with respect to the numerical errors of the employed schemes.  
However, it is still possible to appreciate the absence of numerical errors in the well-balanced case 
and even more so the role of the Lagrangian mesh displacement in the reduction of the convection 
errors. In particular, we would like to emphasize that our moving mesh technique, which makes use of
topology changes and treat them with high order of accuracy, allows to always maintain a high quality mesh 
even on a vortical velocity field studied over a long simulation time.

\subsubsection{Keplerian disk with steep density gradient and Kelvin-Helmholtz instabilities}
\label{ssec.KHinstabilities}

\begin{table}[!tp] 
	\centering
	\numerikNine
	\begin{tabular}{|l|l||cccccc|}
		\hline
		\multicolumn{8}{|c|}{ Keplerian disk with steep density gradient } \\
		\hline
		\hline
		& Time $t = $ & $0.0$ & $0.5$ & $1.0$ & $2.5$ & $5.0$ & $10.0$ \\
		\hline
		\hline
		\multirow{4}{*}{\rotatebox{90}{{DG-$\mathbb{P}_1$}}}
		& $L_2(\rho)$ error & 1.6393E-12 & 1.5792E-12	& 1.5546E-12	& 1.5374E-12	& 1.5260E-12	& 1.5140E-12 \\
		& $L_2(p)$ error  	& 0.0        & 5.3530E-13	& 5.2949E-13	& 5.2403E-13	& 5.1632E-13	& 5.0819E-13 \\
		& No. timestep 	    & 0     	 & 186       	& 349       	& 825       	& 1611	        & 3184       \\
		& No. sliver   	    & 0     	 & 71       	& 225        	& 707       	& 1507	        & 3095       \\
		\cline{2-8}
		\hline
		\hline
		\multirow{4}{*}{\rotatebox{90}{{DG-$\mathbb{P}_2$}}}
		& $L_2(\rho)$ error & 1.5863E-12 & 1.5834E-12	& 1.6313E-12	& 1.7154E-12	& 1.7824E-12	& 1.4623E-12 \\
		& $L_2(p)$ error  	& 0.0        & 7.0141E-13	& 1.0383E-13	& 1.1404E-13	& 7.3145E-13	& 5.2780E-13 \\
		& No. timestep 	    & 0     	 & 344       	& 655       	& 1580       	& 3080	        & 6448       \\
		& No. sliver   	    & 0     	 & 69       	& 205        	& 660       	& 1497	        & 3250      \\
		\cline{2-8}
		\hline
		\hline
		\multirow{4}{*}{\rotatebox{90}{{DG-$\mathbb{P}_3$}}}
		& $L_2(\rho)$ error & 1.5882E-12 & 1.6108E-12	& 1.6991E-12	& 3.5427E-12	& 5.0623E-12	& 3.3474E-12 \\
		& $L_2(p)$ error  	& 0.0        & 7.8564E-13	& 1.2285E-12	& 1.9266E-12	& 3.9819E-12	& 5.3318E-12 \\
		& No. timestep 	    & 0     	 & 576       	& 1104       	& 2580     	    & 5042          & 9940      \\
		& No. sliver   	    & 0     	 & 56       	& 191        	& 667        	& 1457          & 3044       \\
		\cline{2-8}
		\hline
		\hline
		\multirow{4}{*}{\rotatebox{90}{{DG-$\mathbb{P}_4$}}}
		& $L_2(\rho)$ error & 1.5857E-12 & 1.7407E-11	& 2.2113E-10	& 5.8511E-10	    & 4.7131E-10	    & 2.5286E-09 \\
		& $L_2(p)$ error  	& 0.0        & 9.8282E-12	& 1.8642E-09	& 2.9620E-09	    & 9.4133E-08	    & 1.2181E-09 \\
		& No. timestep 	    & 0     	 & 828         	& 1591         	& 3753      		& 7324	   		    & 14805        \\
		& No. sliver   	    & 0     	 & 62      	    & 191        	& 669        		& 1455	      	    & 3092         \\
		\cline{2-8}
		\hline
	\end{tabular}
	\caption{Check of the well-balanced property on a Keplerian disk with a steep density gradient. 
	As in the previous test cases, we can notice that the equilibrium solution, even if initially perturbed, is 
	preserved with machine accuracy for very long simulation times and after handling thousands of sliver elements.}
	\label{table.KeplerianDisk_sharp}
\end{table}

In this section, we consider another equilibrium solution belonging to the family~\eqref{eq.EquilibriaConstraint} but with
a sharp density gradient for $r\rightarrow 1.5$
\begin{equation} 
	\label{eq.equilibriumForKH1}
	\begin{cases} 
		\rho_E = \rho_0 + \rho_1 \text{tanh}\left( \frac{r-r_m}{\sigma} \right),  \\
		u_{r E} = 0, \\
		u_{\phi E} = \sqrt{ \frac{Gm_s}{r} }, \\
		p_E = 1, 
	\end{cases} 
\end{equation} 
with $\rho_0 = 1$, $\rho_1 = 0.25$, $r_m= 1.5$ and $ \sigma = 0.01$. 

%
\begin{figure}[!tb] 
	\centering
	\includegraphics[width=0.31\linewidth, trim=8 8 8 8, clip]{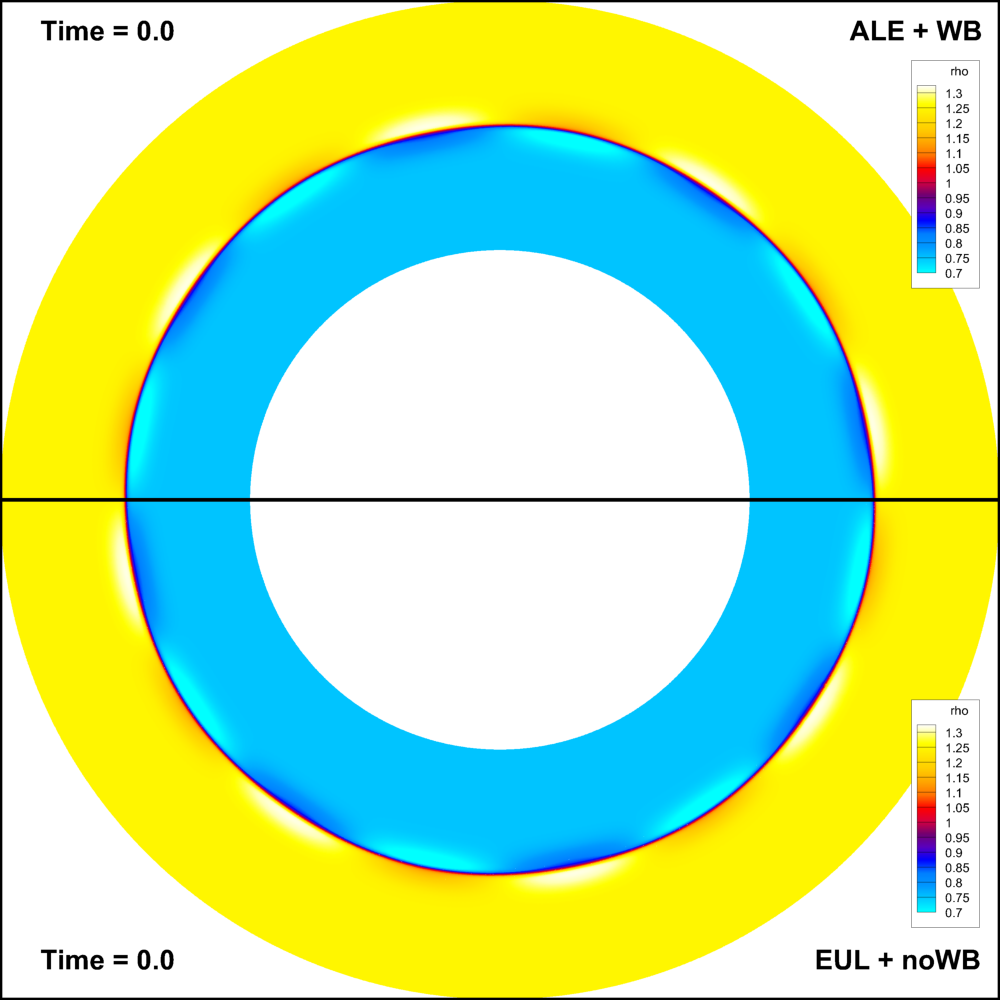}\ \,%
	\includegraphics[width=0.31\linewidth, trim=8 8 8 8, clip]{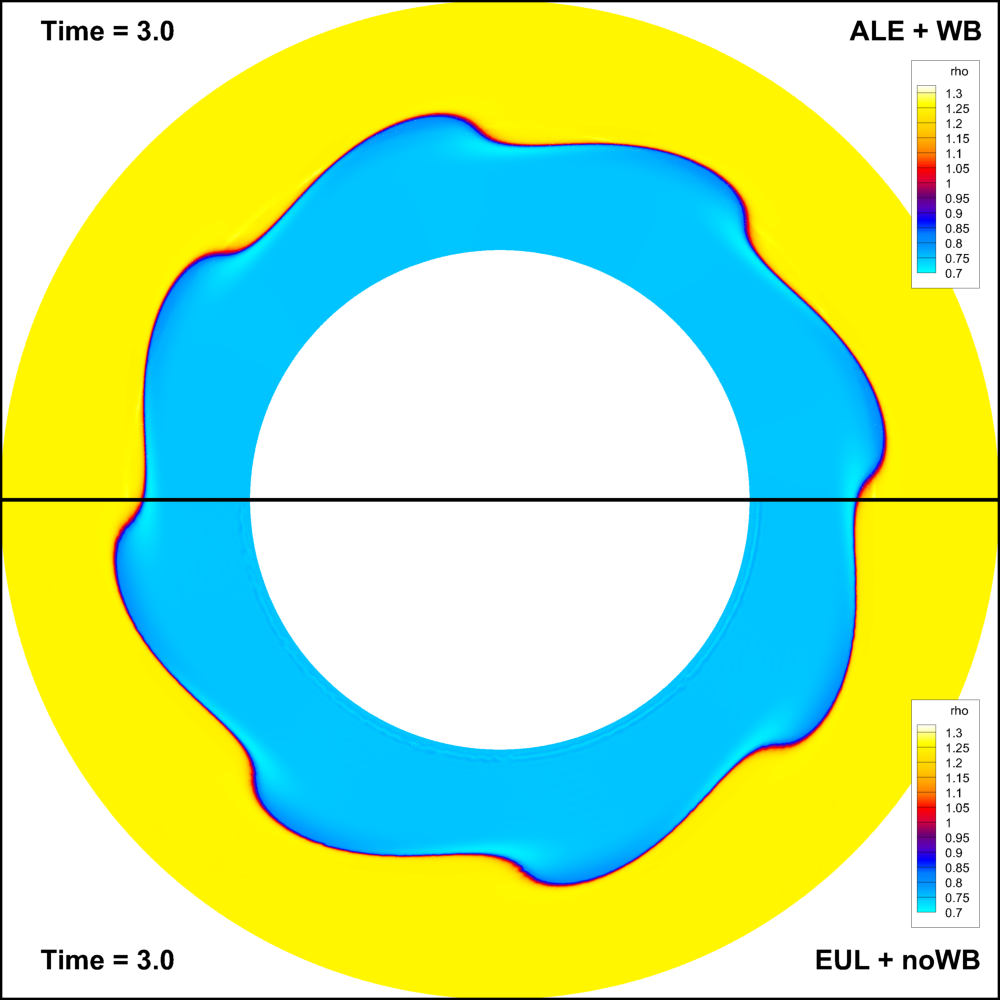}\ \,%
	\includegraphics[width=0.31\linewidth, trim=8 8 8 8, clip]{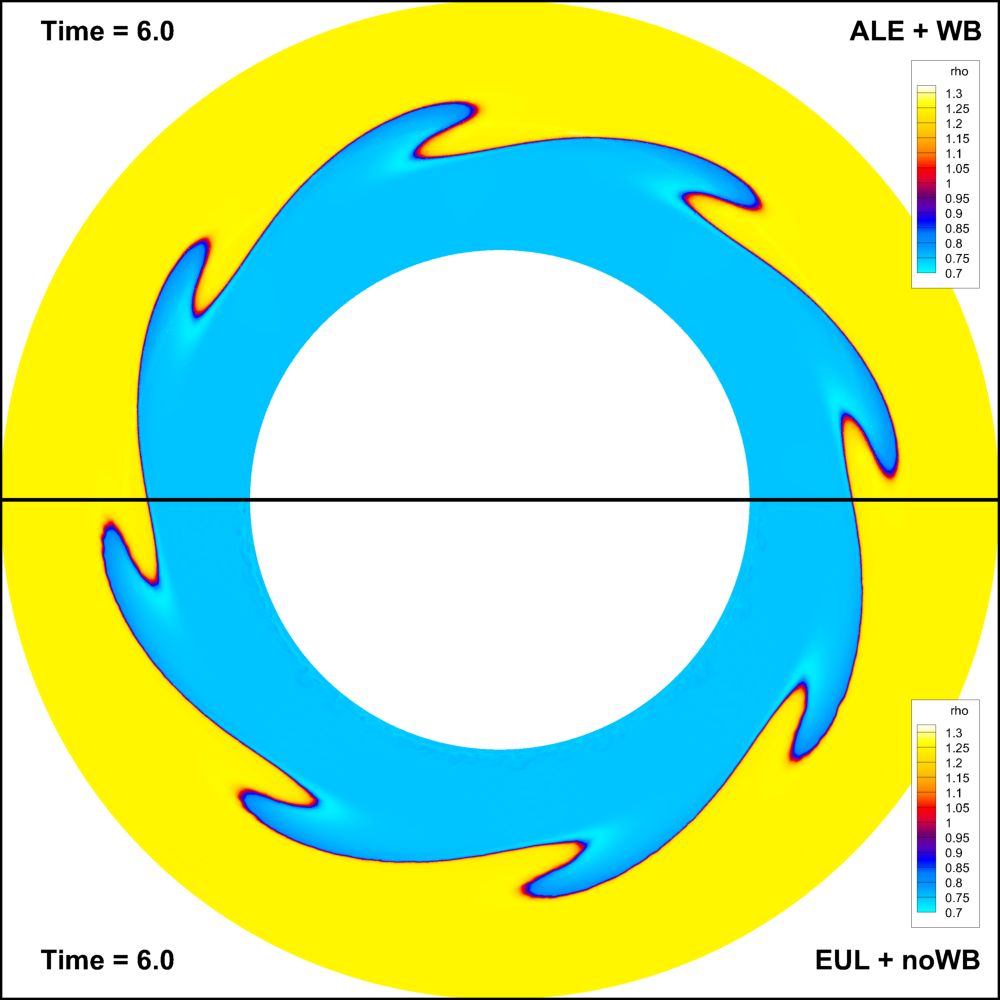}\\[10pt]%
	\includegraphics[width=0.31\linewidth, trim=8 8 8 8, clip]{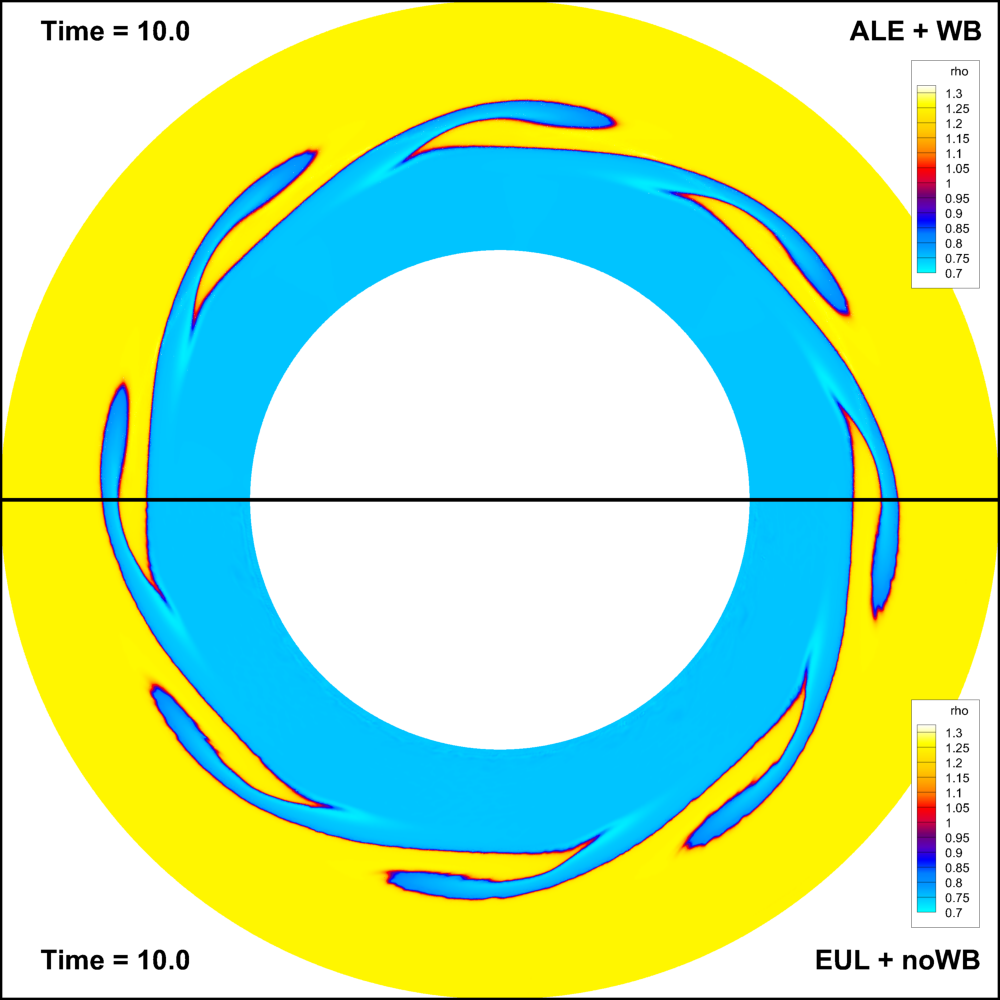}\ \,%
	\includegraphics[width=0.31\linewidth, trim=8 8 8 8, clip]{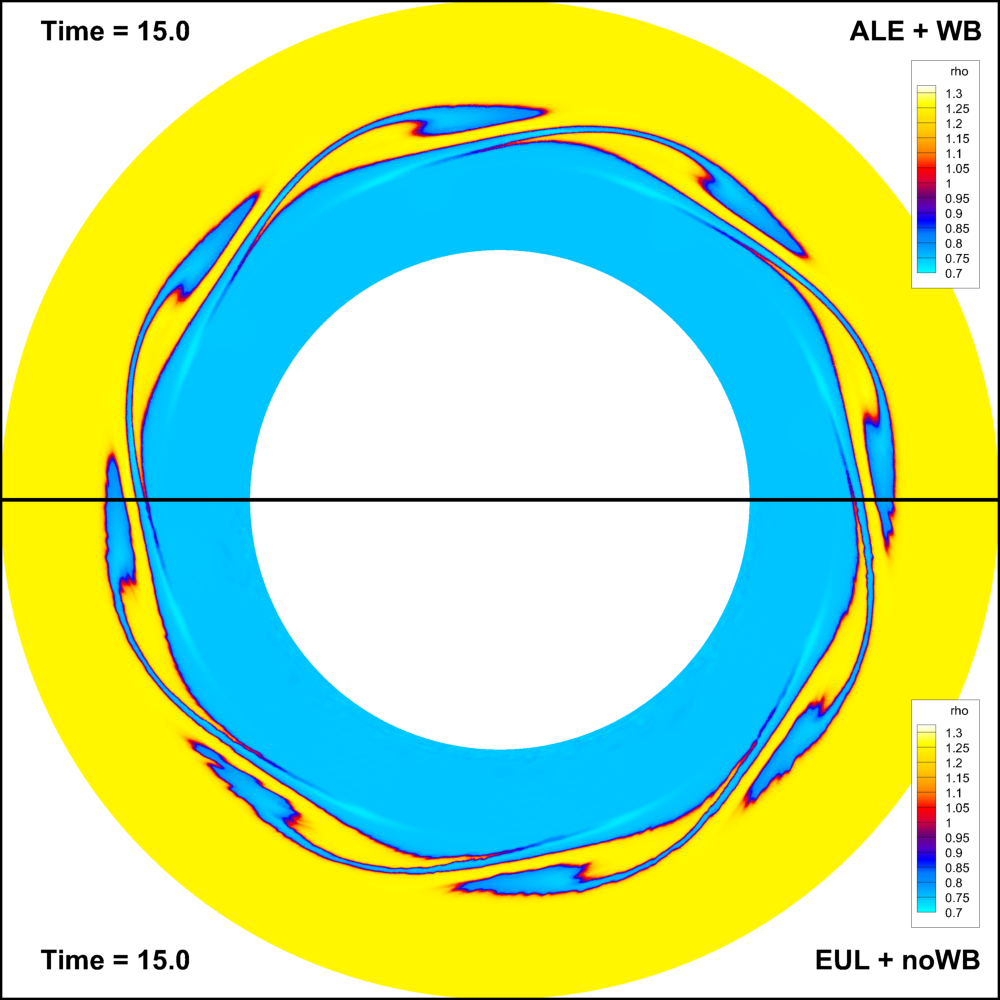}\ \,%
	\includegraphics[width=0.31\linewidth, trim=8 8 8 8, clip]{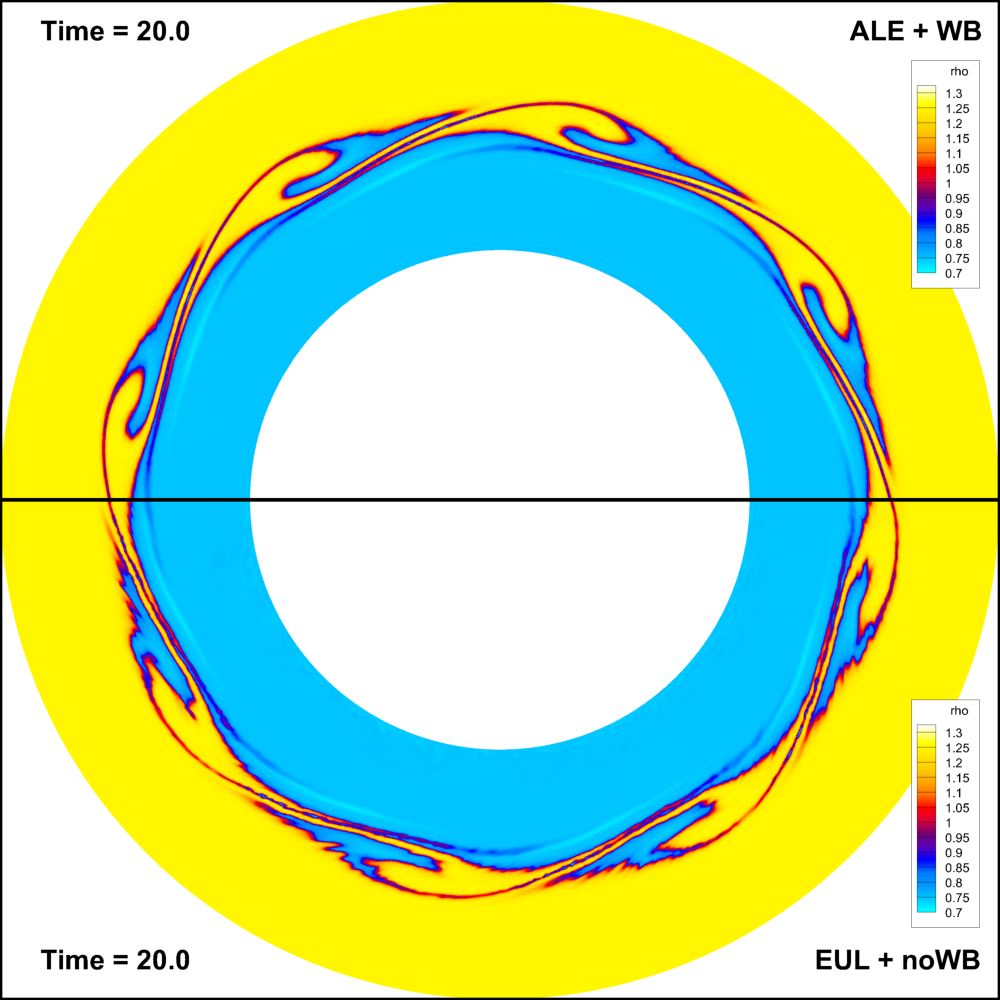}\\[10pt]%
	\includegraphics[width=0.31\linewidth, trim=8 8 8 8, clip]{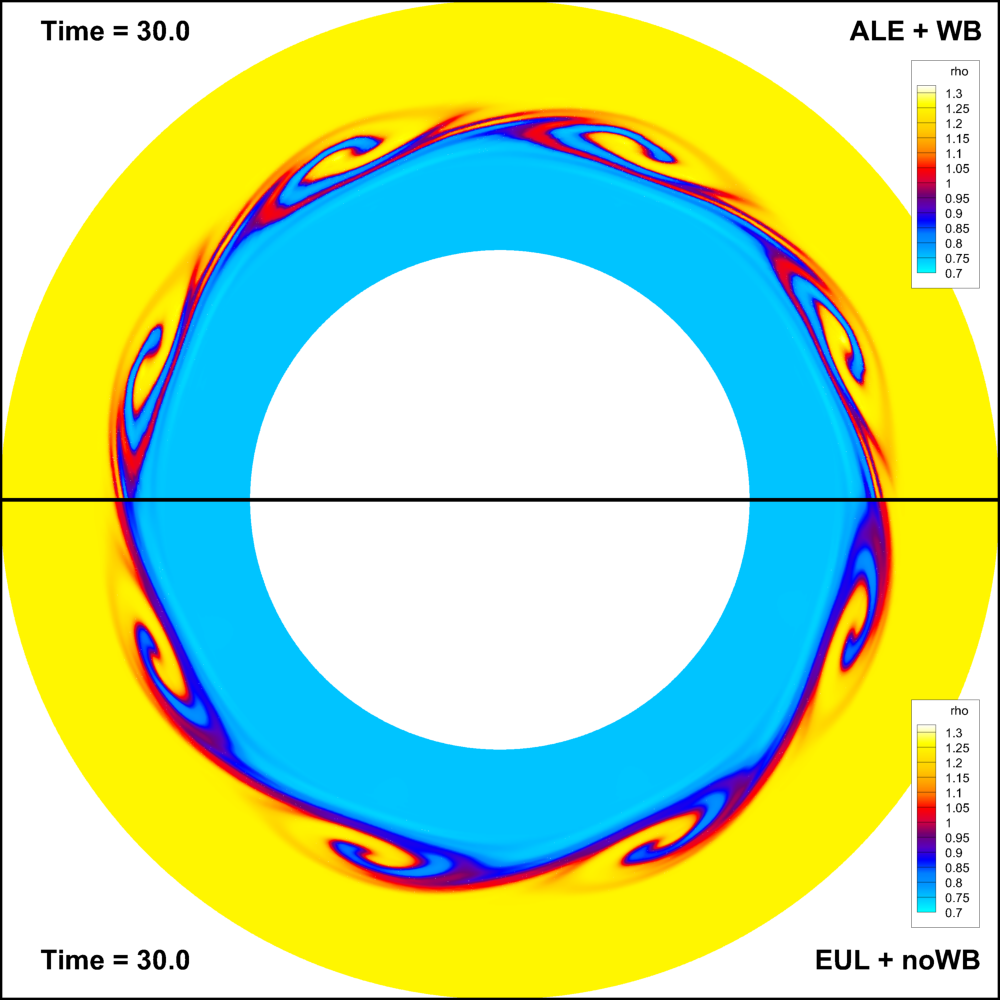}\ \,%
	\includegraphics[width=0.31\linewidth, trim=8 8 8 8, clip]{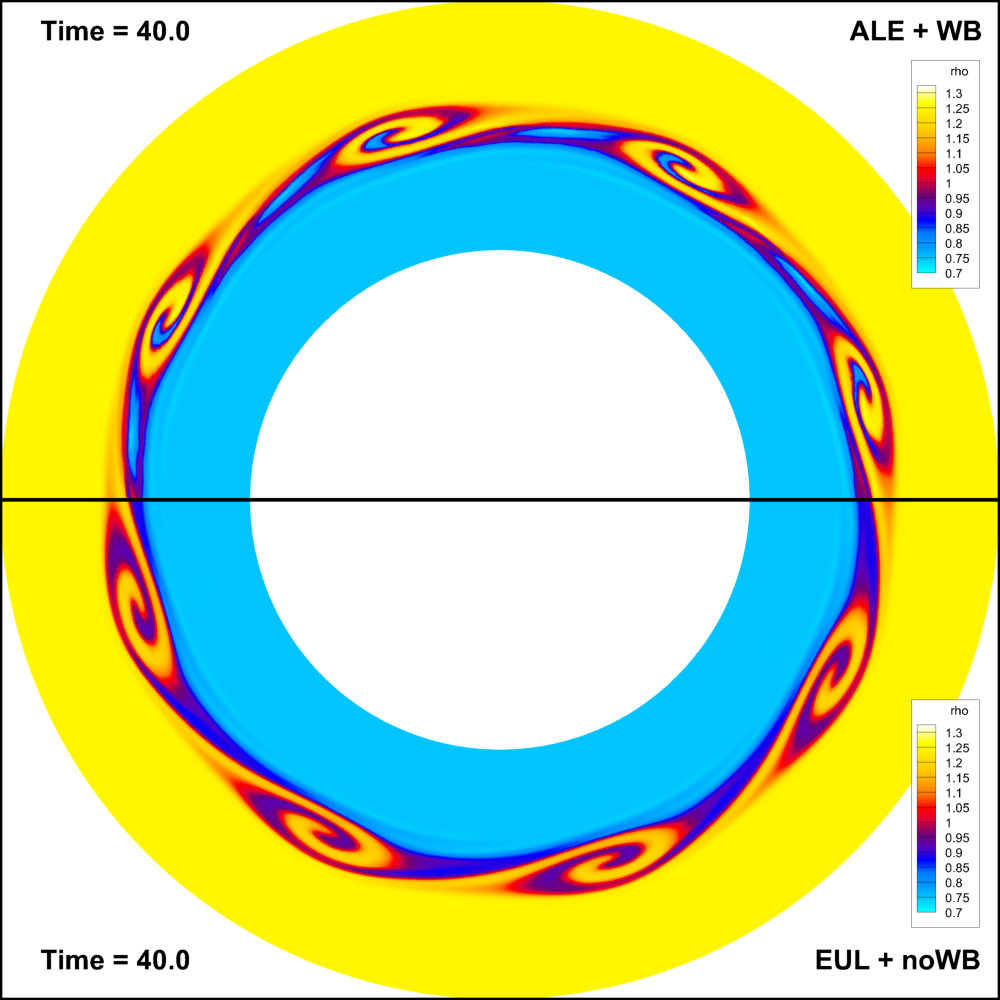}\ \,%
	\includegraphics[width=0.31\linewidth, trim=8 8 8 8, clip]{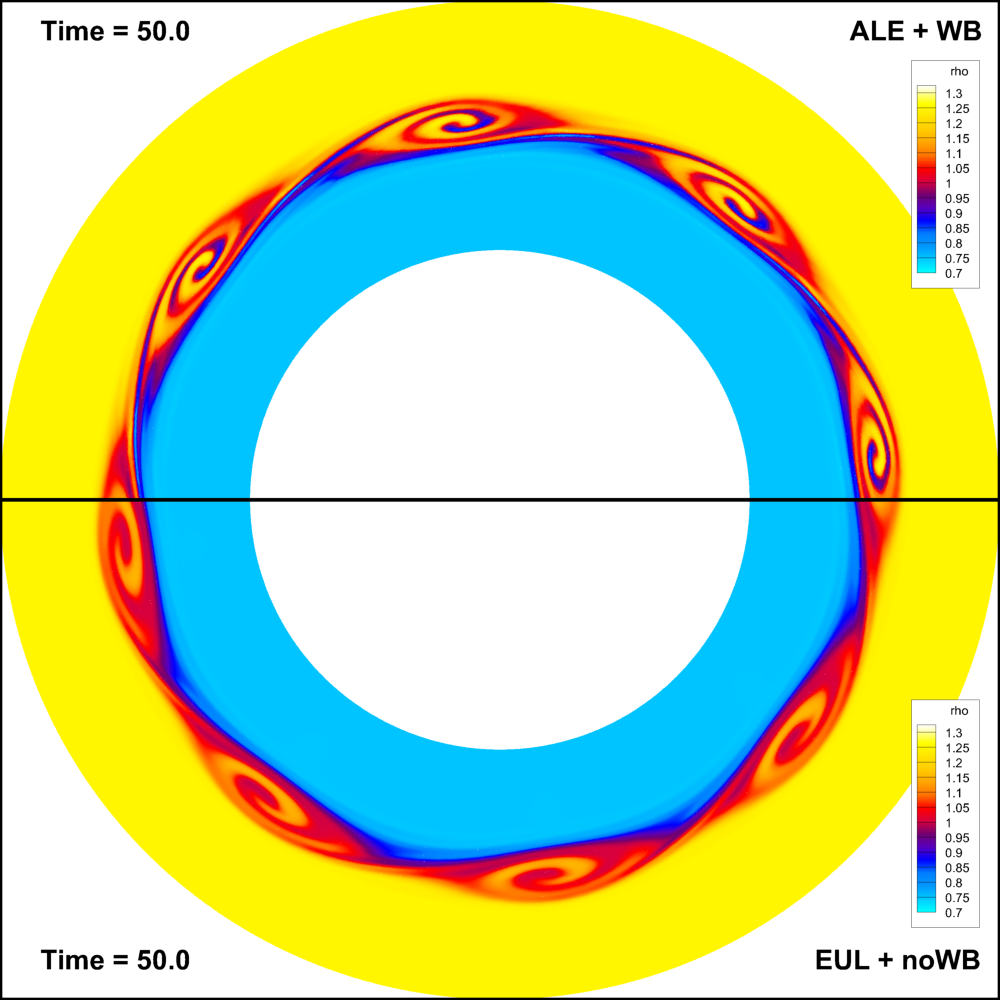}
	\caption{Kelvin-Helmholtz instabilities on a Keplerian disk. We compare here the results obtained with our WB 
	ALE DG scheme of order $3$ with those given by a standard Eulerian scheme of the same order of 
	accuracy on a mesh of $68953$ polygonal elements from time $t=0$ to time $t=50$.}
	\label{fig.KH_film}
\end{figure}
\begin{figure}[!tb] 
	\centering
	\includegraphics[width=0.485\linewidth, trim=8 8 8 8, clip]{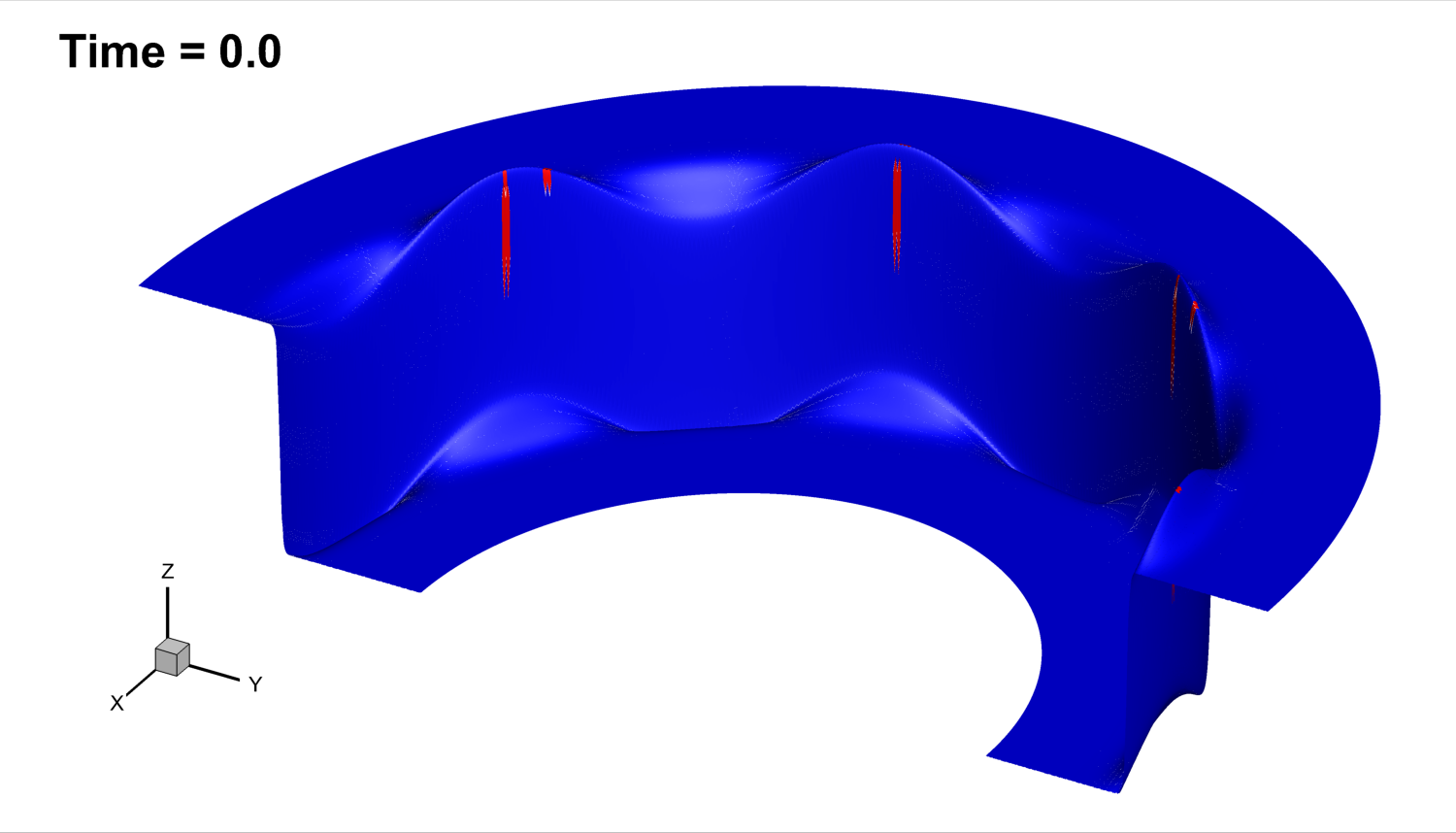}\ 
	\includegraphics[width=0.455\linewidth, trim=8 8 8 8, clip]{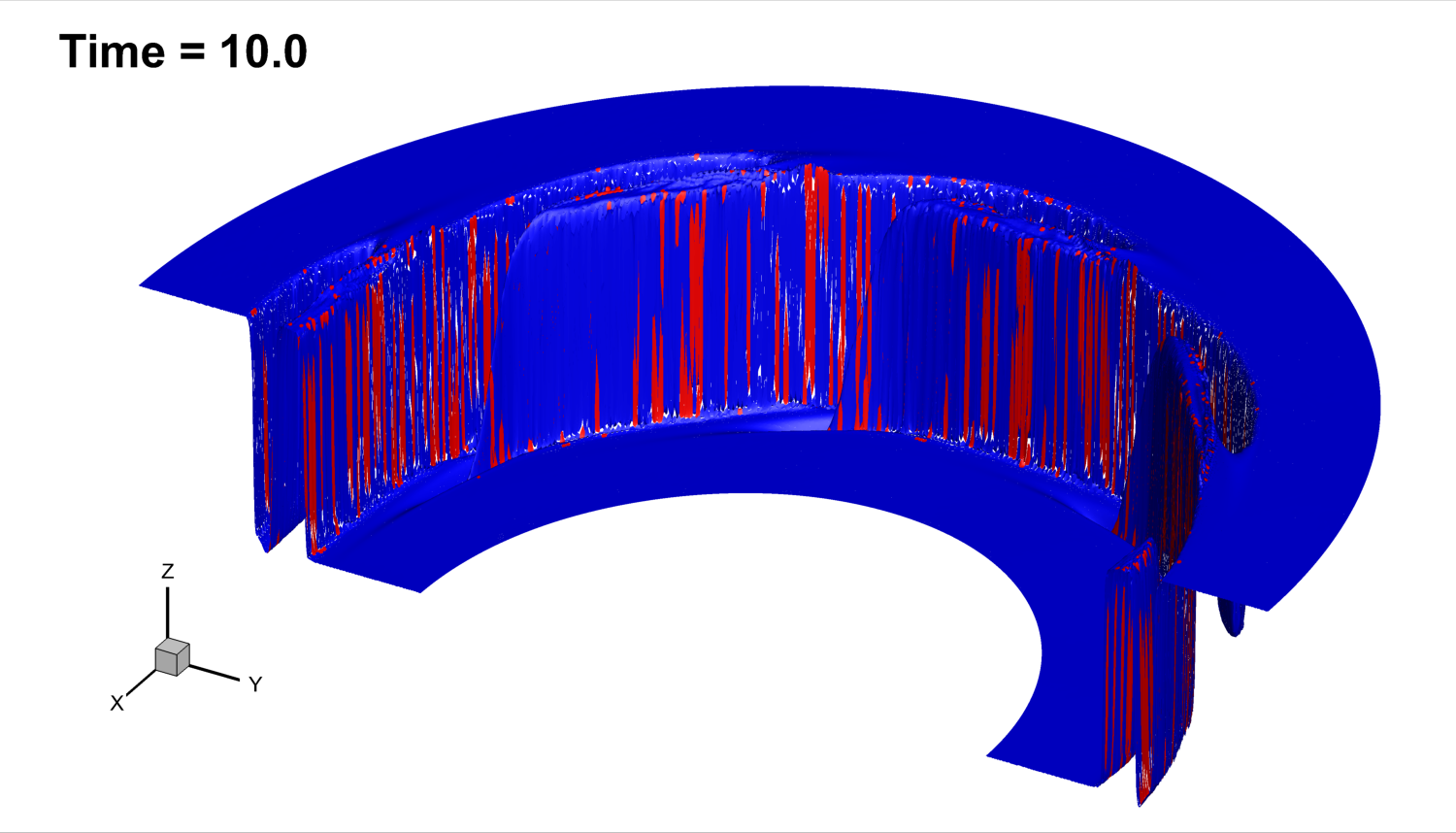}\\[10pt]
	\includegraphics[width=0.485\linewidth, trim=8 8 8 8, clip]{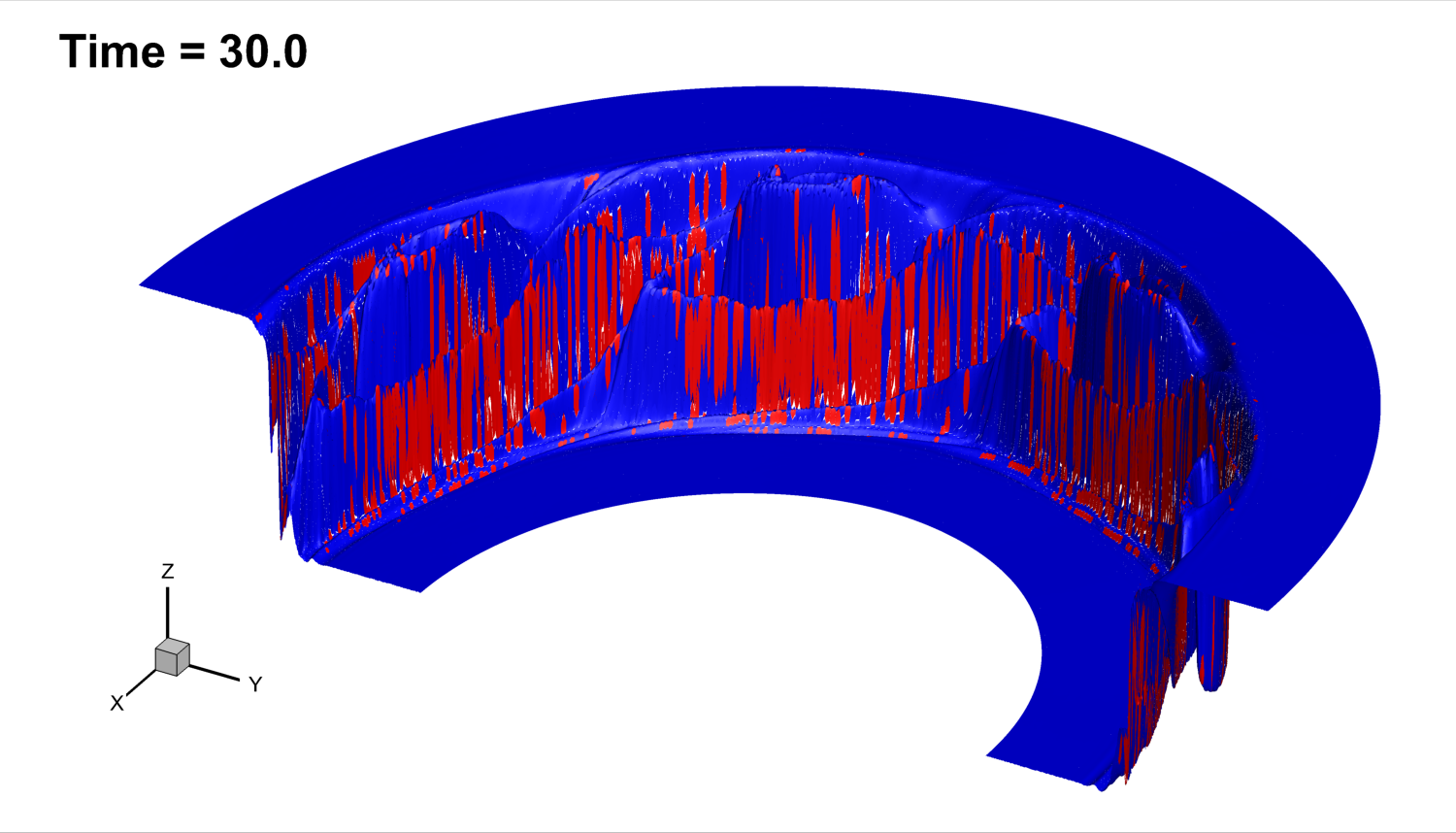}\ 
	\includegraphics[width=0.485\linewidth, trim=8 8 8 8, clip]{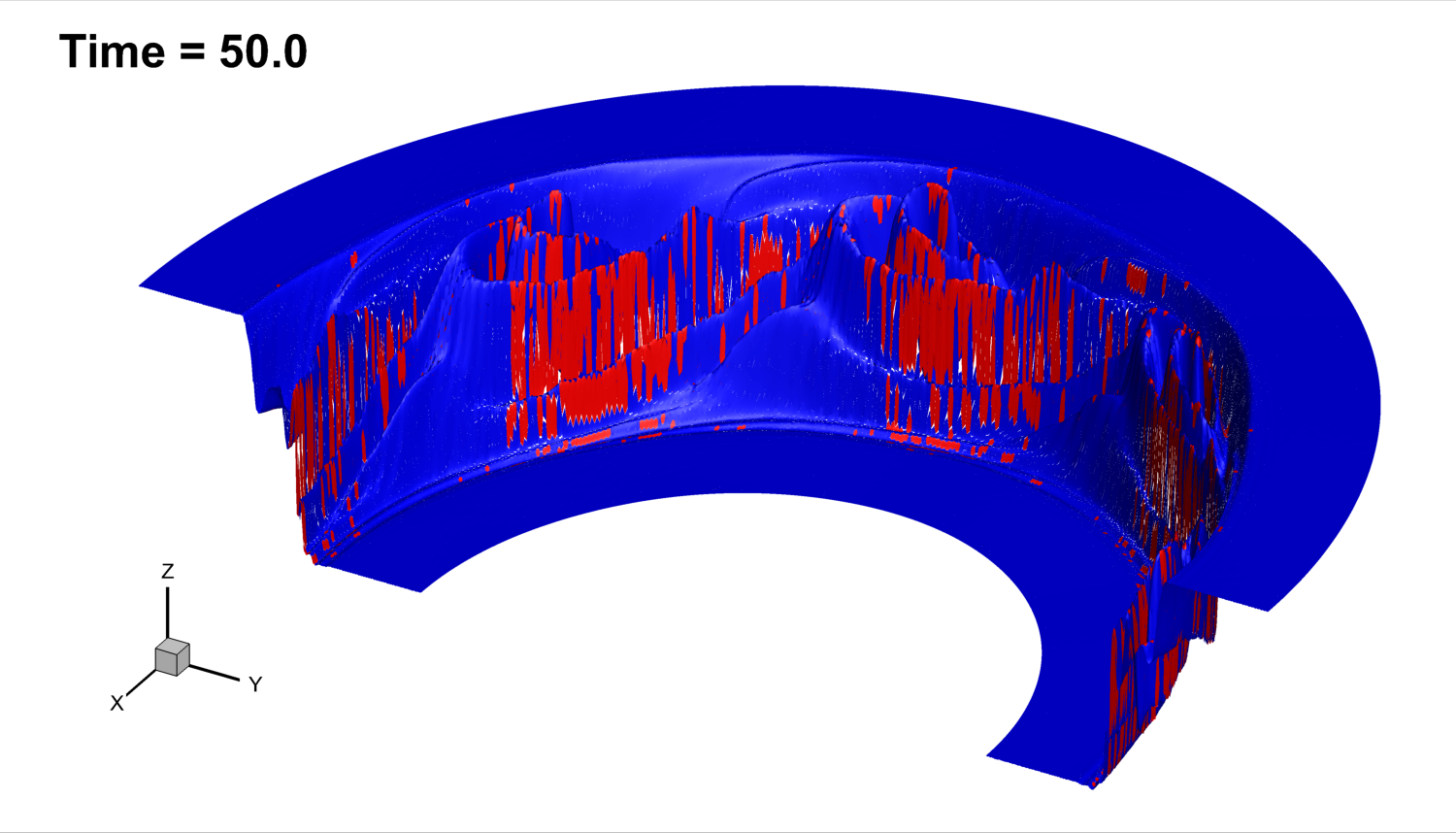}%
	\caption{Kelvin-Helmholtz instabilities on a Keplerian disk. 
		We show here, in red, for the simulation performed with the WB ALE scheme,  
		the cells on which the limiter has been activated at different times. 
		We highlight that  also the finite volume subcell scheme 
		is endowed with well-balanced techniques.}
	\label{fig.KH_limiter}
\end{figure}
%
%

After having verified the well-balanced property, see Table~\ref{table.KeplerianDisk_sharp}, 
we study the evolution of the following initial condition obtained by adding a sinusoidal perturbation 
to the equilibrium profile~\eqref{eq.equilibriumForKH1}
\begin{equation}
	\begin{cases} 
		\rho = \rho_E + A \rho_0  \sin(k \varphi) \text{exp} \left( - \frac{(r - r_m)^2}{ s}   \right ), \\
		u = u_E +  A \sin(k \varphi) \text{exp} \left( - \frac{(r - r_m)^2}{ s }   \right ), \quad  
		v = v_E, \\
		P = P_E + A \sin(k \varphi) \text{exp} \left( - \frac{(r - r_m)^2}{s}   \right ), \\
	\end{cases} 
\end{equation} 
with $ A= 0.1$, $k=8$ and $ s = 0.005$. 
In this flow configuration, with perturbations affecting all the flow quantities, we observe the development of 
Kelvin-Helmholtz instabilities. We report the results obtained with our WB ALE DG scheme and with a 
standard Eulerian scheme, both of order $3$ in Figure~\ref{fig.KH_film}, where we can observe 
the increased resolution of our approach even on a large-amplitude phenomenon (i.e. not due only to a quite 
small initial perturbation, being $A=0.1$).  
Moreover, we make use of this test case to highlight that our scheme is also endowed with 
an \textit{a posteriori} subcell FV limiter, 
whose activations are shown in Figure~\ref{fig.KH_limiter}. We remark that also the subcell limiter
scheme
is equipped with well-balanced techniques, so to not invalidate the benefits given by the 
well-balanced ALE DG method.

\subsection{Frontogenesis on a MHD vortex}
\label{ssec.mhdfrontogenesis}

In this section, we first verify the well-balanced property of our approach on the MHD equations. 
For this purpose, we consider again the vortical stationary solution given in~\eqref{eqn.mhd3d.ic1} and 
we set it up with $\epsilon = \epsilon_E = 5$ both for $\Q_{\text{IC}}$ and $\Q_E$ on a moving 
mesh of $1345$ polygonal elements.
Then, we slightly perturb $\Q_{\text{IC}}$ with a random error of order 1E-12 added everywhere on the 
density profile and we monitor the error evolution in Table~\ref{table.MHDvortex}. 

Once verified that the property is satisfied, we employ our scheme to study the kinematic frontogenesis, 
which is a benchmark arising from the field of computational meteorology~\cite{davies1985comments,cavalcanti2015conservative}, 
usually studied for linearized equations for which also exact solutions are available~\cite{toro2005ader,dumbser2007arbitrary}.
This test case is usually employed also to verify the robustness of moving mesh methods 
because it features a velocity field that strongly stretches the grid and which should be 
well resolved to correctly follow the interface evolution. 

\begin{figure}[!tb]  \centering
	\includegraphics[width=0.33\linewidth, trim=8 8 8 8, clip]{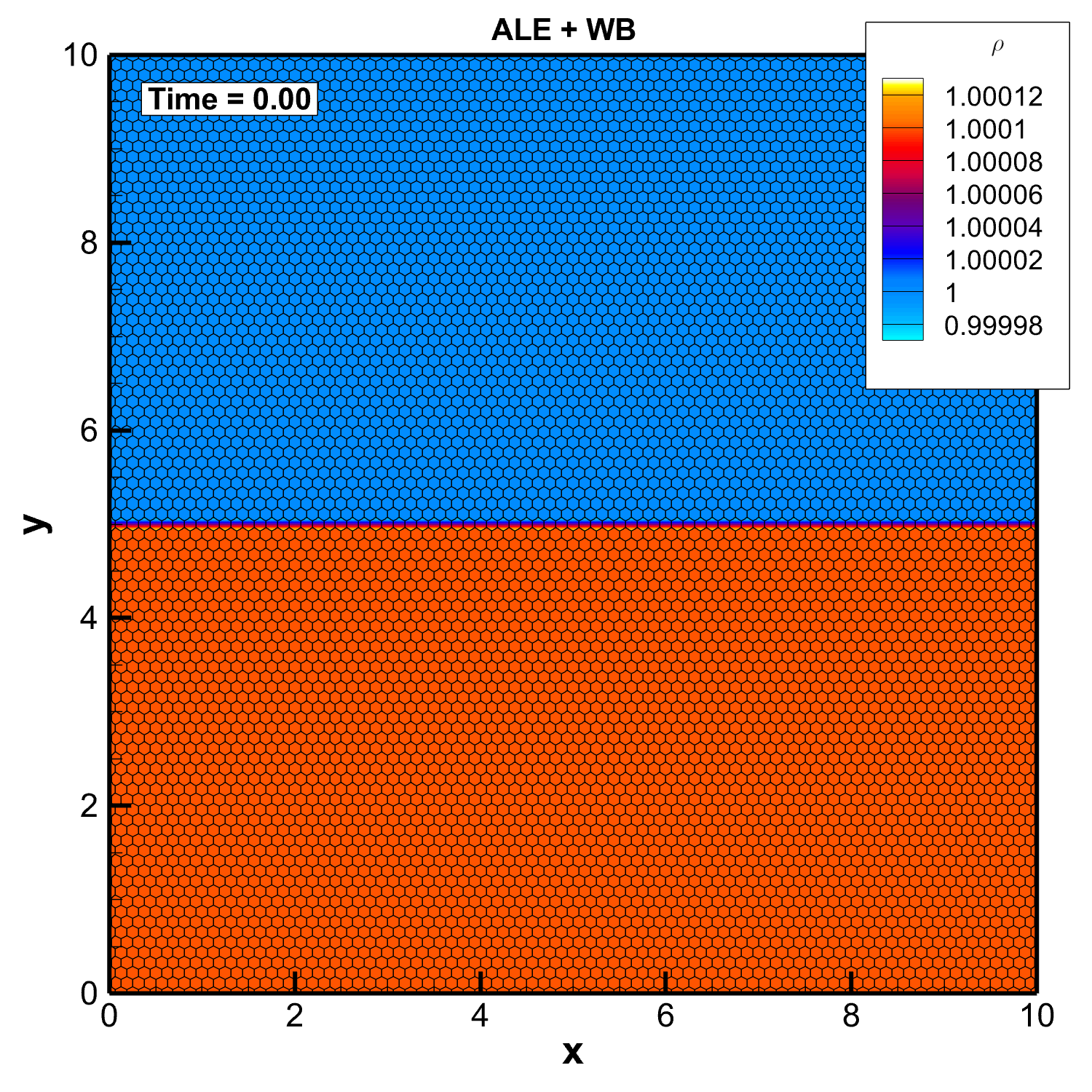}%
	\includegraphics[width=0.33\linewidth, trim=8 8 8 8, clip]{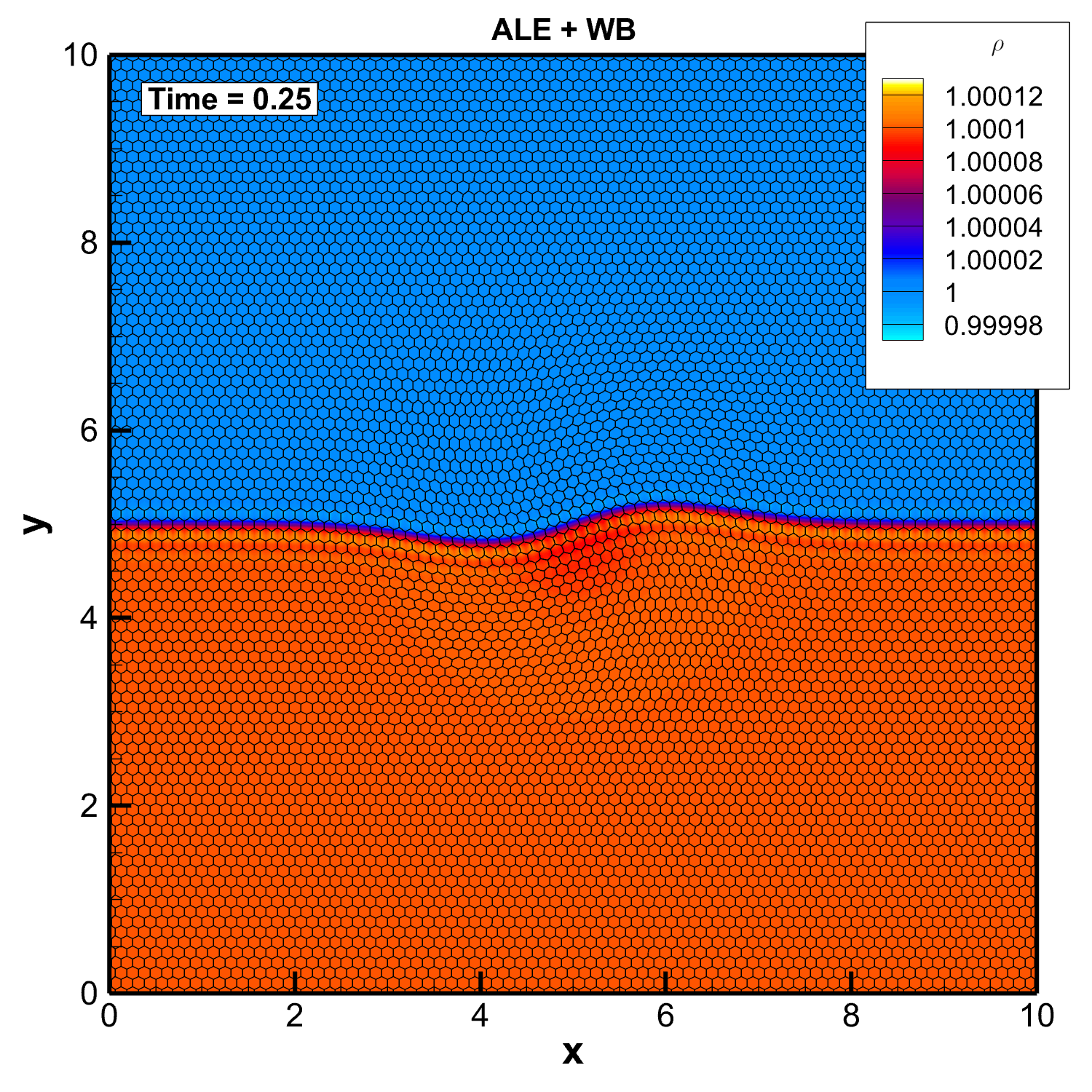}%
	\includegraphics[width=0.33\linewidth, trim=8 8 8 8, clip]{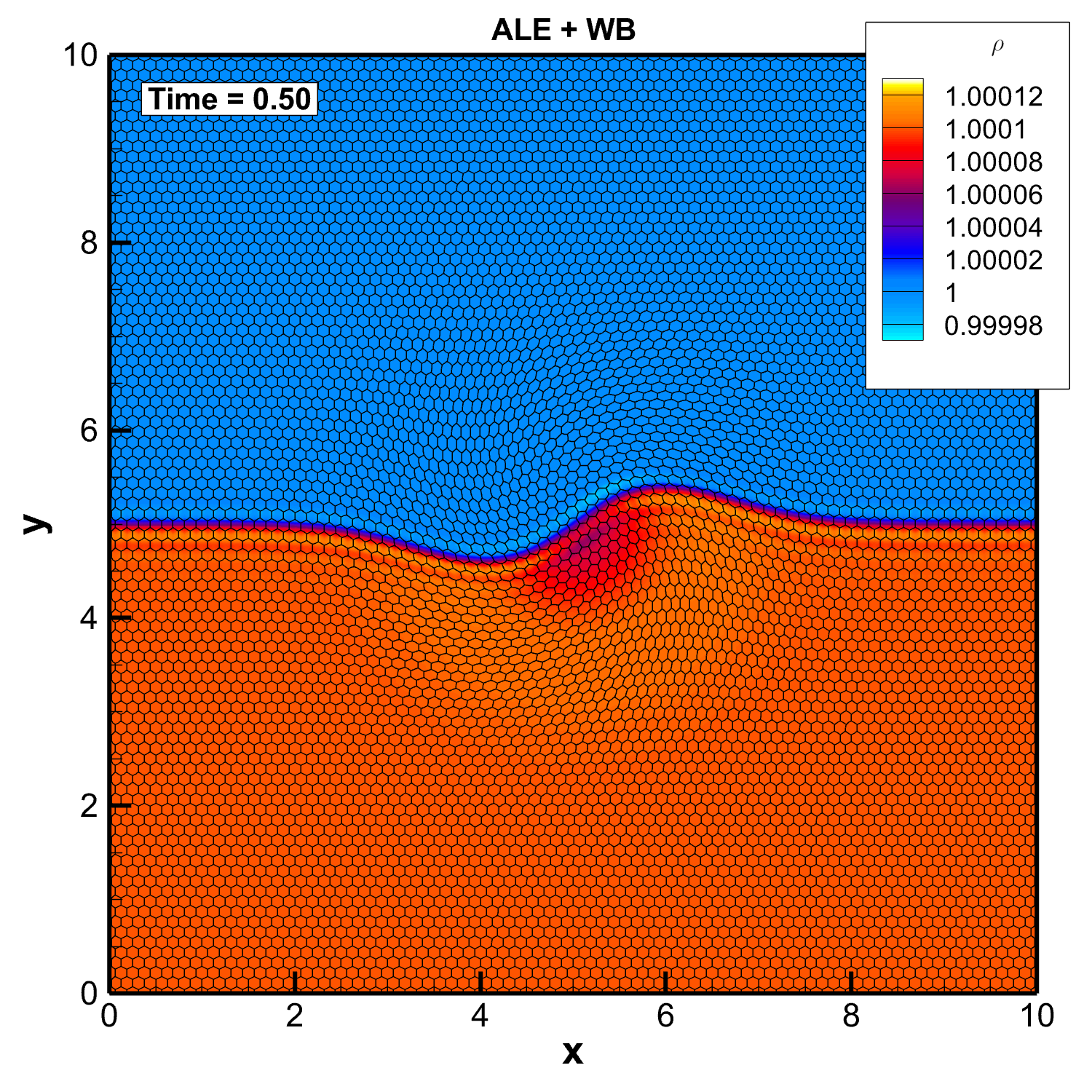}\\[0pt]
	\includegraphics[width=0.33\linewidth, trim=8 8 8 8, clip]{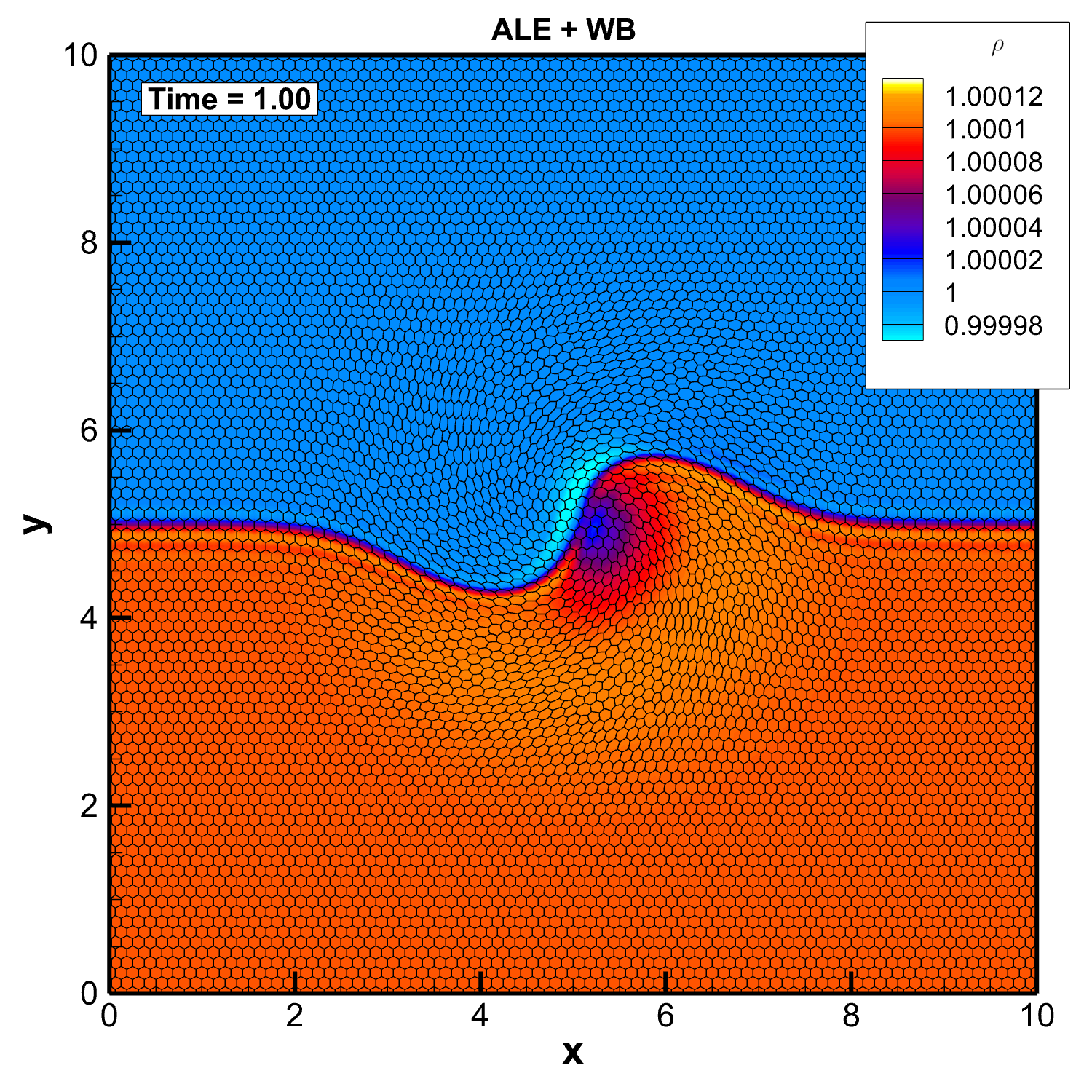}%
	\includegraphics[width=0.33\linewidth, trim=8 8 8 8, clip]{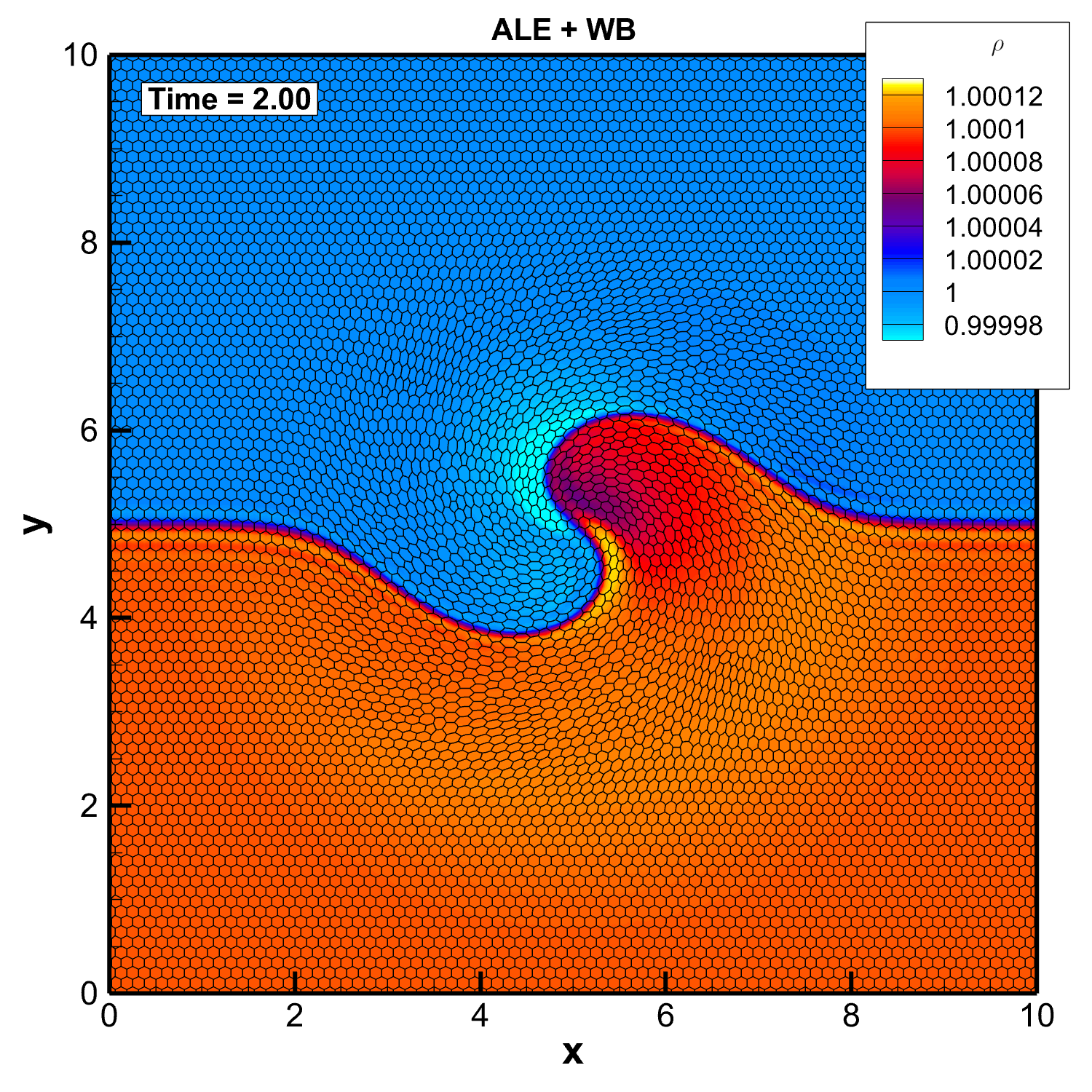}%
	\includegraphics[width=0.33\linewidth, trim=8 8 8 8, clip]{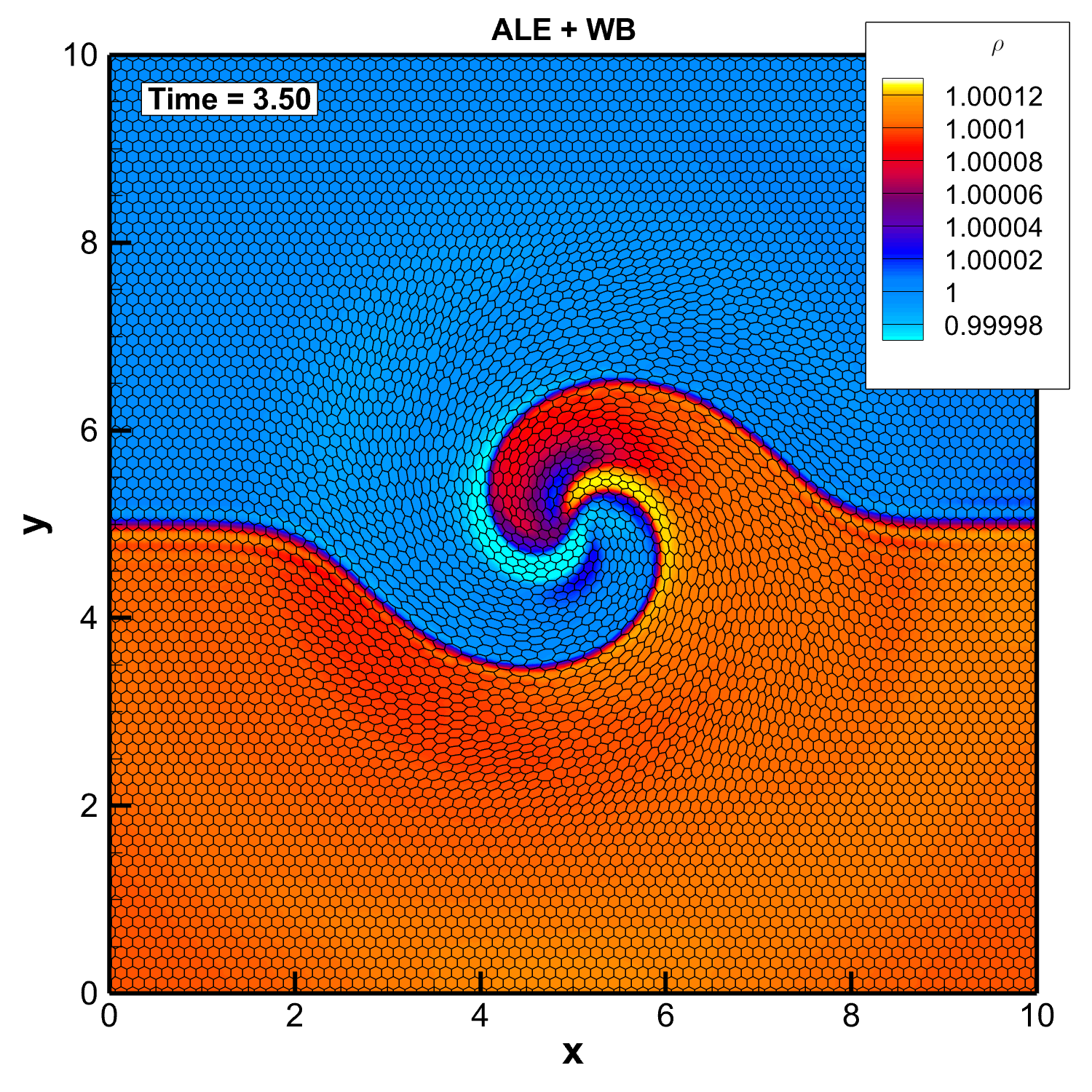}\\[0pt]
	\includegraphics[width=0.33\linewidth, trim=8 8 8 8, clip]{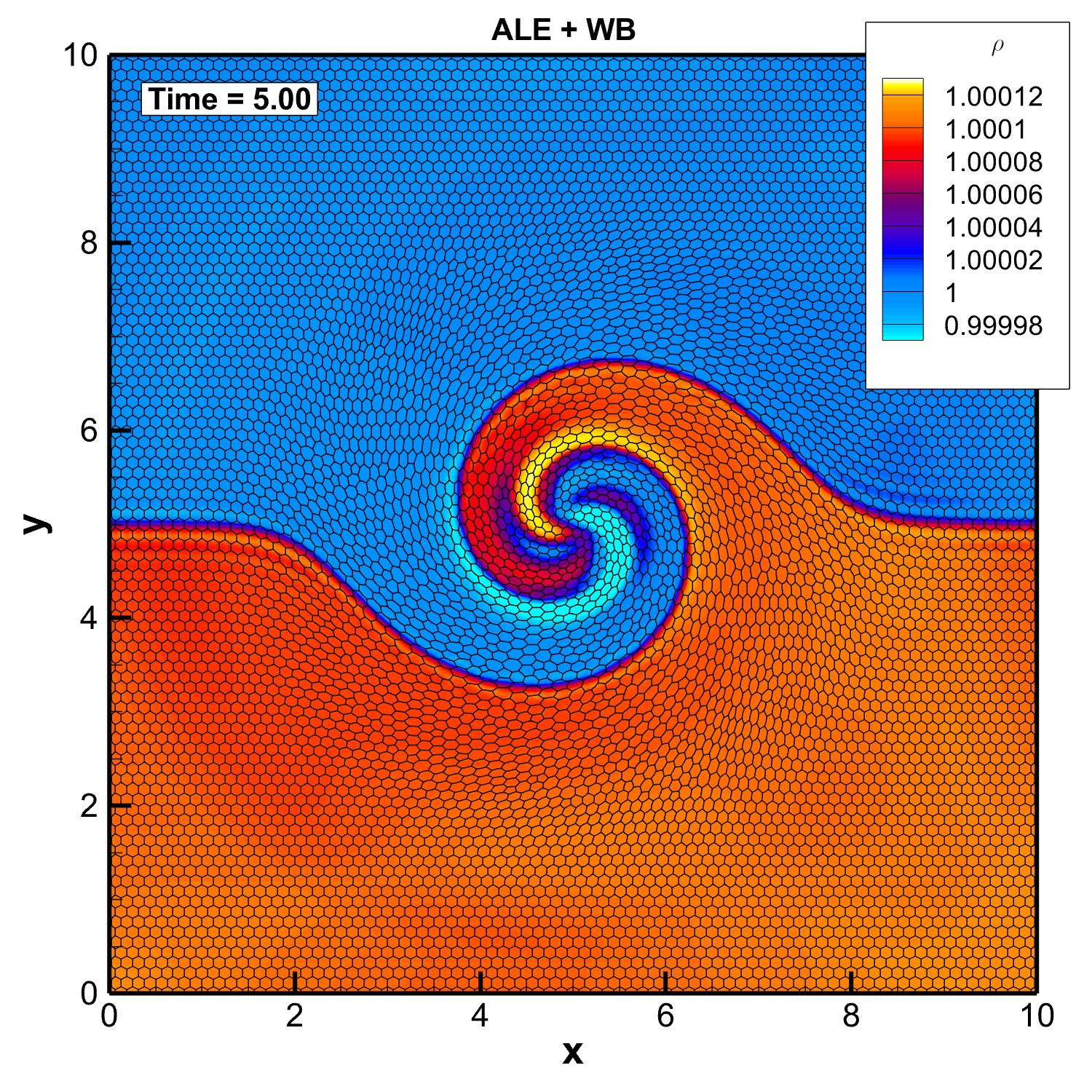}%
	\includegraphics[width=0.33\linewidth, trim=8 8 8 8, clip]{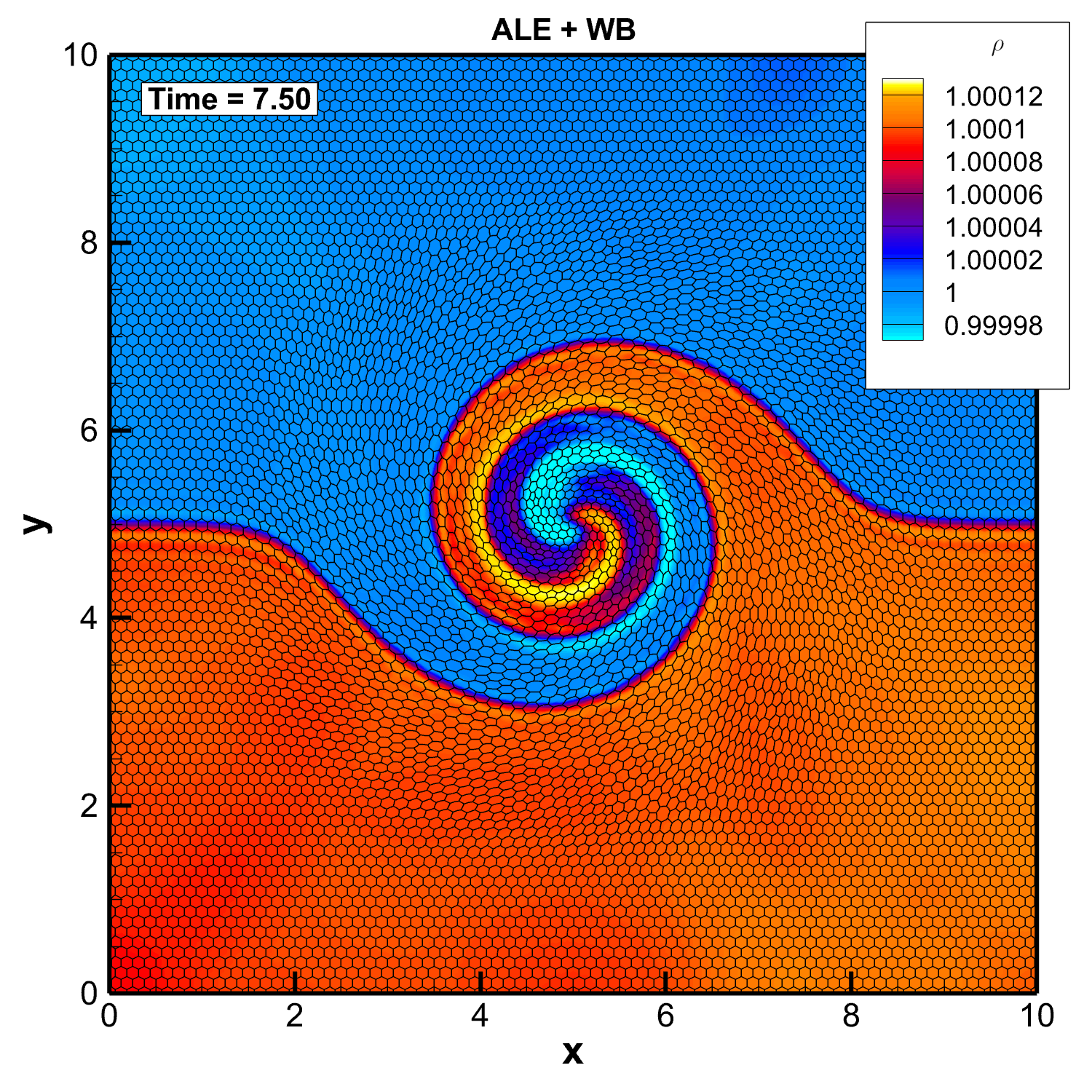}%
	\includegraphics[width=0.33\linewidth, trim=8 8 8 8, clip]{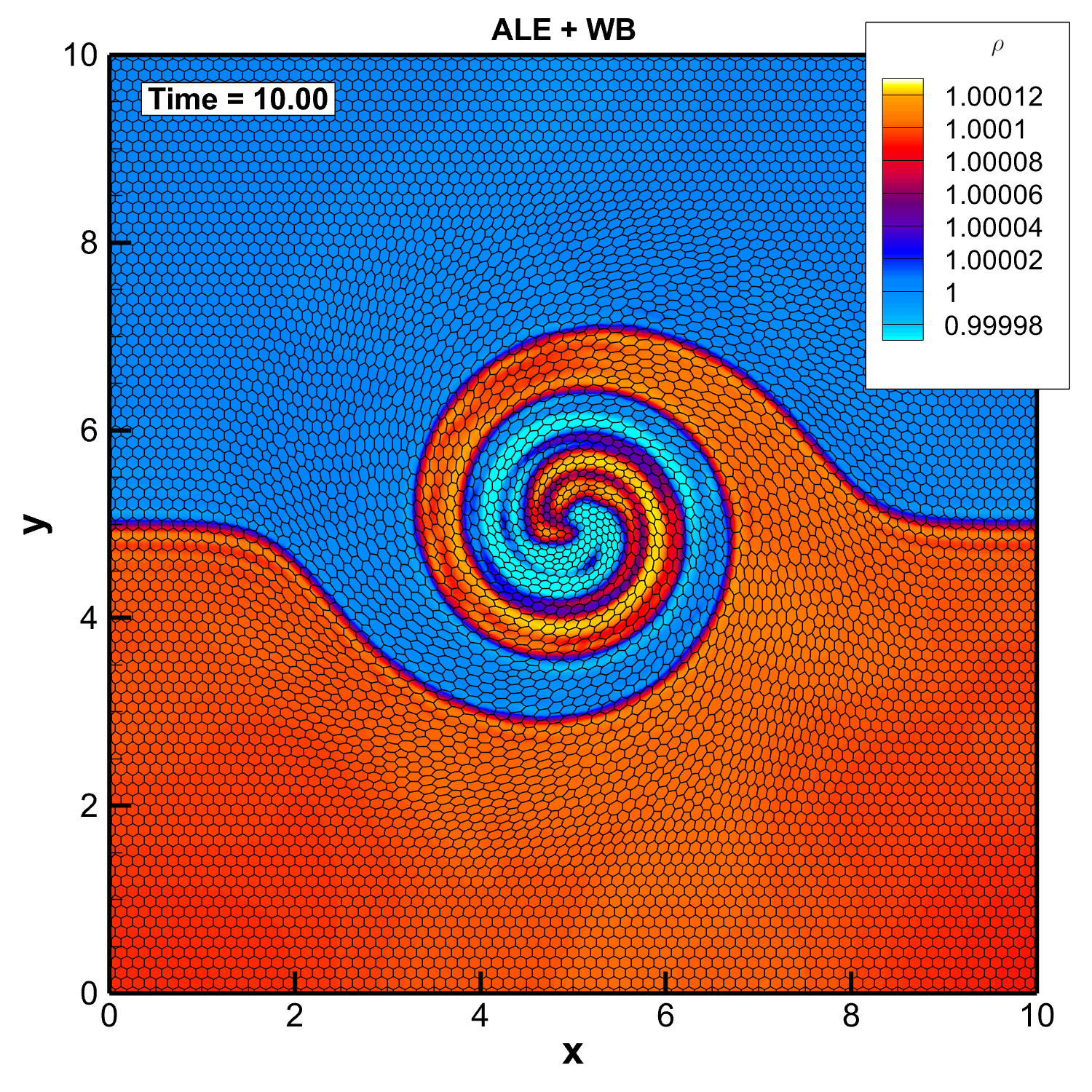}\\[-10pt]
	\caption{Frontogenesis on a MHD vortex solved with our well-balanced ALE DG scheme of order $3$ on 
	a mesh of average size $\bar{h} = 0.12$ (corresponding to 7579 polygonal elements) up to time $t=10$. 
	Thanks to the joint effect of our robust ALE framework and the well-balanced techniques 
	we can clearly model, even on a quite coarse mesh and for a long time, the evolution of a tiny front (of height $h_\epsilon=1e-4$) 
	without seeing the spurious effects of numerical errors.
	}
	\label{fig.MHDvortex_pert_alewb}
\end{figure}

\begin{figure}[!tb] \centering
	\includegraphics[width=0.333\linewidth, trim=8 8 8 8, clip]{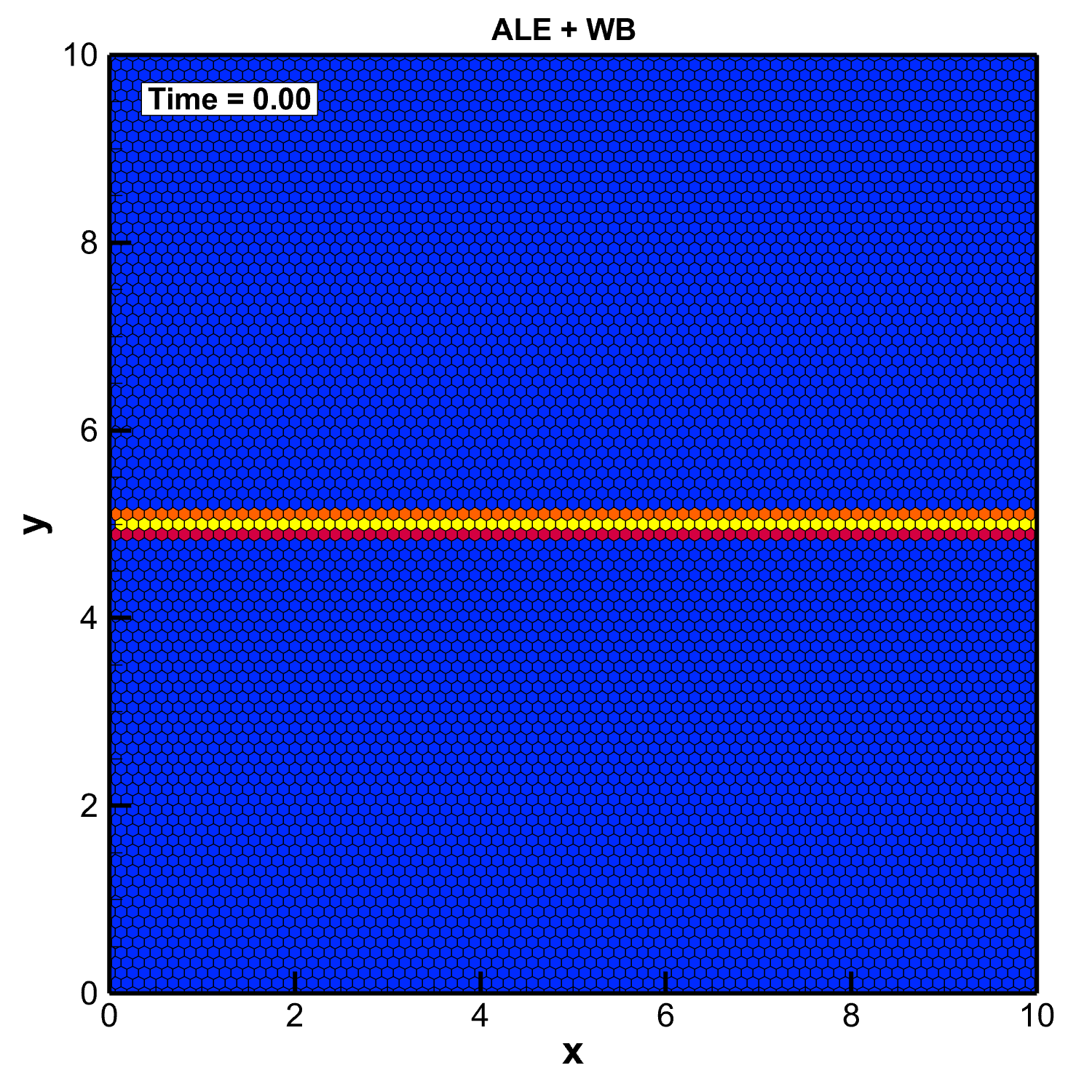}%
	\includegraphics[width=0.333\linewidth, trim=8 8 8 8, clip]{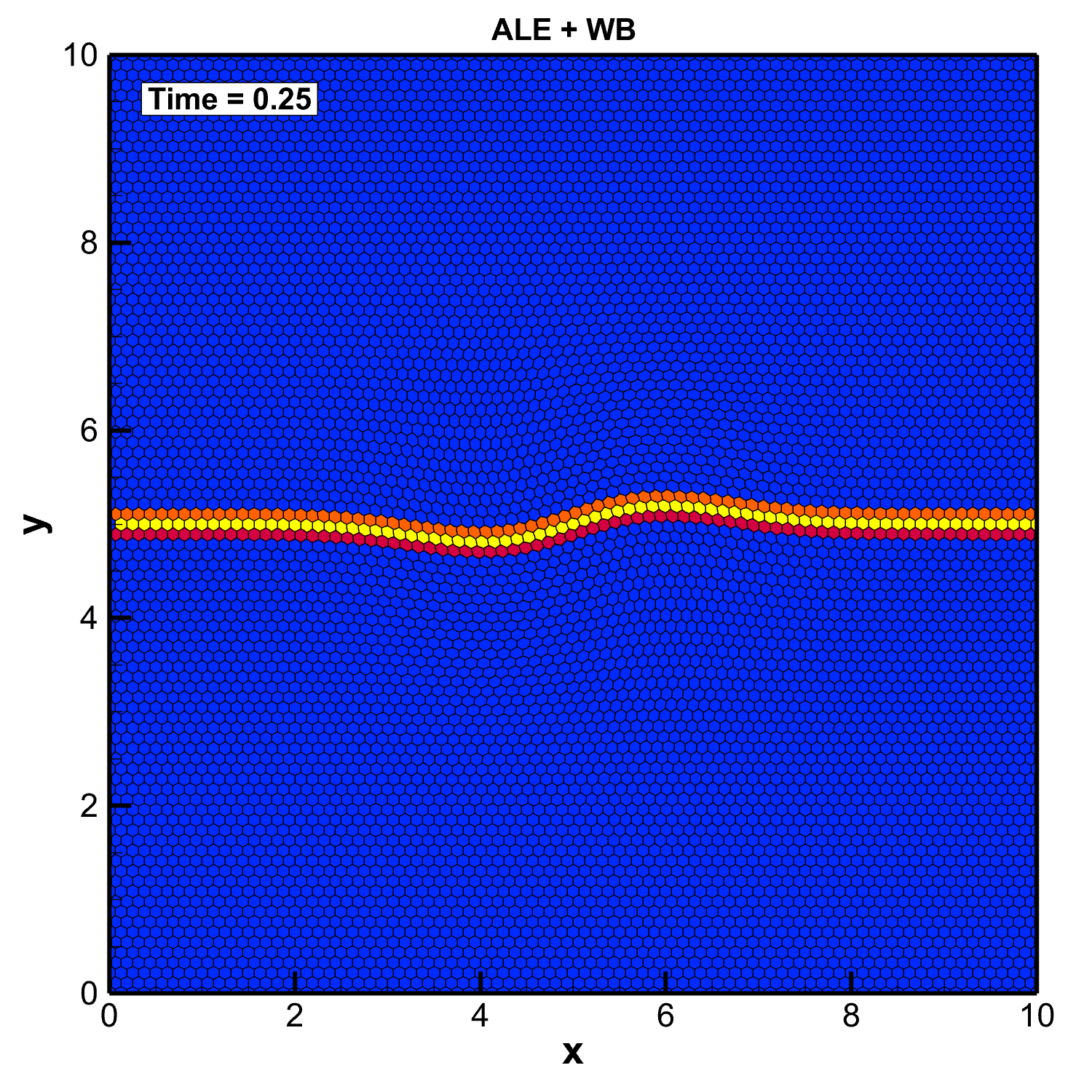}%
	\includegraphics[width=0.333\linewidth, trim=8 8 8 8, clip]{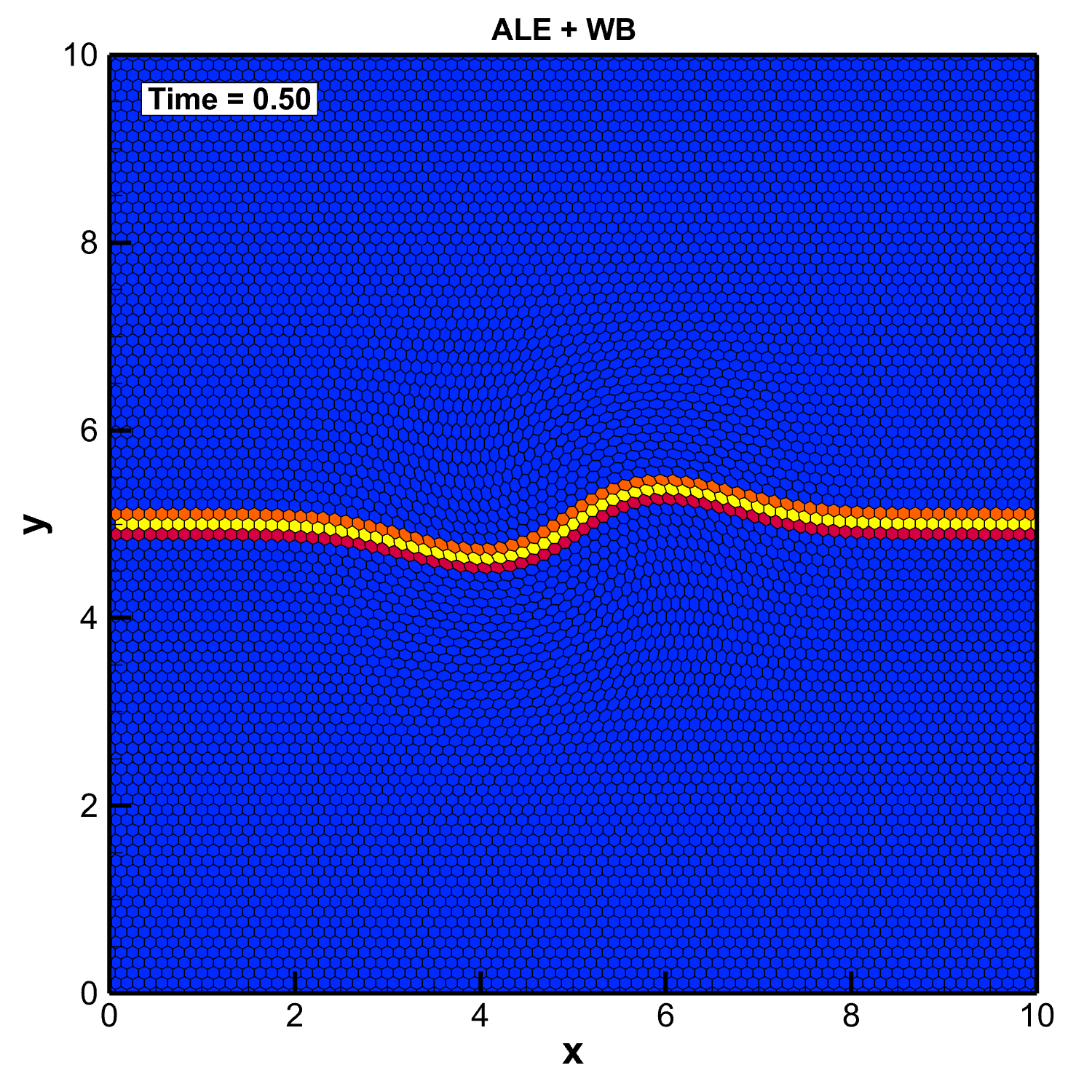}\\[0pt]	
	\includegraphics[width=0.333\linewidth, trim=8 8 8 8, clip]{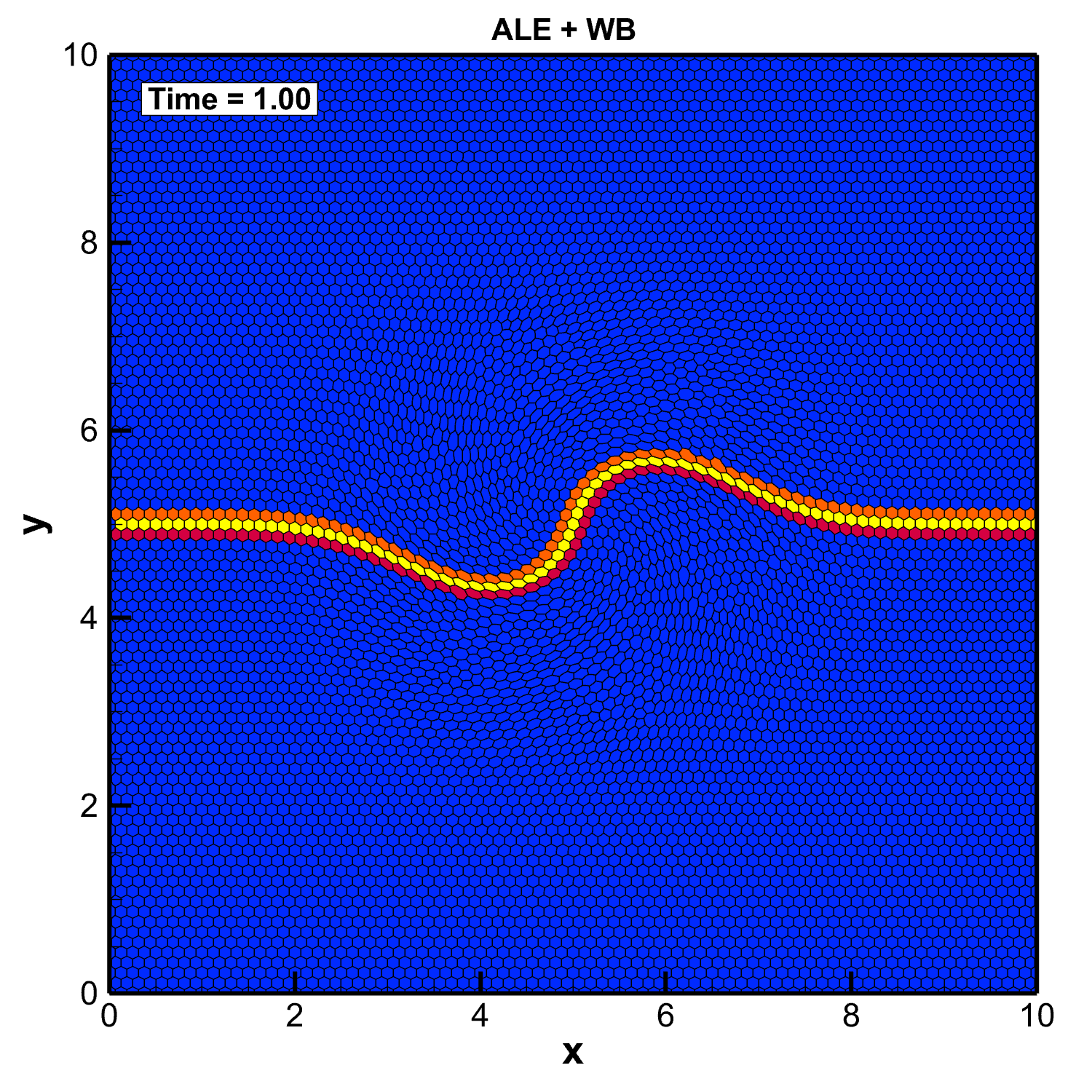}%
	\includegraphics[width=0.333\linewidth, trim=8 8 8 8, clip]{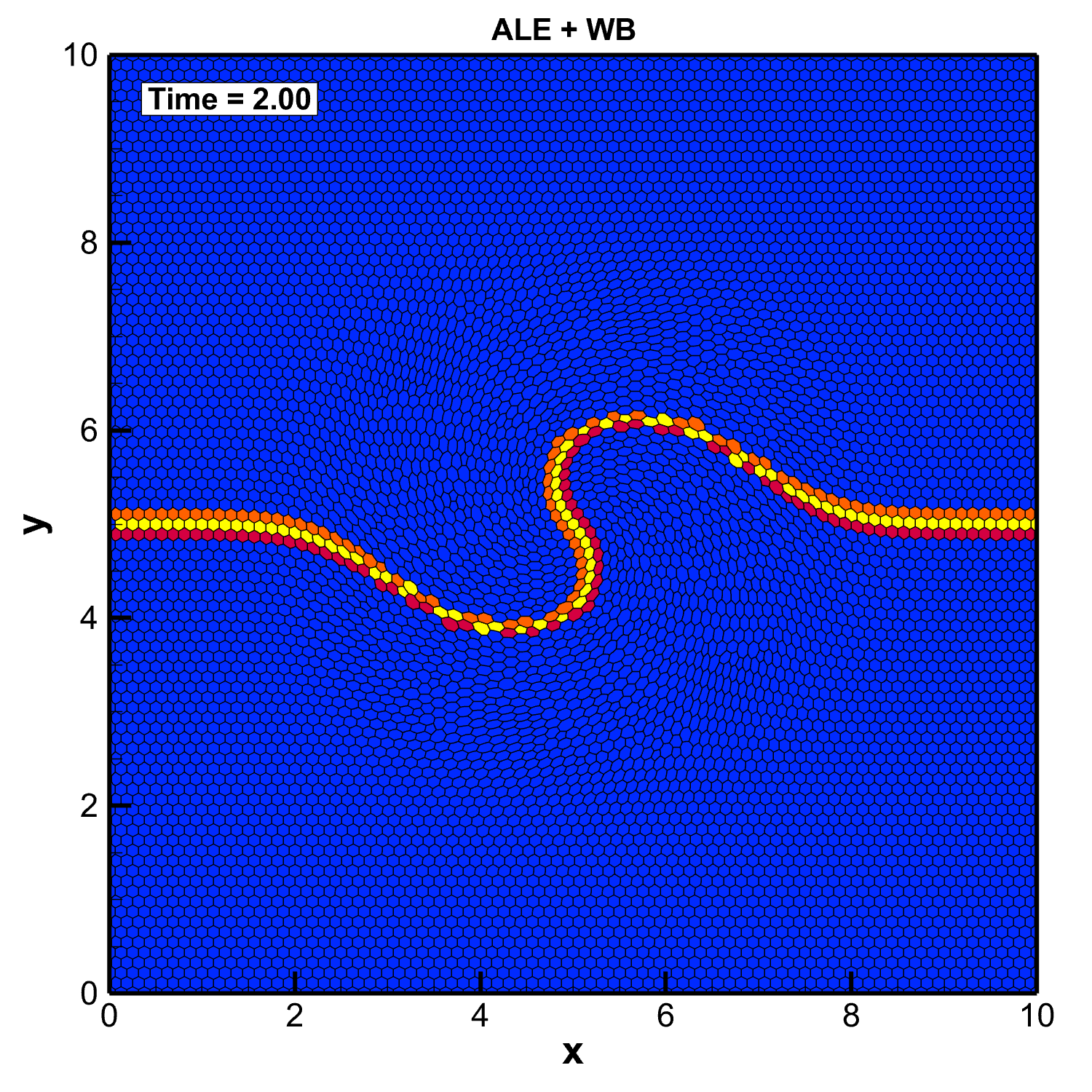}%
	\includegraphics[width=0.333\linewidth, trim=8 8 8 8, clip]{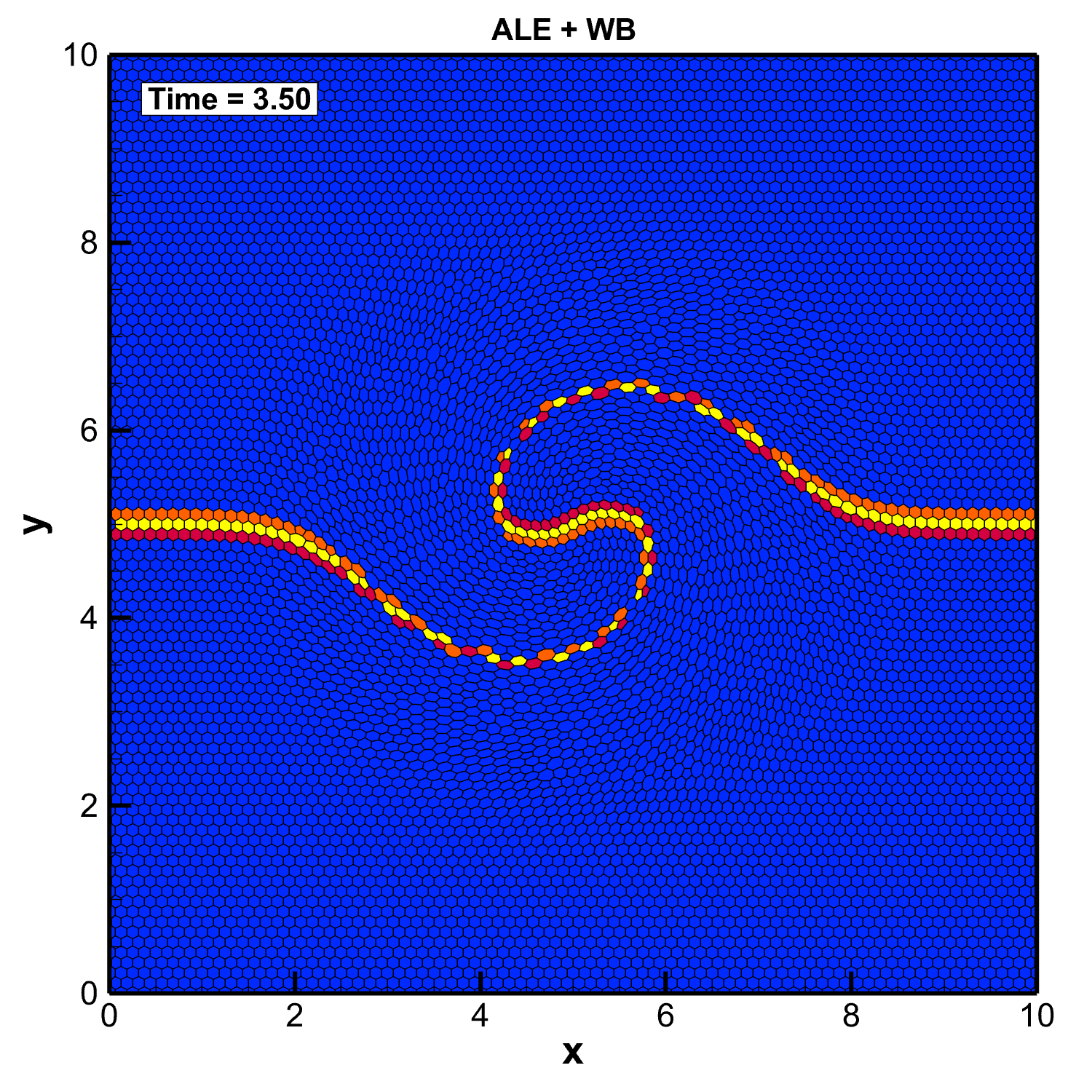}\\[10pt]
	\includegraphics[width=0.33\linewidth, trim=8 8 8 8, clip]{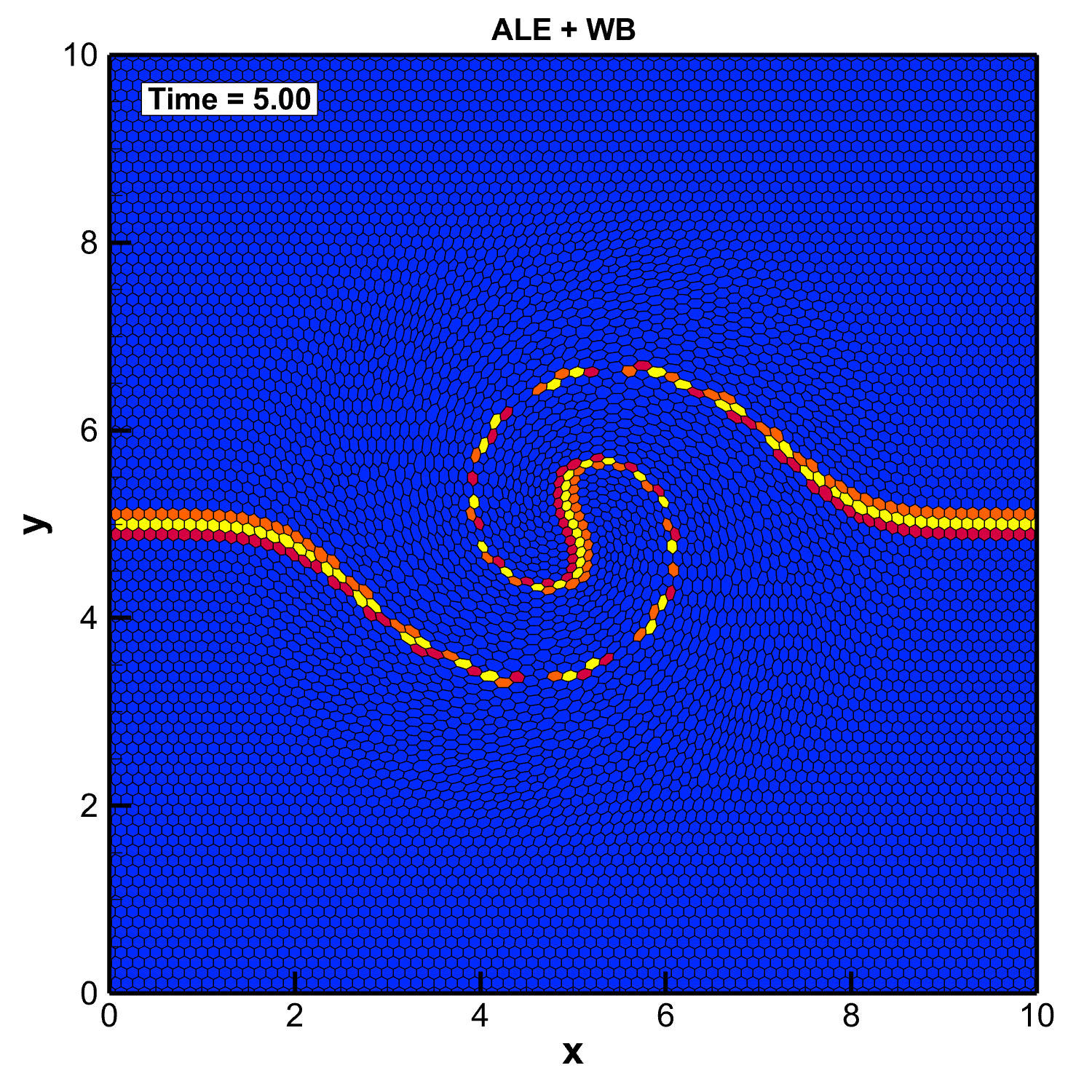}%
	\includegraphics[width=0.33\linewidth, trim=8 8 8 8, clip]{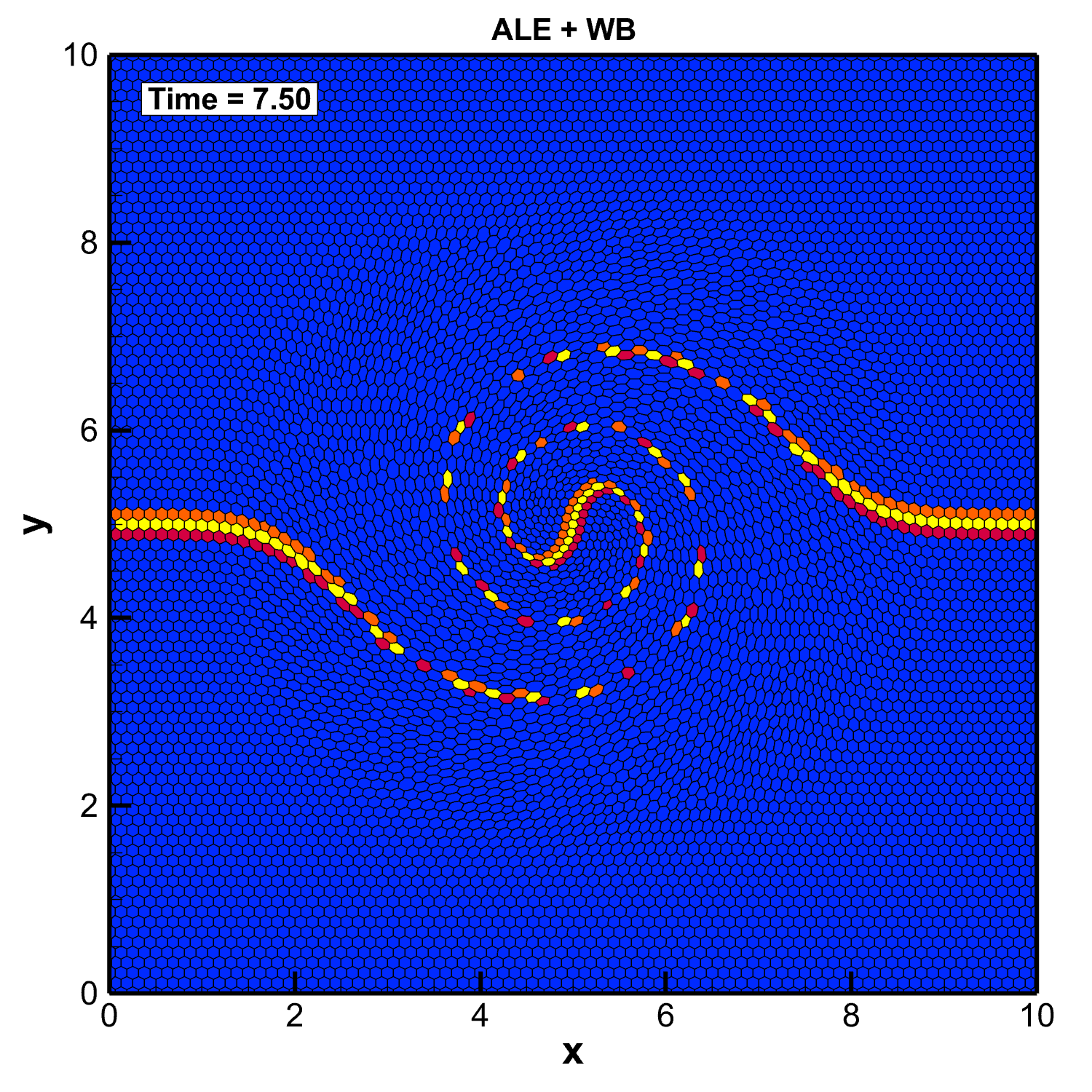}%
	\includegraphics[width=0.33\linewidth, trim=8 8 8 8, clip]{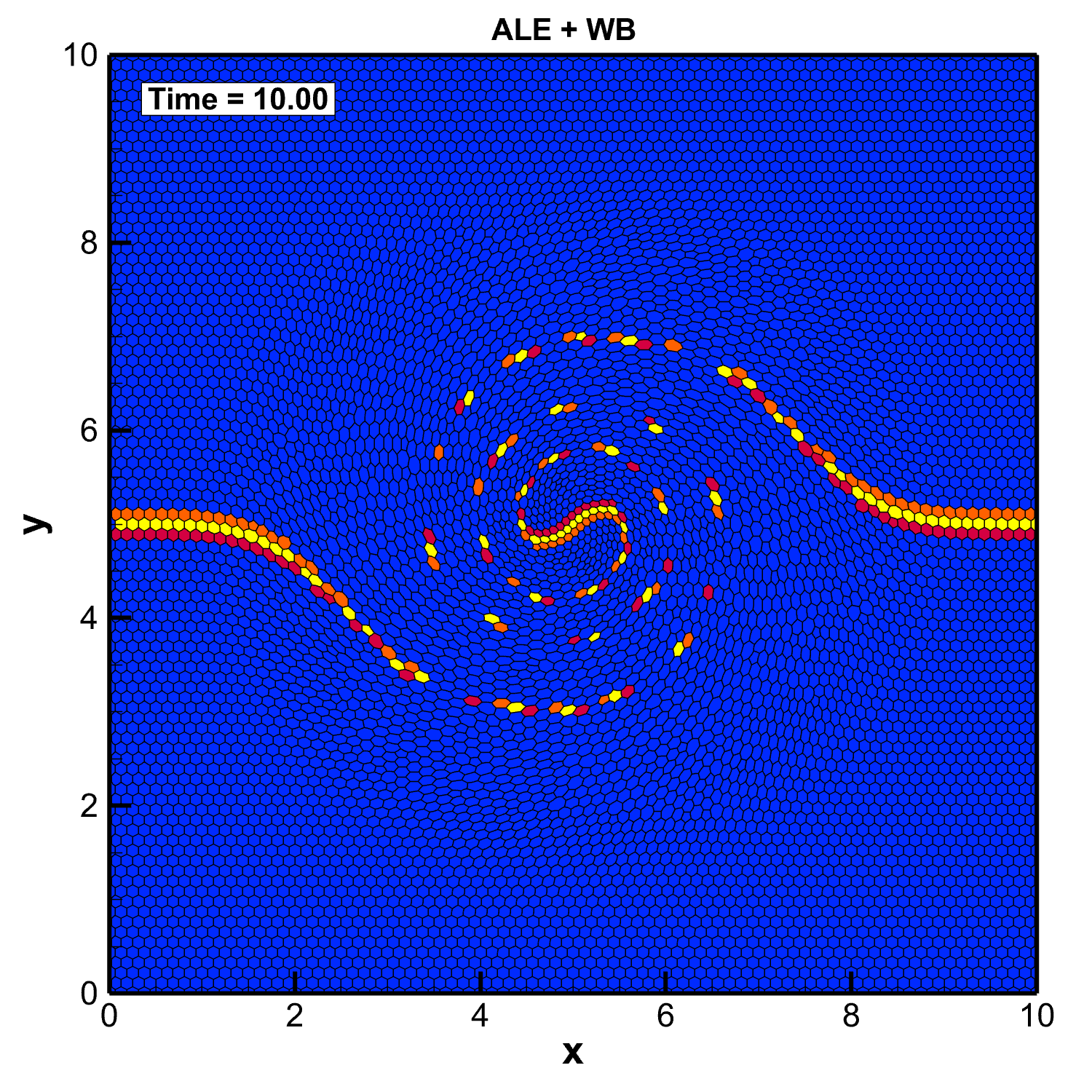}\\[-10pt]
	\caption{Frontogenesis on a MHD vortex. Here, we show the displacement of the elements located at the front interface up to time $t=10$. 
	The mesh closely follows the front interface and in such a way it considerably reduces the convection errors.
	}
	\label{fig.MHDvortex_pert_alewb_numbering}
\end{figure}

\begin{figure}[!tb]  \centering
	\includegraphics[width=0.33\linewidth, trim=8 8 8 8, clip]{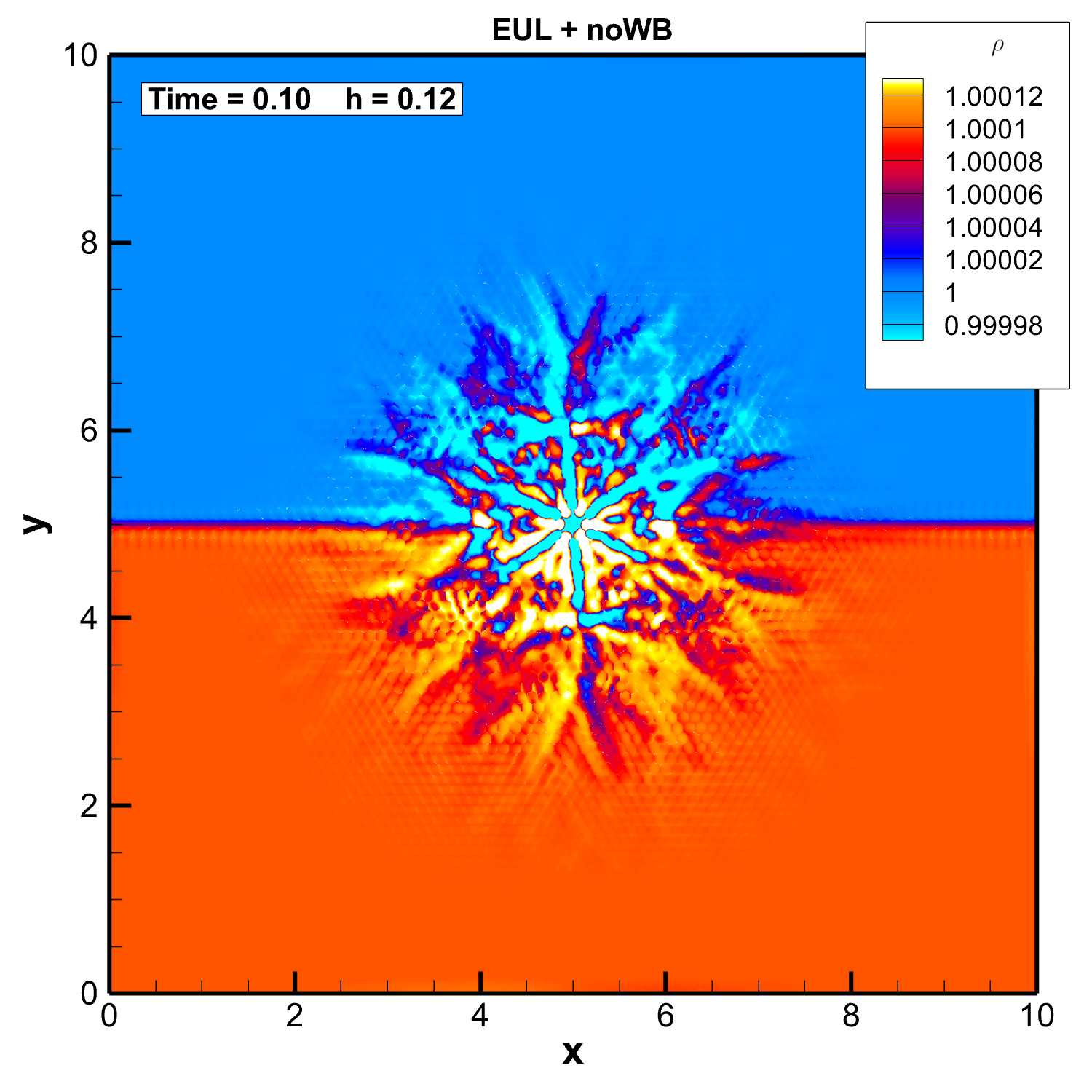}%
	\includegraphics[width=0.33\linewidth, trim=8 8 8 8, clip]{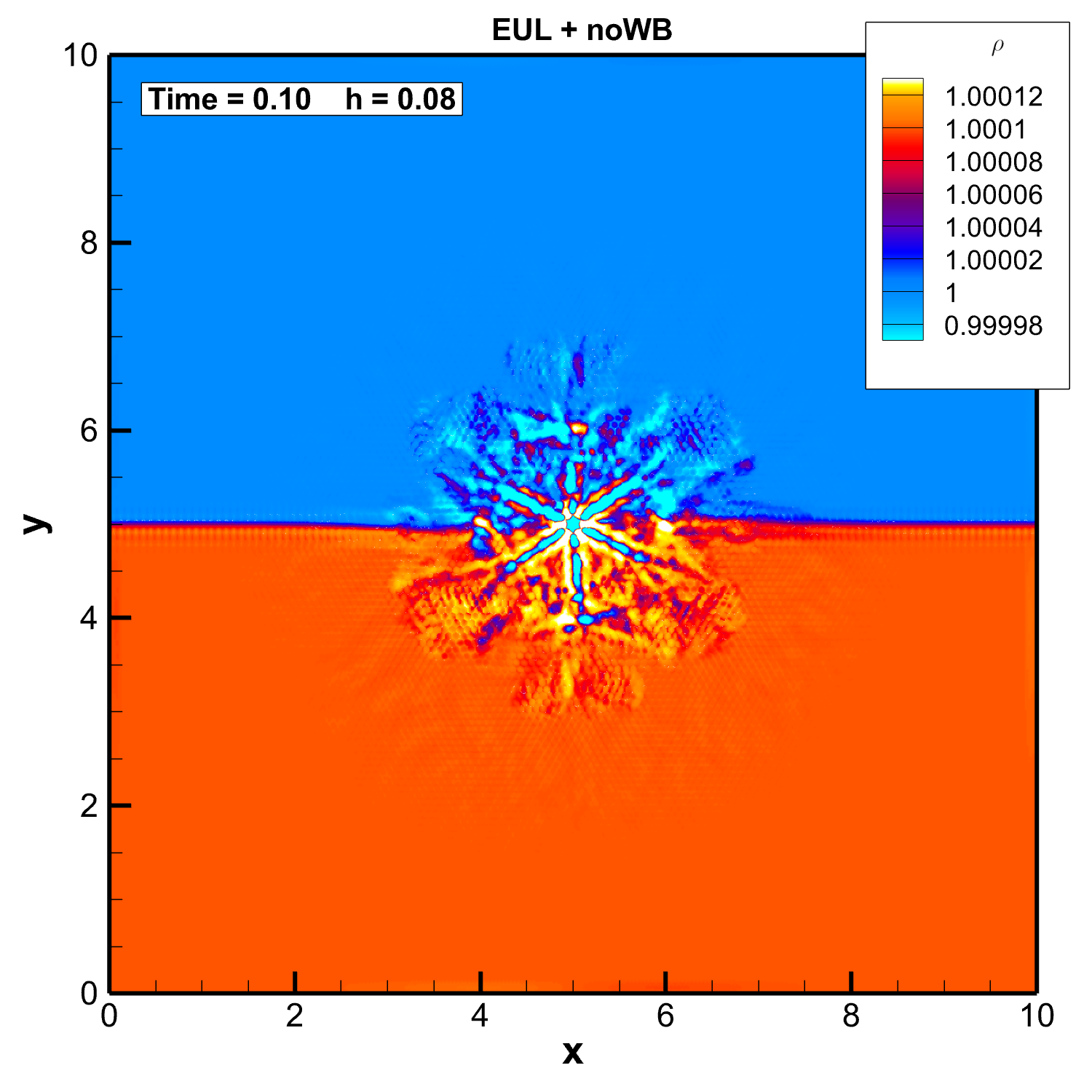}%
	\includegraphics[width=0.33\linewidth, trim=8 8 8 8, clip]{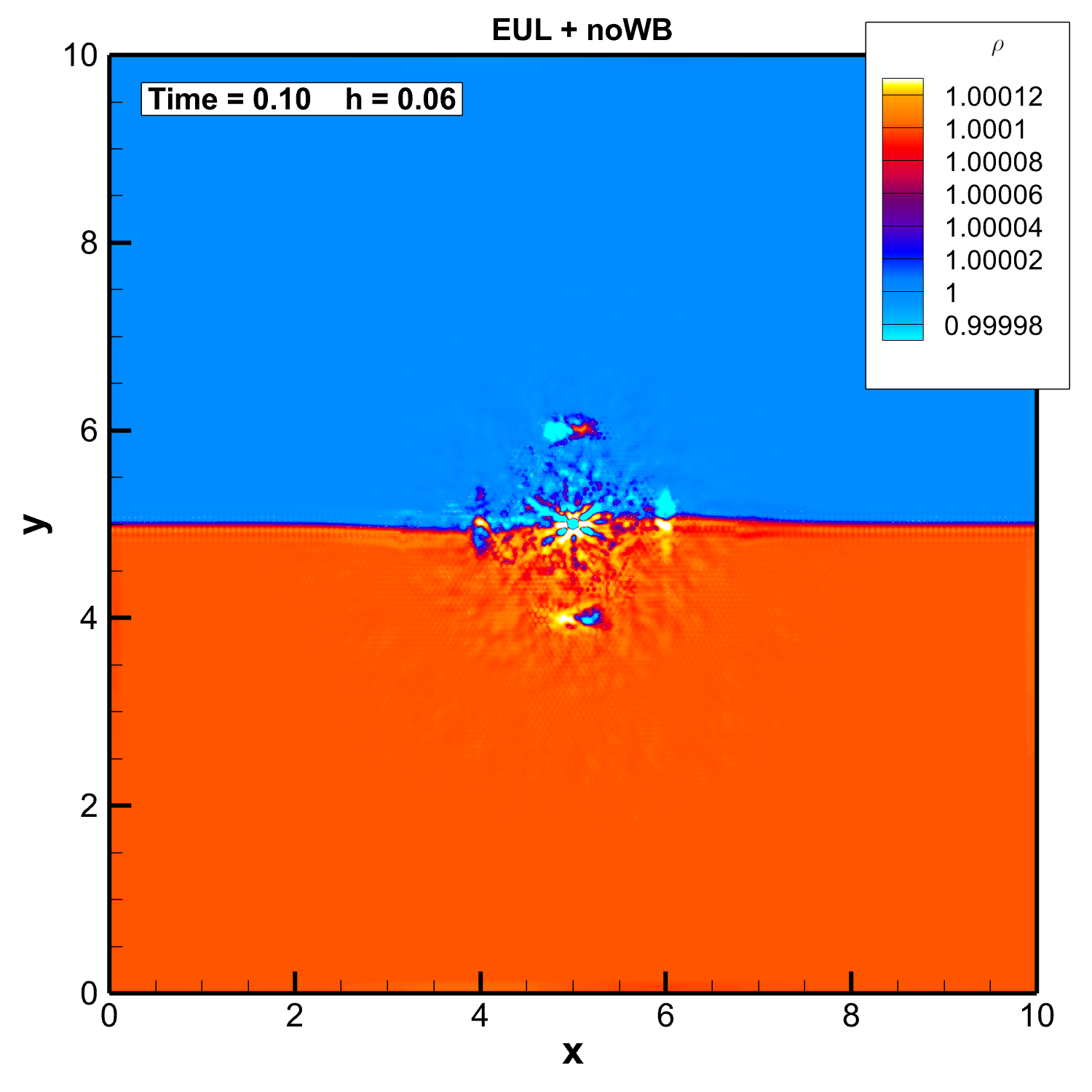}\\[0pt]
	\includegraphics[width=0.33\linewidth, trim=8 8 8 8, clip]{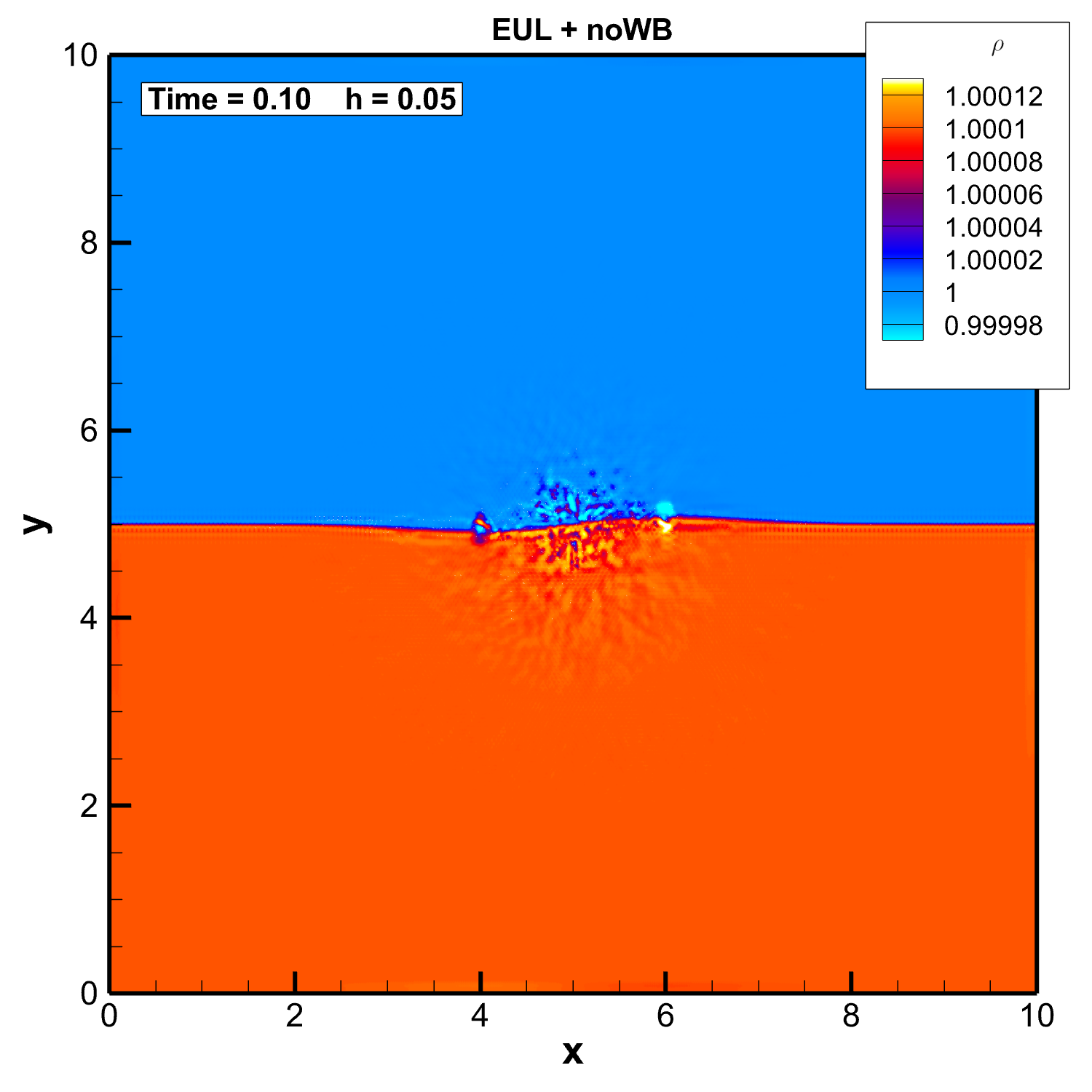}%
	\includegraphics[width=0.33\linewidth, trim=8 8 8 8, clip]{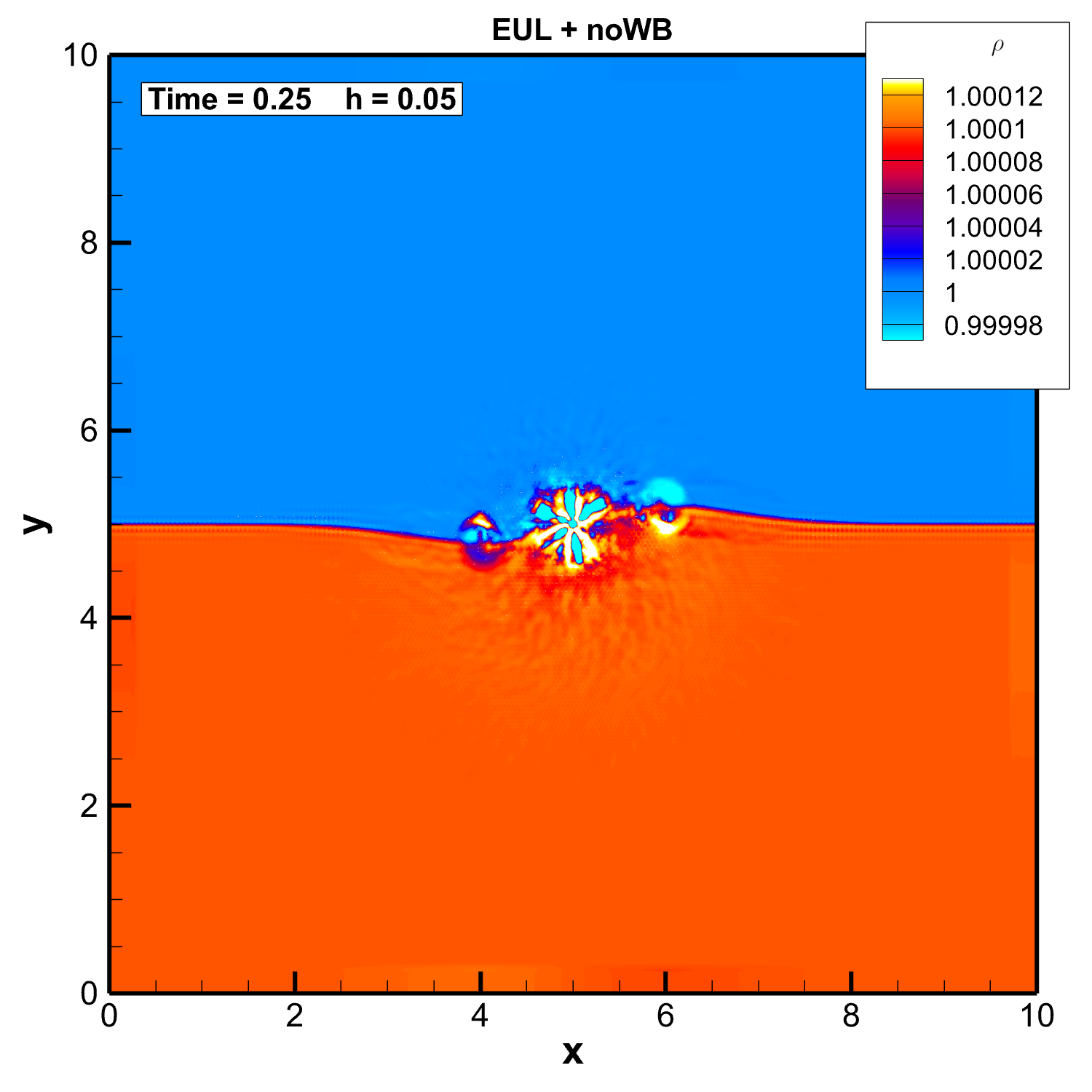}%
	\includegraphics[width=0.33\linewidth, trim=8 8 8 8, clip]{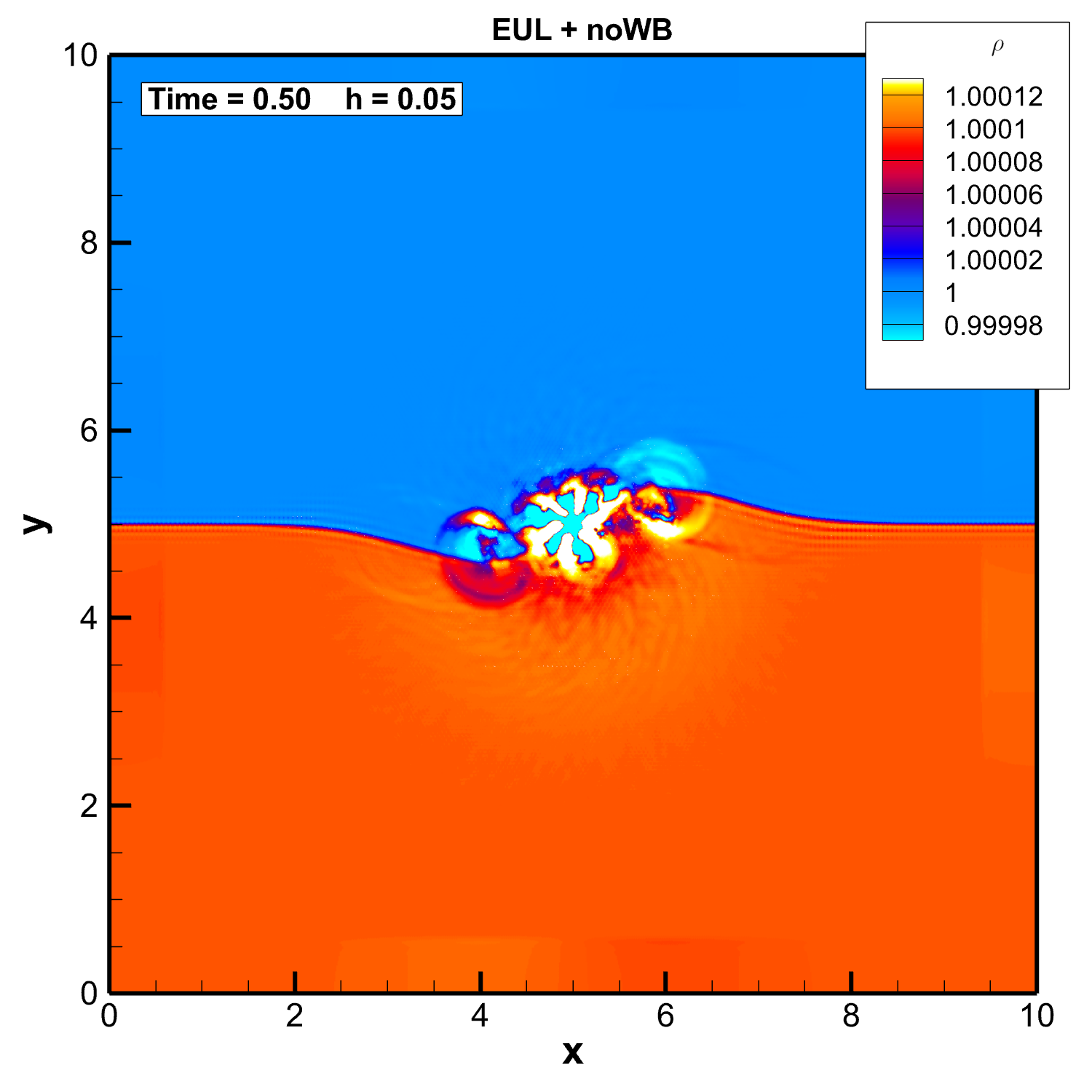}\\[0pt]
	\includegraphics[width=0.33\linewidth, trim=8 8 8 8, clip]{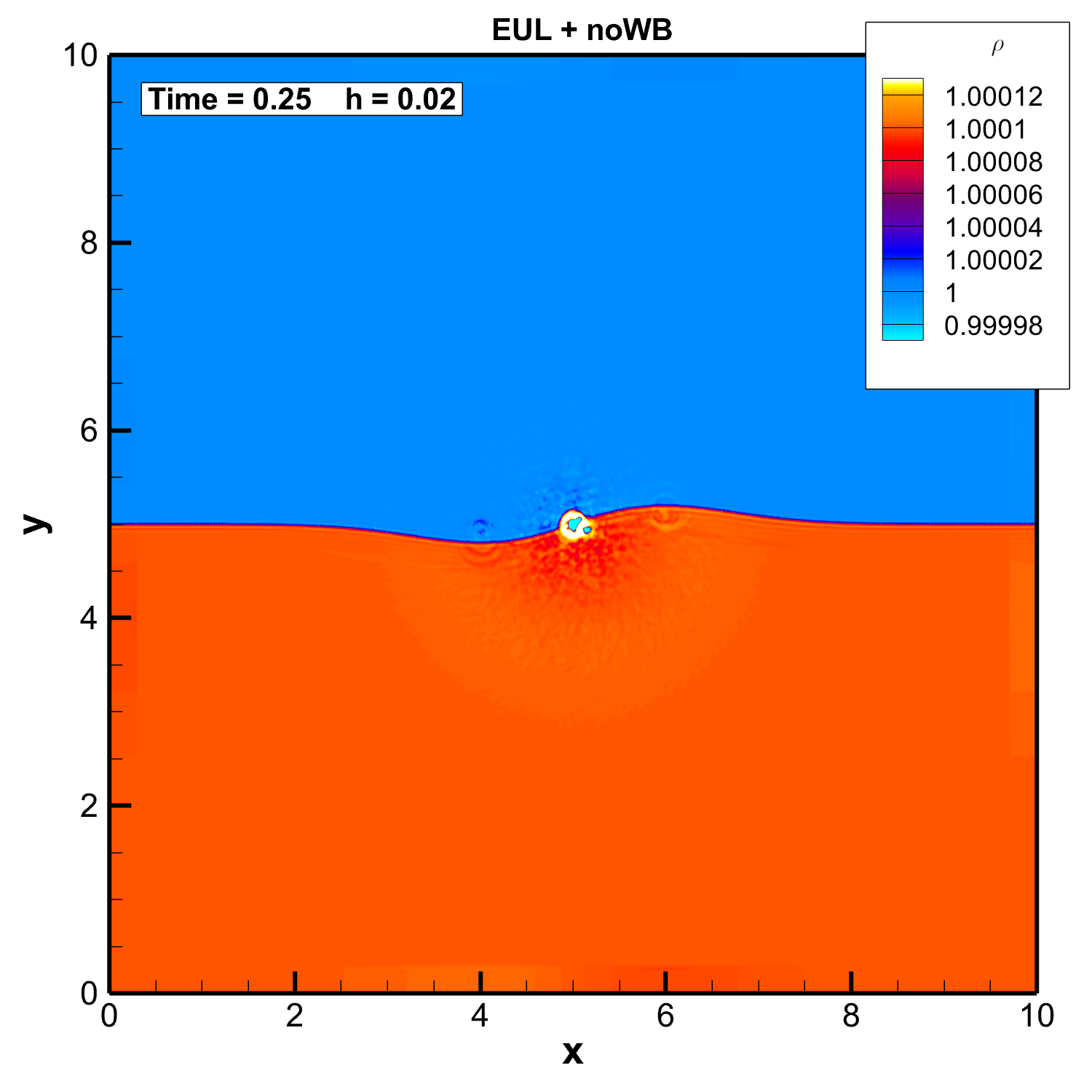}%
	\includegraphics[width=0.33\linewidth, trim=8 8 8 8, clip]{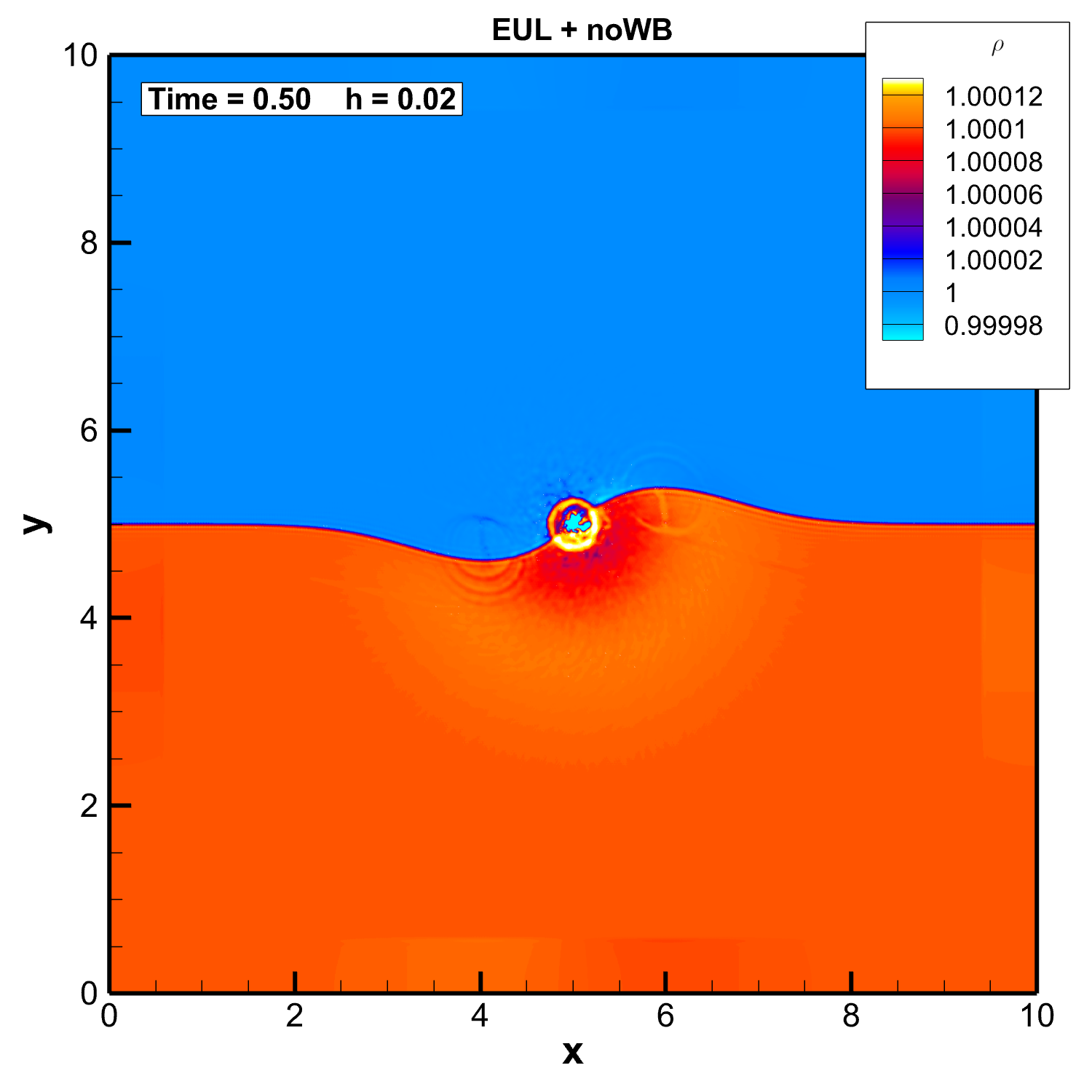}%
	\includegraphics[width=0.33\linewidth, trim=8 8 8 8, clip]{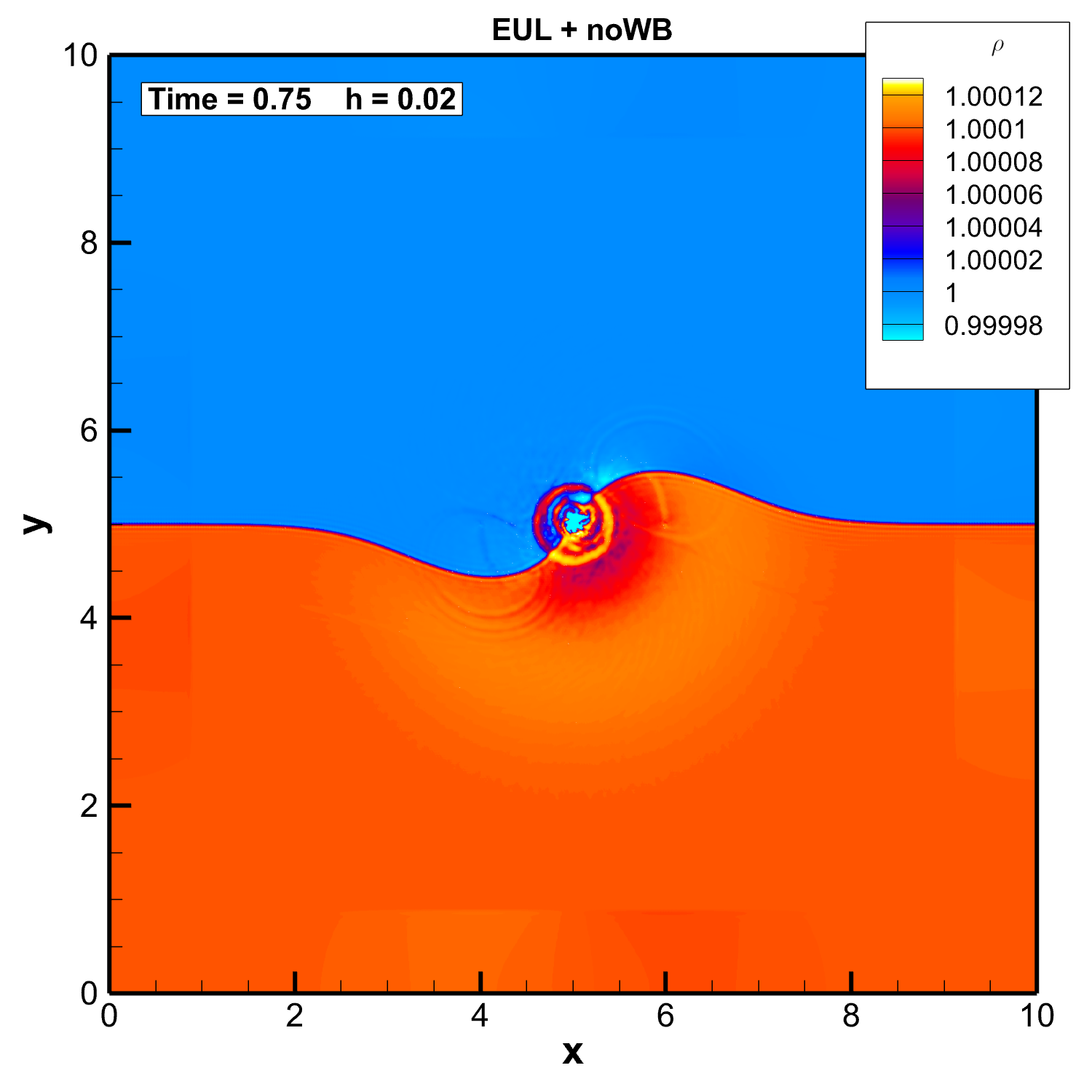}\\[-5pt]
	\caption{Frontogenesis on a MHD vortex solved with a standard Eulerian non well-balanced DG scheme of order $3$.
	Here we show the results obtained in the very initial part of the simulation ($t<1$) on increasingly 
	refined meshes: we start from $h=0.12$ (as in the simulation shown in the previous images and successfully 
	solved with our WB ALE DG scheme up to $t=10$) to $h= 0.02$. 
	We can clearly notice that, even with a very fine mesh (of $142.357$ polygonal elements), 
	the numerical errors accumulate and quickly deteriorate the results.
	}
	\label{fig.MHDvortex_pert_eul}
\end{figure}

Here, we set up the initial condition by taking the MHD vortex given in~\eqref{eqn.mhd3d.ic1} 
to which we add a perturbation of order 1E-4 to the density profile for $y<5$; 
the equilibrium to be preserved $\Q_E$ is taken equal to the MHD vortex~\eqref{eqn.mhd3d.ic1}.
We discretize the domain $\Omega=[0,10]\times[0,10]$ with a mesh made by $7579$ polygonal elements 
and we report the results obtained with our WB ALE DG scheme of order $3$ in 
Figures~\ref{fig.MHDvortex_pert_alewb} and~\ref{fig.MHDvortex_pert_alewb_numbering}. Here we can 
perfectly see the evolution of the density profile and of the mesh elements initially located at 
the interface between the perturbed and not perturbed region for a very long simulation time.  
The study of this tiny perturbation would not be possible with non well-balanced schemes, not
even for short times, since excessive numerical errors would otherwise develop almost immediately. 
This is clearly shown in Figure~\ref{fig.MHDvortex_pert_eul}, where we report the results obtained 
with a non well-balanced Eulerian DG scheme of order $3$ at the very initial moments 
of the simulation with finer and finer meshes. 
In all cases, the numerical errors associated with standard non well-balanced scheme, even if 
of high order of accuracy and on fine meshes, hide 
the frontogenesis evolution highlighting the necessity of well-balanced techniques for this type of
problems. 

\begin{table}[t] 
	\centering
	\numerikNine
	\begin{tabular}{|l|l||cccccc|}
		\hline
		\multicolumn{8}{|c|}{ MHD vortex with $\epsilon=\epsilon_E=5$ } \\
		\hline
		\hline
		& Time $t = $ & $0.0$ & $0.5$ & $1.0$ & $2.5$ & $5.0$ & $10.0$ \\
		\hline
		\hline
		\multirow{4}{*}{\rotatebox{90}{{DG-$\mathbb{P}_1$}}}
		& $L_2(\rho)$ error & 5.3166E-12 & 5.2877E-12	& 5.2655E-12	& 5.2015E-12	& 4.9641E-12	& 4.9723E-12 \\
		& $L_2(p)$ error  	& 0.0        & 5.6353E-12	& 5.7406E-12	& 5.6696E-12	& 4.8272E-12	& 4.9928E-12 \\
		& No. timestep 	    & 0     	 & 116       	& 221       	& 531       	& 1043	        & 2063       \\
		& No. sliver   	    & 0     	 & 0       	    & 13        	& 69        	& 160	        & 369       \\
		\cline{2-8}
		\hline
		\hline
		\multirow{4}{*}{\rotatebox{90}{{DG-$\mathbb{P}_2$}}}
		& $L_2(\rho)$ error & 5.2246E-12 & 5.2074E-12	& 5.1984E-12	& 5.2064E-12	& 5.0257E-12	& 5.0395E-12 \\
		& $L_2(p)$ error  	& 0.0        & 5.5197E-12	& 5.6551E-12	& 5.9331E-12	& 5.2350E-12	& 5.4774E-12 \\
		& No. timestep 	    & 0     	 & 217       	& 419       	& 1017       	& 2008	        & 3985       \\
		& No. sliver   	    & 0     	 & 0       	    & 4        	    & 56       	    & 141	        & 320       \\
		\cline{2-8}
		\hline
		\hline
		\multirow{4}{*}{\rotatebox{90}{{DG-$\mathbb{P}_3$}}}
		& $L_2(\rho)$ error & 5.1698E-12 & 5.1760E-12	& 5.1915E-12	& 5.2039E-12	& 5.2153E-12	& 5.2566E-12 \\
		& $L_2(p)$ error  	& 0.0        & 5.4693E-12	& 5.4849E-12	& 5.4933E-12	& 5.4448E-12	& 5.3596E-11 \\
		& No. timestep 	    & 0     	 & 360       	& 701       	& 1716       	& 3397	        & 6759      \\
		& No. sliver   	    & 0     	 & 0         	& 1        	    & 41        	& 110	        & 259       \\
		\cline{2-8}
		\hline
		\hline
		\multirow{4}{*}{\rotatebox{90}{{DG-$\mathbb{P}_4$}}}
		& $L_2(\rho)$ error & 5.1122E-12 & 5.1284E-12	& 5.1565E-12	& 5.2007E-12	& 5.2894E-12	& 5.3559E-12 \\
		& $L_2(p)$ error  	& 0.0        & 5.4239E-12	& 5.4440E-12	& 5.4853E-12	& 5.5523E-13	& 5.5692E-11 \\
		& No. timestep 	    & 0     	 & 517      	& 1010       	& 2478      	& 4912	        & 9787       \\
		& No. sliver   	    & 0     	 & 0       	    & 0          	& 33        	& 88	        & 232       \\
		\cline{2-8}
		\hline
	\end{tabular}
	\caption{Verification of the well-balanced property on the MHD vortex. 
		As in the previous test cases, we can notice that the equilibrium solution, even if initially 
		perturbed with a random error of 1E-12 distributed everywhere on the domain, is preserved with 
		machine accuracy for very long simulation times and after handling thousands of sliver elements. 
		This holds true for any employed polynomial representation order.}
	\label{table.MHDvortex}
\end{table}



\section{Conclusions and outlook to future works}
\label{sec.conclusions}

In this paper we have developed a novel high order accurate direct Arbitrary-Lagrangian-Eulerian 
discontinuous Galerkin scheme which is \textit{well-balanced} for any \textit{a priori} known 
equilibrium solution of the studied system of hyperbolic PDEs. 
In particular, the chosen equilibrium can be both stationary or time-dependent and it can be 
available in an analytical form or just in a discrete way. 
Moreover the scheme is robust in presence of discontinuities (of the solution, while the equilibrium is always assumed to be smooth)
thanks to the use of an \textit{a posteriori} subcell finite volume limiting strategy, implemented within
the direct ALE context and also made well-balanced.

We have proven the accuracy and robustness of the proposed methodology 
on a wide set of benchmarks, 
clearly demonstrating the capabilities of the novel 
combination of well-balancing and Lagrangian motion in high order DG schemes.
Indeed, after the preliminary work forwarded in~\cite{gaburro2017direct,gaburro2018well}, 
the WB ALE approach has been now generalized to arbitrary high order of accuracy and 
arbitrarily moving polygonal meshes, accounting even for topology changes.
This allows to perform simulations for extremely long times, even when affected by strong 
shear flows or vortical flows, 
while always keeping numerical errors extremely low for near-equilibrium flows. 
We also would like to emphasize the role of the \textit{hole-like sliver} elements, already 
introduced in~\cite{gaburro2020high}, in guaranteeing the {conservativity} and the high 
order of accuracy around an arbitrary topology change, 
whose treatment here has been generalized to the 
well-balanced framework.

\medskip 

The future directions of our work will be found along two main 
axes: one concerning the development of the numerical method, 
the second of applicative nature.

With regards to method development, indeed, 
we will work on the continuous extension of our 
direct ALE framework on moving polygonal meshes with topology changes trying  
\textit{i)}  to improve the mesh optimization techniques applied at each timestep in such a way 
to further ameliorate the Lagrangian character of the algorithm while maintaining a high 
quality of the moving mesh, inspired by~\cite{knupp2002reference,d2016optimization,anderson2018high,dobrev2020simulation}; 
\textit{ii)} to introduce an adaptive mesh refinement strategy
to refine/coarsen the mesh where necessary, 
that will be based on insertion and deletion (see~\cite{guillet2019high}) of generator
points
and the {subsequent formation} of new types of \textit{hole-like} elements;
and \textit{iii)} to improve the employed \textit{a posteriori} finite volume limiter 
considering different subgrids, using for example methodology like~\cite{antonietti2022refinement}.

Concerning applications, 
we plan to employ the present algorithm in more complex practical contexts 
as for example 
\textit{i)} for simulations in the field of continuum mechanics and fluid-structure interaction, 
relying on the general first order hyperbolic model by 
Godunov, Peshkov and Romenski in~\cite{PeshRom2014,HPRmodel,HPRmodelMHD} and 
further studied in~\cite{jackson2019numerical,boscheri2022cell,muller2023sharp,gabriel2021unified,chiocchetti2023gpr},
\textit{ii)} combing the ALE methodology with a diffuse interface approach as those 
forwarded in~\cite{BaerNunziato1986,AbgrallSaurel,DIM2D,DIM3D,RomenskiTwoPhase2010,kemm2020simple, kemm2020simple,chiocchetti2021high}, 
\textit{iii)} and also in challenging astrophysical scenarios.
We remark that the simulations shown in the present work (for example the Keplerian disks) 
can be seen as a Newtonian prototype of interesting applications in the field of astrophysics, 
as the study of neutron star oscillations and black holes characteristics.
Thus, in our future work we plan to apply the proposed methodology to more complex system of 
equations which also include general relativity, like the GRMHD 
model~\cite{del2007echo,bugli2014dynamo,fambri2018ader} and the first order reductions of the 
Einstein field equations presented in~\cite{dumbser2018conformal,dumbser2020glm, gaburro2021well,dumbser2023WBGR}.

\section{Acknowledgments}

E.~Gaburro is member of the INdAM GNCS group in Italy. 
E.~Gaburro gratefully acknowledges the support received from the European Union 
with the ERC Starting Grant \textit{ALcHyMiA} (grant agreement No. 101114995).
Views and opinions expressed are however those of the author only and do not necessarily 
reflect those of the European Union or the European Research Council Executive Agency. 
Neither the European Union nor the granting authority can be held responsible for them.

Finally, E.~Gaburro would like to thank Dr. Simone Chiocchetti for the valuable and fruitful discussions regarding the conception of the method. 

\noindent \section*{In memoriam}

\noindent This paper is dedicated to the memory of Prof. Arturo Hidalgo L\'opez
($^*$July 03\textsuperscript{rd} 1966 - $\dagger$August 26\textsuperscript{th} 2024) of the Universidad Politecnica de Madrid, organizer of HONOM 2019 and active participant in many other editions of HONOM. 
Our thoughts and wishes go to his wife Lourdes and his sister Mar\'ia Jes\'us, whom he left behind.

\bibliographystyle{plain}
\bibliography{./referencesAleVoronoiWB.bib}

\end{document}